\documentclass{amsbook}
\usepackage{amssymb}
\usepackage{mathrsfs}
\usepackage{stmaryrd} 
\usepackage{bbm}
\usepackage{oldgerm}
\usepackage[frenchb]{babel}
\usepackage[T1]{fontenc}
\usepackage[latin1]{inputenc}
\usepackage[all]{xy}
\usepackage{hyperref}
\usepackage{upref}

\newtheorem{teo}[subsection]{Théorème}
\newtheorem{prop}[subsection]{Proposition}
\newtheorem{cor}[subsection]{Corollaire}
\newtheorem{lem}[subsection]{Lemme}

\theoremstyle{definition}

\newtheorem{defi}[subsection]{Définition}
\newtheorem{rema}[subsection]{Remarque}
\newtheorem{remas}[subsection]{Remarques}

\newcommand{\gtimes}{\stackrel{\leftarrow}{\times}}

\mathchardef\mhyphen="2D

\newcommand{\mA}{{\mathbb A}}

\newcommand{\mR}{{\mathbb R}}

\newcommand{\mQ}{{\mathbb Q}}
\newcommand{\mL}{{\mathbb L}}
\newcommand{\mN}{{\mathbb N}}
\newcommand{\mP}{{\mathbb P}}

\newcommand{\mZ}{{\mathbb Z}}
\newcommand{\mF}{{\mathbb F}}
\newcommand{\mG}{{\mathbb G}}

\newcommand{\mU}{{\mathbb U}}
\newcommand{\mV}{{\mathbb V}}

\newcommand{\bA}{{\bf A}}
\newcommand{\bB}{{\bf B}}

\newcommand{\bD}{{\bf D}}

\newcommand{\bL}{{\bf L}}
\newcommand{\bT}{{\bf T}}

\newcommand{\Et}{{\bf \acute{E}t}}
\newcommand{\Sch}{{\bf Sch}}
\newcommand{\Ens}{{\bf Ens}}
\newcommand{\Pt}{{\bf Pt}}
\newcommand{\Mon}{{\bf Mon}}

\newcommand{\Top}{{\bf Top}}
\newcommand{\bHom}{{\bf Hom}}

\newcommand{\bAlg}{{\bf Alg}}
\newcommand{\bMod}{{\bf Mod}}
\newcommand{\aAlg}{{\alpha\mhyphen \bf Alg}}
\newcommand{\aMod}{{\alpha\mhyphen \bf Mod}}

\newcommand{\aI}{{\alpha\mhyphen \bf I}}

\newcommand{\acA}{{\alpha\mhyphen \cA}}

\newcommand{\ahcC}{{\alpha\mhyphen \hcC}}
\newcommand{\atcC}{{\alpha\mhyphen \tcC}}
\newcommand{\ahE}{{\alpha\mhyphen \hE}}
\newcommand{\atE}{{\alpha\mhyphen \tE}}

\newcommand{\et}{{\rm \acute{e}t}}
\newcommand{\fet}{{\rm f\acute{e}t}}

\newcommand{\coim}{{\rm coim}}

\newcommand{\zar}{{\rm zar}}

\newcommand{\coh}{{\rm coh}}
\newcommand{\qcoh}{{\rm qcoh}}
\newcommand{\cart}{{\rm cart}}

\newcommand{\rf}{{\rm f}}
\newcommand{\dR}{{\rm dR}}

\newcommand{\Spec}{{\rm Spec}}

\newcommand{\ob}{{\rm Ob}}

\newcommand{\pr}{{\rm pr}}

\newcommand{\coker}{{\rm coker}}

\newcommand{\im}{{\rm im}}

\newcommand{\p}{{\rm p}}

\newcommand{\gp}{{\rm gp}}
\newcommand{\id}{{\rm id}}
\newcommand{\Tr}{{\rm Tr}}

\newcommand{\Sym}{{\rm Sym}}

\newcommand{\Hom}{{\rm Hom}}

\newcommand{\End}{{\rm End}}

\newcommand{\Gal}{{\rm Gal}}

\newcommand{\AR}{{\rm AR}}

\newcommand{\rE}{{\rm E}}
\newcommand{\rF}{{\rm F}}
\newcommand{\rH}{{\rm H}}
\newcommand{\rT}{{\rm T}}
\newcommand{\rM}{{\rm M}}

\newcommand{\rR}{{\rm R}}
\newcommand{\rS}{{\rm S}}
\newcommand{\rV}{{\rm V}}
\newcommand{\rW}{{\rm W}}

\newcommand{\ra}{{\rm a}}

\newcommand{\rp}{{\rm p}}

\newcommand{\oF}{{\overline{F}}}
\newcommand{\oK}{{\overline{K}}}

\newcommand{\oR}{{\overline{R}}}
\newcommand{\oS}{{\overline{S}}}

\newcommand{\oU}{{\overline{U}}}

\newcommand{\oX}{{\overline{X}}}
\newcommand{\oY}{{\overline{Y}}}

\newcommand{\oa}{{\overline{a}}}
\newcommand{\of}{{\overline{f}}}

\newcommand{\op}{{\overline{p}}}

\newcommand{\ou}{{\overline{u}}}
\newcommand{\ov}{{\overline{v}}}

\newcommand{\ox}{{\overline{x}}}
\newcommand{\oy}{{\overline{y}}}

\newcommand{\oeta}{{\overline{\eta}}}

\newcommand{\oiota}{{\overline{\iota}}}

\newcommand{\ocB}{{\overline{\cB}}}

\newcommand{\ofp}{{\overline{\fp}}}
\newcommand{\otta}{{\overline{\tta}}}

\newcommand{\uE}{{\underline{E}}}

\newcommand{\uX}{{\underline{X}}}
\newcommand{\uY}{{\underline{Y}}}
\newcommand{\uf}{{\underline{f}}}

\newcommand{\urp}{{\underline{\rp}}}

\newcommand{\urho}{{\underline{\rho}}}
\newcommand{\uphi}{{\underline{\phi}}}
\newcommand{\uPhi}{{\underline{\Phi}}}
\newcommand{\uXi}{{\underline{\Xi}}}

\newcommand{\hA}{{\widehat{A}}}
\newcommand{\hB}{{\widehat{B}}}

\newcommand{\hE}{{\widehat{E}}}
\newcommand{\hP}{{\widehat{P}}}

\newcommand{\hRun}{{\widehat{R_1}}}

\newcommand{\hmZ}{{\widehat{\mZ}}}

\newcommand{\halpha}{\widehat{\alpha}}
\newcommand{\hsigma}{\widehat{\sigma}}

\newcommand{\cA}{{\mathscr A}}
\newcommand{\cB}{{\mathscr B}}
\newcommand{\cC}{{\mathscr C}}
\newcommand{\cD}{{\mathscr D}}
\newcommand{\cE}{{\mathscr E}}
\newcommand{\cF}{{\mathscr F}}
\newcommand{\cG}{{\mathscr G}}

\newcommand{\cK}{{\mathscr K}}
\newcommand{\cL}{{\mathscr L}}
\newcommand{\cP}{{\mathscr P}}
\newcommand{\co}{{\mathscr O}}
\newcommand{\cR}{{\mathscr R}}
\newcommand{\cS}{{\mathscr S}}
\newcommand{\cT}{{\mathscr T}}
\newcommand{\cH}{{\mathscr H}}
\newcommand{\cM}{{\mathscr M}}
\newcommand{\cN}{{\mathscr N}}
\newcommand{\cQ}{{\mathscr Q}}

\newcommand{\cZ}{{\mathscr  Z}}

\newcommand{\cHom}{{\mathscr Hom}}

\newcommand{\fA}{{\mathfrak A}}

\newcommand{\fC}{{\mathfrak C}}
\newcommand{\fD}{{\mathfrak D}}

\newcommand{\fF}{{\mathfrak F}}
\newcommand{\fG}{{\mathfrak G}}

\newcommand{\fM}{{\mathfrak M}}

\newcommand{\fS}{{\mathfrak S}}
\newcommand{\fT}{{\mathfrak T}}
\newcommand{\fV}{{\mathfrak V}}
\newcommand{\fW}{{\mathfrak W}}

\newcommand{\fc}{{\mathfrak c}}

\newcommand{\fm}{{\mathfrak m}}

\newcommand{\fp}{{\mathfrak p}}

\newcommand{\tta}{{\tt a}}

\newcommand{\hK}{{\widehat{K}}}
\newcommand{\hM}{{\widehat{M}}}

\newcommand{\hR}{{\widehat{R}}}
\newcommand{\hX}{{\widehat{X}}}
\newcommand{\hY}{{\widehat{Y}}}

\newcommand{\hoR}{{\widehat{\oR}}}

\newcommand{\hcC}{{\widehat{\cC}}}
\newcommand{\hcD}{{\widehat{\cD}}}

\newcommand{\bvR}{{\breve{R}}}

\newcommand{\bvvarphi}{{\breve{\varphi}}}

\newcommand{\coS}{{\check{\oS}}}

\newcommand{\coX}{{\check{\oX}}}
\newcommand{\coY}{{\check{\oY}}}

\newcommand{\tE}{{\widetilde{E}}}
\newcommand{\tF}{{\widetilde{F}}}
\newcommand{\tuE}{{\widetilde{\uE}}}

\newcommand{\tQ}{{\widetilde{Q}}}

\newcommand{\tT}{{\widetilde{T}}}
\newcommand{\tU}{{\widetilde{U}}}

\newcommand{\tX}{{\widetilde{X}}}
\newcommand{\tY}{{\widetilde{Y}}}

\newcommand{\tx}{{\widetilde{x}}}
\newcommand{\ty}{{\widetilde{y}}}

\newcommand{\tcD}{{\widetilde{\cD}}}
\newcommand{\trT}{{\widetilde{\rT}}}
\newcommand{\tOmega}{{\widetilde{\Omega}}}

\newcommand{\tbeta}{{\widetilde{\beta}}}
\newcommand{\tsigma}{{\widetilde{\sigma}}}

\newcommand{\ttau}{{\widetilde{\tau}}}

\newcommand{\talpha}{{\widetilde{\alpha}}}

\newcommand{\tcC}{{\widetilde{\cC}}}

\newcommand{\trH}{{\widetilde{\rH}}}
\newcommand{\trM}{{\widetilde{\rM}}}

\newcommand{\tfD}{{\widetilde{\fD}}}
\newcommand{\trF}{{\widetilde{\rF}}}

\begin{document}

\title{La suite spectrale de Hodge-Tate}
\author{Ahmed Abbes et Michel Gros}
\address{A.A. Laboratoire Alexander Grothendieck, ERL 9216 du CNRS, Institut des Hautes \'Etudes Scientifiques, 35 route de Chartres, 91440 Bures-sur-Yvette, France}
\address{M.G. CNRS UMR 6625, IRMAR, Université de Rennes 1,
Campus de Beaulieu, 35042 Rennes cedex, France}
\email{abbes@ihes.fr}
\email{michel.gros@univ-rennes1.fr}

\begin{abstract}
La construction de la suite spectrale de Hodge-Tate d'une variété propre et lisse sur un corps $p$-adique nous sert de fil conducteur pour revenir, 
en suppléant nombre de détails, sur l'approche de Faltings en théorie de Hodge $p$-adique. 
Cette suite spectrale tire son origine de la suite spectrale de Cartan-Leray pour 
la projection canonique du topos de Faltings sur le topos étale d'un modèle entier de la variété. 
L'aboutissement de cette dernière se calcule grâce au principal théorème de comparaison de Faltings
d'où dérivent tous les autres théorèmes de comparaison entre la cohomologie étale $p$-adique et d'autres cohomologies $p$-adiques.
Son terme initial est lié aux faisceaux des formes différentielles par un morphisme évoquant l'isomorphisme de Cartier.

\vspace{5mm}

\noindent{\sc Abstract.} 
The Hodge-Tate spectral sequence for a proper smooth variety over a $p$-adic field provides a framework for us to  revisit Faltings' approach to $p$-adic Hodge 
theory and to fill in many details. The spectral sequence is obtained from the Cartan-Leray spectral sequence for the canonical projection from the Faltings 
topos to the étale topos of an integral model of the variety.  
Its abutment is computed by Faltings' main comparison theorem from which derive all comparison theorems between $p$-adic étale cohomology
and other $p$-adic cohomologies, and its initial term is related to the sheaf of differential forms by a construction reminiscent of the Cartier isomorphism.
\end{abstract}

\maketitle

\setcounter{tocdepth}{1}
\tableofcontents

\renewcommand{\thesubsection}{\arabic{subsection}}

\numberwithin{equation}{subsection}

\chapter*{Introduction}

\subsection{}\label{intro1}
Le but de ce travail est de consolider et de compléter les fondations de la théorie de Hodge $p$-adique suivant l'approche initiée par Faltings dans \cite{faltings1,faltings2}. 
Parmi celles-ci, Gabber et Ramero se sont concentrés sur  la théorie des presque-anneaux et sur le {\em théorème de presque-pureté} \cite{gr,gr1} 
et Tsuji a donné dans (\cite{agt} V) une présentation concise de la théorie des {\em extensions presque-étales}, plus proche de l'approche originale de Faltings. 
Ce sont à d'autres aspects de la théorie de Faltings que nous nous intéressons ici, à savoir les applications de son théorème de presque-pureté à des calculs de  cohomologie galoisienne, le recollement de ceux-ci grâce au {\em topos de Faltings} et les théorèmes de comparaison entre divers groupes de cohomologie $p$-adique.
Nous avons entamé ce projet dans \cite{agt}, où nous avons repris les aspects de la théorie nécessaires au développement de la correspondance de Simpson $p$-adique,
en particulier ceux concernant le topos de Faltings (\cite{agt} VI). Nous poursuivons le travail dans ce volume en nous focalisant sur deux aspects plus avancés, à savoir
\begin{itemize}
\item[(i)] le principal théorème de comparaison de Faltings;
\item[(ii)] la suite spectrale de Hodge-Tate.
\end{itemize}

\subsection{}\label{intro2}
Soient $K$ un corps de valuation discrète complet de caractéristique
$0$, à corps résiduel parfait de caractéristique $p>0$, $\co_K$ l'anneau de valuation de $K$, 
$\oK$ une clôture algébrique de $K$, $\co_\oK$ la clôture intégrale de $\co_K$ dans $\oK$, $G_K$ le groupe de Galois de $\oK$ sur $K$.
On note $\co_C$ le séparé complété $p$-adique de $\co_\oK$, $\fm_C$ son idéal maximal et $C$ son corps des fractions.
Rappelons l'énoncé du théorème de {\em décomposition de Hodge-Tate} conjecturé par Tate (\cite{tate} Remark page 180)
et démontré indépendamment par Faltings \cite{faltings1,faltings2} et Tsuji \cite{tsuji1}. 
Pour tout $K$-schéma propre et lisse $X$ et tout entier $n\geq 0$, il existe un isomorphisme canonique fonctoriel et $G_K$-équivariant 
\begin{equation}\label{intro2a}
\rH^n_\et(X_\oK,\mQ_p)\otimes_{\mQ_p}C\stackrel{\sim}{\rightarrow}\oplus_{i=0}^n\rH^i(X_C,\Omega^{n-i}_{X_C/C})(i-n).
\end{equation}
Cette énoncé équivaut à l'existence d'une suite spectrale 
canonique fonctorielle et $G_K$-équivariante, la {\em suite spectrale de Hodge-Tate},
\begin{equation}\label{intro2b}
E_2^{i,j}=\rH^i(X_C,\Omega^j_{X_C/C})(-j)\Rightarrow \rH^{i+j}_\et(X_C,\mQ_p)\otimes_{\mQ_p}C.
\end{equation}
L'équivalence des deux énoncés est une conséquence du théorème de Tate sur la cohomologie galoisienne des tordus de $C$ (\cite{tate} Theorem~2 page 176): 
la nullité de $\rH^0(G_K,C(j))$ pour tout entier $j\not=0$ implique la dégénérescence de la suite spectrale en $E_2$ et la nullité  
de $\rH^1(G_K,C(j))$ pour tout entier $j\not=0$ implique sa décomposition. 

Lorsque $X$ admet un modèle semi-stable sur $\co_K$, nous déduirons la suite spectrale \eqref{intro2b} 
d'une suite spectrale de Cartan-Leray pour la projection canonique du topos de Faltings du modèle dans son topos étale.
Nous suivrons la preuve donnée dans \cite{faltings2} à quelques améliorations et simplifications près. 
Notre approche permet entre autres de dégager une analogie entre la suite spectrale de Hodge-Tate et 
la suite spectrale conjuguée en caractéristique $p$. 

\subsection{}\label{intro3}
Faltings a explicitement introduit son topos dans \cite{faltings2} comme le cadre naturel de son approche pour la théorie de Hodge $p$-adique. 
Ce topos répond à deux objectifs. D'une part, il permet de recoller des calculs locaux de cohomologie galoisienne, 
et d'autre part, il permet de comparer diverses cohomologies $p$-adiques. Ce dernier aspect ressort plus clairement de notre interprétation 
du topos de Faltings comme un {\em topos co-évanescent généralisé}, généralisant les produits orientés de topos de Deligne (\cite{agt} VI).  
Rappelons ici sa définition (cf. \ref{tf1}). 

On pose $S=\Spec(\co_K)$ et $\oS=\Spec(\co_\oK)$ et 
on note $s$ (resp.  $\eta$, resp. $\oeta$) le point fermé de $S$ (resp.  générique de $S$, resp. générique de $\oS$).
Pour tout entier $n\geq 1$, on pose $S_n=\Spec(\co_K/p^n\co_K)$. Pour tout $S$-schéma $Y$, on pose 
\begin{equation}\label{intro3a}
\oY=Y\times_S\oS \ \ \ {\rm et}\ \ \ Y_n=Y\times_SS_n.
\end{equation} 

Changeant de notation, $X$ désigne dans la suite de l'introduction un $S$-schéma lisse de type fini. 
Nous considérerons dans ce travail une situation logarithmique lisse plus générale (cf. \ref{TFA3}).
Toutefois, pour simplifier la présentation, nous nous limitons dans cette introduction au cas lisse au sens usuel.   
On désigne par $\Et_{/X}$ le site étale de $X$ et par $X_\et$ son topos étale. On note $E$ la catégorie des morphismes $(V\rightarrow U)$
au-dessus du morphisme canonique $X_\oeta\rightarrow X$, autrement dit des diagrammes commutatifs 
\begin{equation}\label{intro3b}
\xymatrix{V\ar[r]\ar[d]&U\ar[d]\\
X_\oeta\ar[r]&X}
\end{equation}
tel que $U$ soit étale sur $X$ et que le morphisme canonique $V\rightarrow U_\oeta$ soit {\em fini étale}. 
Il est utile de considérer la catégorie $E$ comme fibrée par le foncteur
\begin{equation}\label{intro3c}
\pi\colon E\rightarrow \Et_{/X}, \ \ \ (V\rightarrow U)\mapsto U.
\end{equation}
La fibre de $\pi$ au-dessus d'un objet $U$ de $\Et_{/X}$ est canoniquement équivalente à la catégorie $\Et_{\rf/U_\oeta}$ des revêtements étales de $U_\oeta$. 
On l'équipe de la topologie étale et on note $U_{\oeta,\fet}$ le topos associé (cf. \ref{notconv10}).

On équipe $E$ de la topologie {\em co-évanescente}, c'est-à-dire de la topologie engendrée par les recouvrements 
$\{(V_i\rightarrow U_i)\rightarrow (V\rightarrow U)\}_{i\in I}$ des deux types suivants~:
\begin{itemize}
\item[(v)] $U_i=U$ pour tout $i\in I$ et $(V_i\rightarrow V)_{i\in I}$ est un recouvrement;
\item[(c)] $(U_i\rightarrow U)_{i\in I}$ est un recouvrement et $V_i=V\times_UU_i$ pour tout $i\in I$. 
\end{itemize}
On montre alors que se donner un faisceau $F$ sur $E$ est équivalent à se donner~:
\begin{itemize}
\item[(i)] pour tout objet $U$ de $\Et_{/X}$, un faisceau $F_U$ de $U_{\oeta,\fet}$, à savoir la restriction de $F$ à la fibre
de $\pi$ au-dessus de $U$;
\item[(ii)] pour tout morphisme $f\colon U'\rightarrow U$ de $\Et_{/X}$, un morphisme $\gamma_f\colon F_U\rightarrow f_{\oeta *}(F_{U'})$. 
\end{itemize}
Ces données sont sujettes à des conditions de cocycle (pour la composition des morphismes) et 
à des conditions de recollement (pour les recouvrements de $\Et_{/X}$). 

\subsection{}\label{intro4}
Les foncteurs 
\begin{eqnarray}
\sigma^+\colon \Et_{/X}\rightarrow E, && U\mapsto (U_\oeta\rightarrow U),\label{intro4a}\\
\beta^+\colon \Et_{\rf/X_\oeta}\rightarrow E,&& V\mapsto (V\rightarrow X),\label{intro4b}\\
\psi^+\colon E\rightarrow\Et_{/X_\oeta},&&(V\rightarrow U)\mapsto V,\label{intro4c}
\end{eqnarray}
sont exacts à gauche et continus \eqref{tf1}. Il définissent donc trois morphismes de topos 
\begin{eqnarray}
\sigma\colon \tE&\rightarrow& X_\et,\label{intro4d}\\
\beta\colon \tE&\rightarrow& X_{\oeta,\fet},\label{intro4e}\\
\psi\colon X_{\oeta,\et}&\rightarrow&\tE.\label{intro4f}
\end{eqnarray}
Les deux premiers sont des analogues des {\em projections canoniques} pour les produits orientés de topos de Deligne, et le troisième est un analogue du 
{\em morphisme des cycles (co-)proches}. On rappelle que Deligne a introduit les produits orientés de topos pour généraliser les cycles évanescents 
pour les faisceaux étales relativement à une base de dimension $>1$. 

Pour tout groupe abélien fini $G$, notant encore $G$ le faisceau constant de valeur $G$ de $X_{\oeta,\et}$ ou de $\tE$, 
on a un isomorphisme canonique $G\stackrel{\sim}{\rightarrow} \psi_*(G)$ (\cite{agt} VI.10.9(iii)). 

\begin{prop}[cf. \ref{acycloc2}]\label{intro5}
Pour tout faisceau abélien de torsion localement constant constructible $F$ de $X_{\oeta,\et}$, on a $\rR^i\psi_*(F)=0$ pour tout $i\geq 1$.
\end{prop}
Cet énoncé est une conséquence du fait que tout point géométrique $\ox$ de $X$ admet un voisinage étale $U$ dans $X$ tel que $U_\oeta$ soit un schéma $K(\pi,1)$ 
dans le sens de \ref{Kpun2}. Cette propriété a été démontrée par Faltings (\cite{faltings1} Lemme 2.3 page 281), 
généralisant des résultats d'Artin (\cite{sga4} XI), sous l'hypothèse en vigueur  dans cette introduction que $X$ est lisse. 
Elle a été récemment étendue au cas semi-stable par Achinger \cite{achinger}. 
Ne disposant pas de cette généralisation dans \cite{faltings2}, pour traiter le cas semi-stable Faltings est contraint
de changer de stratégie par rapport à celle de son premier article \cite{faltings1}. Le détour qu'il a proposé est assez technique et compliqué.
Nous reprenons dans ce travail, y compris dans le cas semi-stable, la stratégie initiale. 

La proposition \ref{intro5} implique que pour tout entier $i\geq 0$, on a un isomorphisme canonique 
\begin{equation}\label{intro5a}
\rH^i(X_{\oeta,\et},F)\stackrel{\sim}{\rightarrow}\rH^i(\tE,\psi_*(F)).
\end{equation}
Nous sommes donc ramenés à un calcul de cohomologie dans $\tE$. Pour ce faire, Faltings introduit un anneau de $\tE$ qui joue le rôle d'anneau structural.

\subsection{}\label{intro6} 
Pour tout objet $(V\rightarrow U)$ de $E$, on désigne par $\oU^V$ la fermeture intégrale de $\oU=U\times_S\oS$ dans $V$ et on pose
\begin{equation}
\ocB(V\rightarrow U)=\Gamma(\oU^V,\co_{\oU^V}).
\end{equation} 
Le préfaisceau $\ocB$ sur $E$ ainsi défini est en fait un faisceau. Pour tout entier $n\geq 1$, on pose $\ocB_n=\ocB/p^n\ocB$.

\begin{teo}[cf. \ref{TPCF14}]\label{intro7} 
Supposons $X$ propre, et soient $n$ un entier $\geq 1$, $F$ un faisceau localement constant constructible de $(\mZ/p^n\mZ)$-modules de $X_{\oeta,\et}$.
Alors, pour tout entier $i\geq 0$, le noyau et le conoyau du morphisme canonique 
\begin{equation}\label{intro7a} 
\rH^i(X_{\oeta,\et},F)\otimes_{\mZ_p}\co_C\rightarrow \rH^i(\tE,\psi_*(F)\otimes_{\mZ_p}\ocB)
\end{equation}
sont annulés par $\fm_C$. 
\end{teo}
On dit que le morphisme \eqref{intro7a} est un {\em $\alpha$-isomorphisme} (ou un {\em presque-isomorphisme} dans la terminologie classique) (cf. \ref{alpha2}). 

C'est le principal théorème de comparaison de Faltings d'où dérivent tous les autres théorèmes de comparaison entre la cohomologie étale $p$-adique
et d'autres cohomologies $p$-adiques.  
Nous le démontrerons suivant la stratégie de Faltings dans \cite{faltings2} basée sur la suite exacte d'Artin-Schreier pour le ``perfectisé'' de $\ocB_1$. 
L'un des principaux ingrédients de la preuve \eqref{mptf11}, dû à Faltings et dégagé par Scholze, 
est un énoncé de structure pour les $\varphi$-modules sur le ``perfectisé'' de $\co_\oK/p\co_\oK$ 
vérifiant certains conditions, dont une condition de presque-finitude dans le sens de Faltings. 
Dans notre application à la cohomologie du topos de Faltings annelé par le ``perfectisé'' de $\ocB_1$, 
la preuve de cette dernière condition résulte de la conjugaison de trois ingrédients: 
(i) des calculs locaux de cohomologie galoisienne utilisant le théorème de presque-pureté de Faltings \eqref{cg37}; (ii) une étude 
fine des conditions de presque-finitude pour les faisceaux de modules sur les schémas \eqref{afini6}; (iii) le résultat de Kiehl sur la finitude de la cohomologie 
d'un morphisme propre \eqref{afini8}. La même stratégie de Faltings a été reprise par Scholze pour généraliser l'énoncé aux variétés rigides \cite{scholze1}.

\subsection{}\label{intro8} 
Compte tenu de la proposition \ref{intro5} et du théorème \ref{intro7}, pour établir \eqref{intro2a}, 
nous sommes amenés à calculer les groupes de cohomologie $\rH^*(\tE,\ocB_n)$ pour tout $n\geq 1$. 
La suite spectrale \eqref{intro2b} s'obtient alors à partir de la suite spectrale de Cartan-Leray pour le morphisme $\sigma$ \eqref{intro4d}. Pour ce faire, nous avons
besoin de calculer les faisceaux de cohomologie $\rR^i\sigma_*(\ocB_n)$ pour tout $i\geq 0$. Une suite exacte de Kummer permet de relier ces derniers 
aux faisceaux des formes différentielles $\Omega^i_{\oX_n/\oS_n}$ \eqref{intro3a}. 

Notons $\hbar\colon \oX\rightarrow X$ le morphisme canonique \eqref{intro3a}. On a un homomorphisme canonique 
\begin{equation}\label{intro8a} 
\hbar_*(\co_\oX)\rightarrow \sigma_*(\ocB)
\end{equation}
qui correspond pour chaque objet $U$ de $\Et_{/X}$ au morphisme canonique $\oU^{U_\oeta}\rightarrow \oU$. Pour tout entier $n\geq 1$, 
on désigne par $\cQ_n$ le monoïde de $\tE$ défini par le diagramme cartésien d'homomorphismes de monoïdes
\begin{equation}\label{intro8b} 
\xymatrix{
\cQ_n\ar[r]\ar[d]&{\sigma^*(\hbar_*(\co_\oX^\times))}\ar[d]\\
{\ocB}\ar[r]^{\nu^n}&\ocB}
\end{equation}
où $\nu$ désigne l'endomorphisme de monoïdes de $\ocB$ d'élévation à la puissance $p$-ième et la flèche verticale de droite est induite par \eqref{intro8a}.
Notons $\mu_{p^n,\tE}$ le groupe constant de $\tE$ de valeur le groupe $\mu_{p^n}(\co_\oK)$ des racines $p^n$-ièmes de l'unité dans $\co_\oK$.
Il existe une suite canonique  
\begin{equation}\label{intro8c} 
0\rightarrow \mu_{p^n,\tE}\rightarrow \cQ_n \rightarrow \sigma^*(\hbar_*(\co_\oX^\times))\rightarrow 0,
\end{equation}
qui n'est pas exacte en général. Nous montrons toutefois qu'elle est exacte au-dessus de la fibre spéciale $\tE_s$ de $\tE$, 
c'est-à-dire au-dessus du sous-topos fermé de $\tE$ complémentaire de l'ouvert $\sigma^*(X_\eta)$ (cf. \ref{TFA66}).  
Celui-ci s'insère dans un diagramme commutatif à isomorphisme canonique près
\begin{equation}\label{intro8d} 
\xymatrix{
{\tE_s}\ar[r]^\delta\ar[d]_{\sigma_s}&\tE\ar[d]^\sigma\\
{X_{s,\et}}\ar[r]^a&{X_\et}}
\end{equation}
où $a$ est l'injection canonique, $\delta$ est le plongement canonique et $\sigma_s$ est induit par $\sigma$.
Nous n'utiliserons du topos de Faltings que sa fibre spéciale. On observera que $\ocB_n$ est un objet de $\tE_s$. 

\subsection{}\label{intro9} 
Supposons dans la suite de cette introduction que $k$ soit algébriquement clos et notons $\oa\colon X_s\rightarrow \oX$ l'injection canonique. 
Soit $n$ un entier $\geq 1$. On considère $\co_{\oX_n}$ et $\Omega^1_{\oX_n/\oS_n}$ comme des faisceaux de $X_{s,\et}$. 
Le morphisme $\sigma_s$ et l'homomorphisme \eqref{intro8a} induisent un morphisme de topos annelés
\begin{equation}\label{intro9a} 
\sigma_n\colon (\tE_s,\ocB_n)\rightarrow (X_{s,\et},\co_{\oX_n}).
\end{equation}
La suite \eqref{intro8c} induit une suite exacte de monoïdes de $\tE_s$
\begin{equation}\label{intro9b} 
0\rightarrow \mu_{p^n,\tE_s}\rightarrow \delta^*(\cQ_n)\rightarrow \sigma^*_s(\oa^*(\co_\oX^\times))\rightarrow 0.
\end{equation}
En particulier, $\delta^*(\cQ_n)$ est un faisceau abélien. 
Le morphisme bord de la suite exacte longue de cohomologie obtenue en appliquant le foncteur $\sigma_{n*}$ 
induit un morphisme $\co_{\oX_n}$-linéaire
\begin{equation}\label{intro9c}
(\oa^*\co_\oX^\times)\otimes_\mZ\co_{\oX_n}\rightarrow \rR^1\sigma_{n*}(\ocB_n(1)).
\end{equation}
Par ailleurs, le morphisme $d\log$ induit un morphisme $\co_{\oX_n}$-linéaire surjectif 
\begin{equation}\label{intro9d}
(\oa^*\co_\oX^\times)\otimes_\mZ\co_{\oX_n}\twoheadrightarrow \Omega^1_{\oX_n/\oS_n}.
\end{equation}
Après une légère modification de la torsion à la Tate dans $\ocB_n(1)$, nous montrons que le morphisme \eqref{intro9c} se factorise
à travers \eqref{intro9d}, établissant un lien entre formes différentielles et cohomologie étale.  
La nouvelle torsion est définie grâce aux {\em épaississements infinitésimaux $p$-adiques de Fontaine}, 
qui servent aussi à définir l'outil principal pour établir la factorisation, à savoir le {\em torseur des déformations}.

\subsection{}\label{intro10}
Rappelons que l'anneau 
\begin{equation}\label{intro10a}
(\co_\oK)^\flat= \underset{\underset{x\mapsto x^p}{\longleftarrow}}{\lim}\co_\oK/p\co_\oK
\end{equation}
est intègre et parfait de caractéristique $p$. L'application 
\begin{equation}\label{intro10b}
\{(x^{(n)})_{n\geq 0}\in (\co_C)^{\mN} \ |\ (x^{(n+1)})^p=x^{(n)},\ \forall n\in \mN\}\rightarrow (\co_\oK)^\flat, \ \ \  
(x^{(n)})_{n\geq 0}\mapsto (\ox^{(n)})_{n\geq 0},
\end{equation}
où $\ox^{(n)}$ désigne la classe de $x^{(n)}$ dans $\co_\oK/p\co_\oK$, est un isomorphisme de monoïdes multiplicatifs. 
Pour tout $x\in (\co_\oK)^\flat$, on note $(x^{(n)})_{n\geq 0}$ la suite de $\co_C$ qui lui est associée par l'isomorphisme inverse de 
\eqref{intro10b}. 
On désigne par $\rW((\co_\oK)^\flat)$ l'anneau des vecteurs de Witt à coefficients dans $(\co_\oK)^\flat$ relatif à $p$ et par
\begin{equation}\label{intro10c}
\theta\colon \rW((\co_\oK)^\flat)\rightarrow \co_C, \ \ \ (x_0,x_1,\dots)\mapsto \sum_{n\geq 0} p^nx_n^{(n)}
\end{equation}
l'homomorphisme de Fontaine. 
On fixe une suite $(p_n)_{n\geq 0}$ d'éléments de $\co_\oK$ telle que $p_0=p$ et $p_{n+1}^p=p_n$ pour tout $n\geq 0$ 
et on note $\varpi$ l'élément associé de $(\co_\oK)^\flat$. On pose 
\begin{equation}\label{intro10d}
\xi=[\varpi]-p\in \rW((\co_\oK)^\flat),
\end{equation}
où $[\ ]$ est le représentant multiplicatif. La suite 
\begin{equation}\label{intro10e}
0\longrightarrow \rW((\co_\oK)^\flat) \stackrel{\cdot \xi}{\longrightarrow} \rW((\co_\oK)^\flat)\stackrel{\theta}{\longrightarrow} \co_C\longrightarrow 0
\end{equation}
est exacte. Considérons l'anneau 
\begin{equation}\label{intro10f}
\cA_2(\co_\oK)=\rW((\co_\oK)^\flat)/\ker(\theta)^2
\end{equation}
qui s'insère dans une suite exacte
\begin{equation}\label{intro10g}
0\longrightarrow \co_C\stackrel{\cdot \xi}{\longrightarrow} \cA_2(\co_\oK)
\stackrel{\theta}{\longrightarrow} \co_C \longrightarrow 0,
\end{equation}
où on a encore noté $\theta$ l'homomorphisme induit par $\theta$. L'idéal $\ker(\theta)$ de $\cA_2(\co_\oK)$ est de carré nul. 
C'est un $\co_C$-module libre de base $\xi$ qui sera noté $\xi\co_C$. 
On observera que contrairement à $\xi$, ce module ne dépend pas du choix de la suite $(p_n)_{n\geq 0}$. 
On note $\xi^{-1}\co_C$ le $\co_C$-module dual de $\xi\co_C$. Pour tout $\co_C$-module $M$, on désigne les $\co_C$-modules $M\otimes_{\co_C}(\xi \co_C)$ 
et $M\otimes_{\co_C}(\xi^{-1} \co_C)$ simplement par $\xi M$ et $\xi^{-1} M$, respectivement. 

On a un homomorphisme canonique $\mZ_p(1)\rightarrow ((\co_\oK)^\flat)^\times$. 
Pour tout $\zeta\in \mZ_p(1)$, on note encore $\zeta$ son image dans $((\co_\oK)^\flat)^\times$.
Comme $\theta([\zeta]-1)=0$, on obtient un homomorphisme de groupes
\begin{equation}\label{intro10h}
\mZ_p(1)\rightarrow \cA_2(\co_\oK),\ \ \ 
\zeta\mapsto\log([\zeta])=[\zeta]-1,
\end{equation}
dont l'image est contenue dans $\ker(\theta)=\xi\co_C$. Celui-ci est clairement $\mZ_p$-linéaire. 
Son image engendre l'idéal $p^{\frac{1}{p-1}}\xi\co_C$ de $\cA_2(\co_\oK)$, et le morphisme $\co_C$-linéaire induit
\begin{equation}\label{intro10i}
\co_C(1)\rightarrow p^{\frac{1}{p-1}}\xi \co_C
\end{equation}
est un isomorphisme. 

On déduit de \eqref{intro9c} et \eqref{intro10i} un morphisme $\co_{\oX_n}$-linéaire canonique
\begin{equation}\label{intro10j}
(\oa^*\co_\oX^\times)\otimes_\mZ\co_{\oX_n}\rightarrow \rR^1\sigma_{n*}(\xi\ocB_n).
\end{equation}

\begin{teo}[cf. \ref{TFSLA24}]\label{intro11}
Soit $n$ un entier $\geq 1$. 
\begin{itemize}
\item[{\rm (i)}] Le morphisme \eqref{intro10j} se factorise à travers \eqref{intro9d} et il induit un morphisme 
\begin{equation}\label{intro11a}
\xi^{-1}\Omega^1_{\oX_n/\oS_n}\rightarrow \rR^1\sigma_{n*}(\ocB_n).
\end{equation}
\item[{\rm (ii)}] Il existe un unique homomorphisme de $\co_{\oX_n}$-algèbres graduées 
\begin{equation}\label{intro11b}
\wedge(\xi^{-1}\Omega^1_{\oX_n/\oS_n})\rightarrow \oplus_{i\geq 0}\rR^i\sigma_{n*}(\ocB_n)
\end{equation}
dont la composante de degré un est le morphisme \eqref{intro11a}. Son noyau est annulé par $p^{\frac{2d}{p-1}}\fm$ et son conoyau 
est annulé par $p^{\frac{2d+1}{p-1}}\fm$, où $d=\dim(X/S)$.
\end{itemize}
\end{teo}

La suite spectrale de Hodge-Tate \eqref{intro2b} s'obtient à partir de la suite spectrale de Cartan-Leray pour les morphismes $\sigma_n$ et de l'énoncé ci-dessus. 

Il ressort de notre preuve du théorème \ref{intro11} que le morphisme \eqref{intro11b} est un analogue de l'isomorphisme de Cartier en caractéristique $p$, 
ce qui suggère une analogie entre la suite spectrale de Hodge-Tate et la suite spectrale conjuguée en caractéristique $p$. 

\subsection{}\label{intro12}
Pour établir le théorème \ref{intro11}, on peut supposer $X=\Spec(R)$ affine et {\em petit} dans le sens de Faltings, autrement dit qu'il existe un $S$-morphisme 
étale de $X$ dans le tore $\mG_{m,S}^d$ de dimension $d$ au-dessus de $S$,
et que $X_s$ est non-vide. Pour alléger les notations, on suppose, de plus, $X_\oeta$ connexe. Soient $\ox$
un point géométrique de $X_\oeta$, $(V_i)_{i\in I}$ un revêtement universel de $X_\oeta$ en $\ox$, $\Delta=\pi_1(X_\oeta,\ox)$ le groupe fondamental
de $X_\oeta$ en $\ox$. Pour tout $i\in I$, notons $X_i=\Spec(R_i)$ la fermeture intégrale de $\oX$ dans $V_i$ et posons 
\begin{equation}\label{intro12a}
\oR=\underset{\underset{i\in I}{\longrightarrow}}{\lim}\ R_i,
\end{equation}
qui est une représentation discrète de $\Delta$. On note $\hoR$ le séparé complété $p$-adique de $\oR$ et on pose
\begin{equation}\label{intro12b}
\oR^\flat=\underset{\underset{x\mapsto x^p}{\longleftarrow}}{\lim}\oR/p\oR.
\end{equation}
On a une suite exacte 
\begin{equation}\label{intro12c}
0\longrightarrow \rW(\oR^\flat)\stackrel{\cdot \xi}{\longrightarrow} \rW(\oR^\flat)
\stackrel{\theta}{\longrightarrow} \hoR \longrightarrow 0,
\end{equation}
où $\theta$ est l'homomorphisme canonique analogue à \eqref{intro10c}  (cf. \ref{eipo3})
et $\xi$ est l'élément \eqref{intro10d}. On pose 
\begin{equation}\label{intro12d}
\cA_2(\oR)=\rW(\oR^\flat)/\ker(\theta)^2,
\end{equation} 
qui est donc une extension de $\hoR$ par un idéal de carré nul isomorphe à $\xi\hoR$. 

Comme $X$ est affine, il existe une déformation lisse $\tX$ de $X\otimes_{\co_K}\co_C$ au-dessus de $\cA_2(\co_\oK)$
\begin{equation}\label{intro12e}
\xymatrix{
{X\otimes_{\co_K}\co_C}\ar[r]\ar[d]\ar@{}[rd]|\Box&\tX\ar[d]\\
{\Spec(\co_C)}\ar[r]&{\Spec(\cA_2(\co_\oK))}}
\end{equation}
Fixons une telle déformation dans la suite de cette introduction et considérons le diagramme commutatif privé de la flèche pointillée 
\begin{equation}\label{intro12f}
\xymatrix{
{\Spec(\hoR)}\ar[r]\ar[d]&{\Spec(\cA_2(\oR))}\ar@/^2pc/[dd]\ar@{.>}[d]\\
{X\otimes_{\co_K}\co_C}\ar[r]\ar[d]\ar@{}[rd]|\Box&\tX\ar[d]\\
{\Spec(\co_C)}\ar[r]&{\Spec(\cA_2(\co_\oK))}}
\end{equation}
Les morphismes représentés par des flèches pointillées qui complètent ce diagramme de façon à le laisser commutatif forment un torseur $\cL$
sur $\Spec(\hoR)$ pour la topologie de Zariski sous le $\hoR$-module 
\begin{equation}\label{intro12g}
\Hom_{\hoR}(\Omega^1_{R/\co_K}\otimes_R\hoR,\xi\hoR).
\end{equation} 
On désigne par $\cF$ le $\hoR$-module des fonctions affines sur $\cL$ (cf. \cite{agt} III.4.9). 
Celui-ci s'insère dans une suite exacte canonique
\begin{equation}\label{intro12h}
0\rightarrow \hoR\rightarrow \cF\rightarrow \xi^{-1}\Omega^1_{R/\co_K} \otimes_R \hoR\rightarrow 0.
\end{equation} 
Cette suite induit pour tout entier $n\geq 1$ une suite exacte 
\begin{equation}\label{intro12i}
0\rightarrow \Sym^{n-1}_{\hoR}(\cF)\rightarrow \Sym^{n}_{\hoR}(\cF)\rightarrow \Sym^n_{\hoR}(\xi^{-1}\Omega^1_{R/\co_K}
\otimes_R\hoR)\rightarrow 0.
\end{equation}
Les $\hoR$-modules $(\Sym^{n}_{\hoR}(\cF))_{n\in \mN}$ forment donc un système inductif filtrant 
dont la limite inductive 
\begin{equation}\label{intro12j}
\cC=\underset{\underset{n\geq 0}{\longrightarrow}}\lim\ \Sym^n_{\hoR}(\cF)
\end{equation}
est naturellement munie d'une structure de $\hoR$-algèbre. Le $\hoR$-schéma $\Spec(\cC)$ représente alors canoniquement le torseur $\cL$ (\cite{agt} III.4.10). 

L'action naturelle de $\Delta$ sur le schéma $\cA_2(\oR)$ induit  
une action $\hoR$-semi-linéaire de $\Delta$ sur $\cF$, telle que les morphismes de la suite \eqref{intro12h} soient 
$\Delta$-équivariants. 

\begin{prop}\label{intro13}
Sous les hypothèses de \ref{intro12}, on a un diagramme commutatif
\begin{equation}\label{intro13a}
\xymatrix{
R^\times\ar[r]\ar[d]_{d\log}&{\rH^1(\Delta,\mu_{p^n}(\co_\oK))}\ar[d]^{-d\log([\ ])}\\
{\Omega^1_{R/\co_K}\otimes_{\co_K}(\co_\oK/p^n\co_\oK)}\ar[r]&{\rH^1(\Delta,\xi\oR/p^n\xi\oR)}}
\end{equation}
où la flèche verticale de droite est induite par l'homomorphisme \eqref{intro10h},
la flèche horizontale inférieure est induite par le morphisme bord de la suite exacte \eqref{intro12h}, 
et la flèche horizontale supérieure associe à $t\in R^\times$
la classe du $\mu_{p^n}(\co_\oK)$-torseur $X_\oeta[T]/(T^{p^n}-t)$ au-dessus de $X_\oeta$. 
\end{prop}

On en déduit que l'extension $\cF$ \eqref{intro12h}, dont la classe d'isomorphisme ne dépend pas de la déformation $\tX$, 
fournit la factorisation \eqref{intro11a} recherchée. 

La preuve de la proposition \ref{intro11}(ii) est classique dans la théorie de Faltings. 
L'assertion étant locale, elle résulte d'un calcul de cohomologie galoisienne utilisant le théorème de presque-pureté de Faltings.

\subsection{}\label{intro14}
Pour expliquer l'analogie entre le morphisme \eqref{intro11b} et l'isomorphisme de Cartier en caractéristique $p$, 
nous avons besoin de l'interprétation de ce dernier due à Ogus et Vologodsky \cite{ov} généralisant des résultats de Deligne et Illusie \cite{deligne-illusie}.

Soient $T$ le spectre d'un corps parfait $k$  de caractéristique $p$, $Y$ un $T$-schéma propre et lisse.  
Considérons les deux suites spectrales d'hypercohomologie
\begin{eqnarray}
E^{i,j}_1=\rH^j(Y,\Omega^i_{Y/T})&\Rightarrow& \rH^{i+j}_{\dR}(Y/T)=\rH^{i+j}(Y,\Omega^\bullet_{Y/T}),\label{intro14a}\\
E^{i,j}_2=\rH^i(Y,\cH^j(\Omega^\bullet_{Y/T}))&\Rightarrow& \rH^{i+j}_{\dR}(Y/T).\label{intro14b}
\end{eqnarray} 
La première suite spectrale, dite {\em de Hodge vers de Rham}, ne dégénère pas nécessairement en $E_1$; mais quand elle dégénère en $E_1$, 
la seconde suite spectrale, dite {\em conjuguée}, dégénère en $E_2$. 
C'est une conséquence de l'isomorphisme de Cartier. En effet, considérons le diagramme commutatif
\begin{equation}\label{intro14c}
\xymatrix{
Y\ar[r]^\rF\ar[rd]\ar@/^2pc/[rr]^{\rF_Y}&Y'\ar[r]^{\sigma'}\ar[d]\ar@{}[rd]|\Box&Y\ar[d]\\
&T\ar[r]^{\sigma}&T}
\end{equation}
où $\sigma$ (resp. $\rF_Y$) est l'endomorphisme de Frobenius de $T$ (resp. $Y$), $Y'$ est le changement de base de $Y$ par $\sigma$ et 
$\rF$ est le morphisme de Frobenius relatif de $Y$ sur $T$. Pour tout entier $i\geq 0$, il existe un isomorphisme $\co_{Y'}$-linéaire canonique
\begin{equation}\label{intro14d}
\gamma^i\colon \Omega^i_{Y'/T}\stackrel{\sim}{\rightarrow} \cH^i(\rF_*(\Omega^\bullet_{Y/T})),
\end{equation}
dit isomorphisme de Cartier (inverse). Pour toute section locale $f$ de $\co_Y$, on a 
\begin{eqnarray}
\gamma^0(\sigma'^*f)&=&f^p,\label{intro14e}\\
\gamma^1(d(\sigma'^*f))&=&f^{p-1}df.\label{intro14f}
\end{eqnarray}
On en déduit un isomorphisme $k$-linéaire canonique 
\begin{equation}\label{intro14g}
\sigma^*(E_1^{i,j})\stackrel{\sim}{\rightarrow}E_2^{j,i}.
\end{equation}

\subsection{}\label{intro15}
Posons $\tT=\Spec(\rW_2(k))$ et supposons dans ce numéro qu'il existe un relèvement $\tY$ de $Y$ au-dessus de $\tT$. On note $\tY'$ le changement 
de base de $\tY$ par l'endomorphisme de Frobenius $\sigma_{\tT}$ de $\tT$. Si $\trF\colon \tY\rightarrow \tY'$ est un relèvement de $\rF\colon Y\rightarrow Y'$, 
le morphisme $\trF^*\colon \Omega^1_{\tY'/\tT}\rightarrow \Omega^1_{\tY/\tT}$  est divisible par $p$ et on obtient un diagramme commutatif
\begin{equation}\label{intro15a}
\xymatrix{
{\Omega^1_{\tY'/\tT}}\ar[r]^-(0.5){p^{-1}\trF^*}\ar@{->>}[d]&{\rF_*(\cZ^1_{Y/T})}\ar[d]\\
{\Omega^1_{Y'/T}}\ar[ru]^{\zeta_{\trF}}\ar[r]_-(0.5){\gamma^1}&{\rF_*(\cH^1(\Omega^\bullet_{Y/T}))}}
\end{equation}
où $\cZ^1_{Y/T}$ désigne le faisceau des cocycles du complexe de de Rham $\Omega^\bullet_{Y/T}$ et $\zeta_\trF$ est le morphisme $\co_{Y'}$-linéaire 
induit par $p^{-1}\trF^*$. 

En général, il existe une extension canonique de modules à connexion intégrable sur $Y$
\begin{equation}\label{intro15b}
0\rightarrow (\co_Y,d)\rightarrow (\cE_\tY,\nabla)\rightarrow (\rF^*\Omega^1_{Y'/T},d)\rightarrow 0,
\end{equation}
où $\rF^*\Omega^1_{Y'/T}$ est muni de la connexion définie par descente par Frobenius, vérifiant les deux propriétés suivantes~:
\begin{itemize}
\item[(i)] Le morphisme bord de la suite exacte longue de cohomologie de de Rham
\begin{equation}\label{intro15c}
\rH^0(Y',\Omega^1_{Y'/T})=\rH^0_\dR(Y',\rF^*\Omega^1_{Y'/T})\rightarrow \rH^1_\dR(Y/T)
\end{equation}
est le morphisme défini par Deligne et Illusie à un signe près.
\item[(ii)] Le morphisme bord de la suite exacte longue de faisceaux de cohomologie de de Rham
\begin{equation}\label{intro15d}
\Omega^1_{Y'/T}=\cH^0_\dR(\rF^*\Omega^1_{Y'/T})\rightarrow \cH^1_\dR(\co_Y)
\end{equation}
est l'opposé de l'isomorphisme de Cartier $\gamma^1$. 
\end{itemize}

\subsection{}\label{intro16}
L'extension \eqref{intro15b} est construite comme suit. Posons $\cT_{Y'/T}=\cHom_{\co_{Y'}}(\Omega^1_{Y'/T},\co_{Y'})$.
Les relèvements du morphisme de Frobenius relatif $\rF$ forment un torseur $\cL_\tY$ sur $Y$
pour la topologie de Zariski sous le $\co_Y$-module $\rF^*(\cT_{Y'/T})$.
Notons $\cE_\tY$ le $\co_Y$-module des fonctions affines sur $\cL_\tY$. Celui-ci s'insère dans une suite exacte canonique
\begin{equation}\label{intro16a}
0\rightarrow \co_Y\rightarrow \cE_\tY\rightarrow \rF^*\Omega^1_{Y'/T}\rightarrow 0.
\end{equation}
Cette suite induit pour tout entier $n\geq 1$ une suite exacte 
\begin{equation}\label{intro16b}
0\rightarrow \Sym^{n-1}_{\co_Y}(\cE_\tY)\rightarrow \Sym^{n}_{\co_Y}(\cE_\tY)\rightarrow \Sym^n_{\co_Y}(\rF^*\Omega^1_{Y'/T})\rightarrow 0.
\end{equation}
Les $\co_Y$-modules $(\Sym^{n}_{\hoR}(\cE_\tY))_{n\in \mN}$ forment donc un système inductif filtrant 
dont la limite inductive 
\begin{equation}\label{intro16c}
\cC_\tY=\underset{\underset{n\geq 0}{\longrightarrow}}\lim\ \Sym^n_{\co_Y}(\cE_\tY)
\end{equation}
est naturellement munie d'une structure de $\co_Y$-algèbre. Le $Y$-schéma $\Spec(\cC_\tY)$ représente alors canoniquement le torseur $\cL_\tY$. 
Le torseur $\cL_\tY$ est naturellement muni d'une connexion 
\begin{equation}\label{intro16d}
\nabla\colon \cL_\tY\rightarrow \rF^*(\cT_{Y'/T})\otimes_{\co_Y}\Omega^1_{Y/T}
\end{equation}
compatible avec la connexion sur $\rF^*(\cT_{Y'/T})$ définie par descente de Frobenius. Si $\trF$ est une section locale de $\cL_\tY$, alors 
\begin{equation}\label{intro16e}
\nabla(\tF)\colon \rF^*\Omega^1_{Y'/T}\rightarrow \Omega^1_{Y/T}
\end{equation}
est le morphisme induit par $\zeta_\trF$ \eqref{intro15a}. 

La connexion $\nabla$ sur $\cL_\tY$ induit une connexion sur $\cE_\tY$ dont la $p$-courbure 
\begin{equation}\label{intro16f}
\psi\colon \cE_\tY\rightarrow \cE_\tY\otimes\rF^*\Omega^1_{Y'/T}
\end{equation}
est induite par le morphisme identique de $\rF^*\Omega^1_{Y'/T}$ \eqref{intro16a}.

\subsection{}\label{intro17}
L'article est organisé en deux chapitres, le premier consacré à  des préliminaires techniques, le second à la théorie de 
Hodge $p$-adique proprement dite. Après un rappel  de diverses conventions et 
notations dans la section \ref{notconv}, nous définissons dans la section \ref{Kpun} la notion  de schéma 
$K(\pi,1)$  \eqref{Kpun2}  et en établissons quelques propriétés. La section \ref{EIPF} porte sur les   épaississements $p$-adiques 
infinitésimaux de Fontaine afin de préciser les structures logarithmiques dont nous les munissons \eqref{eip5}. 
Les cinq sections suivantes sont consacrés à établir divers résultats, pour lesquels des références font défaut,  
de \emph{presque-algèbre} omniprésente  dans l'approche de Faltings de la théorie de Hodge $p$-adique.  
Pour alléger, nous avons choisi dans cet article de substituer au préfixe \emph{presque} le préfixe $\alpha$ suggéré 
par le \emph{almost} anglais.  On reprend dans la section \ref{alpha} le formalisme  
des catégories abéliennes localisées dans un cadre légèrement plus général que ce qu'on trouve habituellement dans la 
littérature afin de développer la notion de faisceau sur un site à valeurs dans une catégorie de $\alpha$-modules \eqref{alpha18}.   
La section \ref{finita} aborde l'examen des questions de $\alpha$-finitude pour les modules d'un  topos annelé, 
\emph{type $\alpha$-fini}, \emph{présentation $\alpha$-finie} et 
\emph{$\alpha$-cohérence} \eqref{finita9}, et de leurs  propriétés usuelles de stabilité. 
La section \ref{afini} est consacrée à l'étude de ces propriétés pour les modules sur des schémas. 
Nous y introduisons la notion de module \emph{$\alpha$-quasi-cohérent} \eqref{afini2} et nous établissons pour ceux-ci un critère 
simple pour la présentation $\alpha$-finie, similaire à celui pour le topos ponctuel ({\em i.e.}, pour les modules usuels) \eqref{afini6}. 
Combiné avec un théorème fondamental de Kiehl, nous en déduisons un énoncé de $\alpha$-finitude pour les morphismes propres \eqref{afini8}.
Après la section \ref{aet} constituée de brefs rappels sur les extensions {\em $\alpha$-étales}, on rappelle dans la section \ref{mptf} la 
structure des modules de type $\alpha$-fini sur un anneau de valuation non-discrète de hauteur $1$ \eqref{mptf7}. 
On en déduit en caractéristique $p$ une description de certains $\varphi$-modules à $\alpha$-isomorphisme près \eqref{mptf11}. 
Les trois dernières sections de ce chapitre sont consacrées à diverses questions dans les  topos intervenant 
au cours de cet article. Dans les sections \ref{rec} et \ref{lptce}, nous établissons quelques résultats sur les 
topos co-évanescents de Deligne qui seront utiles pour les développements ultérieurs de ce travail. 
Pour clore ces préliminaires, la section \ref{tf} revient sur le topos de Faltings.

\subsection{}\label{intro18}
Le chapitre \ref{thpadique} débute en précisant, dans sa section \ref{TFA},  
notre cadre géométrique, celui des  morphismes {\em adéquats} $f\colon (X,\cM_X)\rightarrow (S,\cM_S)$ de 
schémas logarithmiques (\cite{agt} III.4.7), et les notations relatives au topos de Faltings associé. 
La donnée d'une carte convenable  pour un tel morphisme 
$f$ nous permet dans la section \ref{cad} de définir et d'étudier les propriétés géométriques de la tour de revêtements qui lui est associée.  
Nous utilisons cette dernière en particulier pour améliorer des énoncés locaux sur la propriété $K(\pi,1)$ \eqref{cad10}. 
Lorsque, de plus, $X$ est affine, nous rappelons dans la section \ref{cg} l'énoncé du théorème de {\em presque-pureté} de Faltings \eqref{cg22},  
et nous en déduisons des énoncés de $\alpha$-finitude de certains groupes de cohomologie galoisienne \eqref{cg37}. 
Nous en donnons aussi quelques conséquences cohomologiques relatives au topos de Faltings (\ref{cg33} et \ref{cg27}).  La section \ref{TPCF}
exploite ces résultats locaux et établit  lorsque $X$ est propre sur $S$  la $\alpha$-finitude de certains groupes de cohomologie 
du topos de Faltings (\ref{TPCF4}). Nous démontrons ensuite le principal théorème de comparaison de Faltings \eqref{TPCF14} 
en utilisant la suite exacte d'Artin-Schreier \eqref{TPCF11} pour le ``perfectisé'' de $\ocB_1$. 
La courte section \ref{acycloc} traite des conséquences cohomologiques pour le topos de Faltings des propriétés $K(\pi,1)$. 
La section \ref{deform} rappelle les définitions du {\em torseur des déformations} et de {\em l'algèbre de Higgs-Tate} relatifs à une déformation logarithmique de $X$ 
au-dessus de l'épaississement infinitésimal universel $p$-adique de Fontaine (\cite{agt} II.10).  
Enfin, dans la section \ref{TFSLA}, nous développons la théorie de Kummer du topos annelé de Faltings et nous en déduisons 
grâce à l'extension de Higgs-Tate le lien recherché entre les faisceaux de cohomologie $\rR^i\sigma_*(\ocB_n)$ et 
les faisceaux des formes différentielles logarithmiques $\tOmega^i_{\oX_n/\oS_n}$ \eqref{TFSLA24}.

\subsection{}\label{intro19}
Cette rédaction laisse temporairement en suspens divers pans de la théorie : cohomologie à supports compacts, 
dualité de Poincaré, généralisation de la suite spectrale de Hodge-Tate à une situation relative, ainsi que les variations de structures de Hodge $p$-adiques
dont la théorie a été très sommairement esquissée par Hyodo dans le cas affine \cite{hyodo}. Nous reviendrons à ces questions dans une diffusion ultérieure de ce travail.   

\subsection*{Remerciements}\label{intro20} 
Les auteurs tiennent  à exprimer toute leur reconnaissance à G. Faltings pour l'inspiration constante suscitée par ses travaux. 
Ils  remercient également chaleureusement T. Tsuji pour les différents échanges qu'ils ont pu avoir sur ces questions.  
Ce travail a bénéficié du soutien du programme ANR Théorie de Hodge $p$-adique et Développements (ThéHopadD) ANR-11-BS01-005.

\renewcommand{\thesection}{\thechapter.\arabic{section}}
\renewcommand{\thesubsection}{\thesection.\arabic{subsection}}

\numberwithin{equation}{subsection}

\chapter{Préliminaires}\label{prelim}

\section{Notations et conventions}\label{notconv}

{\em Tous les anneaux considérés dans cet article possèdent un élément unité~;
les homomorphismes d'anneaux sont toujours supposés transformer l'élément unité en l'élément unité.
Nous considérons surtout des anneaux commutatifs, et lorsque nous parlons d'anneau 
sans préciser, il est sous-entendu qu'il s'agit d'un anneau commutatif~; en particulier, 
il est sous-entendu, lorsque nous parlons d'un topos annelé $(X,A)$ sans préciser, que $A$ est commutatif.

On sous-entend par {\em monoïde}  un monoïde commutatif et unitaire. 
Les homomorphismes de monoïdes sont toujours supposés transformer l'élément unité en l'élément unité.}

\subsection{}\label{notconv1}
Soient $A$ un anneau, $p$ un nombre premier, $n$ un entier $\geq1$.  
On désigne par $\rW(A)$ (resp. $\rW_n(A)$) l'anneau des vecteurs de Witt 
(resp. vecteurs de Witt de longueur $n$) à coefficients dans $A$ relatif à $p$. 
On a un homomorphisme d'anneaux
\begin{equation}\label{notconv1a}
\Phi_n\colon 
\begin{array}[t]{clcr}
\rW_n(A)&\rightarrow& A,\\ 
(x_1,\dots,x_n)&\mapsto&x_1^{p^{n-1}}+p x_2^{p^{n-2}}+\dots+p^{n-1}x_n
\end{array}
\end{equation}
appelé $n$-ième composante fantôme. 
On dispose aussi des morphismes de restriction, de décalage et de Frobenius
\begin{eqnarray}
\rR\colon \rW_{n+1}(A)&\rightarrow& \rW_n(A),\label{notconv1b}\\
\rV\colon \rW_n(A)&\rightarrow& \rW_{n+1}(A),\label{notconv1c}\\
\rF\colon \rW_{n+1}(A)&\rightarrow& \rW_n(A).\label{notconv1d}
\end{eqnarray}
Lorsque $A$ est de caractéristique $p$, $\rF$ induit un endomorphisme de $\rW_n(A)$, encore noté $\rF$.

\subsection{} \label{notconv2}
Pour tout anneau $R$ et tout monoïde $M$, on désigne par $R[M]$ la $R$-algèbre de $M$ et par 
$e\colon M\rightarrow R[M]$ l'homomorphisme canonique, où $R[M]$ est considéré comme un
monoïde multiplicatif. Pour tout $x\in M$, on notera $e^x$ au lieu de $e(x)$. 

On désigne par $\bA_M$ le schéma $\Spec(\mZ[M])$ muni de la structure logarithmique associée à 
la structure pré-logarithmique définie par $e\colon M\rightarrow \mZ[M]$ (\cite{agt} II.5.9). 
Pour tout homomorphisme de monoïdes $\vartheta\colon M\rightarrow N$, 
on note $\bA_\vartheta\colon \bA_N\rightarrow \bA_M$ le morphisme de schémas logarithmiques associé.

\subsection{}\label{notconv3}
Dans tout cet article, on fixe un univers $\mU$ possédant un élément de cardinal infini, 
et un univers $\mV$ tel que $\mU\in \mV$. 
On appelle catégorie des $\mU$-ensembles et l'on note $\Ens$, 
la catégorie des ensembles qui se trouvent dans $\mU$. 
C'est un $\mU$-topos ponctuel que l'on note encore $\Pt$ (\cite{sga4} IV 2.2). 
On désigne par $\Sch$ la catégorie des schémas éléments de $\mU$. 
Sauf mention explicite du contraire, il sera sous-entendu que les anneaux et 
les schémas logarithmiques (et en particulier les schémas)
envisagés dans cet article sont éléments de l'univers~$\mU$.

\subsection{}\label{notconv4}
Pour une catégorie $\cC$, nous notons $\ob(\cC)$ l'ensemble de ses objets,
$\cC^\circ$ la catégorie opposée, et pour $X,Y\in \ob(\cC)$, 
$\Hom_\cC(X,Y)$ (ou $\Hom(X,Y)$ lorsqu'il n'y a aucune ambiguïté) 
l'ensemble des morphismes de $X$ dans $Y$. 

Si $\cC$ et $\cC'$ sont deux catégories, nous désignons par $\Hom(\cC,\cC')$ 
l'ensemble des foncteurs de $\cC$ dans $\cC'$, et  
par $\bHom(\cC,\cC')$ la catégorie des foncteurs de $\cC$ dans $\cC'$. 

Soient $\cE$ une catégorie, $\cC$ et $\cC'$ deux catégories sur $\cE$ (\cite{sga1} VI 2). 
Nous notons $\Hom_{\cE}(\cC,\cC')$ l'ensemble des $\cE$-foncteurs de $\cC$ dans $\cC'$ 
et $\Hom_{\cart/\cE}(\cC,\cC')$ l'ensemble des foncteurs cartésiens (\cite{sga1} VI 5.2).
Nous désignons par $\bHom_{\cE}(\cC,\cC')$ la catégorie des $\cE$-foncteurs de $\cC$ dans $\cC'$ et 
par $\bHom_{\cart/\cE}(\cC,\cC')$ la sous-catégorie pleine formée des foncteurs cartésiens.

\subsection{}\label{notconv5}
Soit $\cC$ une catégorie. On désigne par $\hcC$ la catégorie des préfaisceaux 
de $\mU$-ensembles sur $\cC$, c'est-à-dire la catégorie des foncteurs
contravariants sur $\cC$ à valeurs dans $\Ens$ (\cite{sga4} I 1.2). 
Si $\cC$ est munie d'une topologie (\cite{sga4} II 1.1), on désigne par $\tcC$ le topos 
des faisceaux de $\mU$-ensembles sur $\cC$ (\cite{sga4} II 2.1). 

Pour $F$ un objet de $\hcC$, on note $\cC_{/F}$ la catégorie 
suivante (\cite{sga4} I 3.4.0). Les objets de $\cC_{/F}$
sont les couples formés d'un objet $X$ de $\cC$ 
et d'un morphisme $u$ de $X$ dans $F$. Si $(X,u)$ et $(Y,v)$ sont deux objets, 
un morphisme de $(X,u)$ vers $(Y,v)$ est un morphisme $g\colon X\rightarrow Y$
tel que $u=v\circ g$.

\subsection{}\label{notconv6}
Soient $\cC$ un $\mU$-site (\cite{sga4} II 3.0.2), $\tcC$ le topos des faisceaux de $\mU$-ensembles sur $\cC$,
$X$ un préfaisceau de $\mU$-ensembles sur $\cC$. On munit la catégorie $\cC_{/X}$ de la topologie induite par la topologie de $\cC$ au moyen 
du foncteur ``source'' $j_X\colon \cC_{/X}\rightarrow \cC$ (\cite{sga4} III 3.1).
D'après (\cite{sga4} III 5.2), $j_X$ est un foncteur continu et cocontinu. 
Il définit donc une suite de trois foncteurs adjoints~: 
\begin{equation}\label{notconv6a}
j_{X!}\colon (\cC_{/X})^\sim \rightarrow \tcC,\ \ \ \
j_X^*\colon \tcC \rightarrow (\cC_{/X})^\sim,\ \ \ \
j_{X*}\colon  (\cC_{/X})^\sim\rightarrow \tcC,
\end{equation}
dans le sens que pour deux foncteurs consécutifs de la suite, celui
de droite est adjoint à droite de l'autre. 
Le foncteur $j_{X!}$ se factorise par la catégorie $\tcC_{/X^\ra}$, 
où $X^\ra$ est le faisceau associé à $X$, 
et le foncteur induit $(\cC_{/X})^\sim \rightarrow \tcC_{/X^\ra}$
est une équivalence de catégories (\cite{sga4} III 5.4). 
Le couple de foncteurs $(j_X^*,j_{X*})$ définit un morphisme de topos qu'on notera aussi abusivement
$j_{X}\colon \tcC_{/X^\ra} \rightarrow \tcC$, dit morphisme de localisation de $\tcC$ en $X^\ra$ (cf. \cite{sga4} IV 5.2).
Pour tout objet $F$ de $\tcC$, on pose 
\begin{equation}\label{notconv6b}
F|X^\tta=j_X^*(F).
\end{equation}

Supposons, de plus, que les limites projectives finies soient représentables dans $\cC$ et que $X$ soit un objet de $\cC$.
Soit $e_\cC$ un objet final de $\cC$ (qui existe par hypothèse). 
On note $j_X^+\colon \cC\rightarrow \cC_{/X}$ le foncteur de changement de base par le morphisme canonique
$X\rightarrow e_\cC$. Alors, $j_X^+$ est exact à gauche et continu (\cite{sga4} III 1.6 et 3.3).
Comme $j_X^*$ prolonge $j_X^+$ (\cite{sga4} III 5.4), le morphisme de localisation
$j_X\colon \tcC_{/X^\ra} \rightarrow \tcC$ s'identifie canoniquement au morphisme de topos associé à $j_X^+$.

\subsection{}\label{notconv7} 
Soient $\cC$ un $\mU$-site, $\hcC$ la catégorie des préfaisceaux de $\mU$-ensembles sur $\cC$,
$\tcC$ le topos des faisceaux de $\mU$-ensembles sur $\cC$.
On désigne par $\Mon_\hcC$ (resp. $\Mon_\tcC$) la catégorie des monoïdes (commutatifs et unitaires) 
de $\hcC$ (resp. $\tcC$) (cf. \cite{sga4} I 3.2, resp. II 6.3.1). 
Les $\mU$-limites inductives et projectives dans $\Mon_\hcC$ sont représentables et se calculent terme à terme.
Les $\mU$-limites inductives et projectives dans $\Mon_\tcC$ sont représentables.  
Le foncteur d'injection canonique 
\begin{equation}
i\colon \Mon_\tcC\rightarrow \Mon_\hcC
\end{equation} 
admet un adjoint à gauche
\begin{equation}
\ra\colon \Mon_\hcC\rightarrow \Mon_\tcC, \ \ \ M\mapsto M^\ra,
\end{equation}
qui est exact à gauche (\cite{sga4} II 6.4). 
D'après (\cite{sga4} II 6.5), le foncteur ``faisceau d'ensembles sous-jacent'' $\Mon_\tcC\rightarrow \tcC$ admet un adjoint à gauche 
et commute donc aux limites projectives. De plus, pour tout préfaisceau de monoïdes $M$ sur $\cC$, 
le faisceau d'ensembles sous-jacent à $M^\ra$ est canoniquement isomorphe au faisceau d'ensembles associé
au préfaisceau d'ensembles sous-jacent à $M$. 

On note $1_\tcC$ le monoïde unité de $\tcC$, c'est-à-dire le faisceau de monoïdes associé au préfaisceau de monoïdes 
sur $\cC$ de valeur $1$. C'est un objet initial et final de la catégorie $\Mon_\tcC$.
On appelle {\em conoyau} d'un morphisme $u\colon N\rightarrow M$ de $\Mon_\tcC$ 
la somme amalgamée de $u$ et du morphisme canonique $M\rightarrow 1_\tcC$. 
On a un isomorphisme canonique 
\begin{equation}
\coker(u)\stackrel{\sim}{\rightarrow} \ra(\coker(i(u))).
\end{equation}
Si $u$ est un monomorphisme, on appelle $\coker(u)$ le {\em quotient} de $M$ par $N$ et on le note $M/N$.

\subsection{}\label{notconv8} 
Pour tout morphisme de topos $f\colon Y\rightarrow X$, les foncteurs $f^*$ et $f_*$ induisent un couple de foncteurs adjoints
$f^*\colon \Mon_X\rightarrow \Mon_Y$ et $f_*\colon \Mon_Y\rightarrow \Mon_X$. Le foncteur $f^*$ est exact à gauche
d'après \ref{notconv7}.

\subsection{}\label{notconv9}
Soit $(X,A)$ un $\mU$-topos annelé (\cite{sga4} 11.1.1). On note $\bMod(A)$ ou $\bMod(A,X)$ 
la catégorie des $A$-modules de $X$. Si $M$ est un $A$-module, on désigne par $\rS_A(M)$ 
(resp. $\wedge_A(M)$, resp. $\Gamma_A(M)$) l'algèbre symétrique (resp. extérieure, resp. à puissances divisées) 
de $M$ et pour tout entier $n\geq 0$, par $\rS_A^n(M)$ (resp. $\wedge_A^n(M)$,
resp. $\Gamma_A^n(M)$) sa partie homogène de degré $n$. 
Les formations de ces algèbres commutent à la localisation au-dessus d'un objet de $X$.

\begin{defi}[\cite{sga6} I 1.3.1]\label{notconv14}
Soit $(X,A)$ un topos annelé. 
On dit qu'un $A$-module $M$ de $X$ est {\em localement projectif de type fini} si les conditions 
équivalentes suivantes sont satisfaites :
\begin{itemize}
\item[{\rm (i)}] $M$ est de type fini et le foncteur $\cHom_A(M,\cdot)$ est exact~;
\item[{\rm (ii)}] $M$ est de type fini et tout épimorphisme de $A$-modules $N\rightarrow M$ admet localement une section~;
\item[{\rm (iii)}] $M$ est localement facteur direct d'un $A$-module libre de type fini. 
\end{itemize}
\end{defi}

Lorsque $X$ a suffisamment de points et que pour tout point $x$ de $X$, la fibre de $A$ en $x$
est un anneau local, les $A$-modules localement projectifs de type fini sont les $A$-modules 
localement libres de type fini (\cite{sga6} I 2.15.1).

\subsection{}\label{notconv10}
Pour tout schéma $X$, on désigne par $\Et_{/X}$ le {\em site étale} de $X$, 
c'est-à-dire, la sous-catégorie pleine de $\Sch_{/X}$ \eqref{notconv3} formée des schémas étales sur $X$,
munie de la topologie étale; c'est un $\mU$-site. 
On note $X_\et$ le {\em topos étale} de $X$, c'est-à-dire le topos des faisceaux de $\mU$-ensembles sur $\Et_{/X}$.

On désigne par $\Et_{\coh/X}$ la sous-catégorie pleine de $\Et_{/X}$ formée des schémas étales 
de présentation finie sur $X$, munie de la topologie induite par celle de $\Et_{/X}$; c'est un site $\mU$-petit. 
Si $X$ est quasi-séparé, le foncteur de restriction de $X_\et$ dans le topos des faisceaux de $\mU$-ensembles sur 
$\Et_{\coh/X}$ est une équivalence de catégories (\cite{sga4} VII 3.1 et 3.2). 

On désigne par $\Et_{\rf/X}$ la sous-catégorie pleine de $\Et_{/X}$ formée des schémas étales finis sur $X$, 
munie de la topologie induite par celle de $\Et_{/X}$; c'est un site $\mU$-petit. 
On appelle {\em topos fini étale} de $X$ et on note  $X_\fet$, le topos des faisceaux de $\mU$-ensembles sur $\Et_{\rf/X}$ 
(cf. \cite{agt} VI.9.2). L'injection canonique $\Et_{\rf/X}\rightarrow \Et_{/X}$ induit un morphisme de topos 
\begin{equation}\label{notconv10a}
\rho_X\colon X_\et\rightarrow X_\fet.
\end{equation}

\subsection{}\label{notconv12}
Soit $X$ un schéma. On désigne par $X_\zar$ le topos de Zariski de $X$ et par 
\begin{equation}\label{notconv12b}
u_X\colon X_{\et}\rightarrow X_{\zar}
\end{equation}
le morphisme canonique (\cite{sga4} VII 4.2.2). Si $F$ est un $\co_X$-module quasi-cohérent de $X_\zar$, 
on désigne par $\iota(F)$ le faisceau de $X_\et$ défini pour tout $X$-schéma étale $U$ par (\cite{sga4} VII 2 c))
\begin{equation}\label{notconv12c}
\iota(F)(U)=\Gamma(U,F\otimes_{\co_X}\co_U).
\end{equation}
Il est commode, lorsqu'il n'y aucune risque de confusion, de désigner $\iota(F)$ abusivement par $F$. 
On notera que $\iota(\co_X)$ est un anneau de $X_\et$ et que $\iota(F)$ est un $\iota(\co_X)$-module. 

Notant $\bMod^\qcoh(\co_X,X_\zar)$ la sous-catégorie pleine de $\bMod(\co_X,X_\zar)$ formée des $\co_X$-modules quasi-cohérents \eqref{notconv9},
la correspondance $F\mapsto \iota(F)$ définit un foncteur
\begin{equation}\label{notconv12a}
\iota\colon \bMod^\qcoh(\co_X,X_\zar)\rightarrow \bMod(\co_X,X_\et).
\end{equation}

Pour tout $\co_X$-module quasi-cohérent $F$ de $X_\zar$, on a un isomorphisme canonique 
\begin{equation}\label{notconv12d}
F\stackrel{\sim}{\rightarrow}u_{X*}(\iota(F)).
\end{equation}
Nous considérons donc $u_X$ comme un morphisme de topos annelés 
\begin{equation}\label{notconv12j}
u_X\colon (X_{\et},\co_X)\rightarrow (X_{\zar},\co_X). 
\end{equation}
Nous utilisons pour les modules la notation $u^{-1}_X$ pour désigner l'image
inverse au sens des faisceaux abéliens et nous réservons la notation 
$u^*_X$ pour l'image inverse au sens des modules. L'isomorphisme \eqref{notconv12d} induit par adjonction un morphisme
\begin{equation}\label{notconv12e}
u^*_X(F)\rightarrow \iota(F).
\end{equation}
Celui-ci est un isomorphisme. En effet, pour tout point géométrique $\ox$ de $X$, on a un isomorphisme canonique
\begin{equation}\label{afini12f}
\iota(F)_\ox\stackrel{\sim}{\rightarrow}\underset{\underset{\ox\rightarrow U}{\longleftarrow}}{\lim}\ \Gamma(U,F\otimes_{\co_X}\co_U),
\end{equation}
où la limite est prise sur les voisinages de $\ox$ dans $\Et_{/X}$ (on peut évidemment se limiter au voisinages affines au-dessus d'un ouvert affine de $X$ 
contenant l'image $x$ de $\ox$). Notant $X'$ le localisé strict de $X$ en $\ox$, on en déduit, d'après (\cite{ega4} 8.5.2(i)),
un isomorphisme canonique 
\begin{equation}\label{afini12g}
\iota(F)_\ox\stackrel{\sim}{\rightarrow} \Gamma(X',F\otimes_{\co_X}\co_{X'}).
\end{equation}
Par ailleurs, on a un isomorphisme canonique
\begin{equation}\label{afini12h}
u_X^*(F)_\ox\stackrel{\sim}{\rightarrow} F_x\otimes_{\co_{X,x}} \iota(\co_X)_{\ox}.
\end{equation}
La fibre du morphisme \eqref{notconv12e} en $\ox$ s'identifie au morphisme canonique
\begin{equation}\label{afini12i}
F_x\otimes_{\co_{X,x}} \Gamma(X',\co_{X'}) \rightarrow  \Gamma(X',F\otimes_{\co_X}\co_{X'}),
\end{equation}
qui est un isomorphisme. Par suite, \eqref{notconv12e} est un isomorphisme.  On en déduit que $\iota$ est un adjoint à gauche du foncteur 
$u_{X*}$. La flèche d'adjonction $\id\rightarrow u_{X*}\circ \iota$ est un isomorphisme d'après \eqref{notconv12d}. En particulier, $\iota$ est pleinement fidèle. 

Soient $f\colon X'\rightarrow X$ un morphisme, $F$ un $\co_X$-module quasi-cohérent. 
Pour tout $\co_{X'}$-module quasi-cohérent $G$, on note $\iota'(G)$ le $\co_{X'}$-module associé de $X'_\et$.
Il existe un morphisme canonique fonctoriel
\begin{equation}\label{afini12j}
\iota(F)\rightarrow f_*(\iota'(f^*(F))).
\end{equation}
Le morphisme $f$ induit donc un morphisme de topos annelés que l'on note encore  
\begin{equation}\label{afini12k}
f\colon (X'_\et,\co_{X'})\rightarrow (X_\et,\co_X).
\end{equation}
Nous utilisons pour les modules la notation $f^{-1}$ pour désigner l'image
inverse au sens des faisceaux abéliens et nous réservons la notation $f^*$ pour l'image inverse au sens des modules.

Compte tenu de \eqref{afini12g}, le morphisme \eqref{afini12j} induit par adjonction un isomorphisme
\begin{equation}\label{afini12l}
f^*(\iota(F))\stackrel{\sim}{\rightarrow} \iota'(f^*(F)).
\end{equation}

\subsection{}\label{notconv11}
Soient $X$ un schéma connexe, $\ox$ un point géométrique de $X$.
On désigne par 
\begin{equation}\label{notconv11a}
\omega_\ox\colon \Et_{\rf/X}\rightarrow \Ens
\end{equation}
le foncteur fibre en $\ox$, qui à tout revêtement étale $Y$ de $X$ associe l'ensemble des points géométriques de 
$Y$ au-dessus de $\ox$, par  $\pi_1(X,\ox)$ le groupe fondamental de $X$ en $\ox$ (c'est-à-dire
le groupe des automorphismes du foncteur $\omega_\ox$) et par $\bB_{\pi_1(X,\ox)}$ 
le topos classifiant du groupe profini $\pi_1(X,\ox)$, 
c'est-à-dire la catégorie des $\mU$-ensembles discrets munis d'une action continue à gauche de $\pi_1(X,\ox)$ 
(\cite{sga4} IV 2.7). Alors $\omega_\ox$ induit un foncteur pleinement fidèle 
\begin{equation}\label{notconv11b}
\mu_\ox^+\colon \Et_{\rf/X}\rightarrow \bB_{\pi_1(X,\ox)}
\end{equation}
d'image essentielle la sous-catégorie pleine de $\bB_{\pi_1(X,\ox)}$ formée des ensembles finis
(\cite{sga1} V  §~4 et §~7).
Soit $(X_i)_{i\in I}$ un système projectif sur un ensemble ordonné filtrant $I$ dans $\Et_{\rf/X}$
qui pro-représente $\omega_\ox$, normalisé par le fait que les morphismes de transition $X_i\rightarrow X_j$
$(i\geq j)$ sont des épimorphismes et que tout épimorphisme $X_i\rightarrow X'$ de $\Et_{\rf/X}$ 
est équivalent à un épimorphisme $X_i\rightarrow X_j$ $(j\leq i)$ convenable. 
Un tel pro-objet est essentiellement unique. Il est appelé le {\em revêtement universel normalisé de $X$ en $\ox$} 
ou le {\em pro-objet fondamental normalisé de $\Et_{\rf/X}$ en $\ox$}. 
On notera que l'ensemble $I$ est $\mU$-petit. Le foncteur 
\begin{equation}\label{notconv11c}
\nu_\ox\colon X_\fet\rightarrow \bB_{\pi_1(X,\ox)},
\ \ \ F\mapsto \underset{\underset{i\in I}{\longrightarrow}}{\lim}\ F(X_i)
\end{equation}
est une équivalence de catégories qui prolonge le foncteur $\mu_\ox^+$ (cf. \cite{agt} VI.9.8). 
On l'appelle le {\em foncteur fibre} de $X_\fet$ en $\ox$.

\section{\texorpdfstring{Schémas $K(\pi,1)$}{Schémas K(pi,1)}}\label{Kpun}

\begin{prop}\label{Kpun1}
Soient $X$ un schéma cohérent, n'ayant qu'un nombre fini de composantes connexes, 
$\mP$ un ensemble de nombres premiers, $\rho_X\colon X_\et\rightarrow X_\fet$ le morphisme canonique \eqref{notconv10a}.
Alors, les conditions suivantes sont équivalentes~:
\begin{itemize}
\item[{\rm (i)}] pour tout faisceau abélien de $\mP$-torsion $F$ de $X_\fet$ {\rm (\cite{sga4} IX 1.1)}, 
le morphisme d'adjonction $F\rightarrow \rR\rho_{X*}(\rho_X^*F)$ est un isomorphisme~;
\item[{\rm (ii)}] pour tout faisceau abélien de $\mP$-torsion, localement constant et constructible $F$ de $X_\et$ et tout entier $i\geq 1$, $\rR^i\rho_{X*}(F)=0$~;
\item[{\rm (iii)}] pour tout entier $n$ dont les diviseurs premiers appartiennent à $\mP$, 
tout $(\mZ/n\mZ)$-module localement constant et constructible $F$ de $X_\et$ et tout entier $i\geq 1$, $\rR^i\rho_{X*}(F)=0$~;
\item[{\rm (iv)}] pour tout faisceau abélien de $\mP$-torsion, localement constant et constructible $F$ de $X_\et$, tout entier $i\geq 1$
et tout $\xi\in \rH^i(X,F)$, il existe un revêtement étale surjectif $X'\rightarrow X$ tel que l'image canonique de $\xi$ dans $\rH^i(X',F)$ soit nulle~;
\item[{\rm (v)}] pour tout revêtement étale $Y\rightarrow X$, tout faisceau abélien de $\mP$-torsion, localement constant et constructible $F$ de $Y_\et$, 
tout entier $i\geq 1$ et tout $\xi\in \rH^i(Y,F)$, il existe un revêtement étale surjectif $Y'\rightarrow Y$ tel que l'image canonique de $\xi$ dans 
$\rH^i(Y',F)$ soit nulle~;
\item[{\rm (vi)}] pour tout entier $n$ dont les diviseurs premiers appartiennent à $\mP$, tout $(\mZ/n\mZ)$-module localement constant et constructible $F$ de $X_\et$,
tout entier $i\geq 1$ et tout $\xi\in \rH^i(X,F)$, il existe un revêtement étale surjectif $X'\rightarrow X$ tel que l'image canonique de $\xi$ dans $\rH^i(X',F)$ soit nulle. 
\end{itemize}
\end{prop}

En effet, le morphisme d'adjonction $\id\rightarrow \rho_{X*}\rho_X^*$ est un isomorphisme d'après (\cite{agt} VI.9.18). 
Les conditions (i) et (ii) sont donc équivalentes en vertu de (\cite{agt} VI.9.12, VI.9.14 et VI.9.20) et (\cite{sga4} VI 5.1).  
L'équivalence des conditions (ii) et (iii) est aussi une conséquence de (\cite{agt} VI.9.12 et VI.9.14) et (\cite{sga4} VI 5.1).

Pour tout faisceau abélien $G$ de $X_\et$ et tout entier $i\geq 0$, le faisceau $\rR^i\rho_{X*}(G)$ est le faisceau de $X_\fet$
associé au préfaisceau qui à tout $Y\in \ob(\Et_{\rf/X})$ associe le groupe $\rH^i(Y,G)$ (\cite{sga4} V 5.1). On a donc (iii)$\Rightarrow$(iv) et (v)$\Rightarrow$(iii). 

Montrons (iv)$\Rightarrow$(v). Supposons la condition (iv) satisfaite. Soient $f\colon Y\rightarrow X$ un revêtement étale, $F$ un faisceau abélien de $\mP$-torsion 
localement constant constructible de $Y_\et$, $i$ un entier $\geq 1$, $\xi\in \rH^i(Y,F)$. D'après (\cite{sga4} V 5.3, VIII 5.5 et XVI 2.2),
$f_*(F)$ est un faisceau abélien de $\mP$-torsion localement constant constructible de $X_\et$,  et le morphisme canonique
\begin{equation}
\rH^i(X,f_*(F))\rightarrow \rH^i(Y,F)
\end{equation}
est un isomorphisme.  Notons $\zeta$ l'image inverse de $\xi$. D'après (iv), il existe un revêtement étale surjectif $g\colon X'\rightarrow X$ tel que l'image canonique de 
$\zeta$ dans $\rH^i(X',f_*(F))$ soit nulle. Posons $Y'=Y\times_XX'$ et notons $f'\colon Y'\rightarrow X'$ la projection canonique. 
On a un diagramme commutatif
\begin{equation}
\xymatrix{
{\rH^i(X,f_*(F))}\ar[rr]^-(0.5)u\ar[d]&&{\rH^i(Y,F)}\ar[d]\\
{\rH^i(X',f_*(F))}\ar[r]^-(0.5)v&{\rH^i(X',f'_*(F|Y'))}\ar[r]^-(0.5){u'}&{\rH^i(Y',F)}}
\end{equation}
où les flèches verticales sont les morphismes canoniques, 
les morphismes $u$ et $u'$ sont induits par la suite spectrale de Cartan-Leray et le morphisme $v$ est induit par  
le morphisme de changement de base $f_*(F)|X'\rightarrow f'_*(F|Y')$.
On en déduit que l'image canonique de $\xi$ dans $\rH^i(Y',F)$ est nulle; d'où la condition (v). 

Enfin, les conditions (iv) et (vi) sont équivalentes en vertu de (\cite{sga4} VI 5.2).

\begin{defi}\label{Kpun2}
Soient $X$ un schéma cohérent, n'ayant qu'un nombre fini de composantes connexes, $\mP$ un ensemble de nombres premiers. 
On dit que $X$ est un {\em schéma $K(\pi,1)$ pour les faisceaux abéliens de $\mP$-torsion} (\cite{sga4} IX 1.1) si les conditions équivalentes
de \ref{Kpun1} sont remplies. 
Si $\mP$ est l'ensemble des nombres premiers inversibles dans $\co_X$, on dit simplement que $X$ est un {\em schéma $K(\pi,1)$}.
\end{defi}

\begin{rema}\label{Kpun3}
Soient $X$ un schéma, $F$ un faisceau abélien de torsion, localement constant et constructible de $X_\et$. 
Par descente (\cite{sga1} VIII 2.1, \cite{ega4} 2.7.1 et 17.7.3), $F$ est représentable par un schéma en groupes étale et fini $G$ au-dessus de $X$.  
Le groupe $\rH^1(X,F)$ classifie alors les $G$-fibrés principaux homogènes au-dessus de $X$ (\cite{sga4} VII 2 a)). 
Par descente, tout $G$-fibré principal homogène au-dessus de $X$ est un revêtement étale de $X$. 
Par suite, pour tout $\xi\in \rH^1(X,F)$, il existe un revêtement étale surjectif $X'\rightarrow X$ tel que l'image canonique de $\xi$ dans $\rH^1(X',F)$ soit nulle.
Notant $\rho_X\colon X_\et\rightarrow X_\fet$ le morphisme canonique \eqref{notconv10a}, on en déduit que 
$\rR^1\rho_{X*}(F)=0$ (\cite{sga4} V 5.1). Il en résulte, par la suite spectrale de Cartan-Leray, que le morphisme canonique 
\begin{equation}
\rH^1(X_\fet,\rho_{X*}(F)) \rightarrow \rH^1(X_\et,F)
\end{equation}
est un isomorphisme. 

Supposons, de plus, que $X$ soit cohérent et qu'il n'ait qu'un nombre fini de composantes connexes. Pour tout faisceau abélien de torsion $G$ de $X_\fet$, 
le morphisme d'adjonction $G\rightarrow \rho_{X*}(\rho_X^*G)$ est un isomorphisme d'après (\cite{agt} VI.9.18). On en déduit que le morphisme canonique 
\begin{equation}
\rH^1(X_\fet,G) \rightarrow \rH^1(X_\et,\rho_X^*(G))
\end{equation}
est un isomorphisme. 
\end{rema}

\begin{lem}\label{Kpun4}
Soient $X$ un schéma cohérent et étale localement connexe {\rm (\cite{agt} VI.9.7)}, $f\colon Y\rightarrow X$ un revêtement étale, 
$\mP$ un ensemble de nombres premiers. Alors, 
\begin{itemize}
\item[{\rm (i)}] Si $X$ est un schéma $K(\pi,1)$ pour les faisceaux abéliens de $\mP$-torsion,
il en est de même de $Y$.
\item[{\rm (ii)}] Supposons $f$ surjectif. Pour que $X$ soit un schéma $K(\pi,1)$ pour les faisceaux abéliens de $\mP$-torsion,
il faut et il suffit qu'il en soit de même de $Y$.
\end{itemize}
\end{lem}

Cela résulte aussitôt de \ref{Kpun1}.

\begin{prop}\label{Kpun10} 
Soient $I$ une catégorie filtrante et essentiellement $\mU$-petite {\rm (\cite{sga4} I 8.1.8)}, 
$\varphi\colon I^\circ \rightarrow \Sch$ \eqref{notconv3} un foncteur qui transforme les objets de $I$ en des schémas cohérents n'ayant qu'un nombre
fini de composantes connexes, et les morphismes de $I$ en des morphismes affines, $X$ la limite projective de $\varphi$ {\rm (\cite{ega4} 8.2.3)},
$\mP$ un ensemble de nombres premiers. Supposons que $X$ soit cohérent et qu'il n'ait qu'un nombre fini de composantes connexes. Alors,
\begin{itemize}
\item[{\rm (i)}] Si pour tout $i\in \ob(I)$, $X_i$ est un schéma $K(\pi,1)$ pour les faisceaux abéliens de $\mP$-torsion, il en est de même de $X$. 
\item[{\rm (ii)}] Supposons que $\varphi$ transforme les morphismes de $I$ en des revêtements étales surjectifs. 
Alors, pour que $X$ soit un schéma $K(\pi,1)$ pour les faisceaux abéliens de $\mP$-torsion, il faut et il suffit qu'il en soit de même de $X_i$ pour tout $i\in \ob(I)$,
et il suffit qu'il en soit de même de $X_i$ pour un seul objet $i$ de $I$. 
\end{itemize}
\end{prop}

(i) En effet, soit $F$ un faisceau abélien de $\mP$-torsion, localement constant et constructible de $X_\et$. Par descente, $F$ est représentable par un
schéma en groupes abélien fini et étale au-dessus de $X$. 
D'après (\cite{ega4} 8.8.2, 8.10.5 et 17.7.8), il existe $\iota\in \ob(I)$ et un faisceau abélien de $\mP$-torsion localement constant constructible $F_{\iota}$
de $X_{\iota,\et}$ tels que $F$ soit l'image inverse de $F_\iota$. Quitte à remplacer $I^\circ$ par $I^\circ_{/\iota}$, 
on peut supposer que $\iota$ est un objet initial de $I$. Pour tout $i\in \ob(I)$, on note $F_i$ l'image inverse canonique de $F_\iota$ sur $X_i$. 
En vertu de (\cite{sga4} VII 5.8), pour tout entier $q\geq 0$, on a un isomorphisme canonique 
\begin{equation}\label{Kpun10a} 
\rH^q(X,F)\stackrel{\sim}{\rightarrow}\underset{\underset{i\in I}{\longrightarrow}}{\lim}\ \rH^q(X_i,F_i).
\end{equation}
La proposition s'ensuit aussitôt. 

(ii) Supposons que $X$ soit un schéma $K(\pi,1)$ pour les faisceaux abéliens de $\mP$-torsion. Soit $\iota\in \ob(I)$, $F_\iota$ un 
faisceau abélien de $\mP$-torsion, localement constant et constructible de $X_{\iota,\et}$. Quitte à remplacer $I^\circ$ par $I^\circ_{/\iota}$, on peut supposer 
que $\iota$ est un objet initial de $I$. Pour tout $i\in \ob(I)$, on note $F_i$ (resp. $F$) 
l'image inverse canonique de $F_\iota$ sur $X_i$ (resp. $X$). Soient $n\geq 1$, $\xi\in \rH^n(X_\iota,F_\iota)$. Par hypothèses, 
il existe un revêtement étale surjectif $X'\rightarrow X$ tel que l'image canonique de $\xi$ dans $\rH^n(X',F)$ soit nulle. 
D'après (\cite{ega4} 8.8.2, 8.10.5, 17.7.8 et 8.3.11), il existe un $i\in \ob(I)$, un revêtement étale surjectif $X'_i\rightarrow X_i$
et un $X$-isomorphisme $X'_i\times_{X_i}X \stackrel{\sim}{\rightarrow} X'$. Par un isomorphisme analogue à \eqref{Kpun10a}, on voit qu'il existe 
un morphisme $i\rightarrow j$ de $I$ tel que l'image canonique de $\xi$ dans $\rH^n(X'_i\times_{X_i}X_j,F_j)$ soit nulle. Comme le morphisme $X_j\rightarrow X_\iota$
est un revêtement étale surjectif, on en déduit que $X_\iota$ est $K(\pi,1)$ pour les faisceaux abéliens de $\mP$-torsion. 
La proposition s'ensuit compte tenu de (i) et \ref{Kpun4}(ii).

\begin{cor}\label{Kpun11}
Soient $K$ un corps, $L$ une extension algébrique de $K$, $X$ un $K$-schéma de type-fini,
$\mP$ un ensemble de nombres premiers. Alors, pour que $X$ soit $K(\pi,1)$ pour les faisceaux abéliens de $\mP$-torsion,  
il faut et il suffit qu'il en soit de même de $X\otimes_KL$. 
\end{cor}

Quitte à remplacer $L$ par la clôture séparable de $K$ dans $L$ (\cite{sga4} VIII 1.1), on peut supposer $L$ séparable sur $K$. 
La proposition résulte alors de \ref{Kpun4}(ii) et \ref{Kpun10}(ii). 

\begin{prop}\label{Kpun18}
Soient $X$ un schéma, $X^\circ$ un ouvert de $X$, $\ox$ un point géométrique de $X$, $X'$ le localisé strict de $X$ en $\ox$, 
$\mP$ un ensemble de nombres premiers.  
Pour tout $X$-schéma $Y$, posons $Y^\circ=Y\times_XX^\circ$. Supposons que 
le schéma $X'^\circ$ soit cohérent et qu'il n'ait qu'un nombre fini de composantes connexes. 
Alors, les conditions suivantes sont équivalentes~:
\begin{itemize}
\item[{\rm (i)}] Le schéma $X'^\circ$ est $K(\pi,1)$ pour les faisceaux abéliens de $\mP$-torsion. 
\item[{\rm (ii)}] Pour tout $X$-schéma étale $\ox$-pointé $Y$, tout faisceau abélien de $\mP$-torsion, localement cons\-tant et constructible $\cF$ de $Y^\circ_\et$, 
tout entier $q\geq 1$ et tout $\xi\in \rH^q(Y^\circ,\cF)$, il existe un $Y$-schéma étale $\ox$-pointé $U$ et un revêtement étale surjectif $V\rightarrow U^\circ$
tel que l'image canonique de $\xi$ dans $\rH^q(V,\cF)$ soit nulle.
\end{itemize}
\end{prop}

On peut supposer $X$ affine. On désigne par $\fV_\ox$ la catégorie des $X$-schémas affines étales et $\ox$-pointés.
C'est une catégorie filtrante essentiellement $\mU$-petite (\cite{sga4} IV 6.8.2).
D'après (\cite{ega4} 8.2.3), $X'^\circ$ est la limite projective du foncteur 
\begin{equation}\label{Kpun18a}
\varphi \colon \fV_\ox\rightarrow \Sch, \ \ \  Y\mapsto Y^\circ.
\end{equation} 
Supposons d'abord la condition (ii) satisfaite. 
Soient $F$ un faisceau abélien de $\mP$-torsion localement constant constructible de $X'^\circ_\et$,
$q$ un entier $\geq 1$, $\xi\in \rH^q(X'^\circ,F)$. 
Montrons qu'il existe un revêtement étale surjectif $W\rightarrow X'^\circ$ tel que l'image canonique de $\xi$ dans $\rH^q(W,F)$ soit nulle. 
Par descente, $F$ est représentable par un schéma en groupes abéliens fini et étale au-dessus de $X'^\circ$. 
D'après (\cite{ega4} 8.8.2, 8.10.5 et 17.7.8), quitte à remplacer $X$ par un objet de $\fV_\ox$, on peut supposer qu'il existe un 
faisceau abélien de torsion localement constant constructible $\cF$ de $X^\circ_\et$ tel que $F$ soit isomorphe à l'image inverse de $\cF$. 
En vertu de (\cite{sga4} VII 5.8), on a un isomorphisme canonique 
\begin{equation}\label{Kpun18b}
\rH^q(X'^\circ,F)\stackrel{\sim}{\rightarrow}\underset{\underset{Y\in \fV^\circ_\ox}{\longrightarrow}}{\lim}\ \rH^q(Y^\circ,\cF).
\end{equation}
Il existe donc un objet $Y$ de $\fV_\ox$ tel que $\xi$ soit l'image canonique d'une classe $\zeta\in \rH^q(Y^\circ,\cF)$. 
D'après (ii), il existe un morphisme $U\rightarrow Y$ de $\fV_\ox$ et un revêtement étale surjectif $V\rightarrow U^\circ$ tels que l'image 
canonique de $\zeta$ dans $\rH^q(V,\cF)$ soit nulle. Le revêtement étale surjectif $V\times_{U^\circ}X'^\circ\rightarrow X'^\circ$ 
répond alors à la question~; d'où la condition (i).  

Inversement, supposons la condition (i) satisfaite. Soient $\cF$ un faisceau abélien de $\mP$-torsion, localement constant et constructible de $X^\circ_\et$, 
$q$ un entier $\geq 1$, $\xi\in \rH^q(X^\circ,\cF)$. D'après (i), il existe un revêtement étale surjectif $V\rightarrow X'^\circ$ tel que l'image canonique de $\xi$
dans $\rH^q(V,\cF|X'^\circ)$ soit nulle. D'après (\cite{ega4} 8.8.2, 8.10.5 et 17.7.8), il existe un objet $Y$ de $\fV_\ox$, un revêtement étale surjectif 
$W\rightarrow Y^\circ$ et un $X'^\circ$-isomorphisme $V\stackrel{\sim}{\rightarrow}W\times_{Y^\circ}X'^\circ$. 
En vertu de (\cite{sga4} VII 5.8), on a un isomorphisme canonique 
\begin{equation}\label{Kpun18c}
\rH^q(V,\cF|X'^\circ)\stackrel{\sim}{\rightarrow}\underset{\underset{U\in (\fV^\circ_\ox)_{/Y}}{\longrightarrow}}{\lim}\ \rH^q(W\times_{Y^\circ}U^\circ,\cF).
\end{equation}
Il existe donc un morphisme $U\rightarrow Y$ de $\fV^\circ_\ox$ tel que l'image canonique de $\xi$ dans $\rH^q(W\times_{Y^\circ}U^\circ,\cF)$ soit nulle~;
d'où la condition (ii).

\begin{defi}[\cite{sga4} XI 3.1, \cite{mo} 1.3.1]\label{Kpun5}
On appelle {\em courbe élémentaire} un morphisme de schémas $f\colon X\rightarrow S$ qui s'insère dans un diagramme commutatif
\begin{equation}
\xymatrix{
X\ar[r]^j\ar[rd]_f&\oX\ar[d]_-(0.4)\of&Y\ar[l]_i\ar[ld]^g\\
&S&}
\end{equation}
satisfaisant aux conditions suivantes~:
\begin{itemize}
\item[{\rm (i)}] $\of$ est une courbe relative projective et lisse, à fibres géométriquement irréductibles~; 
\item[{\rm (ii)}] $j$ est une immersion ouverte, $i$ est une immersion fermée, et $X$ est l'ouvert complémentaire de $i(Y)$ dans $\oX$;
\item[{\rm (iii)}] $g$ est un revêtement étale surjectif. 
\end{itemize}
On dit alors aussi que $X$ est une {\em courbe élémentaire au-dessus de $S$}. 
\end{defi}

\begin{defi}[\cite{sga4} XI 3.2, \cite{mo} 1.3.2]\label{Kpun6}
On appelle {\em polycourbe élémentaire} un morphisme de schémas $f\colon X\rightarrow S$ qui admet une factorisation en courbes élémentaires.
On dit alors aussi que $X$ est une {\em polycourbe élémentaire au-dessus de $S$}. 
\end{defi}

\begin{prop}[\cite{sga4} XI 3.3]\label{Kpun7} 
Soient $k$ un corps algébriquement clos, $X$ un $k$-schéma lisse, $x\in X(k)$. Alors, il existe un ouvert de Zariski de $X$ contenant $x$ qui est une 
polycourbe élémentaire au-dessus de $\Spec(k)$. 
\end{prop}

\begin{lem}[d'Abhyankar]\label{Kpun8} 
Soient $f\colon X\rightarrow S$ un morphisme lisse de schémas localement noethériens, 
$D$ un diviseur à croisements normaux sur $X$ relativement à $S$ {\rm (\cite{sga1} XIII 2.1)},
$U$ l'ouvert complémentaire de $D$ dans $X$, $U'$ un revêtement étale surjectif de $U$, modérément ramifié au-dessus de $X$. 
On note $X'$ la clôture intégrale de $X$ dans $U'$ et $D'$ le sous-schéma fermé réduit de $X'$ de même support que $D\times_XX'$. 
Alors, $X'$ est lisse sur $S$, $D'$ est un diviseur à croisements normaux sur $X'$ relativement à $S$,
et le morphisme canonique $D'\rightarrow D$ est un revêtement étale surjectif. 
\end{lem}

Cela résulte de (\cite{sga1} XIII 5.5). 

\begin{prop}\label{Kpun9} 
Soient $S$ un $\mQ$-schéma noethérien, régulier et $K(\pi,1)$, $f\colon X\rightarrow S$ une courbe élémentaire. 
Alors, $X$ est un schéma $K(\pi,1)$.
\end{prop}

On notera d'abord que $X$ est cohérent et n'a qu'un nombre fini de composantes connexes. Considérons un diagramme commutatif 
\begin{equation}\label{Kpun9a} 
\xymatrix{
X\ar[r]^j\ar[rd]_f&\oX\ar[d]_-(0.4)\of&Y\ar[l]_i\ar[ld]^g\\
&S&}
\end{equation}
vérifiant les conditions de \ref{Kpun5}. Soient $F$ un faisceau abélien de torsion, localement constant et constructible sur $X$, 
$n$ un entier $\geq 1$, $\xi\in \rH^n(X,F)$. Montrons qu'il existe un revêtement étale surjectif $X'\rightarrow X$ tel que l'image canonique
de $\xi$ dans $ \rH^n(X',F)$ soit nulle. On peut se borner au cas où $n\geq 2$ \eqref{Kpun3}. 
Il existe un revêtement étale surjectif $Z\rightarrow X$ tel que $F|Z$ soit constant. 
D'après \ref{Kpun7}, le morphisme $Z\rightarrow S$ induit par $f$ est une courbe élémentaire. On peut donc se borner au cas où $F$ est constant,
de valeur un groupe cyclique fini $\Lambda$. 

En vertu du théorème de pureté (\cite{travaux-gabber} XVI 3.1.1), 
on a un isomorphisme canonique $ i^!\Lambda_\oX\stackrel{\sim}{\rightarrow}\Lambda_Y(-1)[-2]$. Le triangle distingué de localisation induit alors 
un triangle de $\bD^+(\oX_\et,\Lambda)$
\begin{equation}\label{Kpun9b} 
\Lambda_\oX\longrightarrow \rR j_*(\Lambda_X)\longrightarrow i_*(\Lambda_Y)(-1)[-1]\stackrel{+1}{\longrightarrow}
\end{equation}
et par suite, $g$ étant fini, un triangle de $\bD^+(S_\et,\Lambda)$
\begin{equation}\label{Kpun9c} 
\rR\of_*(\Lambda_\oX)\longrightarrow \rR f_*(\Lambda_X)\longrightarrow g_*(\Lambda_Y)(-1)[-1]\stackrel{+1}{\longrightarrow}
\end{equation}
Montrons que le morphisme induit 
\begin{equation}\label{Kpun9d} 
g_*(\Lambda_Y)(-1)\rightarrow \rR^2 \of_*(\Lambda_\oX)
\end{equation}
est surjectif. En effet, la formation de ce morphisme commute à tout changement de base $S'\rightarrow S$, d'après (\cite{travaux-gabber} XVI 2.3.2). 
On peut donc se réduire au cas où $S$ est le spectre d'un corps algébriquement clos et où $Y$ est un point fermé de $\oX$. La classe de cycle définie par $Y$
engendre alors $\rH^2(\oX,\Lambda(1))$ (\cite{sga45} Cycle 2.1.5); d'où l'assertion. 

Il résulte de ce qui précède qu'on a un isomorphisme canonique $\of_*(\Lambda_\oX) \stackrel{\sim}{\rightarrow} f_*(\Lambda_X)$ et une suite exacte canonique
\begin{equation}\label{Kpun9e} 
0\rightarrow \rR^1 \of_*(\Lambda_\oX)\rightarrow \rR^1 f_*(\Lambda_X)\rightarrow g_*(\Lambda_Y)(-1) \rightarrow \rR^2\of_*(\Lambda_\oX)\rightarrow 0, 
\end{equation}
et que $\rR^q f_*(\Lambda_X)=0$ pour tout $q\geq 2$. Comme $\of$ et $g$ sont propres et lisses, pour tout $q\geq 0$, 
$\rR^q f_*(\Lambda_X)$ est localement constant constructible (\cite{sga4} XVI 2.2 et IX 2.1).

Considérons la suite spectrale de Cartan-Leray
\begin{equation}
E_2^{a,b}=\rH^a(S,\rR^bf_*(\Lambda_X))\Rightarrow \rH^{a+b}(X,\Lambda_X),
\end{equation}
et notons $(\rE^n_q)_{0\leq q\leq n}$ la filtration aboutissement sur $\rH^n(X,\Lambda_X)$ (on rappelle que n$\geq 2$).
Comme on a $E_\infty^{q,n-q}=E^n_q/E^n_{q+1}$ pour tout $0\leq q\leq n$, on en déduit une suite exacte canonique
\begin{equation}
\rH^n(S,f_*(\Lambda_X))\stackrel{u}{\rightarrow} \rH^n(X,\Lambda_X)\stackrel{v}{\rightarrow} \rH^{n-1}(S,\rR^1 f_*(\Lambda_X))\rightarrow \rH^{n+1}(S,f_*(\Lambda_X)).
\end{equation}
Celle-ci est compatible à tout changement de base $S'\rightarrow S$ dans un sens évident que nous n'explicitons pas. 
Comme le schéma $S$ est $K(\pi,1)$, il existe un revêtement étale surjectif $S'\rightarrow S$ tel que l'image canonique de $v(\xi)$ dans 
$\rH^{n-1}(S',\rR^1 f_*(\Lambda_X))$ soit nulle. On peut donc supposer qu'il existe $\zeta\in \rH^n(S,f_*(\Lambda_X))$ tel que $\xi=v(\zeta)$. 
De même, il existe un revêtement étale surjectif $S''\rightarrow S$ tel que l'image canonique de $\zeta$ dans $\rH^n(S'',f_*(\Lambda_X))$ soit nulle. 
L'image canonique de $\xi$ dans $\rH^n(X\times_SS'',\Lambda_X)$ est donc nulle~; d'où l'assertion recherchée.

\begin{cor}\label{Kpun12}
Soit $k$ un corps algébriquement clos de caractéristique $0$, $X$ un schéma lisse sur $k$, $x\in X$. 
Alors, il existe un ouvert $U$ de $X$ contenant $x$, qui est un schéma $K(\pi,1)$.
\end{cor}

En effet, quitte à remplacer $x$ par une spécialisation, on peut se réduire au cas où $x\in X(k)$. 
La proposition résulte alors de \ref{Kpun7} et \ref{Kpun9}.

\begin{lem}\label{Kpun14}
Soient $S$ un schéma noethérien régulier, $d$ un entier $\geq1$, 
$f\colon X\rightarrow \mA_S^d$ un morphisme lisse d'un schéma $X$ dans l'espace affine de dimension $d$ au-dessus de $S$. 
Notons $\mG^d_{m,S}$ l'ouvert de $\mA_S^d$ où les coordonnées ne s'annulent pas, $U=f^{-1}(\mG^d_{m,S})$ 
et $j\colon U\rightarrow X$ l'injection canonique. Soient $n$ un entier $\geq 1$ inversible dans $\co_S$,  $\Lambda=\mZ/n\mZ$,
$F$ un faisceau de $\Lambda$-modules localement constant et constructible de $X_\et$. 
Notons $\pi\colon \mA_S^d\rightarrow \mA_S^d$ le morphisme défini par l'élévation à la puissance $n$-ième des coordonnées de $\mA_S^d$,
$X'$ le changement de base de $X$ par $\pi$, $\varpi\colon X'\rightarrow X$ la projection canonique, $U'=\varpi^{-1}(U)$,
$j'\colon U'\rightarrow X'$ l'injection canonique et $F'=\varpi^*(F)$. Alors, pour tout entier $q\geq 1$, le morphisme de changement de base
\begin{equation}\label{Kpun14a}
\varpi^*(\rR^qj_*(j^*F))\rightarrow \rR^qj'_*(j'^*F')
\end{equation}
est nul. 
\end{lem}

En effet, la question étant locale pour la topologie étale de $X'$ et donc aussi pour celle de $X$, on peut se borner au cas où $F$ est constant de valeur $\Lambda$. 
Notons $(A_\alpha)_{1\leq \alpha\leq d}$ les axes de coordonnées de $\mA^d_S$ et pour tout $1\leq \alpha\leq d$, posons 
$D_\alpha=X\times_{\mA_S^d} A_\alpha$ et notons $i_\alpha\colon D_\alpha\rightarrow X$ l'injection canonique. 
En vertu de (\cite{travaux-gabber} XVI 3.1.4), on a un isomorphisme canonique
\begin{equation}
\rR^1 j_*(\Lambda_U)\stackrel{\sim}{\rightarrow}\oplus_{1\leq \alpha\leq d}\ i_{\alpha*}(\Lambda_{D_\alpha}(-1)).
\end{equation}
De plus, pour tout entier $q\geq 1$, le morphisme 
\begin{equation}
\wedge^q (\rR^1j_*(\Lambda_U))\rightarrow \rR^qj_*(\Lambda_U)
\end{equation}
défini par cup-produit est un isomorphisme. 
De même, considérons $X'$ comme un $\mA^d_S$-schéma lisse via la projection canonique 
$f'\colon X'\rightarrow \mA_S^d$ et pour tout $1\leq \alpha\leq d$, posons $D'_\alpha=X'\times_{\mA_S^d} A_\alpha$ et notons 
$i'_\alpha\colon D'_\alpha\rightarrow X'$ l'injection canonique.
On a alors des isomorphismes canoniques
\begin{eqnarray}
\rR^1 j'_*(\Lambda_{U'})&\stackrel{\sim}{\rightarrow}&\oplus_{1\leq \alpha\leq d}\ i'_{\alpha*}(\Lambda_{D'_\alpha}(-1)),\\
\wedge^q (\rR^1j'_*(\Lambda_{U'}))&\stackrel{\sim}{\rightarrow}& \rR^qj'_*(\Lambda_{U'}).
\end{eqnarray}
Pour tout $1\leq \alpha\leq d$, $D'_\alpha$ s'identifie canoniquement au sous-schéma réduit de $X'$ sous-jacent à $D_\alpha\times_XX'$. De plus, on a 
$\varpi^*(D_\alpha)=nD'_\alpha$ en tant que diviseurs de Cartier. Il résulte de (\cite{travaux-gabber} XVI 3.4.8) que le diagramme 
\begin{equation}
\xymatrix{
{\varpi^*(\rR^1j_*(\Lambda_U))}\ar[r]^-(0.5)\sim\ar[d]_\lambda&{\oplus_{1\leq \alpha\leq d}\ \varpi^*(i_{\alpha*}(\Lambda_{D_\alpha}(-1)))}\ar[d]^\gamma\\
{\rR^1j'_*(\Lambda_{U'})}\ar[r]^-(0.5)\sim&{\oplus_{1\leq \alpha\leq d}\ i'_{\alpha*}(\Lambda_{D'_\alpha}(-1))}}
\end{equation}
où $\lambda$ est le morphisme de changement de base et $\gamma$ est induit pour chaque $1\leq \alpha \leq d$ par $n$ fois le morphisme de changement de base,
est commutatif. La proposition s'ensuit puisque $\gamma$ est nul.

\begin{lem}\label{Kpun15}
Soient $f\colon Y\rightarrow X$, $g\colon Y'\rightarrow Y$ deux morphismes de schémas
tels que $f$ soit fini et que $g$ soit étale, $\ox$ un point géométrique de $X$ tels que le morphisme
\begin{equation}\label{Kpun15a}
Y'\otimes_X\kappa(\ox)\rightarrow Y\otimes_X\kappa(\ox)
\end{equation}
induit par $g$ soit surjectif. Alors, il existe un $X$-schéma étale $\ox$-pointé $X'$ et un $Y$-morphisme $Y\times_XX'\rightarrow Y'$. 
\end{lem}

Il suffit de montrer que si $X$ est le spectre d'un anneau local strictement hensélien de point fermé $\ox$, alors $g$ admet une
section (\cite{ega4} 8.8.2(i)). On a un $X$-isomorphisme $Y\simeq\amalg_{\oy\in Y\otimes_X\kappa(\ox)}Y_{(\oy)}$,
où $Y_{(\oy)}$ est le localisé strict de $Y$ en le point géométrique $\oy$. 
On peut évidemment se réduire au cas où $Y$ est connexe et non vide. Donc $Y$ est isomorphe au spectre
d'un anneau local strictement hensélien. 
D'après \eqref{Kpun15a} et (\cite{ega4} 18.5.11), $Y'$ est la somme de deux schémas, dont l'un est isomorphe à $Y$; d'où la proposition.

\begin{prop}[\cite{achinger} 8.1]\label{Kpun16}
Soient $S$ un $\mQ$-schéma noethérien régulier, $X$ un $S$-schéma lisse, $D$ un diviseur à croisements normaux sur $X$ 
relativement à $S$ {\rm (\cite{sga1} XIII 2.1)}, 
$U$ l'ouvert complémentaire de $D$ dans $X$, $\ox$ un point géométrique de $X$, 
$F$ un faisceau abélien de torsion, localement constant et constructible sur $U$,
$q$ un entier $\geq 1$, $\xi\in \rH^q(U,F)$. Alors, il existe un $X$-schéma étale $\ox$-pointé $Y$, et posant $V=U\times_XY$, un revêtement étale 
surjectif $V'\rightarrow V$ tels que l'image canonique de $\xi$ dans $\rH^q(V',F)$ soit nulle.
\end{prop}

On notera d'abord que l'assertion est immédiate si $\ox$ est à support dans $U$, auquel cas on peut prendre $Y=V=V'$ (\cite{sga4} V 3.1). 
On peut donc se borner au cas où $\ox$ est à support dans $D$. 
Soient $U'\rightarrow U$ un revêtement étale surjectif tel que le faisceau $F|U'$ soit constant, $X'$ la fermeture intégrale de $X$ dans $U'$. 
En vertu de \ref{Kpun8}, $X'$ est lisse sur $S$, et $U'$ est le complémentaire dans $X'$ d'un diviseur à croisements normaux sur $X'$ relativement à $S$.
Supposons que la proposition soit démontrée pour $(X',U',F|U')$ et pour tout point géométrique de $X'$ au-dessus de $\ox$. 
Il existe donc un morphisme étale $X''\rightarrow X'$ tel que  le morphisme
\begin{equation}\label{Kpun16a}
X''\otimes_X\kappa(\ox)\rightarrow X'\otimes_X\kappa(\ox)
\end{equation}
soit surjectif, et posant $U''=X''\times_XU$, un revêtement étale surjectif $W\rightarrow U''$ tel que l'image canonique de $\xi$ dans $\rH^q(W,F)$ soit nulle.
D'après \ref{Kpun15}, il existe un $X$-schéma étale $\ox$-pointé $Y$ et un $X'$-morphisme $X'\times_XY\rightarrow X''$. Posons $V=Y\times_XU$
et considérons le diagramme commutatif à carrés cartésiens 
\begin{equation}\label{Kpun16b}
\xymatrix{
{W}\ar[d]&{W\times_{U''}(U'\times_UV)}\ar[d]^v\ar[l]&\\
U''\ar[d]&{U'\times_UV}\ar[l]\ar[d]\ar[r]^-(0.5)u&{V}\ar[d]\\
X''&{X'\times_XY}\ar[l]\ar[r]&Y}
\end{equation}
Le couple formé de $Y$ et du revêtement étale surjectif $u\circ v$ répond alors à la question.
On peut donc se borner au cas où $F$ est constant de valeur $\mZ/n\mZ$ et $n$ est un entier $\geq1$. 

La question étant locale pour la topologie étale sur $X$, on peut supposer qu'il existe un entier $d\geq1$ et un morphisme lisse $f\colon X\rightarrow \mA_S^d$  
tels que $D$ soit l'image inverse du diviseur des coordonnées de $\mA_S^d$. Notons $j\colon U\rightarrow X$ l'injection canonique 
et considérons le morphisme composé 
\begin{equation}\label{Kpun16c}
\rH^q(U,F)\rightarrow \rH^0(X,\rR^qj_*F) \rightarrow (\rR^qj_*F)_\ox,
\end{equation}
où la première flèche est induite par la suite spectrale de Cartan-Leray et la seconde flèche est le morphisme canonique. 
Notons $\pi\colon \mA_S^d\rightarrow \mA_S^d$ le morphisme défini par l'élévation à la puissance $n$-ième des coordonnées de $\mA_S^d$,
$X'$ le changement de base de $X$ par $\pi$, $\varpi\colon X'\rightarrow X$ la projection canonique, $U'=\varpi^{-1}(U)$, $F'=F|U'$ et 
$j'\colon U'\rightarrow X'$ l'injection canonique. 
Il résulte de  (\cite{sga4} VIII 5.5) et de la fonctorialité de la suite spectrale de Cartan-Leray que le diagramme 
\begin{equation}\label{Kpun16d}
\xymatrix{
{\rH^q(U,F)}\ar[r]\ar[d]_\lambda&{(\rR^qj_*F)_\ox}\ar[d]^\gamma\\
{\rH^q(U',F')}\ar[r]&{\oplus_{\ox'\in X'\otimes_X\kappa(\ox)}(\rR^qj'_*F')_{\ox'}}}
\end{equation}
où les flèches horizontales sont les morphismes composés \eqref{Kpun16c}, 
$\lambda$ est induit par $\varpi^*$ et $\gamma$ est induit par le morphisme de changement de base \eqref{Kpun14a}, est commutatif. 
En vertu de \ref{Kpun14}, $\gamma$ est nul. Il existe donc un morphisme étale $X''\rightarrow X'$ tel que le morphisme
\begin{equation}\label{Kpun16e}
X''\otimes_X\kappa(\ox)\rightarrow X'\otimes_X\kappa(\ox)
\end{equation}
soit surjectif, et posant $U''=X''\times_XU$, l'image canonique de $\xi$ dans $\rH^q(U'',F)$ soit nulle (\cite{sga4} V 5.1(1) et IV (6.8.4)). 
D'après \ref{Kpun15}, il existe un $X$-schéma étale $\ox$-pointé $Y$ et un $X'$-morphisme $X'\times_XY\rightarrow X''$. Posons $V=Y\times_XU$
et considérons le diagramme commutatif à carrés cartésiens 
\begin{equation}\label{Kpun16f}
\xymatrix{
U''\ar[d]&{U'\times_UV}\ar[l]\ar[d]\ar[r]^-(0.5)w&{V}\ar[d]\\
X''&{X'\times_XY}\ar[l]\ar[r]&Y}
\end{equation}
Le couple formé de $Y$ et du revêtement étale surjectif $w$ répond alors à la question.

\begin{cor}\label{Kpun19}
Soient $k$ un corps algébriquement clos de caractéristique $0$, $X$ un schéma lisse sur $k$, 
$D$ un diviseur à croisements normaux sur $X$, $X^\circ$ l'ouvert complémentaire de $D$ dans $X$, 
$\ox$ un point géométrique de $X$, $X'$ le localisé strict de $X$ en $\ox$. 
Alors, $X'\times_XX^\circ$ est un schéma $K(\pi,1)$. 
\end{cor}

Cela résulte de \ref{Kpun18} et  \ref{Kpun16} puisque $X'\times_XX^\circ$ est noethérien et intègre. 

\subsection{}\label{Kpun20}
Soient $V$ un anneau de valuation discrète, $S=\Spec(V)$,  $s$ (resp.  $\eta$, resp. $\oeta$) le point fermé 
(resp.  le point générique, resp. un point géométrique générique) de $S$. 
On suppose que le corps des fractions de $V$ est de caractéristique $0$, 
et que son corps résiduel est parfait de caractéristique $p>0$.
On munit $S$ de la structure logarithmique $\cM_S$ définie par son point fermé, 
autrement dit, $\cM_S=u_*(\co_\eta^\times)\cap \co_S$, où $u\colon \eta\rightarrow S$ est l'injection canonique.
On renvoie à (\cite{agt} II.5) pour un lexique de géométrie logarithmique. 

\begin{prop}[\cite{achinger} 6.1]\label{Kpun21}
Conservons les hypothèses de \ref{Kpun20}, soient de plus $(X,\cM_X)$ un schéma logarithmique fin,
$f\colon (X,\cM_X)\rightarrow (S,\cM_S)$ un morphisme lisse tel que le schéma usuel $X_\eta$ soit lisse sur $\eta$, 
$\ox$ un point géométrique de $X$. 
Alors, il existe un voisinage étale $U$ de $\ox$ dans $X$ tel que $U_\eta$ soit un schéma $K(\pi,1)$. 
\end{prop}

\begin{cor}[\cite{achinger} 9.5]\label{Kpun22}
Conservons les hypothèses de \ref{Kpun20}, supposons de plus que $S$ soit strictement local. Soient $(X,\cM_X)$ un schéma logarithmique fin,
$f\colon (X,\cM_X)\rightarrow (S,\cM_S)$ un morphisme lisse tel que $X_\eta$ soit lisse sur $\eta$, 
$\ox$ un point géométrique de $X$ au-dessus de $s$, 
$X'$ le localisé strict de $X$ en $\ox$. Alors, $X'_\oeta$ est un schéma $K(\pi,1)$.
\end{cor}

Montrons d'abord que le schéma $X'_\oeta$ n'a qu'un nombre fini de composantes connexes. Soit $\ell$ un nombre premier différent de $p$. 
D'après (\cite{sga45} Th. finitude 3.2), le $\mF_\ell$-espace vectoriel $\rH^0(X'_\oeta,\mF_\ell)$ est de dimension finie. 
Il s'ensuit que l'ensemble $\fC$ des sous-schémas ouverts et fermés de $X'_\oeta$ est fini (\cite{sga4} VIII 6.1). 
Pour tout point $z$ de $X'_\oeta$, notons $U_z$ l'intersection de tous les sous-schémas ouverts et fermés de $X'_\oeta$ contenant $z$. 
Comme $\fC$ est fini, $U_z$ est ouvert et fermé dans $X'_\oeta$, autrement dit, c'est un objet de $\fC$. Par suite, $U_z$ est connexe,  
et est donc égal à la composante connexe de  $X'_\oeta$ contenant $z$. On en déduit que $X'_\oeta$ n'a qu'un nombre fini de composantes connexes.
La proposition résulte alors de \ref{Kpun18} et \ref{Kpun21}.

\begin{rema}
On établira dans \ref{cad10} l'énoncé le plus général de (\cite{achinger} 9.5).
\end{rema}

\section{\texorpdfstring{\'Epaississements infinitésimaux $p$-adiques de Fontaine}
{\'Epaississements infinitésimaux p-adiques de Fontaine}}\label{EIPF}

\subsection{}\label{eip1}
Dans cette section, $p$ désigne un nombre premier. 
Commençons par rappeler la cons\-truction suivante due à Grothendieck (\cite{grot1} IV 3.3).
Soient $A$ une $\mZ_{(p)}$-algèbre, $n$ entier $\geq 1$. L'homomorphisme d'anneaux \eqref{notconv1a}
\begin{equation}\label{eip1a}
\Phi_{n+1}\colon 
\begin{array}[t]{clcr}
\rW_{n+1}(A/p^nA)&\rightarrow& A/p^nA\\
(x_1,\dots,x_{n+1})&\mapsto& x_1^{p^n}+p x_2^{p^{n-1}}+\dots+p^{n}x_{n+1}
\end{array}
\end{equation}
s'annule sur $\rV^n(A/p^nA)$ et induit donc par passage au quotient un homomorphisme d'anneaux
\begin{equation}\label{eip1b}
\Phi'_{n+1}\colon 
\begin{array}[t]{clcr}
\rW_{n}(A/p^nA)&\rightarrow& A/p^nA,\\
(x_1,\dots,x_n)&\mapsto&x_1^{p^n}+p x_2^{p^{n-1}}+\dots+p^{n-1}x_n^p.
\end{array}
\end{equation}
Ce dernier s'annule sur 
\begin{equation}\label{eip1c}
\rW_n(pA/p^nA)=\ker(\rW_n(A/p^nA)\rightarrow \rW_n(A/pA))
\end{equation}
et se factorise à son tour en un homomorphisme d'anneaux
\begin{equation}\label{eip1d}
\theta_n\colon \rW_{n}(A/pA)\rightarrow A/p^nA.
\end{equation}
Il résulte aussitôt de la définition que le diagramme 
\begin{equation}\label{eip1e}
\xymatrix{
{\rW_{n+1}(A/pA)}\ar[r]^-(0.4){\theta_{n+1}}\ar[d]_{\rR\rF}&{A/p^{n+1}A}\ar[d]\\
{\rW_{n}(A/pA)}\ar[r]^-(0.4){\theta_n}&{A/p^nA}}
\end{equation}
où $\rR$ est le morphisme de restriction \eqref{notconv1b}, 
$\rF$ est le Frobenius \eqref{notconv1d} et la flèche non libellée est l'homomorphisme canonique, est commutatif. 

Pour tout homomorphisme de $\mZ_{(p)}$-algèbres commutatives $\varphi\colon A\rightarrow B$, 
le diagramme  
\begin{equation}\label{eip1f}
\xymatrix{
{\rW_n(A/pA)}\ar[r]\ar[d]_{\theta_n}&{\rW_n(B/pB)}\ar[d]^{\theta_n}\\
{A/p^nA}\ar[r]&{B/p^nB}}
\end{equation}
où les flèches horizontales sont les morphismes induits par $\varphi$, est commutatif. 

\subsection{}\label{eipo3}
Soient $A$ une $\mZ_{(p)}$-algèbre, $\hA$ le séparé complété $p$-adique de $A$.
On désigne par $A^\flat$ la limite projective du système projectif $(A/pA)_{\mN}$ 
dont les morphismes de transition sont les itérés de l'endomorphisme de Frobenius de $A/pA$.
\begin{equation}\label{eipo3a}
A^\flat=\underset{\underset{x\mapsto x^p}{\longleftarrow}}{\lim}A/pA.
\end{equation} 
C'est un anneau parfait de caractéristique $p$. 
Pour tout entier $n\geq 1$, la projection canonique $A^\flat\rightarrow A/pA$ sur la $(n+1)$-ième composante
du système projectif $(A/pA)_{\mN}$ ({\em i.e.}, la composante d'indice $n$) induit un homomorphisme
\begin{equation}\label{eipo3b}
\nu_n\colon \rW(A^\flat)\rightarrow \rW_n(A/pA).
\end{equation}
Comme $\nu_n=\rF\circ\rR\circ \nu_{n+1}$, on obtient par passage à la limite projective un homomorphisme 
\begin{equation}\label{eipo3c}
\nu\colon \rW(A^\flat)\rightarrow \underset{\underset{n\geq 0}{\longleftarrow}}{\lim}\ \rW_n(A/pA),
\end{equation}
où les morphismes de transition de la limite projective sont les morphismes $\rF \rR$. On vérifie aussitôt  
qu'il est bijectif. Compte tenu de \eqref{eip1e},
les homomorphismes $\theta_n$ induisent par passage à la limite projective un 
homomorphisme 
\begin{equation}\label{eipo3d}
\theta\colon \rW(A^\flat)\rightarrow \hA.
\end{equation}
On retrouve l'homomorphisme défini par Fontaine 
(\cite{fontaine1} 2.2). On pose
\begin{equation}\label{eipo3e}
\cA_2(A)=\rW(A^\flat)/\ker(\theta)^2,
\end{equation}
et on note encore $\theta\colon \cA_2(A)\rightarrow \hA$ l'homomorphisme induit par $\theta$ (cf. \cite{fontaine3} 1.2.2). 

Pour tout homomorphisme de $\mZ_{(p)}$-algèbres commutatives $\varphi\colon A\rightarrow B$, 
le diagramme  
\begin{equation}\label{eip3f}
\xymatrix{
{\rW(A^\flat)}\ar[r]\ar[d]_{\theta}&{\rW(B^\flat)}\ar[d]^{\theta}\\
{\hA}\ar[r]&{\hB}}
\end{equation}
où les flèches horizontales sont les morphismes induits par $\varphi$, est commutatif \eqref{eip1f}.  
La correspondance $A\mapsto \cA_2(A)$ est donc fonctorielle. 

\begin{rema}\label{eip2}
Soit $k$ un corps parfait.  
La projection canonique $\rW(k)^\flat\rightarrow k$ sur la première composante ({\em i.e.}, d'indice 0)
est un isomorphisme. Elle induit donc un isomorphisme $\rW(\rW(k)^\flat)\stackrel{\sim}{\rightarrow}\rW(k)$,
gs'identifie alors à l'endomorphisme de Frobenius de $\rW(k)$. 
\end{rema}

\begin{prop}[\cite{agt} II.9.5; \cite{tsuji1} A.1.1 et A.2.2]\label{eip4}
Soit $A$ une $\mZ_{(p)}$-algèbre vérifiant les conditions suivantes~:
\begin{itemize}
\item[{\rm (i)}] $A$ est $\mZ_{(p)}$-plat.
\item[{\rm (ii)}] $A$ est intégralement clos dans $A[\frac 1 p]$.
\item[{\rm (iii)}] L'endomorphisme de Frobenius absolu  de $A/pA$ est surjectif. 
\item[{\rm (iv)}] Il existe une suite $(p_n)_{n\geq 0}$ d'éléments de $A$
tels que $p_0=p$ et $p_{n+1}^p=p_n$ pour tout $n\geq 0$.  
\end{itemize}
On désigne par $\varpi$ l'élément de $A^\flat$ induit par la suite $(p_n)_{n\geq 0}$ et on pose  
\begin{equation}\label{eip4a}
\xi=[\varpi]-p \in \rW(A^\flat),
\end{equation}
où $[\ ]$ est le représentant multiplicatif.
Alors la suite
\begin{equation}\label{eip4b}
0\longrightarrow \rW(A^\flat)\stackrel{\cdot \xi}{\longrightarrow} \rW(A^\flat)\stackrel{\theta}{\longrightarrow} 
\hA \longrightarrow 0
\end{equation}
est exacte.
\end{prop}

\subsection{}\label{eip5}
Soient $A$ une $\mZ_{(p)}$-algèbre, $X=\Spec(A)$, $\cM_X$ une structure logarithmique sur $X$ (\cite{agt} II.5.9), 
$M$ un monoïde, $u\colon M\rightarrow \Gamma(X,\cM_X)$ un homomorphisme.  
Considérons le système projectif de monoïdes multiplicatifs $(A)_{n\in \mN}$, 
où les morphismes de transition sont tous égaux à l'élévation à la puissance $p$-ième.
On désigne par $Q$ le produit fibré du diagramme d'homomorphismes de monoïdes 
\begin{equation}\label{eip5a}
\xymatrix{
&{M}\ar[d]\\
{\underset{\underset{x\mapsto x^p}{\longleftarrow}}{\lim}\ A}\ar[r]&A}
\end{equation}
où  la flèche horizontale est la projection sur la première composante ({\em i.e.}, d'indice $0$)
et la flèche verticale est le composé de $u$ et de l'homomorphisme canonique $\Gamma(X,\cM_X)\rightarrow A$, 
et par $q$ l'homomorphisme composé
\begin{equation}\label{eip5b}
Q\longrightarrow
\underset{\underset{x\mapsto x^p}{\longleftarrow}}{\lim}\ A \longrightarrow A^\flat \stackrel{[\ ]}{\longrightarrow} \rW(A^\flat),
\end{equation} 
où la première et la deuxième flèches sont les homomorphismes canoniques \eqref{eipo3a}
et $[\ ]$ est le représentant multiplicatif. Il résulte aussitôt des définitions que le diagramme 
\begin{equation}\label{eip5c}
\xymatrix{
Q\ar[r]\ar[d]_{q}&{M}\ar[d]\\
{\rW(A^\flat)}\ar[r]^-(0.4)\theta&{\hA}}
\end{equation}
où les flèches non libellées sont les morphismes canoniques, est commutatif. 

On pose $\hX=\Spec(\hA)$ que l'on munit de la structure logarithmique $\cM_{\hX}$ image inverse de $\cM_X$.  
On munit $\Spec(\rW(A^\flat))$ de la structure logarithmique $\cQ$ associée à la structure pré-logarithmique définie par $q$ 
\eqref{eip5b}. D'après \eqref{eip5c}, $\theta$ induit un morphisme 
\begin{equation}\label{eip5d}
(\hX,\cM_\hX)\rightarrow (\Spec(\rW(A^\flat)),\cQ).
\end{equation}

\begin{prop}[\cite{agt} II.9.7]\label{eip6}
Conservons les hypothèses de \ref{eip5}, notons, de plus, $X^\circ$ l'ouvert maximal de $X$ où la structure logarithmique 
$\cM_X$ est triviale et supposons les conditions suivantes remplies~:
\begin{itemize}
\item[{\rm (a)}] $A$ est intègre et normal. 
\item[{\rm (b)}] $X^\circ$ est un $\mQ$-schéma non-vide et simplement connexe. 
\item[{\rm (c)}] $M$ est intègre et il existe un monoïde fin et saturé $M'$ et un homomorphisme $v\colon M'\rightarrow M$
tels que l'homomorphisme induit  $M'\rightarrow M/M^\times$ soit un isomorphisme.  
\end{itemize}
Alors~:
\begin{itemize}
\item[{\rm (i)}] Le monoïde $Q$ est intègre et le groupe $M'^\gp$ est libre. 
\item[{\rm (ii)}] On peut compléter le diagramme \eqref{eip5a} en un diagramme commutatif 
\begin{equation}\label{eip6a}
\xymatrix{
{M'}\ar[r]^v\ar[d]_-(0.4)w&{M}\ar[d]\\
{\underset{\underset{x\mapsto x^p}{\longleftarrow}}{\lim}\ A}\ar[r]&A}
\end{equation}
Notons $\beta\colon M'\rightarrow Q$ l'homomorphisme induit. 
\item[{\rm (iii)}] La structure logarithmique $\cQ$ sur $\Spec(\rW(A^\flat))$ est associée à la structure pré-logarithmi\-que 
définie par l'homomorphisme composé
\begin{equation}\label{eip6b}
M'\stackrel{\beta}{\rightarrow} Q \stackrel{q}{\rightarrow} \rW(A^\flat).
\end{equation}
En particulier, le schéma logarithmique $(\Spec(\rW(A^\flat)),\cQ)$ est fin et saturé.
\item[{\rm (iv)}] Si de plus, l'homomorphisme composé $u\circ v\colon M'\rightarrow \Gamma(X,\cM_X)$ est 
une carte pour $X$, alors  le morphisme \eqref{eip5d} est strict.
\end{itemize}
\end{prop}

\section{\texorpdfstring{Faisceaux de $\alpha$-modules}{Faisceaux de alpha-modules}}\label{alpha}

\subsection{}\label{alpha1}
Dans cette section, $\Lambda$ désigne un sous-groupe ordonné et dense de $\mR$, 
$\Lambda^+$ l'ensemble des éléments strictement positifs de $\Lambda$, et 
$R$ un anneau muni d'une suite d'idéaux principaux $\fm_\varepsilon$, indexée par $\Lambda^+$. 
Pour tout $\varepsilon\in \Lambda^+$, on choisit un générateur $\pi^\varepsilon\in R$ de $\fm_\varepsilon$. 
On suppose, de plus, que pour tout $\varepsilon\in \Lambda^+$, $\pi^\varepsilon$ n'est pas un diviseur de zéro dans $R$ 
et que pour tous $\varepsilon,\delta\in \Lambda^+$, 
$\pi^\varepsilon\cdot \pi^\delta=u_{\varepsilon,\delta}\cdot \pi^{\varepsilon+\delta}$, 
où $u_{\varepsilon,\delta}$ est une unité de $R$. On pose $\fm=\cup_{\varepsilon\in \Lambda^+}\fm_\varepsilon$. 
Ces hypothèses sont celles fixées dans (\cite{faltings2} page 186 et \cite{agt} V Notations). 
On notera que $\fm$ est $R$-plat et que $\fm^2=\fm$. 

\subsection{}\label{alpha2}
Soient $\cA$ une catégorie abélienne qui est une $\mU$-catégorie \eqref{notconv3} (\cite{sga4} I 1.1), 
$\End(\id_\cA)$ l'anneau des endomorphismes du foncteur identique de $\cA$, $\varphi\colon R\rightarrow \End(\id_\cA)$ un homomorphisme.  
Pour tout objet $M$ de $\cA$ et tout $\gamma\in R$, on note $\mu_\gamma(M)$ l'endomorphisme de $M$ défini par $\varphi(\gamma)$. 
On observera que pour tout morphisme $f\colon M\rightarrow N$ de $\cA$, on a $\mu_\gamma(N)\circ f= f\circ \mu_\gamma(M)$.
En particulier, pour tous objets $M$ et $N$ de $\cA$, $\Hom(M,N)$ est naturellement muni d'une structure de $R$-module. 

Suivant (\cite{faltings2} page 187), on dit qu'un objet $M$ de $\cA$  
est {\em $\alpha$-nul} (ou {\em presque-nul}) s'il est annulé par tout élément de $\fm$, 
{\em i.e.}, si  $\mu_\gamma(M)=0$ pour tout $\gamma\in \fm$.
On vérifie aussitôt que pour toute suite exacte $0\rightarrow M'\rightarrow M\rightarrow M''\rightarrow 0$ de $\cA$, pour que 
$M$ soit $\alpha$-nul, il faut et il suffit que $M'$ et $M''$ soient $\alpha$-nuls. 
On appelle catégorie des {\em $\alpha$-objets} (ou {\em presque-objets}) de $\cA$ et l'on note 
$\acA$ le quotient de la catégorie $\cA$ par la sous-catégorie épaisse formée des objets $\alpha$-nuls (\cite{gabriel} III §~1). On note 
\begin{equation}\label{alpha2a}
\alpha\colon \cA\rightarrow \acA, \ \ \ M\mapsto \alpha(M)
\end{equation}
le foncteur canonique~; on notera aussi $M^\alpha$ au lieu de $\alpha(M)$ lorsqu'il n'y a aucun risque de confusion. 
La catégorie $\acA$ est abélienne et le foncteur $\alpha$ est exact (\cite{gabriel} III §~1 prop.~1).
Si $f\colon M\rightarrow N$ est un morphisme  de $\cA$, pour que $\alpha(f)$ soit nul (resp. un monomorphisme,
resp. un épimorphisme), il faut et il suffit que $\im(f)$ (resp. $\ker(f)$, resp. $\coker(f)$) soit $\alpha$-nul 
(\cite{gabriel} III §~1 lem.~2). On dit que $f$ est {\em $\alpha$-injectif} (resp. {\em $\alpha$-surjectif}, resp. un {\em $\alpha$-isomorphisme})  
si $\alpha(f)$ est injectif (resp. surjectif, resp. un isomorphisme), autrement dit si son noyau (resp. son conoyau, resp. son noyau et son conoyau) sont $\alpha$-nuls. 

La famille des $\alpha$-isomorphismes de $\cA$ permet un calcul de fractions bilatère (\cite{illusie1} I 1.4.2).
La catégorie $\acA$ s'identifie à la catégorie localisée de $\cA$ 
par rapport aux $\alpha$-isomorphismes et $\alpha$ est le foncteur canonique (de localisation). 
On laissera le soin au lecteur de vérifier ces propriétés, valables d'ailleurs pour tout quotient d'une catégorie 
abélienne par une sous-catégorie épaisse (\cite{gz} I 2.5(d)).

\begin{lem}\label{alpha3}
Les hypothèses étant celles de \eqref{alpha2}, soient, de plus, $f\colon M\rightarrow N$ un morphisme de $\cA$, $\gamma \in R$ tels que 
le noyau et le conoyau de $f$ soient annulés par $\gamma$. Alors, il existe un morphisme $g\colon N\rightarrow M$ de $\cA$ tel que 
$g\circ f=\mu_{\gamma^2}(M)$ et $f\circ g=\mu_{\gamma^2}(N)$.
\end{lem}

En effet, comme $\mu_\gamma(\coker(f))=0$, le morphisme $\mu_\gamma(N)\colon N\rightarrow N$ se factorise uniquement 
à travers un morphisme $\varphi_\gamma\colon N\rightarrow \im(f)$. 
Comme $\mu_\gamma(\ker(f))=0$, le morphisme $\mu_\gamma(M)\colon M\rightarrow M$ se factorise uniquement 
à travers un morphisme $\psi_\gamma\colon \im(f)\rightarrow M$. 
Il est clair que $g=\psi_\gamma\circ \varphi_\gamma\colon N\rightarrow M$ répond à la question.

\subsection{}\label{alpha33}
Soit $\cA$ une catégorie abélienne tensorielle qui est une $\mU$-catégorie, autrement dit, $\cA$ est une catégorie abélienne
munie d'un foncteur bi-additif  $\otimes\colon \cA\times \cA\rightarrow \cA$ et d'un objet unité $A$, 
vérifiant certaines conditions (\cite{dm} 1.15), et soit $\phi\colon R\rightarrow \End(A)$ un homomorphisme. 
On a un homomorphisme canonique $\End(A)\rightarrow \End(\id_\cA)$ (cf. \cite{sr} I 1.3.3.3 et 2.3.3). 
On peut donc définir la catégorie quotient $\acA$ suivant \ref{alpha2}. Il résulte aussitôt de \ref{alpha3} que pour tout 
$\alpha$-isomorphisme $f\colon M\rightarrow M'$ de $\cA$ et tout $N\in \ob(\cA)$, $f\otimes \id_N\colon M\otimes N\rightarrow M'\otimes N$ 
est un $\alpha$-isomorphisme. Par suite, le produit tensoriel induit un foncteur 
\begin{equation}\label{alpha33a}
\acA\times \acA\rightarrow \acA, \ \ \ (M,N)\mapsto M\otimes N,
\end{equation}
qui fait de $\acA$ une catégorie abélienne tensorielle, dont $A^\alpha$ est un objet unité. 

Le foncteur $\alpha$ induit un homomorphisme $\End(A)\rightarrow \End(A^\alpha)$. Par suite, pour tous objets $M$ et $N$ de $\acA$, 
$\Hom_{\acA}(M,N)$ est canoniquement muni d'une structure de $\End(A)$-module. 
On voit aussitôt que pour tout $P\in \ob(\cC)$, l'application 
\begin{equation}\label{alpha33b}
\Hom_{\acA}(M,N)\rightarrow \Hom_{\acA}(M\otimes P,N\otimes P)
\end{equation}
définie par fonctorialité est $\End(A)$-linéaire (\cite{sr} I 2.2.6). 
 
\subsection{}\label{alpha4}
Pour toute $R$-algèbre $A$ (appartenant à $\mU$), on désigne par $\bMod(A)$ la catégorie abélienne tensorielle 
des $A$-modules qui se trouvent dans $\mU$. Prenant pour $\phi\colon R\rightarrow A=\End(A)$ l'homomorphisme structural,  
on appelle catégorie des {\em $\alpha$-$A$-modules} et l'on note $\aMod(A)$ 
le quotient de la catégorie abélienne $\bMod(A)$ par la sous-catégorie épaisse des $A$-modules $\alpha$-nuls. 
Nous utiliserons les conventions de notation de \ref{alpha2} et \ref{alpha33}. On observera en particulier que pour tous
$\alpha$-$A$-modules $M$ et $N$, $\Hom_{\aMod(A)}(M,N)$ est naturellement muni d'une structure de $A$-module.  

Pour tout homomorphisme de $R$-algèbres $A\rightarrow B$, le foncteur d'oubli  
$\bMod(B)\rightarrow \bMod(A)$ induit un foncteur exact
\begin{equation}\label{alpha4a}
\aMod(B)\rightarrow \aMod(A).
\end{equation}

\begin{lem}\label{alpha7}
Pour qu'un morphisme de $R$-modules $M\rightarrow N$ soit un $\alpha$-isomorphisme,
il faut et il suffit que le morphisme induit $\fm\otimes_RM\rightarrow \fm\otimes_RN$ soit un isomorphisme. 
\end{lem}
En effet, la condition est suffisante puisque le morphisme canonique 
$\fm\otimes_RM\rightarrow M$ est un $\alpha$-isomorphisme \eqref{alpha3},
et elle est nécessaire car $\fm$ est $R$-plat et pour tout $R$-module $\alpha$-nul $P$, $\fm\otimes_RP=0$. 

\subsection{}\label{alpha8} 
Soit $A$ une $R$-algèbre. 
D'après \ref{alpha7}, la catégorie des $\alpha$-isomorphismes de $\bMod(A)$ de but un $A$-module $M$ donné 
admet un objet initial à savoir le morphisme canonique $\fm\otimes_RM\rightarrow M$. 
Par suite, pour tous $A$-modules $M$ et $N$, on a un isomorphisme canonique
\begin{equation}\label{alpha8a}
\Hom_{\aMod(A)}(M^\alpha,N^\alpha)\stackrel{\sim}{\rightarrow}\Hom_{\bMod(A)}(\fm\otimes_RM,N).
\end{equation}
On voit aussitôt que cet isomorphisme est $A$-linéaire.

On désigne par $\sigma_*$ le foncteur 
\begin{equation}\label{alpha8b}
\sigma_*\colon \aMod(A)\rightarrow \bMod(A),\ \ \ P\mapsto \Hom_{\aMod(A)}(A^\alpha,P),
\end{equation}
et par $\sigma_!$ le foncteur
\begin{equation}\label{alpha8c}
\sigma_!\colon \aMod(A)\rightarrow \bMod(A), \ \ \ P\mapsto \fm\otimes_R\sigma_*(P).
\end{equation}
D'après \eqref{alpha8a}, pour tout $A$-module $M$, on a un isomorphisme canonique fonctoriel
\begin{equation}\label{alpha8e}
\sigma_*(M^\alpha)\stackrel{\sim}{\rightarrow}\Hom_{R}(\fm,M).
\end{equation}
On en déduit que pour tout homomorphisme de $R$-algèbres $A\rightarrow B$, les diagrammes
\begin{equation}\label{alpha8d}
\xymatrix{
{\aMod(B)}\ar[r]^{\sigma_*}\ar[d]&{\bMod(B)}\ar[d]\\
{\aMod(A)}\ar[r]^{\sigma_*}&{\bMod(A)}}\ \ \
\xymatrix{
{\aMod(B)}\ar[r]^{\sigma_!}\ar[d]&{\bMod(B)}\ar[d]\\
{\aMod(A)}\ar[r]^{\sigma_!}&{\bMod(A)}}
\end{equation}
où les flèches verticales sont les foncteurs d'oubli \eqref{alpha4a}, sont commutatifs à isomorphismes canoniques près~;
ce qui justifie l'abus d'omettre $A$ dans les notations $\sigma_*$ et $\sigma_!$.

\begin{prop}\label{alpha9}
Soit $A$ une $R$-algèbre.
\begin{itemize}
\item[{\rm (i)}] Le foncteur $\sigma_*$ \eqref{alpha8b} est un adjoint à droite du foncteur de localisation $\alpha$ \eqref{alpha2a}.
\item[{\rm (ii)}] Le morphisme d'adjonction $\alpha\circ \sigma_* \rightarrow \id$ est un isomorphisme, 
{\em i.e.}, le foncteur $\sigma_*$ est pleinement fidèle.
\item[{\rm (iii)}] Le morphisme d'adjonction $\id\rightarrow \sigma_*\circ\alpha$
induit un isomorphisme $\alpha\stackrel{\sim}{\rightarrow} \alpha\circ\sigma_*\circ\alpha$.
\item[{\rm (iv)}] Le foncteur $\sigma_!$ \eqref{alpha8c} est un adjoint à gauche du foncteur de localisation $\alpha$. 
\item[{\rm (v)}]  Le morphisme d'adjonction $\id \rightarrow \alpha\circ \sigma_!$ 
est un isomorphisme, {\em i.e.}, le foncteur $\sigma_!$ est pleinement fidèle.
\end{itemize}
\end{prop}

Soient $M$, $N$ deux $A$-modules, $P$ un $\alpha$-$A$-module.

(i) On a des isomorphismes canoniques fonctoriels
\begin{eqnarray}
\Hom_{\aMod(A)}(M^\alpha,N^\alpha)&\stackrel{\sim}{\rightarrow}&\Hom_A(\fm\otimes_RM,N)\\
&\stackrel{\sim}{\rightarrow}&\Hom_A(M,\Hom_R(\fm,N))\nonumber\\
&\stackrel{\sim}{\rightarrow}&\Hom_A(M,\sigma_*(N^\alpha)).\nonumber\label{alpha9a}
\end{eqnarray}
L'assertion s'ensuit compte tenu de (\cite{gz} I 1.2).

(ii) Le morphisme d'adjonction
$\alpha(\sigma_*(M^\alpha))\rightarrow M^\alpha$ correspond par \eqref{alpha8a} et \eqref{alpha8e} au morphisme canonique
\begin{equation}\label{alpha9b}
\fm\otimes_R\Hom_R(\fm,M)\rightarrow M,
\end{equation} 
qui est un $\alpha$-isomorphisme car les morphismes canoniques $\fm\otimes_RM\rightarrow M$ et  
$M\rightarrow \Hom_R(\fm,M)$ sont des $\alpha$-isomorphismes \eqref{alpha3}.

(iii) Le morphisme d'adjonction $M\rightarrow \sigma_*(\alpha(M))$  
s'identifie par \eqref{alpha8e} au morphisme canonique $M\rightarrow \Hom_R(\fm,M)$
qui est un $\alpha$-isomorphisme. 

(iv) D'après (ii), on a des isomorphismes canoniques fonctoriels
\begin{eqnarray}
\Hom_{\aMod(A)}(P,\alpha(M))&\stackrel{\sim}{\rightarrow}&\Hom_A(\alpha(\sigma_*(P)),\alpha(M))\\
&\stackrel{\sim}{\rightarrow}&\Hom_A(\fm\otimes_R\sigma_*(P),M)\nonumber\\
&\stackrel{\sim}{\rightarrow}&\Hom_A(\sigma_!(P),M).\nonumber\label{alpha9c}
\end{eqnarray}

(v) Le morphisme d'adjonction $M^\alpha\rightarrow \alpha(\sigma_!(M^\alpha))$ correspond 
par \eqref{alpha8a} et \eqref{alpha8e} au morphisme canonique 
\begin{equation}\label{alpha9d}
\fm\otimes_RM\rightarrow \fm\otimes_R\Hom_R(\fm,M),
\end{equation}
qui est un isomorphisme \eqref{alpha7}.

\begin{cor}\label{alpha10}
Soit $A$ une $R$-algèbre. 
\begin{itemize}
\item[{\rm (i)}] Le foncteur de localisation $\alpha$ \eqref{alpha2a} 
commute aux limites inductives (resp. projectives) représentables et
le foncteur $\sigma_*$ \eqref{alpha8b} (resp. $\sigma_!$ \eqref{alpha8c}) commute aux limites projectives (resp. inductives) 
représentables.
\item[{\rm (ii)}] Les foncteurs $\alpha$ et $\sigma_!$ sont exacts, et le foncteur $\sigma_*$ est exact à gauche. 
\item[{\rm (iii)}] Pour toute catégorie $\mU$-petite $I$ et tout foncteur $\varphi\colon I\rightarrow \aMod(A)$, les
limites projective et inductive de $\varphi$ sont représentables et les morphismes canoniques 
\begin{eqnarray}
\underset{\underset{I}{\rightarrow}}{\lim}\ \varphi \rightarrow 
\alpha(\underset{\underset{I}{\rightarrow}}{\lim}\ \sigma_*\circ \varphi)\label{alpha10a}\\
\alpha(\underset{\underset{I}{\leftarrow}}{\lim}\ \sigma_!\circ \varphi)
\rightarrow \underset{\underset{I}{\leftarrow}}{\lim}\ \varphi\label{alpha10b}
\end{eqnarray}
sont des isomorphismes. 
\end{itemize}
\end{cor}

(i) Cela résulte de \ref{alpha9}(i)-(iv) et (\cite{sga4} I 2.11). 

(ii) En effet, il résulte de (i) que $\alpha$ est exact, $\sigma_*$ est exact à gauche et $\sigma_!$ est exact à droite. 
Comme $\fm$ est $R$-plat, $\sigma_!$ est aussi exact à gauche \eqref{alpha8c}. 

(iii) Cela résulte de (i), \ref{alpha9}(ii)-(v) et du fait que les $\mU$-limites inductives et projectives sont représentables dans 
la catégorie $\bMod(A)$.  

\begin{cor}\label{alpha100}
Pour toute $R$-algèbre $A$, la catégorie abélienne $\aMod(A)$ est une $\mU$-catégorie vérifiant la propriété {\rm (AB 5)} de 
{\rm (\cite{tohoku} §~1.5)} et admettant un générateur, à savoir $A^\alpha$. 
\end{cor}

En effet, comme $\alpha$ admet un adjoint à droite, à savoir $\sigma_*$ \eqref{alpha9}, la sous-catégorie épaisse de $\bMod(A)$ formée 
des $A$-modules $\alpha$-nuls est localisante dans le sens de (\cite{gabriel} p. 372). 
Par suite, d'après ({\em loc. cit.}, III §~2 lem.~4 et §~4 prop.~9), la catégorie $\aMod(A)$ est une $\mU$-catégorie avec générateur, 
à savoir $A^\alpha$, et limites inductives exactes~; d'où la proposition compte tenu de ({\em loc. cit.}, I §~6 prop.~6).

\subsection{}\label{alpha5}
Soit $A$ une $R$-algèbre. On appelle {\em $\alpha$-$A$-algèbre} (ou {\em $A^\alpha$-algèbre})
un monoïde unitaire commutatif de $\aMod(A)$. 
On désigne par $\bAlg(A)$ la catégorie des $A$-algèbres qui se trouvent dans $\mU$ et par $\aAlg(A)$ la catégorie des 
$\alpha$-$A$-algèbres. Le foncteur de localisation $\alpha$ \eqref{alpha2a}
étant monoïdal, il induit un foncteur que l'on note encore 
\begin{equation}\label{alpha5a}
\alpha\colon \bAlg(A)\rightarrow \aAlg(A).
\end{equation}

Compte tenu de l'isomorphisme canonique $A^\alpha\stackrel{\sim}{\rightarrow}A^\alpha\otimes_{A^\alpha}A^\alpha$,
le foncteur $\sigma_*$ \eqref{alpha8b} induit un foncteur que l'on note encore
\begin{equation}\label{alpha5b}
\sigma_*\colon \aAlg(A)\rightarrow \bAlg(A), \ \ \ P\mapsto \Hom_{\aMod(A)}(A^\alpha,P). 
\end{equation} 

Pour toute $A$-algèbre $B$, l'isomorphisme $\fm\stackrel{\sim}{\rightarrow} \fm\otimes_R\fm$ induit sur
$\Hom_R(\fm,B)$ une structure canonique de $A$-algèbre. On a un isomorphisme canonique fonctoriel de $A$-algèbres
\begin{equation}\label{alpha5c}
\sigma_*(B^\alpha)\stackrel{\sim}{\rightarrow}\Hom_{R}(\fm,B).
\end{equation}

\begin{prop}\label{alpha6}
Soit $A$ une $R$-algèbre.
\begin{itemize}
\item[{\rm (i)}] Le foncteur $\sigma_*$ \eqref{alpha5b} est un adjoint à droite du foncteur de localisation $\alpha$ \eqref{alpha5a}.
\item[{\rm (ii)}] Le morphisme d'adjonction $\alpha\circ \sigma_* \rightarrow \id$ est un isomorphisme, 
{\em i.e.}, le foncteur $\sigma_*$ est pleinement fidèle.
\item[{\rm (iii)}] Le morphisme d'adjonction $\id\rightarrow \sigma_*\circ\alpha$
induit un isomorphisme $\alpha\stackrel{\sim}{\rightarrow} \alpha\circ\sigma_*\circ\alpha$.
\end{itemize}
\end{prop}

Soient $B$, $C$ deux $A$-algèbres, $u\colon B^\alpha\rightarrow C^\alpha$ un morphisme de $\aMod(A)$,
$v\colon \fm\otimes_RB\rightarrow C$ et $w\colon B\rightarrow \Hom_R(\fm,C)$
les morphismes $A$-linéaires associés \eqref{alpha8a}. On vérifie aussitôt que les conditions suivantes sont équivalentes~:
\begin{itemize}
\item[(a)] $u$ est un morphisme de $A^\alpha$-algèbres.
\item[(b)]  Le diagramme 
\begin{equation}\label{alpha6b}
\xymatrix{
{(\fm\otimes_RB)\otimes_A(\fm\otimes_RB)}\ar[r]^-(0.5){v\otimes v}\ar@{=}[d]&{C\otimes_AC}\ar[r]^-(0.5){\mu_C}&C\\
{(\fm\otimes_R\fm)\otimes_R(B\otimes_AB)}\ar[rr]^-(0.5){\mu_\fm\otimes_R\mu_B}&&{\fm\otimes_RB}\ar[u]_v}
\end{equation}
où $\mu_\fm$, $\mu_B$ et $\mu_C$ désignent les morphismes de multiplication de $\fm$, $B$ et $C$, 
respectivement, est commutatif. 
\item[(c)] $w$ est  un morphisme de $A$-algèbres.
\end{itemize}
On en déduit des isomorphismes canoniques fonctoriels
\begin{eqnarray}
\Hom_{\aAlg(A)}(B^\alpha,C^\alpha)&\stackrel{\sim}{\rightarrow}&\Hom_{\bAlg(A)}(B,\Hom_R(\fm,C))\nonumber\\
&\stackrel{\sim}{\rightarrow}&\Hom_{\bAlg(A)}(B,\sigma_*(C^\alpha)).\label{alpha6a}
\end{eqnarray}
La proposition (i) s'ensuit compte tenu de (\cite{gz} I 1.2). On notera que les morphismes d'adjonction 
$\alpha\circ \sigma_* \rightarrow \id$ et $\id\rightarrow \sigma_*\circ\alpha$ s'identifient aux morphismes d'adjonction 
pour les $A$-modules \eqref{alpha9}. Les propositions (ii) et (iii) résultent alors de \ref{alpha9}(ii)-(iii).

\subsection{}\label{alpha39}
Dans ce numéro, si $D$ est une $R$-algèbre, nous affecterons d'un indice $D$ les foncteurs 
$\alpha$ \eqref{alpha2a} et $\sigma_*$ \eqref{alpha8b} pour les $D$-modules 
ainsi que leurs variantes \eqref{alpha5a} et \eqref{alpha5b} pour les $D$-algèbres.  

Soit $A$ une $R$-algèbre. On pose $B=\alpha_R(A)$ \eqref{alpha5a}  
et on désigne par $\bMod(B)$ la catégorie des $B$-modules unitaires de $\aMod(R)$. Le foncteur $\alpha_R$ étant monoïdal, 
il induit un foncteur 
\begin{equation}\label{alpha39a}
\beta\colon \bMod(A)\rightarrow \bMod(B).
\end{equation} 
Celui-ci transforme clairement les $\alpha$-isomorphismes en des isomorphismes. Il induit donc un foncteur 
\begin{equation}\label{alpha39b}
b\colon \aMod(A)\rightarrow \bMod(B).
\end{equation} 

On pose $A'=\sigma_{R*}(B)$ \eqref{alpha5b}. 
D'après \ref{alpha6}, on a un homomorphisme canonique de $R$-algèbres $\lambda\colon A\rightarrow A'$,
qui induit un isomorphisme $\alpha_R(A)\stackrel{\sim}{\rightarrow}\alpha_R(A')$ \eqref{alpha5a}. 
Pour tous $R^\alpha$-modules $P$ et $Q$, on a morphisme $R$-linéaire
\begin{equation}\label{alpha39c}
\Hom_{\aMod(R)}(R^\alpha,P)\otimes_R \Hom_{\aMod(R)}(R^\alpha,Q)\rightarrow \Hom_{\aMod(R)}(R^\alpha,P\otimes_{R^\alpha}Q),
\end{equation}
défini par fonctorialité et composition \eqref{alpha33b}.
Par suite, le foncteur $\sigma_{R*}$ \eqref{alpha8b} induit un foncteur 
\begin{equation}\label{alpha39d}
\tau'_*\colon \bMod(B)\rightarrow \bMod(A').
\end{equation}
Composant avec le foncteur induit par $\lambda$, on obtient un foncteur
\begin{equation}\label{alpha39e}
\tau_*\colon \bMod(B)\rightarrow \bMod(A)
\end{equation}

Soient $M$ et $N$ deux $A$-modules, $u\colon \beta(M)\rightarrow \beta(N)$ un morphisme de $\aMod(R)$,
$v\colon \fm\otimes_RM\rightarrow N$ et $w\colon M\rightarrow \Hom_R(\fm,N)$
les morphismes $R$-linéaires associés à $u$ \eqref{alpha8a}. On vérifie aussitôt que $u$ est $B$-linéaire si et seulement si $v$
(ou ce qui revient au même $w$) est $A$-linéaire.  
L'isomorphisme \eqref{alpha8a} induit donc un isomorphisme canonique 
\begin{equation}
\Hom_{\bMod(B)}(\beta(M),\beta(N))\stackrel{\sim}{\rightarrow} \Hom_{\bMod(A)}(\fm\otimes_RM,N).
\end{equation}
Calquant alors la preuve de \ref{alpha9}, on en déduit que $\tau_*$ est un adjoint à droite de $\beta$; 
le morphisme d'adjonction $\beta\circ \tau_* \rightarrow \id$ est un isomorphisme~; 
et le morphisme d'adjonction $\id\rightarrow \tau_*\circ\beta$
induit un isomorphisme $\beta\stackrel{\sim}{\rightarrow} \beta\circ\tau_*\circ\beta$. 

On pose
\begin{equation}\label{alpha39f}
t=\alpha_A\circ \tau_*\colon \bMod(B)\rightarrow \aMod(A).
\end{equation}
L'isomorphisme $\beta\circ \tau_* \stackrel{\sim}{\rightarrow} \id$ induit un isomorphisme 
\begin{equation}\label{alpha39g}
b\circ t \stackrel{\sim}{\rightarrow} \id
\end{equation}
Par ailleurs, le morphisme d'adjonction $\id\rightarrow \tau_*\circ\beta$
induit un isomorphisme $\alpha_A\stackrel{\sim}{\rightarrow} \alpha_A\circ\tau_*\circ\beta$. 
Compte tenu de (\cite{gz} I 1.2), on en déduit un isomorphisme 
\begin{equation}\label{alpha39h}
\id\stackrel{\sim}{\rightarrow} t\circ b.
\end{equation} 
On vérifie aussitôt que les isomorphismes \eqref{alpha39g} et \eqref{alpha39h} font de $t$ un adjoint à droite de $b$. 
Par suite, $b$ et $t$ sont des équivalences de catégories. 
On a des isomorphismes canoniques
\begin{equation}
\beta\stackrel{\sim}{\rightarrow} b\circ \alpha_A \ \ \ {\rm et} \ \ \ \alpha_A\stackrel{\sim}{\rightarrow} t \circ \beta,
\end{equation}
compatibles aux isomorphismes \eqref{alpha39g} et \eqref{alpha39h}. 
On en déduit aussitôt un isomorphisme
\begin{equation}
\sigma_A \stackrel{\sim}{\rightarrow} \tau_*\circ b.
\end{equation}

\subsection{}\label{alpha13}
Soit $\mV$ un univers tel que $\mU\subset \mV$. 
On affectera d'un indice $\mU$ ou $\mV$ les catégories et foncteurs dépendants de l'univers. 
Soit $A$ une $R$-algèbre appartenant à $\mU$. 
On a un foncteur pleinement fidèle canonique
\begin{equation}\label{alpha13a}
\bMod_\mU(A)\rightarrow \bMod_\mV(A).
\end{equation}
Celui-ci induit un foncteur
\begin{equation}\label{alpha13b}
\aMod_\mU(A)\rightarrow \aMod_\mV(A)
\end{equation}
qui s'insère dans un diagramme commutatif 
\begin{equation}\label{alpha13c}
\xymatrix{
{\bMod_\mU(A)}\ar[r]\ar[d]_{\alpha_\mU}&{\bMod_\mV(A)}\ar[d]^{\alpha_\mV}\\
{\aMod_\mU(A)}\ar[r]&{\aMod_\mV(A)}}
\end{equation}
Le foncteur \eqref{alpha13b} est pleinement fidèle. En effet, comme le foncteur \eqref{alpha13a} est pleinement fidèle,
il suffit de montrer que pour tout objet $M$ de $\bMod_\mU(A)$, notant $\aI_\mU(M)$ (resp. $\aI_\mV(M)$)
la catégorie des $\alpha$-isomorphismes de $\bMod_\mU(A)$ (resp. $\bMod_\mV(A)$) de but $M$, 
le foncteur d'inclusion $\phi_M\colon \aI^\circ_\mU(M)\rightarrow \aI^\circ_\mV(M)$ est cofinal. 
Soit $f\colon N\rightarrow M$ un objet de $\aI_\mV(M)$.  D'après \ref{alpha3},
pour tout $\varepsilon \in \Lambda^+$, il existe un morphisme $A$-linéaire $g_\varepsilon\colon M\rightarrow N$ 
tel que $f\circ g_\varepsilon=\pi^\varepsilon \id_M$ et $g_\varepsilon\circ f=\pi^\varepsilon \id_N$. 
Comme l'ensemble $\Lambda^+$ est $\mU$-petit, 
$\cup_{\varepsilon\in \Lambda^+}\im(g_\varepsilon)$
est représentable par un sous-objet $N'$ de $N$ appartenant à $\bMod_\mU(A)$. 
Comme l'injection canonique $g\colon N'\rightarrow N$ est clairement un $\alpha$-isomorphisme, 
$f\circ g\colon N'\rightarrow M$ est un objet de $\aI_\mU(M)$, d'où l'assertion recherchée en vertu de (\cite{sga4} I 8.1.3(c)).

\subsection{}\label{alpha11}
Dans la suite de cette section, $\cC$ désigne une $\mU$-catégorie \eqref{notconv3} (\cite{sga4} I 1.1) et
$\hcC$ la catégorie des préfaisceaux de $\mU$-ensembles sur $\cC$ (\cite{sga4} I 1.2).
On note $R_\hcC$ le préfaisceau d'anneaux constant sur $\cC$ de valeur $R$.
On désigne par $\bMod(R_\hcC)$ la catégorie abélienne tensorielle des  $\mU$-préfaisceaux de $R$-modules sur $\cC$,
que l'on identifie à la catégorie des $(R_\hcC)$-modules de $\hcC$ (\cite{sga4} I 3.2). 
Prenant pour $\phi\colon R\rightarrow \End(R_\hcC)=R$ l'homomorphisme identique, 
on appelle catégorie des {\em $\alpha$-$(R_\hcC)$-modules} et l'on note $\aMod(R_\hcC)$ 
le quotient de la catégorie abélienne $\bMod(R_\hcC)$ par la sous-catégorie épaisse des $(R_\hcC)$-modules 
$\alpha$-nuls. Nous utiliserons les conventions de notation de \ref{alpha2} et \ref{alpha33}.

\begin{defi}\label{alpha14}
On appelle catégorie des {\em préfaisceaux de $\alpha$-$R$-modules} sur $\cC$ et l'on note $\ahcC$ 
la catégorie des préfaisceaux sur $\cC$ à valeurs dans la catégorie $\aMod(R)$, {\em i.e.}, 
la catégorie des foncteurs de $\cC^\circ$ à valeurs dans $\aMod(R)$ (\cite{sga4} II 6.0).
\end{defi}

On note $R^\alpha_\hcC$ le préfaisceau constant sur $\cC$ de valeur $R^\alpha$. 
Pour deux préfaisceaux de $\alpha$-$R$-modules $M$ et $N$, 
on définit le produit tensoriel $M\otimes_{R^\alpha_\hcC}N$  par la formule
suivante~: pour tout $X\in \ob(\cC)$, 
\begin{equation}\label{alpha14a}
(M\otimes_{R^\alpha_\hcC}N)(X)=M(X)\otimes_{R^\alpha}N(X), 
\end{equation}
où le terme de droite désigne le produit tensoriel dans la catégorie $\aMod(R)$ \eqref{alpha33a}. 
Munie de ce produit, $\ahcC$ est une catégorie abélienne tensorielle, ayant $R^\alpha_{\hcC}$ pour objet unité.

\begin{prop}\label{alpha15}
Les $\mU$-limites inductives et projectives dans $\ahcC$ sont représentables. Pour tout $X\in\ob(\cC)$, le foncteur 
\begin{equation}\label{alpha15a}
\ahcC\rightarrow \aMod(R),\ \ \ M\mapsto M(X)
\end{equation}
commute aux limites inductives et projectives.
\end{prop}
Cela résulte immédiatement de \ref{alpha10}. 

\begin{cor}\label{alpha150}
La catégorie $\ahcC$ est une catégorie abélienne vérifiant l'axiome {\rm (AB 5)} de {\rm (\cite{tohoku} §~1.5)}. 
\end{cor}

C'est une conséquence de \ref{alpha100} et \ref{alpha15}.

\subsection{}\label{alpha16}
Le foncteur $\alpha\colon \bMod(R)\rightarrow \aMod(R)$ définit un foncteur exact et monoïdal 
\begin{equation}\label{alpha16a}
\halpha\colon \bMod(R_\hcC)\rightarrow \ahcC. 
\end{equation}
Celui-ci transforme les modules $\alpha$-nuls en dans le préfaisceau de $\alpha$-$R$-modules nul. 
Il induit donc un foncteur exact et monoïdal
\begin{equation}\label{alpha16b}
u\colon \aMod(R_\hcC)\rightarrow \ahcC. 
\end{equation}

Les foncteurs $\sigma_*$ \eqref{alpha8b} et $\sigma_!$ \eqref{alpha8c} (pour les $R$-modules) induisent des foncteurs 
que l'on note respectivement
\begin{eqnarray}
\hsigma_*\colon \ahcC&\rightarrow& \bMod(R_\hcC),\label{alpha16c}\\
\hsigma_!\colon \ahcC&\rightarrow& \bMod(R_\hcC).\label{alpha16d}
\end{eqnarray}
D'après \ref{alpha9}, $\hsigma_*$ (resp. $\hsigma_!$) est un adjoint à droite (resp. à gauche) de $\halpha$; 
les morphismes d'adjonction $\halpha\circ \hsigma_* \rightarrow \id$ et $\id\rightarrow \halpha \circ \hsigma_!$ 
sont des isomorphismes~; et le morphisme d'adjonction $\id\rightarrow \hsigma_*\circ\halpha$
induit un isomorphisme $\halpha\stackrel{\sim}{\rightarrow} \halpha\circ\hsigma_*\circ\halpha$. 

On note $v$ le foncteur composé
\begin{equation}\label{alpha16e}
v=\alpha\circ \hsigma_*\colon  \ahcC\rightarrow \aMod(R_\hcC),
\end{equation}
où $\alpha$ désigne le foncteur de localisation pour les $R$-modules de $\hcC$. 
L'isomorphisme $\halpha\circ \hsigma_*\stackrel{\sim}{\rightarrow} \id$ induit un isomorphisme 
\begin{equation}\label{alpha16f}
u\circ v\stackrel{\sim}{\rightarrow} \id.
\end{equation}
D'après \ref{alpha9}(iii), le morphisme d'adjonction $\id\rightarrow \hsigma_*\circ \halpha$ induit un isomorphisme 
$\alpha\stackrel{\sim}{\rightarrow} v\circ u\circ \alpha$. Compte tenu de (\cite{gz} I 1.2), on en déduit un isomorphisme 
\begin{equation}\label{alpha16g}
\id\stackrel{\sim}{\rightarrow} v\circ u.
\end{equation} 
On vérifie aussitôt que les isomorphismes \eqref{alpha16f} et \eqref{alpha16g} font de $v$ un adjoint à droite de $u$. 
Par suite, $u$ et $v$ sont des équivalences de catégories abéliennes tensorielles, quasi-inverses l'une de l'autre. 
On a des isomorphismes canoniques
\begin{equation}\label{alpha16h}
\halpha\stackrel{\sim}{\rightarrow} u\circ \alpha \ \ \ {\rm et} \ \ \ \alpha\stackrel{\sim}{\rightarrow} v\circ \halpha,
\end{equation}
compatibles aux isomorphismes \eqref{alpha16f} et \eqref{alpha16g}. 
On en déduit aussitôt que le foncteur
\begin{equation}\label{alpha16i}
\hsigma_*\circ u \colon \aMod(R_\hcC)\rightarrow \bMod(R_\hcC)
\end{equation}
est un adjoint à droite du foncteur de localisation $\alpha$, et que le foncteur
\begin{equation}\label{alpha16j}
\hsigma_!\circ u \colon \aMod(R_\hcC)\rightarrow \bMod(R_\hcC)
\end{equation}
est un adjoint à gauche de $\alpha$. 

\subsection{}\label{alpha36}
La donnée d'un monoïde commutatif unitaire de $\bMod(R_\hcC)$ (resp. $\ahcC$) est équivalente à la donnée 
d'un préfaisceau sur $\cC$ à valeurs dans la catégorie $\bAlg(R)$ (resp. $\aAlg(R)$) \eqref{alpha5}. 
On note $\bAlg(R_\hcC)$ (resp. $\bAlg(\ahcC)$) la catégorie des monoïdes commutatifs unitaires de $\bMod(R_\hcC)$ 
(resp. $\ahcC$). Le foncteur $\halpha$ \eqref{alpha16a} étant monoïdal, il induit un foncteur que l'on note encore
\begin{equation}\label{alpha36a}
\halpha\colon \bAlg(R_\hcC)\rightarrow \bAlg(\ahcC).
\end{equation}
Par ailleurs, compte tenu de \ref{alpha5}, le foncteur $\hsigma_*$ \eqref{alpha16c} induit un foncteur que l'on note encore
\begin{equation}\label{alpha36b}
\hsigma_*\colon \bAlg(\ahcC)\rightarrow \bAlg(R_\hcC).
\end{equation}
D'après \ref{alpha6}, $\hsigma_*$ est un adjoint à droite de $\halpha$; 
le morphisme d'adjonction $\halpha\circ \hsigma_* \rightarrow \id$  
est un isomorphisme~; et le morphisme d'adjonction $\id\rightarrow \hsigma_*\circ\halpha$
induit un isomorphisme $\halpha\stackrel{\sim}{\rightarrow} \halpha\circ\hsigma_*\circ\halpha$.

\subsection{}\label{alpha17}
Lorsque la catégorie $\cC$ est $\mU$-petite, $\hcC$ est une $\mU$-catégorie.  
Mais lorsque $\cC$ est une $\mU$-catégorie, $\hcC_\mU$ n'est pas en général une $\mU$-catégorie (\cite{sga4} I 1.2). 
Soit $\mV$ un univers tel que $\cC\in \mV$ et $\mU\subset \mV$. 
On affectera d'un indice $\mU$ ou $\mV$ les catégories et foncteurs dépendants de l'univers. 
Le foncteur pleinement fidèle canonique $\bMod_\mU(R)\rightarrow \bMod_\mV(R)$ induit un foncteur pleinement fidèle 
\begin{equation}\label{alpha17a}
\bMod(R_{\hcC_\mU})\rightarrow \bMod(R_{\hcC_\mV}).
\end{equation}
De même, le foncteur pleinement fidèle canonique $\aMod_\mU(R)\rightarrow \aMod_\mV(R)$ \eqref{alpha13b} induit un foncteur pleinement fidèle 
\begin{equation}\label{alpha17b}
\ahcC_\mU\rightarrow \ahcC_\mV.
\end{equation}

Tout $(R_{\hcC_\mU})$-module $M$ définit un foncteur
\begin{equation}\label{alpha17c}
M\colon (\hcC_\mU)^\circ\rightarrow \bMod_\mV(R), \ \ \ X\mapsto M(X)=\Hom_{\hcC_\mU}(X,M).
\end{equation}
Pour tout objet $P$ de $\ahcC_\mU$, on considère le foncteur 
\begin{equation}\label{alpha17d}
\hP\colon (\hcC_\mU)^\circ\rightarrow \aMod_\mV(R), \ \ \ X\mapsto \alpha(\hsigma_*(P)(X)),
\end{equation}
où $\hsigma_*$ est le foncteur \eqref{alpha16c}. On obtient ainsi un foncteur
\begin{equation}\label{alpha17e}
\ahcC\rightarrow \bHom((\hcC_\mU)^\circ,\aMod_\mV(R)), \ \ \ P\mapsto \hP.
\end{equation}
D'après \ref{alpha9}(ii), on a un isomorphisme canonique fonctoriel $\hP|\cC^\circ\stackrel{\sim}{\rightarrow} P$. 
Il est alors commode et sans risque d'ambiguïté de noter $\hP$ encore $P$.

D'après (\cite{sga4} I 3.5), pour tout $(R_{\hcC_\mU})$-module $M$ et tout $X\in \ob(\hcC_\mU)$, on a un isomorphisme canonique et fonctoriel en $M$,
\begin{equation}\label{alpha17f}
M(X)\stackrel{\sim}{\rightarrow}\underset{\underset{(Y,u)\in (\cC_{/X})^\circ}{\longleftarrow}}{\lim}\ M(Y).
\end{equation}
Par suite, en vertu de \ref{alpha10}(i), pour tout objet $P$ de $\ahcC_\mU$ et tout $X\in \ob(\hcC_\mU)$, on a un isomorphisme canonique fonctoriel en $P$,
\begin{equation}\label{alpha17g}
P(X)\stackrel{\sim}{\rightarrow}\underset{\underset{(Y,u)\in (\cC_{/X})^\circ}{\longleftarrow}}{\lim}\ P(Y).
\end{equation}

\subsection{}\label{alpha20}
Dans la suite de cette section, on se donne une topologie sur $\cC$ qui en fait un $\mU$-site, {\em i.e.}, tel que $\cC$ admette une
$\mU$-petite famille topologiquement génératrice (\cite{sga4} II 3.0.2). 
On désigne par $\tcC$ le topos des faisceaux de $\mU$-ensembles sur $\cC$.
On note $R_\tcC$ le faisceau d'anneaux constant de valeur $R$ sur $\cC$ et $\fm_\tcC$ le faisceau d'idéaux constant de valeur 
$\fm$ sur $\cC$,  {\em i.e.}, les faisceaux associés au préfaisceaux constants sur $\cC$ de valeurs $R$ et $\fm$, respectivement 
(\cite{sga4} II 6.4). 

Soit $A$ une $(R_\tcC)$-algèbre de $\tcC$, {\em i.e.}, un $\mU$-faisceau de $R$-algèbres sur $\cC$ (\cite{sga4} II 6.3.1). 
On désigne par $\bMod(A)$ la catégorie abélienne tensorielle des $A$-modules de $\tcC$. 
Prenant pour $\phi\colon R\rightarrow \Gamma(\tcC,A)=\End(A)$ l'homomorphisme canonique, 
on appelle catégorie des {\em $\alpha$-$A$-modules} et l'on note $\aMod(A)$ 
le quotient de la catégorie abélienne $\bMod(A)$ par la sous-catégorie épaisse des $A$-modules 
$\alpha$-nuls. Nous utiliserons les conventions de notation de \ref{alpha2} et \ref{alpha33}. 
On observera en particulier que pour tous $A$-modules $M$ et $N$, $\Hom_{\aMod(A)}(M,N)$ est 
naturellement muni d'une structure de $\Gamma(\tcC,A)$-module.

\begin{lem}\label{alpha21}
Soit $A$ une $(R_\tcC)$-algèbre de $\tcC$. 
\begin{itemize}
\item[{\rm (i)}] Pour tout $A$-module $M$, les conditions suivantes sont équivalentes~:
\begin{itemize}
\item[{\rm (a)}] $M$ est $\alpha$-nul en tant qu'objet de $\bMod(A)$;
\item[{\rm (b)}] $M$ est $\alpha$-nul en tant qu'objet de $\bMod(R_\tcC)$;
\item[{\rm (c)}] $M$ est $\alpha$-nul en tant qu'objet de $\bMod(R_\hcC)$; 
\item[{\rm (d)}] pour tout $U\in \ob(\cC)$, $\fm M(U)=0$; 
\item[{\rm (e)}] $\fm_\tcC M=0$.
\end{itemize}
\item[{\rm (ii)}] Pour tout $(R_\hcC)$-module $\alpha$-nul $P$ de $\hcC$, le $(R_\tcC)$-module associé $P^\tta$ 
de $\tcC$ est $\alpha$-nul. 
\item[{\rm (iii)}]  Pour tout morphisme de $A$-modules $f\colon M\rightarrow N$, les conditions suivantes sont équivalentes~:
\begin{itemize}
\item[$(1)$] $f$ est un $\alpha$-isomorphisme en tant que morphisme de $\bMod(A)$~; 
\item[$(2)$] $f$ est un $\alpha$-isomorphisme en tant que morphisme de $\bMod(R_\tcC)$~; 
\item[$(3)$] $f$ est un $\alpha$-isomorphisme en tant que morphisme de $\bMod(R_\hcC)$~;  
\item[$(4)$] le morphisme $\fm_\tcC\otimes_{R_\tcC}M\rightarrow \fm_\tcC\otimes_{R_\tcC}N$ induit par $f$ est un isomorphisme.
\end{itemize} 
\item[{\rm (iv)}] Pour tous $A$-modules $M$ et $N$, on a un isomorphisme $\Gamma(\tcC,A)$-linéaire canonique
\begin{equation}\label{alpha21a}
\Hom_{\aMod(A)}(M^\alpha,N^\alpha)\stackrel{\sim}{\rightarrow}\Hom_{\bMod(A)}(\fm_\tcC\otimes_{R_\tcC}M,N).
\end{equation}
\end{itemize}
\end{lem}

(i) En effet, les conditions (a), (b), (c) et (d) sont clairement équivalentes. Comme $\fm_\tcC M$ est le faisceau associé au préfaisceau 
$U\mapsto \fm M(U)$, (d) implique (e). Il est clair que (e) implique (a). 

(ii) Cela résulte aussitôt de (i)

(iii) Les conditions (1) et (2) sont équivalentes en vertu de (i). 
L'implication $(2)\Rightarrow(3)$ résulte de \ref{alpha3} et l'implication $(3)\Rightarrow(2)$ est une conséquence de (ii).
L'implication $(4)\Rightarrow(2)$ est une conséquence du fait que le morphisme canonique $\fm_\tcC\otimes_{R_\tcC} M\rightarrow M$
est un $\alpha$-isomorphisme \eqref{alpha3}. Par ailleurs, le $(R_\tcC)$-module $\fm_\tcC$ est plat (\cite{sga4} V 1.7.1)
et pour tout $(R_\tcC)$-module $\alpha$-nul $F$, on a $\fm_\tcC\otimes_{R_\tcC}F=0$; d'où l'implication $(2)\Rightarrow(4)$.

(iv) D'après (iii), le morphisme canonique $\fm_\tcC\otimes_{R_\tcC}M\rightarrow M$ est un objet initial de
la catégorie des $\alpha$-isomorphismes de $\bMod(A)$ de but $M$; d'où  l'isomorphisme \eqref{alpha21a}.
On vérifie aussitôt qu'il est $\Gamma(\tcC,A)$-linéaire.

\begin{defi}\label{alpha18}
On dit qu'un préfaisceau de $\alpha$-$R$-modules $F$ sur $\cC$ est {\em séparé} 
(resp. est un {\em faisceau}) 
si  pour tout objet $X$ de $\cC$ et tout crible couvrant $\cR$ de $X$, le morphisme canonique 
\begin{equation}\label{alpha18a}
F(X)\rightarrow F(\cR)
\end{equation}
est un monomorphisme (resp. isomorphisme) (cf. \ref{alpha17}). 
On note $\atcC$ la sous-catégorie pleine de $\ahcC$ formée des faisceaux de $\alpha$-$R$-modules.
\end{defi}

D'après \eqref{alpha17g}, pour qu'un préfaisceau de $\alpha$-$R$-modules 
$F$ sur $\cC$ soit séparé (resp. un faisceau), 
il faut et il suffit que pour tout objet $X$ de $\cC$ et tout crible couvrant $\cR$ de $X$, le morphisme canonique
\begin{equation}\label{alpha18b}
F(X)\rightarrow \underset{\underset{(Y,u)\in (\cC_{/\cR})^\circ}{\longleftarrow}}{\lim}\ F(Y)
\end{equation}
soit un monomorphisme (resp. isomorphisme). 
Les faisceaux de $\alpha$-$R$-modules sont donc les faisceaux sur 
$\cC$ à valeurs dans la catégorie $\aMod(R)$ dans le sens de (\cite{sga4} II 6.1).

\begin{prop}\label{alpha182}
Soient $F$ un préfaisceau de $\alpha$-$R$-modules sur $\cC$,  $X\in \ob(\cC)$, $\cR$ un crible de $X$. 
Pour que le morphisme canonique $F(X)\rightarrow F(\cR)$
soit un monomorphisme (resp. isomorphisme), il faut et il suffit qu'il en soit de même 
du morphisme canonique \eqref{alpha16c}
\begin{equation}\label{alpha182a}
\hsigma_*(F)(X)\rightarrow \hsigma_*(F)(\cR).
\end{equation}
\end{prop} 

En effet, d'après \ref{alpha10}(i) et \eqref{alpha17g}, on a un isomorphisme canonique 
\begin{equation}\label{alpha182b}
\sigma_*(F(\cR))\stackrel{\sim}{\rightarrow} \hsigma_*(F)(\cR).
\end{equation}
L'image du morphisme $u\colon F(X)\rightarrow F(\cR)$ par le foncteur $\sigma_*$ s'identifie alors au morphisme \eqref{alpha182a}.
Si $u$ est un monomorphisme (resp. isomorphisme), il en est de même de $\sigma_*(u)$ en vertu de \ref{alpha10}(i).
Inversement, si $\sigma_*(u)$ est un monomorphisme (resp. isomorphisme), il en est de même de $u$ 
d'après \ref{alpha9}(ii) et \ref{alpha10}(i).

\begin{cor}\label{alpha183}
Pour qu'un préfaisceau de $\alpha$-$R$-modules $F$ sur $\cC$ soit séparé (resp. un faisceau), il faut et il suffit 
que le préfaisceau de $R$-modules $\hsigma_*(F)$ sur $\cC$ \eqref{alpha16c} soit séparé (resp. un faisceau).
\end{cor}

\begin{prop}\label{alpha181}
Pour tout $X\in \ob(\cC)$, soit $\cK(X)$ un ensemble de cribles de $X$. Supposons que les $\cK(X)$ soient stables par
changement de base et qu'ils engendrent la topologie de $\cC$. Alors, pour qu'un préfaisceau de $\alpha$-$R$-modules
$F$ sur $\cC$ soit séparé (resp. un faisceau), il faut et il suffit que pour tout $X\in \ob(\cC)$ et tout $\cR \in \cK(X)$, 
le morphisme canonique
\begin{equation}\label{alpha181a}
F(X)\rightarrow F(\cR)
\end{equation}
soit un monomorphisme (resp. isomorphisme). 
\end{prop} 
Cela résulte \ref{alpha182} et (\cite{sga4} II 2.3).

\begin{cor}\label{alpha19}
Si la topologie de $\cC$ est définie par une prétopologie, pour qu'un 
préfaisceau de $\alpha$-$R$-modules $F$ sur $\cC$ soit un faisceau, il faut et il suffit que pour tout objet $X$ de $\cC$
et tout recouvrement $(X_i\rightarrow X)_{i\in I}$, la suite de $\alpha$-$R$-modules
\begin{equation}\label{alpha19a}
0\rightarrow F(X)\rightarrow \prod_{i\in I}F(X_i)\rightarrow \prod_{(i,j)\in I^2}F(X_i\times_XX_j),
\end{equation}
où la dernière flèche est la différence des morphismes induits par les projections de $X_i\times_XX_j$ sur les deux facteurs, 
soit exacte. 
\end{cor}

Cela résulte de \ref{alpha181}, \eqref{alpha17g} et (\cite{sga4} I 2.12).

\subsection{}\label{alpha24}
Rappelons la définition du foncteur $L$ sur la catégorie $\bMod(R_\hcC)$ (\cite{sga4} II 3.0.5). 
Soient $\mV$ un univers tel que $\cC\in \mV$ et $\mU\subset \mV$, $G$ une $\mU$-petite famille topologiquement génératrice de $\cC$ 
(\cite{sga4} II 3.0.1). Pour tout objet $X$ de $\cC$, on désigne par $J(X)$ l'ensemble des cribles couvrants de $X$ et 
par $J_G(X)$ l'ensemble des cribles couvrants de $X$ engendrés par une famille $(X_i\rightarrow X)_{i\in I}$
telle que $X_i\in G$ pour tout $i\in I$. L'ensemble $J(X)$ est $\mV$-petit, et ordonné par l'inclusion, il est cofiltrant. 
Pour tout $\mU$-préfaisceau de $R$-modules $F$ sur $\cC$, 
\begin{equation}\label{alpha24a}
\underset{\underset{\cR\in J(X)^\circ}{\longrightarrow}}{\lim}\ F(\cR)
\end{equation}
est représentable par un $R$-module de $\mV$ \eqref{alpha17}. D'après (\cite{sga4} II 3.0.4), 
pour tout $\cR\in J_G(X)$, $F(\cR)$ est $\mU$-petit, et comme $J_G(X)$ est un $\mU$-petit ensemble cofinal dans $J(X)$ ({\em loc. cit.}), 
il résulte de (\cite{sga4} I 2.3.3) que la limite inductive \eqref{alpha24a} est représentable par un $R$-module $\mU$-petit. 
Choisissons pour tout $F$ et tout $X$ un $R$-module appartenant à $\mU$ qui représente cette limite inductive et posons 
\begin{equation}\label{alpha24b}
LF(X)=\underset{\underset{\cR\in J(X)^\circ}{\longrightarrow}}{\lim}\ F(\cR).
\end{equation}
Pour tout morphisme $f\colon Y\rightarrow X$ de $\cC$, le foncteur de changement de base $f^*\colon J(X)\rightarrow J(Y)$ définit un morphisme
$LF(f)\colon LF(X)\rightarrow LF(Y)$ faisant de $X\mapsto LF(X)$ un $\mU$-préfaisceau de $R$-modules sur $\cC$. 
Pour tout $X\in \ob(\cC)$, le morphisme identique de $X$ étant un objet de $J(X)$, 
on a une application canonique $\ell(F)(X)\colon F(X)\rightarrow LF(X)$. 
On définit ainsi un morphisme $\ell(F)\colon F\rightarrow LF$ de $\bMod(R_\hcC)$.
La correspondance $F\mapsto LF$ est clairement fonctorielle en $F$ et les $\ell(F)$ définissent un morphisme de foncteurs 
$\ell\colon \id\rightarrow L$.

Avec les notations de \ref{alpha16}, on désigne par $\cL$ le foncteur composé 
\begin{equation}\label{alpha24c}
\cL=\halpha\circ L\circ \hsigma_*\colon \ahcC\rightarrow \ahcC.
\end{equation}
Le morphisme $\ell$ et l'isomorphisme $\halpha\circ \hsigma_*\stackrel{\sim}{\rightarrow} \id$ induisent un morphisme de foncteur 
$\lambda\colon \id\rightarrow \cL$. 
D'après \ref{alpha10}(i) et compte tenu de la définition \eqref{alpha17d}, 
pour tout objet $P$ de $\ahcC$ et tout $X\in \ob(\cC)$, on a un isomorphisme canonique fonctoriel 
\begin{equation}\label{alpha24d}
\cL P(X)\stackrel{\sim}{\rightarrow}\underset{\underset{\cR\in J(X)^\circ}{\longrightarrow}}{\lim}\ P(\cR).
\end{equation}
Le composé de cet isomorphisme et du morphisme $\lambda(P)(X)\colon P(X)\rightarrow \cL P(X)$ n'est autre que 
le morphisme induit par l'objet $\id_X$ de $J(X)$.

\begin{prop}\label{alpha26}
{\rm (i)}\  Le foncteur $\cL$ est exact à gauche. 

{\rm (ii)}\ Le foncteur $L$ transforme les $\alpha$-isomorphismes en des $\alpha$-isomorphismes. 

{\rm (iii)}\ Le diagramme 
\begin{equation}\label{alpha26a}
\xymatrix{
{\bMod(R_\hcC)}\ar[r]^L\ar[d]_\halpha&{\bMod(R_\hcC)}\ar[d]^\halpha\\
{\ahcC}\ar[r]^\cL&{\ahcC}}
\end{equation}
est commutatif à un isomorphisme canonique près
\begin{equation}\label{alpha26b}
\halpha\circ L\stackrel{\sim}{\rightarrow} \cL\circ \halpha,
\end{equation} 
induit par le morphisme d'adjonction $\id\rightarrow \hsigma_*\circ \halpha$.

{\rm (iv)}\ Pour tout préfaisceau de $\alpha$-$R$-modules $P$, le préfaisceau $\cL P$ est séparé. 

{\rm (v)}\ Pour qu'un préfaisceau de $\alpha$-$R$-modules $P$ soit séparé, 
il faut et il suffit que $\lambda(P)\colon P\rightarrow \cL P$ soit un monomorphisme.
Le préfaisceau $\cL P$ est alors un faisceau. 

{\rm (vi)}\ Pour qu'un préfaisceau de $\alpha$-$R$-modules $P$ soit un faisceau, 
il faut et il suffit que $\lambda(P)\colon P\rightarrow \cL P$ soit un isomorphisme. 
\end{prop}

(i) En effet, les foncteurs $\halpha$, $L$ et $\hsigma_*$ sont exacts à gauche en vertu de \ref{alpha10}, \ref{alpha15} et (\cite{sga4} II 3.2(i)). 

(ii) En vertu de \ref{alpha10}(i), pour tout $(R_\hcC)$-module $F$ et tout $X\in \ob(\cC)$, on a un isomorphisme canonique
\begin{equation}\label{alpha26c}
\halpha(L F)(X)\stackrel{\sim}{\rightarrow}\underset{\underset{\cR\in J(X)^\circ}{\longrightarrow}}{\lim}\ \alpha(F(\cR)).
\end{equation}
Compte tenu de \eqref{alpha17f}, pour tout crible $\cR$ de $X$, on a un isomorphisme canonique 
\begin{equation}\label{alpha26d}
\alpha(F(\cR))\stackrel{\sim}{\rightarrow}\underset{\underset{(Y,u)\in (\cC_{/\cR})^\circ}{\longleftarrow}}{\lim}\ \alpha(F(Y)).
\end{equation}
La proposition s'ensuit.

(iii) Cela résulte de (ii) et \ref{alpha9}(iii). 

(iv) D'après \ref{alpha10}(i), $\halpha$ transforme les préfaisceaux séparés de $R$-modules en des 
préfaisceaux séparés de $\alpha$-$R$-modules. 
La proposition résulte donc de (\cite{sga4} II 3.2(i)).

(v) Si $P$ est un préfaisceau séparé de $\alpha$-$R$-modules, 
$\lambda(P)$ est un monomorphisme en vertu de \eqref{alpha24d}
car une limite inductive filtrante de monomorphismes est un monomorphisme \eqref{alpha150}. 
Inversement, si $\lambda(P)$ est un monomorphisme, $P$ est un sous-préfaisceau 
d'un préfaisceau séparé de $\alpha$-$R$-modules~; il est donc séparé.  
Dans ce cas, $\hsigma_*(P)$ est un préfaisceau séparé de $R$-modules d'après \ref{alpha10}(i); donc
 $L\circ \hsigma_*(P)$ est un faisceau de $R$-modules en vertu de (\cite{sga4} II 3.2(iii)) et par suite 
$\cL P$ est un faisceau de $\alpha$-$R$-modules d'après \ref{alpha10}(i).

(vi) En effet, la condition est nécessaire compte tenu de \eqref{alpha24d} et elle est suffisante en vertu de (v).

\begin{prop}\label{alpha25}
Le foncteur d'inclusion $\iota\colon \atcC\rightarrow \ahcC$ admet un adjoint à gauche
\begin{equation}\label{alpha25a}
\otta\colon \ahcC\rightarrow \atcC
\end{equation}
tel que le composé $\iota\circ \otta$ soit canoniquement isomorphe au foncteur $\cL\circ \cL$ \eqref{alpha24c}.
Pour tout préfaisceau de $\alpha$-$R$-modules $P$, le morphisme d'adjonction $P\rightarrow \iota\circ \otta (P)$ se déduit par l'isomorphisme 
précédent du morphisme $\lambda(\cL P)\circ \lambda(P)\colon P\rightarrow \cL\circ\cL(P)$. 
\end{prop}

En effet, d'après \ref{alpha26}, il existe un foncteur $\otta\colon \ahcC\rightarrow \atcC$ tel que $\iota\circ \otta=\cL\circ \cL$. 
On a un morphisme de foncteurs 
\begin{equation}\label{alpha25b}
\id\rightarrow \iota\circ \otta,
\end{equation}
défini pour tout préfaisceau de $\alpha$-$R$-modules $P$, par le morphisme 
$\lambda(\cL P)\circ \lambda(P)\colon P\rightarrow \cL\circ\cL(P)$. 
Si $P$ est un faisceau de $\alpha$-$R$-modules, 
$\lambda(\cL P)\circ \lambda(P)$ est un isomorphisme. On en déduit un isomorphisme 
$\iota\circ \otta\circ \iota\stackrel{\sim}{\rightarrow} \iota$ et par suite un isomorphisme  
\begin{equation}\label{alpha25c}
\otta\circ \iota\stackrel{\sim}{\rightarrow} \id.
\end{equation}
On vérifie aussitôt que les morphismes \eqref{alpha25b} et \eqref{alpha25c} font de $\otta$ un adjoint à gauche de $\iota$; d'où la proposition.

\begin{defi}\label{alpha29}
Pour tout préfaisceau de $\alpha$-$R$-modules $F$, 
on appelle $\otta(F)$ le {\em faisceau associé} à $F$ \eqref{alpha25a}. 
\end{defi}

\begin{prop}\label{alpha28}
{\rm (i)} Le foncteur $\otta$ \eqref{alpha25a} commute aux limites inductives et est exact.

{\rm (ii)} Les $\mU$-limites inductives dans $\atcC$ sont représentables. Pour toute catégorie $\mU$-petite $I$ et 
tout foncteur $\varphi\colon I\rightarrow \atcC$, le morphisme canonique
\begin{equation}\label{alpha28a}
\underset{\underset{I}{\rightarrow}}{\lim}\ \varphi \rightarrow 
\otta(\underset{\underset{I}{\rightarrow}}{\lim}\ \iota\circ \varphi)
\end{equation}
est un isomorphisme.

{\rm (iii)} Les $\mU$-limites projectives dans $\atcC$ sont représentables. Pour tout objet $X$ de $\cC$, le foncteur 
$F\mapsto F(X)$ commute aux limites projectives~; {\em i.e.}, le foncteur d'inclusion $\iota\colon \atcC\rightarrow \ahcC$ commute aux limites projectives. 
\end{prop}

En effet, $\otta$ commute aux limites inductives et $\iota$ commute aux limites projectives d'après \ref{alpha25}.
Soit $\mV$ un univers tel que $\cC\in \mV$ et $\mU\subset \mV$.
Comme les limites projectives commutent aux limites projectives, il résulte de \ref{alpha15} et \eqref{alpha17g}
que pour tout $X\in \ob(\cC)$ et tout crible $\cR$ de $X$, le foncteur \eqref{alpha17e}
\begin{equation}
\ahcC_\mU\rightarrow \aMod_\mV(R), \ \ \ F\mapsto F(\cR)
\end{equation} 
commute aux limites projectives.
Les propositions (ii) et (iii) s'ensuivent compte tenu de \ref{alpha15}.  
Comme $\iota\circ \otta\simeq \cL\circ \cL$ d'après \ref{alpha25} 
et que $\cL$ est exact à gauche en vertu de \ref{alpha26}(i), $\otta$ est exact à gauche et donc exact~; d'où la proposition (i).

\begin{prop}\label{alpha30}
La catégorie $\atcC$ est une catégorie abélienne vérifiant l'axiome {\rm (AB 5)} de {\rm (\cite{tohoku} §~1.5)}. 
\end{prop}

Il résulte de \ref{alpha150} et \ref{alpha28} que $\atcC$ est une catégorie additive où les noyaux et les conoyaux sont représentables. 
Plus précisément, soit $u\colon F\rightarrow G$ un morphisme de $\atcC$. D'après \ref{alpha28}, on a des isomorphismes canoniques
\begin{eqnarray}
\iota(\ker(u))&\stackrel{\sim}{\rightarrow}&\ker(\iota(u)),\label{alpha30a}\\
\coker(u)&\stackrel{\sim}{\rightarrow}&\otta(\coker(\iota(u))).\label{alpha30b}
\end{eqnarray}
On en déduit un isomorphisme canonique
\begin{equation}\label{alpha30c}
\coim(u)\stackrel{\sim}{\rightarrow}\otta(\coim(\iota(u))).
\end{equation}
Montrons que le morphisme canonique $\coim(u)\rightarrow \im(u)$ est un isomorphisme. 
D'après \ref{alpha28}(iii), il suffit de montrer que la suite 
\begin{equation}
0\rightarrow \iota(\coim(u))\rightarrow \iota(G)\rightarrow \iota(\coker(u))
\end{equation}
est exacte. Compte tenu de \eqref{alpha30b} et \eqref{alpha30c} et comme le morphisme d'adjonction $\id\rightarrow \otta\circ \iota$ est un isomorphisme, cette suite est isomorphe à l'image par le foncteur $\iota\circ \otta$ de la suite 
\begin{equation}
0\rightarrow \coim(\iota(u))\rightarrow  \iota(G)\rightarrow \coker(\iota(u)).
\end{equation}
Or, cette suite est exacte et le foncteur $\iota\circ \otta$ est exact à gauche d'après \ref{alpha28}. 
Par suite, $\atcC$ est une catégorie abélienne. Comme $\ahcC$ est une catégorie abélienne vérifiant (AB 5) \eqref{alpha150},
il en est de même de $\atcC$ en vertu de \ref{alpha28}.

\subsection{}\label{alpha22}
On identifie la catégorie $\bMod(R_\tcC)$ à la catégorie des $\mU$-faisceaux de $R$-modules sur $\cC$ 
(\cite{sga4} II 6.3.1). 
D'après \ref{alpha10}(i), le foncteur $\halpha$ \eqref{alpha16a} transforme les faisceaux de $R$-modules en des 
faisceaux de $\alpha$-$R$-modules.
Il définit donc un foncteur 
\begin{equation}\label{alpha22a}
\talpha\colon \bMod(R_\tcC)\rightarrow \atcC,
\end{equation} 
qui s'insère dans un diagramme commutatif 
\begin{equation}\label{alpha22g}
\xymatrix{
{\bMod(R_\tcC)}\ar[r]^-(0.5)\talpha\ar[d]_i&{\atcC}\ar[d]^\iota\\
{\bMod(R_\hcC)}\ar[r]^-(0.5)\halpha&{\ahcC}}
\end{equation}
où $i$ est $\iota$ sont les foncteurs d'injection canoniques. 
D'après \ref{alpha21}(iii), $\talpha$ transforme les $\alpha$-isomor\-phismes en des isomorphismes. Il induit donc un foncteur
\begin{equation}\label{alpha22b}
\mu\colon \aMod(R_\tcC)\rightarrow \atcC.
\end{equation} 
D'autre part, le foncteur $\hsigma_*$ \eqref{alpha16c} transforme les faisceaux de 
$\alpha$-$R$-modules en des faisceaux de $R$-modules d'après \ref{alpha182}. Il induit donc un foncteur 
\begin{equation}\label{alpha22c}
\tsigma_*\colon \atcC\rightarrow \bMod(R_\tcC).
\end{equation}
On note $\nu$ le foncteur composé
\begin{equation}\label{alpha22d}
\nu=\alpha\circ \tsigma_*\colon \atcC\rightarrow \aMod(R_\tcC),
\end{equation} 
où $\alpha$ est le foncteur de localisation pour les $(R_\tcC)$-modules \eqref{alpha20}.
D'après \ref{alpha16}, le foncteur $\tsigma_*$ est un adjoint à droite de $\talpha$; 
le morphisme d'adjonction $\talpha\circ \tsigma_* \rightarrow \id$ 
est un isomorphisme~; et le morphisme d'adjonction $\id\rightarrow \tsigma_*\circ\talpha$
induit un isomorphisme $\talpha\stackrel{\sim}{\rightarrow} \talpha\circ\tsigma_*\circ\talpha$. 

L'isomorphisme $\talpha\circ \tsigma_*\stackrel{\sim}{\rightarrow} \id$ induit un isomorphisme 
\begin{equation}\label{alpha22e}
\mu\circ \nu\stackrel{\sim}{\rightarrow} \id.
\end{equation}
D'après \ref{alpha21}(iii), le morphisme d'adjonction $\id\rightarrow \tsigma_*\circ \talpha$ induit un isomorphisme 
$\alpha\stackrel{\sim}{\rightarrow} \nu\circ \mu\circ \alpha$. On en déduit un isomorphisme 
\begin{equation}\label{alpha22f}
\id\stackrel{\sim}{\rightarrow} \nu\circ \mu.
\end{equation} 
On vérifie aussitôt que les isomorphismes \eqref{alpha22e} et \eqref{alpha22f} font de $\nu$ un adjoint à droite de $\mu$. 
Par suite, $\mu$ et $\nu$ sont des équivalences de catégories, quasi-inverses l'une de l'autre. 
Ce sont donc des foncteurs exacts (\cite{gabriel} I §~1 prop.~13). Il s'ensuit que le foncteur $\talpha$ est aussi exact.

\begin{prop}\label{alpha27}
Le diagramme 
\begin{equation}\label{alpha27a}
\xymatrix{
{\bMod(R_\hcC)}\ar[r]^-(0.5){\tta}\ar[d]_{\halpha}&{\bMod(R_\tcC)}\ar[d]^{\talpha}\\
{\ahcC}\ar[r]^-(0.5)\otta&{\atcC}}
\end{equation}
où $\otta$ est le foncteur \eqref{alpha25a} et $\tta$ est le foncteur ``faisceau de $R$-modules associé'' {\rm (\cite{sga4} II 6.4)},
est commutatif à isomorphisme canonique près. 
\end{prop}

En effet, notant $i\colon \bMod(R_\tcC)\rightarrow \bMod(R_\hcC)$ le foncteur d'injection canonique, on a un isomorphisme canonique
$i\circ \tta\stackrel{\sim}{\rightarrow}L\circ L$ (\cite{sga4} II 3.4 et 6.4). 
Par suite, en vertu de \ref{alpha26}(iii), on a un isomorphisme canonique
\begin{equation}
\halpha \circ i\circ \tta\stackrel{\sim}{\rightarrow}\iota\circ \otta\circ \halpha.
\end{equation}
La proposition s'ensuit puisqu'on a $\halpha\circ i=\iota\circ \talpha$ \eqref{alpha22a}.

\begin{prop}\label{alpha32}
La catégorie abélienne $\atcC$ admet une famille de générateurs, indexée par un ensemble appartenant à $\mU$. 
\end{prop}

Cela résulte de \ref{alpha22}, (\cite{sga4} II 6.7) et (\cite{gabriel} III §~2 lem.~4). 

\subsection{}\label{alpha34}
Pour tout $R^\alpha$-module $P$, on appelle {\em faisceau de $\alpha$-$R$-modules constant de valeur $P$ sur $\cC$} 
et l'on note $P_\tcC$ le faisceau associé au 
préfaisceau de $\alpha$-$R$-modules constant de valeur $P$ sur $\cC$ \eqref{alpha14}. D'après \ref{alpha27}, 
pour tout $R$-module $M$,
notant $M_\tcC$ le faisceau constant de valeur $M$ sur $\cC$, on a un isomorphisme canonique 
\begin{equation}\label{alpha34a}
M^\alpha_{\tcC}\stackrel{\sim}{\rightarrow}\talpha(M_\tcC).
\end{equation}

On munit $\atcC$ de la structure de catégorie abélienne tensorielle déduite de celle de $\aMod(R_\tcC)$ \eqref{alpha33a} 
via l'équivalence de catégories $\mu$ \eqref{alpha22b}; on note $\otimes_{R^\alpha_\tcC}$ le produit tensoriel dans $\atcC$. 
Le foncteur $\talpha$ est donc monoïdal~: pour tous $R_\tcC$-modules $M$ et $N$, on a un isomorphisme canonique 
\begin{equation}\label{alpha34b}
\talpha(M\otimes_{R_\tcC}N)\stackrel{\sim}{\rightarrow}M^\alpha\otimes_{R^\alpha_\tcC}N^\alpha.
\end{equation}
D'après \eqref{alpha34a}, $R^\alpha_\tcC$ est un objet unité de $\atcC$.

\begin{prop}\label{alpha35}
Pour tous préfaisceaux de $\alpha$-$R$-modules $M$ et $N$ sur $\cC$, on a un isomorphisme canonique fonctoriel 
\begin{equation}\label{alpha35a}
\otta(M) \otimes_{R^\alpha_\tcC}\otta(N)\stackrel{\sim}{\rightarrow}\otta(M\otimes_{R^\alpha_\hcC}N),
\end{equation}
où $M\otimes_{R^\alpha_\hcC}N$ désigne le produit tensoriel dans $\ahcC$ \eqref{alpha14a}. 
\end{prop}

Cela résulte aussitôt de \ref{alpha27} et (\cite{sga4} IV 12.10) puisque le foncteur $u$ \eqref{alpha16b} est une équivalence 
de catégories tensorielles.

\subsection{}\label{alpha37}
Considérons les foncteurs adjoints \eqref{alpha25}
\begin{equation}\label{alpha37c}
\xymatrix{{\atcC}\ar@<1ex>[r]^\iota&{\ahcC}\ar@<1ex>[l]^\otta}
\end{equation}
D'après \ref{alpha35}, pour tout objet $A$ de $\atcC$, on a un isomorphisme canonique
\begin{equation}\label{alpha37d}
\otta(\iota(A)\otimes_{R^\alpha_\hcC}\iota(A))\stackrel{\sim}{\rightarrow} A\otimes_{R^\alpha_\tcC}A.
\end{equation}
On en déduit par adjonction que la donnée d'une structure de monoïde commutatif unitaire de $\atcC$ sur $A$ 
est équivalente à la donnée sur $\iota(A)$ d'une structure de monoïde commutatif unitaire de $\ahcC$.
Par suite, la donnée d'un monoïde commutatif unitaire de $\atcC$ est équivalente à la donnée d'un 
préfaisceau sur $\cC$ à valeur dans $\aAlg(R)$ \eqref{alpha5} dont le préfaisceau de $\alpha$-$R$-modules  
sous-jacent est un faisceau \eqref{alpha36}. 
On note $\bAlg(\atcC)$ la catégorie des monoïdes commutatifs unitaires de $\atcC$. 
  
De même, la donnée d'un monoïde commutatif unitaire de $\bMod(R_\tcC)$ est équivalente à la donnée 
d'une $(R_\tcC)$-algèbre de $\tcC$ (\cite{sga4} II 6.3.1 et IV 12.10). 
On note $\bAlg(R_\tcC)$ la catégorie des monoïdes commutatifs unitaires de $\bMod(R_\tcC)$.

Le foncteur $\talpha$ \eqref{alpha22a} étant monoïdal, il induit un foncteur que l'on note encore
\begin{equation}\label{alpha37e}
\talpha\colon \bAlg(R_\tcC)\rightarrow \bAlg(\atcC).
\end{equation}
Compte tenu de \ref{alpha5}, le foncteur $\tsigma_*$ \eqref{alpha22c} induit un foncteur que l'on note encore
\begin{equation}\label{alpha37f}
\tsigma_*\colon \bAlg(\atcC)\rightarrow \bAlg(R_\tcC).
\end{equation}
D'après \ref{alpha6}, $\tsigma_*$ est un adjoint à droite de $\talpha$; 
le morphisme d'adjonction $\talpha\circ \tsigma_* \rightarrow \id$  
est un isomorphisme~; et le morphisme d'adjonction $\id\rightarrow \tsigma_*\circ\talpha$
induit un isomorphisme $\talpha\stackrel{\sim}{\rightarrow} \talpha\circ\tsigma_*\circ\talpha$. 

\subsection{}\label{alpha38}
Dans ce numéro, si $D$ est une $(R_\tcC)$-algèbre de $\tcC$ (resp. une $R$-algèbre), nous affecterons d'un indice $D$ le foncteur 
de localisation $\alpha$ \eqref{alpha2a} pour les $D$-modules.

Soit $A$ une $(R_\tcC)$-algèbre de $\tcC$. On pose $B=\talpha(A)$ \eqref{alpha37e}
et on désigne par $\bMod(B)$ la catégorie des $B$-modules unitaires de $\atcC$. 
Il résulte aussitôt de \ref{alpha35} 
que la donnée d'une structure de $B$-module unitaire sur un faisceau de $\alpha$-$R$-modules 
$M$ est équivalente à la donnée
pour tout $X\in \ob(\cC)$ d'une structure de $\alpha_R(A(X))$-module unitaire sur $M(X)$ dans le sens de \ref{alpha39} telles que 
pour tout morphisme $X\rightarrow Y$ de $\cC$, le morphisme $M(Y)\rightarrow M(X)$ 
soir linéaire relativement au morphisme de $R^\alpha$-algèbres $\alpha_R(A(Y))\rightarrow \alpha_R(A(X))$. 

Le foncteur $\talpha$ \eqref{alpha22a} étant monoïdal, il définit un foncteur 
\begin{equation}\label{alpha38a}
\tbeta\colon \bMod(A)\rightarrow \bMod(B).
\end{equation}
Celui-ci transforme les $\alpha$-isomorphismes en des isomorphismes d'après \ref{alpha21}(iii). Il induit donc un foncteur \eqref{alpha20}
\begin{equation}\label{alpha38b}
b\colon \aMod(A)\rightarrow \bMod(B).
\end{equation}
On pose $A'=\tsigma_*(B)$ \eqref{alpha37f}. D'après \ref{alpha37}, 
on a un homomorphisme canonique de $(R_\tcC)$-algèbres $\lambda\colon 
A\rightarrow A'$, qui induit un isomorphisme $\talpha(A)\stackrel{\sim}{\rightarrow}\talpha(A')$. 
Compte tenu de \ref{alpha39}, le foncteur $\tsigma_*$ \eqref{alpha22c} induit un foncteur 
\begin{equation}\label{alpha38c}
\ttau'_*\colon \bMod(B)\rightarrow \bMod(A').
\end{equation}
Composant avec le foncteur induit par $\lambda$, on obtient un foncteur 
\begin{equation}\label{alpha38d}
\ttau_*\colon \bMod(B)\rightarrow \bMod(A).
\end{equation}
On pose 
\begin{equation}\label{alpha38e}
t=\alpha_A\circ \ttau_*\colon \bMod(B)\rightarrow \aMod(A).
\end{equation}

D'après \ref{alpha39} et \ref{alpha16}, le foncteur $\ttau_*$ est un adjoint à droite de $\tbeta$; 
le morphisme d'adjonction $\tbeta\circ \ttau_* \rightarrow \id$ 
est un isomorphisme~; et le morphisme d'adjonction $\id\rightarrow \ttau_*\circ\tbeta$
induit un isomorphisme $\tbeta\stackrel{\sim}{\rightarrow} \tbeta\circ\ttau_*\circ\tbeta$. 

L'isomorphisme $\tbeta\circ \ttau_*\stackrel{\sim}{\rightarrow} \id$ induit un isomorphisme 
\begin{equation}\label{alpha38f}
b\circ t\stackrel{\sim}{\rightarrow} \id.
\end{equation}
D'après \ref{alpha21}(iii), le morphisme d'adjonction $\id\rightarrow \ttau_*\circ \tbeta$ induit un isomorphisme 
$\alpha_A\stackrel{\sim}{\rightarrow} t\circ b\circ \alpha_A$. On en déduit un isomorphisme 
\begin{equation}\label{alpha38g}
\id\stackrel{\sim}{\rightarrow} t\circ b.
\end{equation} 
On vérifie aussitôt que les isomorphismes \eqref{alpha38f} et \eqref{alpha38g} font de $t$ un adjoint à droite de $b$. 
Par suite, $b$ et $t$ sont des équivalences de catégories, quasi-inverses l'une de l'autre.

\section{\texorpdfstring{Conditions de $\alpha$-finitude}{Conditions de alpha-finitude}}\label{finita}

\subsection{}\label{finita1}
Les hypothèses et notations de \ref{alpha1} sont en vigueur dans cette section. 
On se donne, de plus, un $\mU$-topos $\cT$ \eqref{notconv3} annelé par une $R$-algèbre $A$ (\cite{sga4} IV 11.1.1). 
Nous considérons toujours $\cT$ comme muni de sa topologie canonique (\cite{sga4} II 2.5), qui en fait un $\mU$-site.
Pour tout objet $X$ de $\cT$, le topos $\cT_{/X}$ sera annelé par l'anneau $A|X$ \eqref{notconv6}. 

\subsection{}\label{finita3}
Soit $\gamma\in R$. On dit qu'un morphisme de $A$-modules est un {\em $\gamma$-isomorphisme} 
si son noyau et son conoyau sont annulés par $\gamma$. 

On appelle {\em suite de $A$-modules} un complexe de longueur finie de $A$-modules, 
à différentielle de degré $1$, le degré étant écrit en exposant.
On dit qu'une suite de $A$-modules est {\em $\gamma$-exacte} si ses groupes de cohomologie sont annulés par $\gamma$,
et qu'elle est {\em $\alpha$-exacte} si elle est $\gamma$-exacte pour tout $\gamma\in \fm$.

\begin{defi}\label{finita9}
Soient $F$ un $A$-module, $n$ un entier $\geq 0$, $\gamma\in R$.
\begin{itemize}
\item[(i)] On dit que $F$ est {\em de $n$-présentation $\gamma$-finie} si 
la sous-catégorie pleine de $\cT$ formée des objets $X$ tels qu'il existe une suite $\gamma$-exacte \eqref{finita3}
\begin{equation}
E^{-n}\rightarrow E^{-n+1}\rightarrow \dots\rightarrow E^0\rightarrow F|X\rightarrow 0,
\end{equation} 
avec $E^i$ un $(A|X)$-module libre de type fini pour tout $-n\leq i\leq 0$, est un raffinement de l'objet final de $\cT$. 
On dit que $F$ est de {\em type $\gamma$-fini} (resp. {\em présentation $\gamma$-finie}) 
s'il est de $0$-présentation $\gamma$-finie (resp. $1$-présentation $\gamma$-finie).

\item[(ii)] On dit que $F$ est de {\em $n$-présentation $\alpha$-finie} s'il est de $n$-présentation 
$\gamma$-finie pour tout $\gamma\in \fm$. 
On dit que $F$ est de {\em type $\alpha$-fini} (resp. {\em présentation $\alpha$-finie}) 
s'il est de $0$-présentation $\alpha$-finie (resp. $1$-présentation $\alpha$-finie).

\item[(iii)] On dit que $F$ est {\em $\alpha$-cohérent} s'il est de type $\alpha$-fini et si 
pour tout objet $X$ de $\cT$ et tout $(A|X)$-morphisme $u\colon E\rightarrow F|X$, où $E$ est un $(A|X)$-module libre de type fini, 
$\ker(u)$ est un $(A|X)$-module de type $\alpha$-fini. 
\end{itemize}
\end{defi}

Si $\gamma$ est une unité de $R$, la notion de $n$-présentation $\gamma$-finie correspond 
à la notion standard de $n$-présentation finie introduite dans (\cite{sga6} I 2.8). 
La notion de $\alpha$-cohérence est modelée sur la notion standard de cohérence (\cite{sga6} I 3.1). 

\begin{rema}\label{finita2}
(i)\ Supposons que $\cT$ soit le topos ponctuel, annelé par un anneau standard $A$. 
Pour qu'un $A$-module $F$ soit de type $\alpha$-fini (resp. présentation $\alpha$-finie), il faut et il suffit que pour tout $\gamma\in \fm$, 
il existe un $\gamma$-isomorphisme $f\colon G\rightarrow F$ avec $G$ un $A$-module de type fini (resp. de présentation finie) 
\eqref{finita3}.

(ii)\ Nous donnerons dans \ref{afini5} un cas intéressant pour lequel on dispose d'une caractérisation des modules 
de type $\alpha$-fini (resp. présentation $\alpha$-finie) similaire à celle pour le topos ponctuel (i). 
\end{rema}

\begin{lem}\label{finita15}
Soient $\gamma,\gamma'\in R$, $E\rightarrow F$ un $\gamma$-isomorphisme de $A$-modules, $n$ un entier $\geq 0$. 
\begin{itemize}
\item[{\rm (i)}] Si $E$ est de $n$-présentation $\gamma'$-finie, $F$ est de $n$-présentation $\gamma\gamma'$-finie.
\item[{\rm (ii)}] Si $F$ est de $n$-présentation $\gamma'$-finie, $E$ est de $n$-présentation $\gamma^2\gamma'$-finie.
\end{itemize}
\end{lem}

L'assertion (i) est immédiate et l'assertion (ii) résulte aussitôt de \ref{alpha3}.

\begin{lem}\label{finita16}
Soient $f\colon E\rightarrow F$ un morphisme de $A$-modules qui est un $\alpha$-isomorphisme, $n$ un entier $\geq 0$. 
Pour que $E$ soit de $n$-présentation $\alpha$-finie (resp. $\alpha$-cohérent), il faut et il suffit qu'il en soit de même de $F$.
\end{lem}

L'assertion relative à la $n$-présentation $\alpha$-finie résulte aussitôt de \ref{finita15}. Supposons que $F$ soit $\alpha$-cohérent et montrons 
que $E$ est $\alpha$-cohérent. On sait déjà que $E$ est de type $\alpha$-fini. Soient $X$ un objet de $\cT$, $E'$ un $(A|X)$-module libre de type fini,
$u\colon E'\rightarrow E|X$ un $(A|X)$-morphisme. Le morphisme canonique $\ker(u)\rightarrow \ker(f\circ u)$ étant un $\alpha$-isomorphisme, on en déduit 
que $\ker(u)$ est un $(A|X)$-module de type $\alpha$-fini. Par suite, $E$ est $\alpha$-cohérent. 
Inversement, supposons que $E$ soit $\alpha$-cohérent et montrons que $F$ est $\alpha$-cohérent. On sait que $F$ est de type $\alpha$-fini. 
Soient $X$ un objet de $\cT$, $F'$ un $(A|X)$-module libre de type fini, $v\colon F'\rightarrow F|X$ un $(A|X)$-morphisme. Pour tout $\gamma\in \fm$,
il existe un $\gamma^2$-isomorphisme $g\colon F\rightarrow E$ \eqref{alpha3}. Le morphisme canonique
$\ker(v)\rightarrow \ker(g\circ v)$ est donc un $\gamma^2$-isomorphisme. Comme $\ker(g\circ v)$ est de type $\alpha$-fini, 
on déduit que $\ker(v)$ est de type $\gamma^5$-fini d'après \ref{finita15}(ii)
et est donc de type $\alpha$-fini. Par suite, $F$ est $\alpha$-cohérent.

\begin{lem}\label{finita14}
Soient $f\colon (\cT',A')\rightarrow (\cT,A)$ un morphisme de topos annelés, $\gamma\in R$.
\begin{itemize}
\item[{\em (i)}] Si $u\colon E\rightarrow F$ est un $\gamma$-isomorphisme de $A$-modules, alors $f^*(u)\colon f^*(E)\rightarrow f^*(F)$ est un 
$\gamma^2$-isomorphisme. 
\item[{\rm (ii)}] Si $E\rightarrow F\rightarrow G\rightarrow 0$ est une suite $\gamma$-exacte de $A$-modules, 
son image inverse $f^*(E)\rightarrow f^*(F)\rightarrow f^*(G)\rightarrow 0$ est $\gamma^2$-exacte.
\item[{\rm (iii)}] Si $F$ est un $A$-module de type $\gamma$-fini 
(resp. de présentation $\gamma$-finie, resp. de type $\alpha$-fini, resp. de présentation $\alpha$-finie), son image inverse $f^*(F)$ est de type $\gamma$-fini 
(resp. de présentation $\gamma^2$-finie, resp. de type $\alpha$-fini, resp. de présentation $\alpha$-finie). 
\end{itemize}
\end{lem}

(i) Cela résulte aussitôt de \ref{alpha3}.

(ii) Notant $G'$ le conoyau du morphisme $E\rightarrow F$, il revient au même de dire que la suite 
$E\rightarrow F\rightarrow G\rightarrow 0$ est $\gamma$-exacte ou que le morphisme induit $G'\rightarrow G$ 
est un $\gamma$-isomorphisme. La proposition résulte alors de (i). 

(iii) Cela résulte aussitôt de (ii) et des définitions.

\begin{lem}\label{finita10}
Soient $X$ un objet de $\cT$, $F$ un $A$-module. 
Si $F$ est $\alpha$-cohérent, il en est de même du $(A|X)$-module $F|X$. 
\end{lem}

Cela résulte aussitôt des définitions \ref{finita9}. 

\begin{lem}\label{finita11}
Soient $(X_i)_{i\in I}$ un raffinement de l'objet final de $\cT$, $F$ un $A$-module, $n$ un entier $\geq 0$, $\gamma\in R$. 
Si pour tout $i\in I$, le $(A|X_i)$-module $F|X_i$ est de $n$-présentation $\gamma$-finie 
(resp. de $n$-présentation $\alpha$-finie, resp. $\alpha$-cohérent), il en est de même de $F$.
\end{lem}

Cela résulte aussitôt des définitions \ref{finita9}.

\begin{prop}\label{finita8}
Soient $F$ un $A$-module, $n$ un entier $\geq 0$, $\gamma\in R$.
Supposons que $\cT$ soit équivalent au topos des faisceaux de $\mU$-ensembles sur un $\mU$-site $\cC$. 
Pour tout objet $X$ de $\cC$, on note $X^\tta$ le faisceau associé à $X$. Alors,
\begin{itemize} 
\item[{\rm (i)}] Pour que $F$ soit de $n$-présentation $\gamma$-finie, il faut et il suffit que pour tout $U\in \ob(\cC)$, 
la sous-catégorie pleine de $\cC_{/U}$ formée des objets $X\rightarrow U$ tels qu'il existe une suite $\gamma$-exacte 
\begin{equation}
E^{-n}\rightarrow E^{-n+1}\rightarrow \dots\rightarrow E^0\rightarrow F|X^\tta\rightarrow 0,
\end{equation} 
avec $E^i$ un $(A|X^\tta)$-module libre de type fini pour tout $-n\leq i\leq 0$, soit un raffinement de $U$. 
\item[{\rm (ii)}] Pour que $F$ soit $\alpha$-cohérent, il faut et il suffit que $F$ soit de type $\alpha$-fini et que pour tout $X\in \ob(\cC)$
et tout morphisme $u\colon E\rightarrow F|X^\tta$, où $E$ est un $(A|X^\tta)$-module libre de type fini, 
$\ker(u)$ soit un $(A|X^\tta)$-module de type $\alpha$-fini. 
\end{itemize}
\end{prop}

Cela résulte aussitôt de \ref{finita11} et (\cite{sga4} II 4.4 et 4.10).

\subsection{}\label{finita4}
Soient $\gamma\in R$, 
\begin{equation}\label{finita4a}
\xymatrix{
&E\ar[r]^u\ar[d]_e&F\ar[r]^v\ar[d]_f&G\ar[d]_g\ar[r]&0\\
0\ar[r]&E'\ar[r]^{u'}&F'\ar[r]^{v'}&G'&}
\end{equation}
un diagramme commutatif de $A$-modules tel que les lignes soient $\gamma$-exactes. 
Ce dernier induit des morphismes $E'\rightarrow \ker(v')$ et $\coker(u)\rightarrow G$ et un diagramme commutatif 
\begin{equation}\label{finita4c}
\xymatrix{
&E\ar[r]^u\ar[d]_{e'}&F\ar[r]\ar[d]_f&{\coker(u)}\ar[d]_{g'}\ar[r]&0\\
0\ar[r]&{\ker(v')}\ar[r]&F'\ar[r]^{v'}&G'&}
\end{equation}
On en déduit un diagramme commutatif 
\begin{equation}\label{finita4d}
\xymatrix{
{\ker(e)}\ar[rd]^{u_1}\ar@{^(->}[d]_a&&&{\coker(e)}\ar[d]_c\ar[rd]^{u'_1}&&\\
{\ker (e')}\ar[r]&{\ker(f)}\ar[r]\ar[rd]_{v_1}&{\ker(g')}\ar[d]^b\ar[r]&{\coker(e')}\ar[r]&{\coker(f)}\ar[rd]_{v'_1}\ar[r]&{\coker(g')}\ar@{->>}[d]^-(0.5)d\\
&&{\ker(g)}&&&{\coker(g)}}
\end{equation}
où la ligne centrale est exacte, $a$ est injectif et $d$ est surjectif. Il résulte aussitôt des hypothèses que $a$, $b$, 
$c$ et $d$ sont des $\gamma$-isomorphismes. En particulier, les suites $(u_1,v_1)$ et $(u'_1,v'_1)$ sont $\gamma^2$-exactes.

\begin{prop}\label{finita12}
Soient $\gamma,\gamma',\gamma''$ et $\lambda$ des éléments de $R$, 
$0\rightarrow F'\rightarrow F \rightarrow F''\rightarrow 0$ une suite $\lambda$-exacte de $A$-modules. 
\begin{itemize}
\item[{\rm (i)}] Si $F$ est de type $\gamma$-fini, $F''$ est de type $(\lambda\gamma)$-fini. 
\item[{\rm (ii)}] Si $F'$ est de type $\gamma'$-fini (resp. présentation $\gamma'$-finie) et
$F''$ est de type $\gamma''$-fini (resp. présentation $\gamma''$-finie), 
$F$ est de type $(\lambda^2\gamma'\gamma'')$-fini (resp. présentation $(\lambda^5\gamma'^3\gamma'')$-finie). 
\item[{\rm (iii)}] Si $F$ est de type $\gamma$-fini et $F''$ est de présentation $\gamma''$-finie, 
$F'$ est de type $(\lambda^8\gamma\gamma''^4)$-fini. 
\item[{\rm (iv)}] Si $F$ est de présentation $\gamma$-finie et $F'$ est de type $\gamma'$-fini, 
$F''$ est de présentation $(\lambda\gamma^3\gamma')$-finie.
\end{itemize}
\end{prop} 

(i) C'est immédiat. 

(ii) La question étant locale, on peut supposer qu'il existe des morphismes $A$-linéaires    
$f'\colon E'\rightarrow F'$ et $f''\colon E''\rightarrow F''$ avec $E'$ et $E''$ des $A$-modules libres de type fini
dont les conoyaux sont annulés par $\gamma'$ et $\gamma''$ respectivement. Notons $u'\colon E'\rightarrow F$ le morphisme
induit par $f'$. D'autre part, il existe un morphisme $u''\colon E''\rightarrow F$ qui relève $\lambda f''$. On voit aussitôt que 
le conoyau du morphisme $u'+u''\colon E'\oplus E''\rightarrow F$ est annulé par $\lambda^2\gamma'\gamma''$.
Par suite, $F$ est de type $(\lambda^2\gamma'\gamma'')$-fini. D'après \ref{finita4}, on a une suite $(\lambda^2\gamma')$-exacte
\begin{equation}\label{finita12a}
0\rightarrow \ker(f')\rightarrow \ker(u'+u'')\rightarrow \ker(\lambda f'')\rightarrow 0.
\end{equation}
Supposons $F'$ de présentation $\gamma'$-finie et $F''$ de présentation $\gamma''$-finie, et montrons que $F$ est de présentation 
$(\lambda^4\gamma'^3\gamma'')$-finie. La question étant locale, on peut supposer de plus  
$\ker(f')$ de type $\gamma'$-fini et  $\ker(f'')$ de type $\gamma''$-fini. Par suite, $\ker(\lambda f'')$ est de type $(\lambda\gamma'')$-fini.
Il résulte alors de \eqref{finita12a} et de ce qui précède que $\ker(u'+u'')$ est de type 
$(\lambda^5\gamma'^3\gamma'')$-fini~; d'où l'assertion. 

(iii) La question étant locale, on peut se borner au cas où il existe une suite $\gamma''$-exacte 
\begin{equation}
G''\stackrel{g}{\longrightarrow} E''\stackrel{f}{\longrightarrow} F''\longrightarrow 0
\end{equation} 
avec $E''$ et $G''$ des $A$-modules libres de type fini. Il existe alors un diagramme commutatif 
\begin{equation}
\xymatrix{
&G''\ar[r]^{\lambda g}\ar[d]_v&E''\ar[r]^{\lambda f}\ar[d]_{u}&F''\ar@{=}[d]\ar[r]&0\\
0\ar[r]&F'\ar[r]&F\ar[r]&F''\ar[r]&0}
\end{equation}
Les lignes étant $(\lambda^2\gamma'')$-exactes,
le morphisme canonique $\coker(v)\rightarrow \coker(u)$ est un $(\lambda^4\gamma''^2)$-isomorphisme en vertu de \ref{finita4}. 
D'après (i), $\coker(u)$ est de type $\gamma$-fini. Donc $\coker(v)$ est  de type $(\lambda^8\gamma''^4\gamma)$-fini
en vertu de (ii) ou de \ref{alpha3}. Compte tenu de la suite exacte $0\rightarrow \im(v)\rightarrow F'\rightarrow \coker(v)\rightarrow 0$,
on en déduit  par (ii) que $F'$ est de type $(\lambda^8\gamma\gamma''^4)$-fini.

(iv) La question étant locale, on peut se borner au cas où il existe un $A$-module libre de type fini $E$ 
et un morphisme $A$-linéaire $f\colon E\rightarrow F$ dont le conoyau est annulé par $\gamma$ et dont le noyau est de type
$\gamma$-fini. Notons $u\colon F\rightarrow F''$ le morphisme donné et posons $f''=u\circ f\colon E\rightarrow F''$. 
On a alors une suite exacte 
\begin{equation}
0\rightarrow \ker(f)\rightarrow \ker(f'')\rightarrow \ker(u)\rightarrow \coker(f).
\end{equation}
Le conoyau du morphisme canonique $F'\rightarrow \ker(u)$ est annulé par $\lambda$. Par suite, $\ker(u)$ est de type 
$(\lambda\gamma')$-fini. On en déduit que $\ker(f'')$ est de type $(\lambda\gamma^3\gamma')$-fini en vertu de (ii). 
Comme $\coker(f'')$ est annulé par $\lambda\gamma$, $F''$ est de présentation $(\lambda\gamma^3\gamma')$-finie.

\begin{cor}\label{finita13}
Soit $0\rightarrow F'\rightarrow F \rightarrow F''\rightarrow 0$ une suite $\alpha$-exacte de $A$-modules. Alors, 
\begin{itemize}
\item[{\rm (i)}] Si $F$ est de type $\alpha$-fini, il en est de même de $F''$. 
\item[{\rm (ii)}] Si $F'$ et $F''$ sont de type $\alpha$-fini (resp. présentation $\alpha$-finie), il en est de même de $F$. 
\item[{\rm (iii)}] Si $F$ est de type $\alpha$-fini et $F''$ est de présentation $\alpha$-finie, $F'$ est de type $\alpha$-fini. 
\item[{\rm (iv)}] Si $F$ est de présentation $\alpha$-finie et $F'$ est de type $\alpha$-fini, $F''$ est de présentation $\alpha$-finie.
\end{itemize}
\end{cor} 

\begin{prop}\label{finita17}
Soit $0\rightarrow F'\stackrel{u}{\rightarrow} F \stackrel{v}{\rightarrow} F''\rightarrow 0$ une suite $\alpha$-exacte de $A$-modules. 
Si deux $A$-modules sont $\alpha$-cohérents, il en est de même du troisième.
\end{prop}

Supposons d'abord que $F$ et $F''$ soient $\alpha$-cohérents. D'après \ref{finita13}(iii), $F'$ est de type $\alpha$-fini. La question étant locale, 
on peut supposer qu'il existe un morphisme $A$-linéaire et $\alpha$-surjectif $w\colon A^n\rightarrow F$. 
Comme $F''$ est cohérent, $E=\ker(u\circ w)$ est de type $\alpha$-fini. Il en est alors de même de $w(E)\subset F$, 
qui est donc $\alpha$-cohérent. 
L'injection canonique $w(E)\rightarrow \ker(v)$ et le morphisme canonique $F'\rightarrow \ker(v)$ sont des $\alpha$-isomorphismes.
On en déduit que $F'$ est $\alpha$-cohérent d'après \ref{finita16}.

Supposons ensuite que $F$ et $F'$ soient $\alpha$-cohérents. Il est loisible de remplacer $F'$ par $\ker(v)$ \eqref{finita16}.
Le morphisme canonique $\coker(u)\rightarrow F''$ étant un $\alpha$-isomorphisme, 
on peut supposer que c'est un isomorphisme \eqref{finita16}. D'après \ref{finita13}(i), $F''$ est de type $\alpha$-fini. 
Soient $X$ un objet de $\cT$, $f\colon (A|X)^n\rightarrow F''|X$ un morphisme $(A|X)$-linéaire. Montrons que $\ker(f)$ est de type $\alpha$-fini.
La question étant locale, on peut supposer qu'il existe un morphisme $(A|X)$-linéaire $f'\colon (A|X)^n\rightarrow F|X$ tel que $f=v\circ f'$.
Comme $F'$ est de type $\alpha$-fini, on peut supposer qu'il existe un morphisme $(A|X)$-linéaire et $\alpha$-surjectif $h\colon (A|X)^m\rightarrow F'|X$. 
On en déduit un morphisme $g\colon (A|X)^{m+n}\rightarrow F|X$ qui s'insère dans un diagramme commutatif 
\begin{equation}
\xymatrix{
{F'|X}\ar[r]^u&{F|X}\ar[r]^v&{F''|X}\\
(A|X)^m\ar[r]^i\ar[u]^h&(A|X)^{m+n}\ar[r]^\pi\ar[u]^g&(A|X)^n\ar[u]_f}
\end{equation}
où identifiant $(A|X)^{m+n}$ à $(A|X)^m\oplus (A|X)^n$, $i(x)=x+0$ et $\pi(x+x')=x'$. On en déduit une suite exacte
\begin{equation}
\ker(g)\rightarrow \ker(f)\rightarrow \coker(h).
\end{equation} 
Par hypothèse, $\ker(g)$ est de type $\alpha$-fini et $\coker(h)$ est $\alpha$-nul~; donc $\ker(f)$ est de type $\alpha$-fini, ce qui prouve que 
$F''$ est $\alpha$-cohérent. 

Supposons enfin que $F'$ et $F''$ soient $\alpha$-cohérents. Il est loisible de remplacer $F'$ par $\ker(v)$ \eqref{finita16}. 
D'après \ref{finita13}(ii), $F$ est de type $\alpha$-fini. Soient $X$ un objet de $\cT$, $f\colon (A|X)^n\rightarrow F|X$ un morphisme $(A|X)$-linéaire. 
Montrons que $\ker(f)$ est de type $\alpha$-fini. Comme $F''$ est $\alpha$-cohérent, $\ker(v\circ f)$ est de type $\alpha$-fini. 
La question étant locale, on peut supposer qu'il existe un morphisme $(A|X)$-linéaire et $\alpha$-surjectif $w\colon (A|X)^m\rightarrow \ker(v\circ f)$. 
Il existe un morphisme $(A|X)$-linéaire $f'\colon (A|X)^m\rightarrow F'|X$ qui s'insère dans un diagramme commutatif 
\begin{equation}
\xymatrix{
{F'|X}\ar[r]^u&{F|X}\ar[r]^v&{F''|X}\\
{(A|X)^m}\ar[r]^{w'}\ar[u]^{f'}&{(A|X)^n}\ar[u]_f&}
\end{equation}
où $w'$ est induit par $w$. On en déduit un morphisme $(A|X)$-linéaire et $\alpha$-surjectif  $\ker(f')\rightarrow \ker(f)$. 
Comme $F'$ est cohérent, $\ker(f')$ est de type $\alpha$-fini~; donc $\ker(f)$ est de type $\alpha$-fini d'après \ref{finita13}(i), ce qui prouve que 
$F$ est $\alpha$-cohérent.

\begin{cor}\label{finita18}
Soit $F$, $G$ deux $A$-modules $\alpha$-cohérents, $u\colon F\rightarrow G$ un morphisme $A$-linéaire.
Alors, $\ker(u)$, $\im(u)$ et $\coker(u)$ sont $\alpha$-cohérents.
\end{cor} 

En effet, $\im(u)$ est clairement de type $\alpha$-fini~; étant un sous-module d'un $A$-module $\alpha$-cohérent, il est $\alpha$-cohérent.  
La proposition résulte alors de \ref{finita17} appliqué aux suites exactes
$0\rightarrow \ker(u)\rightarrow F\rightarrow \im(u)\rightarrow 0$ et $0\rightarrow \im(u)\rightarrow G\rightarrow \coker(u)\rightarrow 0$.

\section{\texorpdfstring{Modules $\alpha$-cohérents sur un schéma}{Modules alpha-cohérents sur un schéma}}\label{afini}

\subsection{}\label{afini1}
Les hypothèses et notations de \ref{alpha1} sont en vigueur dans cette section. 
Pour tout $R$-schéma $X$, on désigne par $\bMod(\co_X)$ la catégorie des $\co_X$-modules de $X_\zar$ \eqref{notconv12}
et  par $\aMod(\co_X)$ la catégorie des $\alpha$-$\co_X$-modules, c'est-à-dire 
le quotient de $\bMod(\co_X)$ par la sous-catégorie pleine formée des $\co_X$-modules $\alpha$-nuls (cf. \ref{alpha20}).
Lorsque nous parlons de $\co_X$-module sans préciser, il est sous-entendu qu'il s'agit d'un $\co_X$-module de $X_\zar$.

\begin{prop}\label{afini2}
Soient $X$ un $R$-schéma, $F$ un $\co_X$-module. Les conditions suivantes sont équivalentes~:
\begin{itemize}
\item[{\rm (i)}] Il existe un $\co_X$-module quasi-cohérent $G$ et un isomorphisme de $\alpha$-$\co_X$-modules
$\alpha(F)\stackrel{\sim}{\rightarrow} \alpha(G)$. 
\item[{\rm (ii)}] Pour tout $x\in X$, il existe un voisinage ouvert $U$ de $x$ dans $X$, un $\co_U$-module 
quasi-cohérent $G$ et un isomorphisme de $\alpha$-$\co_U$-modules
$\alpha(F|U)\stackrel{\sim}{\rightarrow} \alpha(G)$. 
\item[{\rm (iii)}] Le $\co_X$-module $\fm\otimes_RF$ est quasi-cohérent. 
\end{itemize}
\end{prop}

On a clairement (i)$\Rightarrow$(ii). Montrons (ii)$\Rightarrow$(iii). 
La question étant locale, on peut supposer la condition (i) remplie. D'après \ref{alpha21}(iv),
il existe donc un $\co_X$-module quasi-cohérent $G$ et un $\alpha$-isomorphisme $u\colon\fm\otimes_RF\rightarrow G$. 
Comme $\fm\otimes_R\fm\simeq \fm$, le morphisme $\fm\otimes_RF\rightarrow \fm\otimes_RG$ induit par $u$
est un isomorphisme en vertu de \ref{alpha21}(iii); d'où la condition (iii). 
Enfin, l'implication (iii)$\Rightarrow$(i) est immédiate puisque le morphisme canonique $\fm\otimes_RF\rightarrow F$
est un $\alpha$-isomorphisme. 

\begin{defi}\label{afini3}
Soient $X$ un $R$-schéma, $F$ un $\co_X$-module. 
On dit que $F$ est {\em $\alpha$-quasi-cohérent} s'il remplit les conditions de \ref{afini2}.
\end{defi}

\begin{lem}\label{afini4}
Soient $X$ un $R$-schéma cohérent ({\em i.e.}, quasi-compact et quasi-séparé), $\gamma\in R$, 
$F$ un $\co_X$-module quasi-cohérent de type $\gamma$-fini \eqref{finita9}.
Alors, il existe un $\co_X$-module quasi-cohérent de présentation finie $G$ et un morphisme $\co_X$-linéaire
$f\colon G\rightarrow F$ de conoyau annulé par $\gamma$.
\end{lem}

En effet, d'après (\cite{ega1n} 6.9.12), $F$ est isomorphe à la limite inductive d'un système inductif filtrant 
de $\co_X$-modules de présentation finie $(F_i)_{i\in I}$. 
Pour tout $i\in I$, notons $f_i\colon F_i\rightarrow F$ le morphisme canonique. 
Montrons qu'il existe $i\in I$ tel que le conoyau de $f_i$ soit annulé par $\gamma$. 
Comme $X$ est quasi-compact et que la catégorie $I$ est filtrante, on peut supposer qu'il existe
un morphisme $\co_X$-linéaire $u\colon E\rightarrow F$ avec $E$ un $\co_X$-module libre de type fini, 
dont le conoyau est annulé par $\gamma$. Il existe $i\in I$ et un morphisme $\co_X$-linéaire $u_i\colon E\rightarrow F_i$
tel que $u=f_i\circ u_i$. Par suite, le conoyau de $f_i$ est annulé par $\gamma$; d'où la proposition.

\begin{prop}\label{afini5}
Soient $X$ un $R$-schéma cohérent ({\em i.e.}, quasi-compact et quasi-séparé), $F$ un $\co_X$-module $\alpha$-quasi-cohérent.
Pour que $F$ soit de type $\alpha$-fini (resp. présentation $\alpha$-finie) \eqref{finita9}, 
il faut et il suffit que pour tout $\gamma\in \fm$, 
il existe un $\gamma$-isomorphisme $f\colon G\rightarrow F$ avec $G$ un $\co_X$-module quasi-cohérent de type fini 
(resp. de présentation finie).
\end{prop}

Considérons d'abord l'assertion non-respée.  
La condition est clairement suffisante. Inversement, supposons $F$ de type $\alpha$-fini et 
montrons que la condition est satisfaite. Soit $\gamma\in \fm$. 
On peut se borner au cas où $F$ est quasi-cohérent \eqref{afini2}. D'après \ref{afini4}, 
il existe un $\co_X$-module quasi-cohérent de présentation finie $G'$ et un morphisme $\co_X$-linéaire
$f'\colon G'\rightarrow F$ de conoyau annulé par $\gamma$. 
Le morphisme $f\colon \im(f')\rightarrow F$ induit par $f'$ satisfait alors à la condition pour $\gamma$.  

Considérons ensuite l'assertion respée. La condition est clairement suffisante. 
Inversement, supposons $F$ de présentation $\alpha$-finie et montrons que la condition est satisfaite. Soit $\gamma\in \fm$. 
On peut se borner au cas où $F$ est quasi-cohérent \eqref{afini2}. D'après \ref{afini4}, 
il existe un $\co_X$-module quasi-cohérent de présentation finie $G'$ et un morphisme $\co_X$-linéaire
$f'\colon G'\rightarrow F$ de conoyau annulé par $\gamma$. En vertu de \ref{finita12}(iii), 
$\ker(f')$ est quasi-cohérent de type $\gamma^{12}$-fini. D'après \ref{afini4}, il existe 
un $\co_X$-module quasi-cohérent de présentation finie $G''$ et un morphisme $f''\colon G''\rightarrow G'$ tel que 
$f'\circ f''=0$ et que $\ker(f')/\im(f'')$ soit annulé par $\gamma^{12}$. Le morphisme $f\colon G'/\im(f'')\rightarrow F$ induit par $f'$
répond alors à la condition pour $\gamma^{12}$. 

\begin{cor}\label{afini10}
Soient $X$ un $R$-schéma affine d'anneau $A$, $F$ un $\co_X$-module $\alpha$-quasi-cohérent. 
Pour que $F$ soit de type $\alpha$-fini (resp. de présentation $\alpha$-finie), 
il faut et il suffit que le $A$-module $\Gamma(X,\fm\otimes_RF)$
soit  de type $\alpha$-fini (resp. de présentation $\alpha$-finie).
\end{cor}

Compte tenu de \ref{finita16}, 
on peut se réduire au cas où $F$ est quasi-cohérent~; il faut alors démontrer que $F$ est de type $\alpha$-fini 
(resp. de présentation $\alpha$-finie)
si et seulement si le $A$-module $\Gamma(X,F)$ est de type $\alpha$-fini (resp. de présentation $\alpha$-finie). 
Cette assertion résulte de \ref{afini5}.

\begin{prop}\label{afini11}
Soient $f\colon X'\rightarrow X$ un $R$-morphisme fidèlement plat et quasi-compact
de $R$-schémas, $F$ un $\co_X$-module $\alpha$-quasi-cohérent.
Pour que $F$ soit de type $\alpha$-fini (resp. de présentation $\alpha$-finie), il faut et il suffit
que son image inverse $F'$ sur $X'$ le soit. 
\end{prop}

En effet, la condition est nécessaire en vertu de \ref{finita14}. Montrons qu'elle est suffisante. 
On peut supposer $X$ affine. Remplaçant $X$ par une somme d'ouverts affines qui recouvrent $X'$, 
on est ramené au cas où $X'$ est également affine. L'assertion résulte alors de \ref{afini10} et (\cite{agt} V.8.1 et V.8.5).

\begin{prop}\label{afini6}
Soient $X$ un $R$-schéma cohérent, 
$F$ un $\co_X$-module $\alpha$-quasi-cohérent. Supposons que l'anneau $\co_X$ soit cohérent.
Alors, les conditions suivantes sont équivalentes~:
\begin{itemize}
\item[{\rm (i)}] $F$ est $\alpha$-cohérent \eqref{finita9};
\item[{\rm (ii)}] $F$ est de présentation $\alpha$-finie~;
\item[{\rm (iii)}] pour tout $\gamma\in \fm$, il existe un $\gamma$-isomorphisme 
$G\rightarrow F$ avec $G$ un $\co_X$-module cohérent.
\end{itemize}
\end{prop}

On a clairement (i)$\Rightarrow$(ii). L'implication (ii)$\Rightarrow$(iii) résulte de \ref{afini5}. Montrons (iii)$\Rightarrow$(i). 
Supposons la condition (iii) satisfaite. D'après \ref{afini5}, $F$ est de type $\alpha$-fini. Il suffit donc de montrer que pour tout
morphisme $\co_X$-linéaire $f\colon E\rightarrow F$ où $E$ est un $\co_X$-module libre de type fini, $\ker(f)$ est de type $\alpha$-fini. 
Soient $\gamma\in \fm$, $g\colon G\rightarrow F$ un $\gamma$-isomorphisme tel que $G$ soit un $\co_X$-module cohérent.
Il existe un morphisme $\co_X$-linéaire $f'\colon E\rightarrow G$ tel que $g\circ f'=\gamma f$. 
Le conoyau de l'injection canonique $\ker(f')\rightarrow \ker(\gamma f)$ est annulé par $\gamma$. 
Par ailleurs, la multiplication par $\gamma$ dans $E$ induit un morphisme $\ker(\gamma f)\rightarrow \ker(f)$
dont le conoyau est clairement annulé par $\gamma$. On en déduit un morphisme $\co_X$-linéaire 
$\ker(f')\rightarrow \ker(f)$ dont le conoyau est annulé par $\gamma^2$. Comme $G$ est cohérent, $\ker(f')$ est de type fini. 
Par suite, $\ker(f)$ est de type $\gamma^2$-fini et donc de type $\alpha$-fini.

\begin{cor}\label{afini12}
Soient $f\colon X'\rightarrow X$ un $R$-morphisme fidèlement plat et quasi-compact
de $R$-schémas, $F$ un $\co_X$-module $\alpha$-quasi-cohérent.
Supposons que les anneaux $\co_X$ et $\co_{X'}$ soient cohérents. Pour que $F$ soit $\alpha$-cohérent, 
il faut et il suffit que son image inverse $F'$ sur $X'$ le soit.
\end{cor}

Cela résulte de \ref{afini6} et \ref{afini11}.

\begin{defi}[\cite{egr1} 1.4.1]\label{afini7}
On dit qu'un anneau $A$ est {\em universellement cohérent} 
si l'anneau de polynômes $A[t_1,\dots,t_n]$ est cohérent pour tout entier $n\geq 0$. 
\end{defi}

\begin{teo}\label{afini8}
Soient $Y$ un $R$-schéma, $f\colon X\rightarrow Y$ un morphisme propre de présentation finie,
$F$ un $\co_X$-module $\alpha$-quasi-cohérent et $\alpha$-cohérent.
Supposons que $Y$ admette un recouvrement par des ouverts affines $(U_i)_{i\in I}$
tels que $\Gamma(U_i,\co_Y)$ soit universellement cohérent pour tout $i\in I$. 
Alors, pour tout entier $q\geq 0$, $\rR^qf_*(F)$ est un $\co_Y$-module $\alpha$-quasi-cohérent et $\alpha$-cohérent. 
\end{teo}

On notera d'abord que les anneaux $\co_X$ de $X_\zar$ et $\co_Y$ de $Y_\zar$ sont cohérents d'après (\cite{egr1} 1.4.2 et 1.4.3).
Par ailleurs, la question étant locale sur $Y$, on peut le supposer affine, de sorte que $X$ est cohérent. 
Le morphisme canonique $\rR^qf_*(\fm\otimes_RF)\rightarrow \rR^qf_*(F)$ est un $\alpha$-isomorphisme \eqref{alpha3}. 
Comme $\fm\otimes_RF$ est quasi-cohérent, $\rR^qf_*(\fm\otimes_RF)$ 
est quasi-cohérent et donc $\rR^qf_*(F)$ est $\alpha$-quasi-cohérent.
Par ailleurs, pour tout $\gamma\in \fm$, il existe un $\gamma$-isomorphisme 
$u\colon G\rightarrow F$ avec $G$ un $\co_X$-module cohérent d'après \ref{afini6}. 
Il existe donc un morphisme $\co_X$-linéaire $v\colon F\rightarrow G$ tel que $u\circ v=\gamma^2\cdot \id_F$ et 
$v\circ u=\gamma^2\cdot \id_G$ \eqref{alpha3}. En vertu  de (\cite{egr1} 1.4.8), le $\co_Y$-module
$\rR^qf_*(G)$ est cohérent. Par suite, $\rR^qf_*(F)$ est $\alpha$-cohérent \eqref{afini6}.

\begin{cor}\label{afini9}
Soient $Y$ un $R$-schéma affine d'anneau $A$, $f\colon X\rightarrow Y$ un morphisme propre de présentation finie,
$F$ un $\co_X$-module $\alpha$-quasi-cohérent et $\alpha$-cohérent.
Supposons que $Y$ admette un recouvrement par des ouverts affines $(U_i)_{i\in I}$
tels que $\Gamma(U_i,\co_Y)$ soit universellement cohérent pour tout $i\in I$. 
Alors, pour tout entier $q\geq 0$, le $A$-module $\rH^q(X,F)$ est $\alpha$-cohérent. 
\end{cor}

On notera d'abord que l'anneau $\co_X$ de $X_\zar$ est cohérent (\cite{egr1} 1.4.2 et 1.4.3).
Pour tout $\gamma\in \fm$, il existe un $\gamma$-isomorphisme $u\colon G\rightarrow F$ avec $G$ un $\co_X$-module cohérent d'après \ref{afini6}. 
Il existe donc un morphisme $\co_X$-linéaire $v\colon F\rightarrow G$ tel que $u\circ v=\gamma^2\cdot \id_F$ et 
$v\circ u=\gamma^2\cdot \id_G$ \eqref{alpha3}. En vertu  de (\cite{egr1} 1.4.8 et 1.4.3), le $A$-module
$\rH^q(X,G)$ est cohérent. Par suite, $\rH^q(X,F)$ est $\alpha$-cohérent \eqref{afini6}. 

\begin{rema}\label{afini99}
On notera que (\cite{egr1} 1.4.8), utilisé dans les preuves de \ref{afini8} et \ref{afini9}, 
est un corollaire d'un théorème de Kiehl (\cite{kiehl} 2.9'a) (cf. \cite{egr1} 1.4.7)
\end{rema}

\section{\texorpdfstring{Rappel sur les extensions $\alpha$-étales}{Extensions alpha-étales}}\label{aet}

\subsection{}\label{aet1}
Les hypothèses et notations de \ref{alpha1} sont en vigueur dans cette section. On renvoie à (\cite{agt} V) pour 
les définitions et les principales propriétés des notions utilisées en {\em $\alpha$-algèbre} 
(ou presque-algèbre dans {\em loc. cit.})~: 
modules {\em $\alpha$-projectifs} ({\em loc. cit.} V.3.2), {\em $\alpha$-plats} ({\em loc. cit.} V.6.1) et {\em $\alpha$-fidèlement plats} ({\em loc. cit.} V.6.4) et {\em $\alpha$-rang} d'un module $\alpha$-projectif de type $\alpha$-fini
({\em loc. cit.} V.5.14).

\begin{defi}[\cite{agt} V.7.1]\label{aet2}
Soit $f\colon A\rightarrow B$ un homomorphisme de $R$-algèbres. On dit que $f$ est {\em $\alpha$-étale}
si les conditions suivantes sont satisfaites~:
\begin{itemize}
\item[(i)] $B$ est $\alpha$-projectif de type $\alpha$-fini et de $\alpha$-rang fini en tant que $A$-module;
\item[(ii)] $B$ est $\alpha$-projectif en tant que $B\otimes_AB$-module.
\end{itemize}  
\end{defi}

On renvoie à (\cite{agt} V.7) pour les principales propriétés de cette notion.

\begin{prop}[\cite{agt} V.7.11]\label{aet3}
Soit $f\colon A\rightarrow B$ un homomorphisme $\alpha$-étale de $R$-algèbres tel que pour tout
(ou ce qui revient au même pour un) $\alpha\in \Lambda^+$, $\pi^\alpha$ ne soit un diviseur de zéro dans $A$,  
et que $\Spec(A)$ soit normal et localement irréductible {\rm (\cite{agt} III.3.1)}. 
Notons $A_\pi$ (resp. $B_\pi$) l'anneau des fractions 
de $A$ (resp. $B$) dont les dénominateurs sont de la forme $\pi^\alpha$ pour $\alpha\in \Lambda^+$, et $B'$ la clôture 
intégrale de $A$ dans $B_\pi$. Alors, l'homomorphisme canonique $B\rightarrow B_\pi$ se factorise en
\begin{equation}
B\rightarrow \Hom_R(\fm,B')\rightarrow \Hom_R(\fm,B_\pi) \stackrel{\sim}{\rightarrow} B_\pi,
\end{equation}
où la première flèche est un $\alpha$-isomorphisme et les autres flèches sont les homomorphismes canoniques  \eqref{alpha5c}. En particulier, l'homomorphisme $A\rightarrow B'$ est $\alpha$-étale.  
\end{prop} 

En effet, on peut se borner au cas où $A$ est intègre (\cite{agt} V.7.4(2) et V.7.8), auquel cas l'assertion est démontrée dans
(\cite{agt} V.7.11).

\begin{defi}\label{aet4}
Soient $A$ une $R$-algèbre, $G$ un groupe fini, $B$ une $A[G]$-algèbre.
On dit que $B$ est un {\em $\alpha$-$G$-torseur sur $A$} s'il existe un homomorphisme $\alpha$-fidèlement
plat $A\rightarrow A'$ et un homomorphisme de $A'[G]$-algèbres 
\begin{equation}
B\otimes_AA'\rightarrow \prod_GA',
\end{equation}
où $G$ agit sur le terme de droite par multiplication à droite sur lui même, 
qui est un $\alpha$-isomorphisme.
\end{defi}

On renvoie à (\cite{agt} V.12) pour les principales propriétés de cette notion. 

\begin{prop}[\cite{agt} V.12.5] \label{aet5}
Soient $A$ une $R$-algèbre, $G$ un groupe fini, $B$ un $\alpha$-$G$-torseur sur $A$, $M$ un $B$-module
muni d'une action semi-linéaire de $G$. Alors, le morphisme canonique 
$B\otimes_AM^G\rightarrow M$ est un $\alpha$-isomorphisme.
\end{prop}

\begin{cor}[\cite{agt} V.12.6] \label{aet6}
Soient $A$ une $R$-algèbre, $G$ un groupe fini, $B$ un $\alpha$-$G$-torseur sur $A$. 
Alors, l'homomorphisme naturel $A\rightarrow B^G$ est un $\alpha$-isomorphisme.
\end{cor}

\begin{prop}[\cite{agt} V.12.8]\label{aet7}
Soient $A$ une $R$-algèbre, $G$ un groupe fini, $B$ un $\alpha$-$G$-torseur sur $A$, $M$
un $B$-module muni d'une action semi-linéaire de $G$, $\Tr_G$ le morphisme $A$-linéaire
de $M$ défini par $\sum_{\sigma\in G}\sigma$. Alors, pour tout $q\geq 1$,
$\rH^q(G,M)$ et $M^G/\Tr_G(M)$ sont $\alpha$-nuls. 
\end{prop}

\begin{prop}[\cite{agt} V.12.9]\label{aet8}
Soit $A$ une $R$-algèbre telle que pour tout
(ou ce qui revient au même pour un) $\alpha\in \Lambda^+$, $\pi^\alpha$ ne soit un diviseur de zéro dans $A$,  
et que $\Spec(A)$ soit normal et localement irréductible {\rm (\cite{agt} III.3.1)}. 
Notons $A_\pi$  l'anneau des fractions de $A$ dont les dénominateurs sont de la forme $\pi^\alpha$ pour $\alpha\in \Lambda^+$. Soient $U$ un ouvert de $\Spec(A_\pi)$, $V\rightarrow U$ un morphisme étale fini et galoisien de groupe 
$G$, $Y$ la fermeture intégrale de $\Spec(A)$ dans $V$, $B=\Gamma(Y,\co_Y)$. 
Si $B$ est une extension $\alpha$-étale de $A$, alors $B$ est un $\alpha$-$G$-torseur au-dessus de $A$. 
\end{prop}

On peut se borner au cas où $U$ est irréductible (\cite{agt} V.7.4(3)). 
Il suffit alors de calquer la preuve de (\cite{agt} V.12.9) en remplaçant (\cite{agt} V.7.11) par \ref{aet3}.

\section{\texorpdfstring{Modules de type $\alpha$-fini sur un anneau de valuation non discrète de hauteur $1$}
{Modules de type alpha-fini sur un anneau de valuation non discrète de hauteur 1}}\label{mptf}

\subsection{}\label{mptf1}
Dans cette section, $R$ désigne un anneau muni d'une valuation non discrète de hauteur $1$, 
$v\colon R\rightarrow \mR$. On note $\fm$ l'idéal maximal de $R$ et on pose $\Lambda=v(R)$ et $\Lambda^+=v(\fm)$. 
Pour tout $\varepsilon\in \Lambda^+$, on note $\fm_\varepsilon$ l'idéal de $R$ formé des éléments $x\in R$ tels que $v(x)\geq \varepsilon$.
On observera que les hypothèses de \ref{alpha1} sont satisfaites. Il est donc loisible de considérer dans ce contexte
les notions introduites dans les sections \ref{alpha}, \ref{finita} et \ref{afini}. Nous considérons toujours $R$ 
comme muni de la topologie définie par l'idéal $\fm_\varepsilon$ pour un élément quelconque $\varepsilon\in \Lambda^+$; 
c'est un anneau prévaluatif séparé (\cite{egr1} 1.9.8). 

\begin{prop}[\cite{egr1} 1.12.15]\label{mptf2}
Si $R$ est complet et séparé, toute $R$-algèbre topologiquement de présentation finie est un anneau 
universellement cohérent \eqref{afini7}. 
\end{prop}

\subsection{}\label{mptf3}
Suivant (\cite{scholze1} 2.2), pour tous $R$-modules $M$ et $N$ et tout $\gamma \in \fm$, 
on dit que $M$ et $N$ sont {\em $\gamma$-équivalents} et on note $M\approx_\gamma N$,
s'il existe deux $R$-morphismes $f\colon M\rightarrow N$ et $g\colon N\rightarrow M$ 
tels que $f\circ g=\gamma\cdot \id_N$ et $g\circ f=\gamma\cdot \id_M$. 
On dit que $M$ et $N$ sont {\em $\alpha$-équivalents} et on note $M\approx N$
si $M$ et $N$ sont $\gamma$-équivalents pour tout $\gamma\in \fm$. 
Les relations $\approx_\gamma$ et $\approx$ sont clairement symétriques. 
Pour tous $R$-modules $M$, $N$ et $P$ et tous $\gamma,\delta \in \fm$, si $M\approx_\gamma N$
et $N\approx_\delta P$, alors $M\approx_{\gamma\delta} P$. En particulier, la relation $\approx$ est transitive. 

\subsection{}\label{mptf4}
Pour tout $\varepsilon\in \mR_{\geq 0}$, on désigne par $I_\varepsilon$ l'idéal de $R$ formé des éléments 
$x\in R$ tels que $v(x)>\varepsilon$. Alors $I_\varepsilon\approx R$. 
En effet, pour tout $\gamma \in \fm$, il existe $\gamma_1,\gamma_2\in R$ tels que $\gamma\gamma_1=\gamma_2$
et $v(\gamma_1)\leq \varepsilon <v(\gamma_2)$. Les morphismes 
$R\rightarrow I_\varepsilon$ et $I_\varepsilon\rightarrow R$ définis par les multiplications par $\gamma_2$ et 
$\gamma_1^{-1}$, respectivement, montrent que $I_\varepsilon\approx_\gamma R$. 
Par ailleurs, comme $\fm\cdot I_\varepsilon=I_\varepsilon$,
si $\varepsilon\not\in \Lambda$, $R$ et $I_\varepsilon$ ne sont pas $\alpha$-isomorphes en vertu de \ref{alpha7} et \eqref{alpha8a}. 

\begin{remas}\label{mptf5}  
(i)\ Deux $R$-modules $\alpha$-isomorphes sont $\alpha$-équivalents \eqref{alpha3}. 
La réciproque n'est pas vraie \eqref{mptf4}. 

(ii)\ Pour qu'un $R$-module $M$ soit {\em de type $\alpha$-fini} (resp. {\em de présentation $\alpha$-finie}) 
il faut et il suffit que pour tout $\gamma\in \fm$, il existe un $R$-module de type fini (resp. de présentation finie) 
$N$ tel que $M\approx_\gamma N$, d'après \ref{finita2}(i) et \ref{alpha3}.
\end{remas}

\begin{prop}[\cite{scholze1} 2.6]\label{mptf62}
Tout $R$-module de type $\alpha$-fini est de présentation $\alpha$-finie. 
\end{prop}

\begin{cor}\label{mptf63}
Tout sous-$R$-module d'un $R$-module de type $\alpha$-fini est de type $\alpha$-fini. 
\end{cor}

Cela résulte de \ref{mptf62} et \ref{finita13}.

\begin{teo}[\cite{scholze1} 2.5]\label{mptf7}
Soit $M$ un $R$-module de type $\alpha$-fini. Alors, il existe une unique suite décroissante de nombres réels positifs 
$(\varepsilon_i)_{i\geq 1}$, tendant vers $0$, et un unique entier $n\geq 0$ tels que 
\begin{equation}\label{mptf7a}
M\approx R^n\oplus (\oplus_{i\geq 1}R/I_{\varepsilon_i}). 
\end{equation}
\end{teo}

\subsection{}\label{mptf8}
On désigne par $\fS$ le $\mR$-espace vectoriel des suites de nombres réels 
$(\varepsilon_i)_{i\geq 1}$, tendant vers $0$, et par $\fS_{\geq}$ le sous-ensemble formé des suites décroissantes.
Pour tout $\varepsilon=(\varepsilon_i)_{i\geq 1}\in \fS$, on pose 
\begin{equation}\label{mptf8a}
\|\varepsilon\|=\sup_{i\geq 1}(|\varepsilon_i|).
\end{equation}
On munit $\fS_{\geq}$ de la relation d'ordre définie pour toutes suites
$\varepsilon=(\varepsilon_i)_{i\geq 1}$ et $\varepsilon'=(\varepsilon'_i)_{i\geq 1}$ de $\fS_\geq$, par $\varepsilon\geq \varepsilon'$
si pour tout $i\geq 1$, $\varepsilon_1+\dots+\varepsilon_i\geq \varepsilon'_1+\dots+\varepsilon'_i$.

\subsection{}\label{mptf9}
En vertu de \ref{mptf7}, on peut associer à tout $R$-module de torsion de type $\alpha$-fini $M$ 
une unique suite $\varepsilon_{M}=(\varepsilon_{M,i})_{i\geq 1}$ de $\fS_{\geq}$ telle que 
\begin{equation}\label{mptf9a}
M\approx \oplus_{i\geq 1}R/I_{\varepsilon_{M,i}}. 
\end{equation}
On appelle {\em longueur} de $M$ et l'on note $\lambda(M)$ l'élément 
\begin{equation}\label{mptf9b}
\lambda(M)=\sum_{i\geq 1}\varepsilon_{M,i}\in \mR_{\geq 0}\cup\{\infty\}.
\end{equation}
D'après (\cite{scholze1} 2.11), la suite $\varepsilon_M$ et la longueur $\lambda(M)$ vérifient les propriétés suivantes~:

\begin{itemize}
\item[(i)] Pour tout $R$-module de torsion et de présentation finie $M$, il existe des éléments  
$\gamma_i\in R$ $(1\leq i\leq n)$ tels que $v(\gamma_1)\geq v(\gamma_2)\geq \dots\geq v(\gamma_n)$ et un $R$-isomorphisme
\begin{equation}
M\stackrel{\sim}{\rightarrow} R/\gamma_1R\oplus R/\gamma_2R\oplus \dots \oplus R/\gamma_nR.
\end{equation}
On a alors $\varepsilon_{M,i}=v(\gamma_i)$ pour tout $1\leq i\leq n$ et $\varepsilon_{M,i}=0$ pour tout $i>n$.
\item[(ii)] Pour que deux $R$-modules de torsion de type $\alpha$-fini $M$ et $N$ 
soient $\gamma$-équivalents (pour $\gamma\in \fm$), il faut et il suffit que 
$\|\varepsilon_M-\varepsilon_N\|\leq v(\gamma)$.
En particulier, pour qu'un $R$-module de torsion de type $\alpha$-fini $M$ soit $\alpha$-nul, il faut et il suffit que 
la suite $\varepsilon_M$ soit nulle. 
\item[(iii)] Pour tous $R$-modules de torsion de type $\alpha$-fini $M$ et $M'$, si $M'$ est un sous-quotient de $M$
alors $\varepsilon_{M',i}\leq \varepsilon_{M,i}$ pour tout $i\geq 1$. 
\item[(iv)] Pour toute suite exacte $0\rightarrow M'\rightarrow M\rightarrow M''\rightarrow 0$ de $R$-modules de torsion 
de type $\alpha$-fini, on a $\varepsilon_M\leq \varepsilon_{M'}+\varepsilon_{M''}$ et $\lambda(M)= \lambda(M')+\lambda(M'')$. 
De plus, si $\varepsilon_M=\varepsilon_{M'}$ (resp. $\varepsilon_M=\varepsilon_{M''}$), alors $M''$ (resp. $M'$) est $\alpha$-nul.
\end{itemize}

\begin{lem}\label{mptf12}
Soient $M$ un $R$-module, $\varpi$ un élément non nul de $\fm$.  On pose $\hR=\underset{\longleftarrow}{\lim} \ R/\varpi^nR$ 
et $\hM=\underset{\longleftarrow}{\lim} \ M/\varpi^nM$, et on suppose que le $R$-module $M/\varpi M$ est de type $\alpha$-fini. Alors, 
le $\hR$-module $\hM$ est de type $\alpha$-fini.
\end{lem}

On notera d'abord que $\hR$ et $\hM$ sont complets et séparés pour les topologies $\varpi$-adiques et pour tout entier $n\geq1$, que le morphisme canonique 
$h_n\colon \hM\rightarrow M/\varpi^n M$ est surjectif (\cite{ac} III §~2.11 prop.~11) et que le morphisme canonique $\hR/\varpi^n \hR\rightarrow R/\varpi^n R$
est un isomorphisme (\cite{ac} III §~2.11 cor.~1 à prop.~11). Soient $d$ un entier $\geq 1$, $\varphi\colon \hR^d\rightarrow \hM$
un morphisme $\hR$-linéaire, $\varphi_1\colon (R/\varpi R)^d\rightarrow M/\varpi M$ le morphisme induit par $\varphi$. 
Supposons que $\coker(\varphi_1)$ soit annulé par un élément $\gamma\in \fm$ tel que $v(\gamma)<v(\varpi)$. Alors, la suite 
\begin{equation}\label{mptf12a}
\hR^d\stackrel{\varphi}{\rightarrow} \hM \stackrel{u}{\rightarrow}\coker(\varphi_1)\rightarrow 0,
\end{equation}
où $u$ est le morphisme canonique, est exacte.
En effet, $u$ est surjectif puisque $h_1$ est surjectif. Soit $x\in \hM$ tel que $u(x)=0$.
On montre par récurrence que pour tout entier $n\geq 1$, qu'il existe $y_n\in \hR^d$ et $x_n\in \hM$ tel que $x=\varphi(y_n)+\gamma^{1-n}\varpi^n x_n$
et $y_{n+1}-y_n\in \gamma^{-n}\varpi^n \hR^d$. On en déduit qu'il existe $y\in \hR^d$ tel que $x=\varphi(y)$; d'où l'assertion recherchée.
Par suite, $\coker(\varphi)$ est annulé par $\gamma$, ce qui implique la proposition.

\subsection{}\label{mptf10}
On suppose dans la suite de cette section que $R$ est de caractéristique $p$, et on note $\varphi$ son endomorphisme de Frobenius. 
Soit $\varpi$ un élément non nul de $\fm$. Pour tout entier $n\geq 1$, on pose $R_n=R/\varpi^n R$, $\bvR=(R_n)_{n\geq 1}$
et $\bvR^{(p)}=(R_{pn})_{n\geq 1}$. On considère $\bvR$ et $\bvR^{(p)}$ comme des anneaux du topos $\Ens^{\mN^\circ}$ (\cite{agt} III.7.7). 
L'endomorphisme $\varphi$ induit un homomorphisme d'anneaux que l'on note
\begin{equation}\label{mptf10a}
\bvvarphi\colon \bvR\rightarrow \bvR^{(p)}.
\end{equation}
Pour tout $\bvR$-module $M=(M_n)_{n\geq 1}$, on pose $M_0=0$, $M[-1]=(M_n)_{n\geq 0}$ qui est naturellement un $\bvR$-module, 
et $M^{(p)}=(M_{pn})_{n\geq 1}$ qui est naturellement un $\bvR^{(p)}$-module. Les correspondances $M\mapsto M[-1]$ et $M\mapsto M^{(p)}$ sont clairement
fonctorielles, et on a un morphisme $\bvR$-linéaire canonique fonctoriel $M\rightarrow M[-1]$.

\begin{lem}[\cite{faltings2} page 224, \cite{scholze1} 2.12]\label{mptf11}
Conservons les hypothèses de \ref{mptf10}, soient de plus $M=(M_n)_{n\geq 1}$ un $\bvR$-module, 
$u\colon M[-1]\rightarrow M$ un morphisme $\bvR$-linéaire, 
\begin{equation}\label{mptf11a}
\phi\colon M\otimes_{\bvR}\bvR^{(p)}\stackrel{\sim}{\rightarrow} M^{(p)}
\end{equation}
un isomorphisme $\bvR^{(p)}$-linéaire. On pose $\hR=\underset{\longleftarrow}{\lim} \ R_n$ et $\hM=\underset{\longleftarrow}{\lim} \ M_n$. On note encore $\varphi$ 
l'endomorphisme de Frobenius de $\hR$ et on désigne par
\begin{equation}\label{mptf11b}
\uphi\colon \hM\rightarrow\hM
\end{equation}
le morphisme $\hR$-semi-linéaire induit par $\phi$. On suppose que 
\begin{itemize}
\item[{\rm (a)}] le $R$-module $M_1$ est de type $\alpha$-fini~;
\item[{\rm (b)}] le composé $M[-1]\stackrel{u}{\rightarrow}M\rightarrow M[-1]$, où la seconde flèche est le morphisme canonique, 
est induit par la multiplication par $\varpi$.
\item[{\rm (c)}] pour tout entier $n\geq 0$, la suite de $R$-modules
\begin{equation}\label{mptf11c}
\xymatrix{M_1\ar[r]\ar[r]^-(0.5){u^{\circ n}}&{M_{n+1}}\ar[r]&{M_n}}
\end{equation}
est exacte au centre.  
\end{itemize}
Alors,
\begin{itemize}
\item[{\rm (i)}] Pour tout entier $n\geq 1$, le morphisme canonique $M_{n+1}\rightarrow M_n$ est $\alpha$-surjectif
et le morphisme canonique $\hM/\varpi^n \hM \rightarrow M_n$ est un $\alpha$-isomorphisme. 
\item[{\rm (ii)}] Le $\hR$-module $\hM$ est de type $\alpha$-fini. 
\item[{\rm (iii)}] Si le corps des fractions de $R$ est algébriquement clos, il existe un entier $d\geq 0$ et un morphisme $\hR$-linéaire 
\begin{equation}\label{mptf11d}
\fm\otimes_R\hM\rightarrow \hR^d,
\end{equation}
compatible aux morphismes $\uphi$ et $\varphi$, et qui est un $\alpha$-isomorphisme.
\end{itemize}
\end{lem}

En effet, il résulte de (c), \ref{finita13} et \ref{mptf9}(iii)-(iv) que pour tout entier $n\geq 1$, le $R$-module $M_n$ est de type $\alpha$-fini  et on a 
$\varepsilon_{M_{n+1}}\leq \varepsilon_{M_n}+\varepsilon_{M_1}$. Par suite, $\varepsilon_{M_n}\leq n \varepsilon_{M_1}$. 
Par ailleurs, on a $\varepsilon_{M_{pn}}=p\varepsilon_{M_n}$ d'après \eqref{mptf11a}. On en déduit que $\varepsilon_{M_n}= n \varepsilon_{M_1}$.
Notons $M^*_{n+1}$ et $M'_n$ les images des première et seconde flèches de \eqref{mptf11c}, respectivement, de sorte qu'on a une suite exacte
$0\rightarrow M^*_{n+1}\rightarrow M_{n+1}\rightarrow M'_n\rightarrow 0$. D'après \ref{mptf9}, on a 
\begin{equation}\label{mptf11e}
(n+1) \varepsilon_{M_1}= \varepsilon_{M_{n+1}}\leq \varepsilon_{M^*_{n+1}}+\varepsilon_{M'_n}\leq \varepsilon_{M_1}+\varepsilon_{M_n}=(n+1) \varepsilon_{M_1}.
\end{equation}
Par suite, on a $\varepsilon_{M^*_{n+1}}=\varepsilon_{M_1}$ et $\varepsilon_{M'_n}=\varepsilon_{M_n}$. Il résulte alors de \ref{mptf9}(iv) que 
les morphismes canoniques $M_1\rightarrow M^*_{n+1}$ et $M'_n\rightarrow M_n$ sont des $\alpha$-isomorphismes, en particulier, la suite 
\begin{equation}\label{mptf11f}
0\longrightarrow M_1\stackrel{u^{\circ n}}{\longrightarrow} M_{n+1}\longrightarrow M_n\longrightarrow 0
\end{equation}
est $\alpha$-exacte. Pour tout entier $m\geq 1$, le diagramme suivant 
\begin{equation}\label{mptf11g}
\xymatrix{
M_1\ar[r]^-(0.5){u^{\circ m}}\ar@{=}[d]&M_{m+1}\ar[r]\ar[d]^-(0.5){u^{\circ n}}&{M_m}\ar[d]^-(0.5){u^{\circ n}}\\
M_1\ar[r]^-(0.5){u^{\circ (m+n)}}&M_{m+n+1}\ar[r]\ar[d]&{M_{m+n}}\ar[d]\\
&M_n\ar@{=}[r]&M_n}
\end{equation}
permet de déduire (par récurrence sur $m$) que la suite 
\begin{equation}\label{mptf11h}
0\longrightarrow M_m\stackrel{u^{\circ n}}{\longrightarrow} M_{n+m}\longrightarrow M_n\longrightarrow 0
\end{equation}
est $\alpha$-exacte. Le système projectif $(\alpha(M_n))_{n\geq 1}$ de $\aMod(R)$ vérifie donc la condition de Mittag-Leffler. Compte tenu de \ref{alpha10}
et (\cite{ega3} 0.13.2.2), on en déduit que la suite de $\hR$-modules
\begin{equation}\label{mptf11i}
0\longrightarrow \hM\stackrel{\varpi^n}{\longrightarrow} \hM\longrightarrow M_n\longrightarrow 0
\end{equation}
est $\alpha$-exacte~; d'où la proposition (i). Le morphisme
\begin{equation}
M^\dagger=\underset{\longleftarrow}{\lim} \ \hM/\varpi^n\hM \rightarrow \hM
\end{equation}
déduit par passage à la limite projective, est donc un $\alpha$-isomorphisme d'après \ref{alpha10}.
Il résulte de (i) et \ref{mptf12} que le $\hR$-module $M^\dagger$ est de type $\alpha$-fini~; d'où la proposition (ii).

Supposons enfin le corps des fractions $K$ de $R$ algébriquement clos. Il résulte du lemme de Krasner que le corps des fractions $\hK$ de $\hR$ 
est aussi algébriquement clos. Par suite, $\varphi\colon \hR\rightarrow \hR$ est un isomorphisme~; 
il en est alors de même de $\bvvarphi\colon \bvR\rightarrow \bvR^{(p)}$ \eqref{mptf10a} et donc de $\uphi\colon \hM\rightarrow \hM$ \eqref{mptf11b}. 
On note encore $\varphi$ l'endomorphisme de Frobenius de $\hK$. 
D'après (ii) et (\cite{katz} 4.1.1), il existe un entier $d\geq 0$ et un isomorphisme $\hK$-linéaire
\begin{equation}
f\colon \hM\otimes_\hR\hK\stackrel{\sim}{\rightarrow}\hK^d
\end{equation}
compatible à $\uphi$ et $\varphi$. Notons $N$ l'image canonique de $\hM$ dans $\hM\otimes_\hR\hK$.
Il résulte de la suite $\alpha$-exacte \eqref{mptf11i} que le morphisme canonique $\hM\rightarrow N$ est un $\alpha$-isomorphisme.
Par ailleurs, on a $(\uphi\otimes \varphi)(N)=N$. Comme $N$ est de type $\alpha$-fini d'après (ii), 
il existe $\gamma\in \fm$ tel que l'on ait $\gamma \hR^d\subset f(N)\subset \gamma^{-1} \hR^d$. Appliquant les puissances négatives de $\varphi$ à ces inclusions, 
on en déduit que $f(\fm\otimes_R N)\subset \hR^d$ et que le morphisme 
\begin{equation}
\fm\otimes_R N\rightarrow \hR^d
\end{equation}
induit par $f$ est un $\alpha$-isomorphisme~; d'où la proposition (iii).

\section{Topos co-évanescents}\label{rec}

\subsection{}\label{rec1}
Soient $X$ et $Y$ deux $\mU$-sites dans lesquels les 
limites projectives finies sont représentables, $f^+\colon X\rightarrow Y$ un foncteur continu et exact à gauche. 
On note $\tX$ et $\tY$ les topos des faisceaux de $\mU$-ensembles sur $X$ et $Y$, respectivement, 
$f\colon \tY\rightarrow \tX$ le morphisme de topos associé à $f^+$ et 
$e_X$ et $e_Y$ des objets finaux de $X$ et $Y$, respectivement, qui existent par hypothèse.

On désigne par $\cD_{f^+}$ la catégorie des paires $(U, V\rightarrow f^+(U))$, où $U$ est un objet de $X$ et  
$V\rightarrow f^+(U)$ est un morphisme de $Y$; un tel objet sera noté $(V\rightarrow U)$. 
Soient $(V\rightarrow U)$, $(V'\rightarrow U')$ deux objets de $\cD_{f^+}$. 
Un morphisme de $(V'\rightarrow U')$ dans $(V\rightarrow U)$
est la donnée de deux morphismes $V'\rightarrow V$ de $Y$ et $U'\rightarrow U$ de $X$, tels que le diagramme
\[
\xymatrix{
V'\ar[r]\ar[d]&{f^+(U')}\ar[d]\\
{V}\ar[r]&{f^+(U)}}
\] 
soit commutatif. Il résulte aussitôt des hypothèses que les limites projectives finies dans $\cD_{f^+}$ sont représentables.
On munit $\cD_{f^+}$ de la topologie {\em co-évanescente} (\cite{agt} VI.4.1), c'est-à-dire, 
la topologie engendrée par les recouvrements 
\[
\{(V_i\rightarrow U_i)\rightarrow (V\rightarrow U)\}_{i\in I}
\] 
des deux types suivants~:
\begin{itemize}
\item[$(\alpha)$] $U_i=U$ pour tout $i\in I$, et $(V_i\rightarrow V)_{i\in I}$ est une famille couvrante.
\item[$(\beta)$] $(U_i\rightarrow U)_{i\in I}$ est une famille couvrante, 
et pour tout $i\in I$, le morphisme canonique $V_i\rightarrow V\times_{f^+(U)}f^+(U_i)$ est un isomorphisme. 
\end{itemize}
Le site ainsi défini est appelé site {\em co-évanescent} associé au foncteur $f^+$; c'est un $\mU$-site. 
On désigne par $\hcD_{f^+}$ (resp. $\tcD_{f^+}$) la catégorie des préfaisceaux (resp. le topos des faisceaux) de 
$\mU$-ensembles sur $\cD_{f^+}$. Pour tout préfaisceau $F$ sur $\cD_{f^+}$, on note $F^\ra$ le faisceau associé.

Les foncteurs 
\begin{eqnarray}
\rp_1^+\colon X&\rightarrow& \cD_{f^+},\ \ \ U\mapsto (f^+(U)\rightarrow U),\label{rec1a}\\
\rp_2^+\colon Y&\rightarrow& \cD_{f^+},\ \ \ V\mapsto (V\rightarrow e_X),\label{rec1b}
\end{eqnarray}
sont exacts à gauche et continus. Ils définissent donc deux morphismes de topos 
\begin{eqnarray}
\rp_1\colon \tcD_{f^+}&\rightarrow& \tX,\label{rec1c}\\
\rp_2\colon \tcD_{f^+}&\rightarrow& \tY.\label{rec1d}
\end{eqnarray}
On a un $2$-morphisme 
\begin{equation}\label{rec1e}
\tau\colon f\rp_2\rightarrow \rp_1
\end{equation}
défini par le morphisme de foncteurs $(f\rp_2)_*\rightarrow \rp_{1*}$ suivant~: pour tout faisceau $F$ sur $\cD_{f^+}$ et 
tout $U\in \ob(X)$, 
\begin{equation}\label{rec1f}
f_*(\rp_{2*}(F))(U)\rightarrow \rp_{1*}(F)(U)
\end{equation}
est l'application canonique
\[
F(f^+(U)\rightarrow e_X)\rightarrow F(f^+(U)\rightarrow U).
\]

D'après (\cite{agt} VI.3.7 et VI.4.10), le quadruplet $(\tcD_{f^+},\rp_1,\rp_2,\tau)$ est universel dans le sens suivant~:
pour tout $\mU$-topos $T$ muni de deux morphismes de topos $a\colon T\rightarrow \tX$ et $b\colon T\rightarrow \tY$ 
et d'un $2$-morphisme $t\colon fb\rightarrow a$, il existe un triplet 
\[
(h\colon T\rightarrow \tcD_{f^+}, \alpha\colon \rp_1h\stackrel{\sim}{\rightarrow}a,  
\beta\colon \rp_2h\stackrel{\sim}{\rightarrow}b),
\] 
unique à isomorphisme unique près, formé d'un morphisme de topos $h$ et de deux isomorphismes 
de morphismes de topos $\alpha$ et $\beta$, tel que le diagramme 
\begin{equation}\label{rec1g}
\xymatrix{
{f\rp_2h}\ar[r]^{\tau*h}\ar[d]_{f*\beta}&{\rp_1h}\ar[d]^{\alpha}\\
{fb}\ar[r]^t&{a}}
\end{equation}
soit commutatif. Par suite, $\tcD_{f^+}$ ne dépend que du morphisme $f\colon \tY\rightarrow \tX$ (cf. \cite{agt} VI.3.9). 
On l'appelle le topos {\em co-évanescent} de $f$ et on le note aussi
$\tX\gtimes_{\tX}\tY$ (\cite{agt} VI.2.12). 

Le foncteur 
\begin{equation}\label{rec1h}
\Psi^+\colon \cD_{f^+}\rightarrow Y, \ \ \ (V\rightarrow U)\mapsto V
\end{equation}
est continu et exact à gauche (\cite{agt} VI.4.13). Il définit donc un morphisme de topos
\begin{equation}\label{rec1i}
\Psi\colon \tY\rightarrow \tX\gtimes_{\tX}\tY
\end{equation}
tel que $\rp_1\Psi=f$, $\rp_2\Psi=\id_\tY$ et $\tau*\Psi=\id_f$. Ces relations déterminent $\Psi $
compte tenu de la propriété universelle du topos co-évanescent.
\begin{equation}\label{rec1j}
\xymatrix{
&\tY\ar[dl]_f\ar[rd]^{\id_\tY}\ar[d]_{\Psi}&\\
{\tX}\ar[dr]&{\tX\gtimes_\tX\tY}\ar[r]^-(0.4){\rp_2}\ar[l]_-(0.4){\rp_1}&{\tY}\ar[ld]^f\\
&{\tX}&}
\end{equation}
Le morphisme $\Psi$ est appelé morphisme des {\em cycles co-proches}.
De la relation $\rp_{2*}\Psi_*=\id_{\tY}$, on obtient par adjonction un morphisme 
\begin{equation}\label{rec1k}
\rp_2^*\rightarrow \Psi_*.
\end{equation}
Celui-ci est un isomorphisme d'après (\cite{agt} VI.4.14), en particulier, le foncteur $\Psi_*$ est exact.

\subsection{}\label{rec4}
Conservons les hypothèses et notations de \ref{rec1}. Le foncteur 
\begin{equation}\label{rec4a}
\pi\colon \cD_{f^+}\rightarrow X, \ \ \ (V\rightarrow U)\mapsto U
\end{equation}
est fibrant, canoniquement clivé et normalisé (\cite{sga4} VI 7.1). Pour tout $U\in \ob(X)$, la catégorie fibre au-dessus de 
$U$ est canoniquement équivalente à la catégorie $Y_{/f^+(U)}$, et pour tout morphisme
$u\colon U'\rightarrow U$ de $X$, le foncteur image inverse $j_u^+\colon Y_{/f^+(U)}\rightarrow Y_{/f^+(U')}$ 
n'est autre que le foncteur de changement de base par $f^+(u)$. 
Munissant chaque fibre $Y_{/f^+(U)}$ de la topologie induite par celle de $Y$ \eqref{notconv6}, 
$\pi$ devient un $\mU$-site fibré, vérifiant les conditions de (\cite{agt} VI.5.1). On désigne par 
\begin{equation}\label{rec4b}
\cF_{f^+}\rightarrow X
\end{equation}
le $\mU$-topos fibré associé à $\pi$ (\cite{sga4} VI 7.2.6). La catégorie fibre de $\cF_{f^+}$ au-dessus de tout $U\in \ob(X)$
est canoniquement équivalente au topos $\tY_{/f^*(U)}$ \eqref{notconv6}, et le foncteur image inverse par tout  
morphisme $u\colon U'\rightarrow U$ de $X$ s'identifie au foncteur image inverse $j_u^*$ du morphisme 
de localisation 
\begin{equation}\label{rec4i}
j_{u}\colon \tY_{/f^*(U')}\rightarrow \tY_{/f^*(U)}
\end{equation} 
associé à $f^*(u)$ (\cite{sga4} IV 5.5). On note
\begin{equation}\label{rec4c}
\cF^\vee_{f^+}\rightarrow X^\circ
\end{equation}
la catégorie fibrée obtenue en associant à tout $U\in \ob(X)$ la catégorie $\tY_{/f^*(U)}$, et à tout morphisme 
$u\colon U'\rightarrow U$ de $X$ le foncteur image directe $j_{u*}$ du morphisme $j_u$ \eqref{rec4i}. 
On désigne par $\hY$ la catégorie des préfaisceaux de $\mU$-ensembles sur $Y$, par 
\begin{equation}\label{rec4d}
\cP^\vee_{f^+}\rightarrow X^\circ
\end{equation}
la catégorie fibrée obtenue en associant à tout $U\in \ob(X)$ la catégorie $\hY_{/f^+(U)}$ des préfaisceaux de $\mU$-ensembles
sur $Y_{/f^+(U)}$, et à tout morphisme $u\colon U'\rightarrow U$ de $X$ le foncteur 
``restriction de Weil'' $j_{u*}\colon \hY_{/f^+(U')}\rightarrow \hY_{/f^+(U)}$ défini, 
pour tous $F\in \ob(\hY_{/f^+(U')})$ et $V\in \ob(Y_{/f^+(U)})$, par  
\begin{equation}\label{rec4e}
j_{u*}(F)(V)=F(j_u^+(V))=F(V\times_{f^+(U)}f^+(U')).
\end{equation}

On a une équivalence de catégories 
\begin{eqnarray}\label{rec4f}
\hcD_{f^+}&\stackrel{\sim}{\rightarrow} &\bHom_{X^\circ}(X^\circ,\cP^\vee_{f^+})\\
F&\mapsto&\{U\mapsto F_U\},\nonumber
\end{eqnarray}
définie, pour tout $(V\rightarrow U)\in \ob(\cD_{f^+})$, par la relation  
\begin{equation}\label{rec4g}
F_U(V)=F(V\rightarrow U).
\end{equation}
On identifiera dans la suite $F$ à la section $\{U\mapsto F_U\}$ qui lui est associée par cette équivalence.

D'après (\cite{agt} VI.4.4), le foncteur \eqref{rec4f} induit un foncteur pleinement fidèle 
\begin{equation}\label{rec4h}
\tcD_{f^+}\rightarrow \bHom_{X^\circ}(X^\circ,\cF^\vee_{f^+})
\end{equation}
d'image essentielle les sections $\{U\mapsto F_U\}$ vérifiant une condition de recollement.

\subsection{}\label{rec3}
Soient $X$, $Y$, $X'$ et $Y'$ des $\mU$-sites dans lesquels les limites projectives sont représentables,
\begin{equation}\label{rec3a}
\xymatrix{
X\ar[r]^{f^+}\ar[d]_{u^+}&Y\ar[d]^{v^+}\\
X'\ar[r]^{f'^+}&Y'}
\end{equation}
un diagramme de foncteurs continus et exacts à gauche, commutatif à isomorphisme près~:
\begin{equation}\label{rec3b}
a\colon v^+\circ f^+\stackrel{\sim}{\rightarrow} f'^+\circ u^+.
\end{equation}
On note $\tX$, $\tY$, $\tX'$ et $\tY'$ les topos des faisceaux de $\mU$-ensembles sur $X$, $Y$, $X'$ et $Y'$,
respectivement,  
\begin{equation}\label{rec3c}
\xymatrix{
{\tY'}\ar[r]^{f'}\ar[d]_v&{\tX'}\ar[d]^u\\
{\tY}\ar[r]^f&{\tX}}
\end{equation}
le diagramme de morphismes de topos déduit de \eqref{rec3a} et 
\begin{equation}\label{rec3d}
\alpha\colon uf'\stackrel{\sim}{\rightarrow} fv
\end{equation}
l'isomorphisme induit par $a$. On a un foncteur 
\begin{equation}\label{rec3e}
\Phi^+\colon \cD_{f^+}\rightarrow \cD_{f'^+}, 
\end{equation}
qui à tout objet $(V\rightarrow U)$ de $\cD_{f^+}$ associe l'objet $(v^+(V)\rightarrow u^+(U))$
de $\cD_{f'^+}$ défini par le composé
\begin{equation}\label{rec3f}
\xymatrix{
{v^+(V)}\ar[r]&{v^+(f^+(U))}\ar[r]^{a(U)}&{f'^+(u^+(U))}}.
\end{equation}
Il est clairement exact à gauche. 

Pour tout $U\in \ob(X)$, on désigne par $v^+_U$ le foncteur composé
\begin{equation}\label{rec3g}
\xymatrix{
{Y_{/f^+(U)}}\ar[rr]^-(0.5){v^+_{/f^+(U)}}&&{Y'_{/v^+(f^+(U))}}\ar[r]^\sim&{Y'_{/f'^+(u^+(U))}}},
\end{equation}
où la seconde flèche est l'équivalence de catégories induite par l'isomorphisme $a(U)$.
Celui-ci est exact à gauche et continu lorsque l'on munit 
$Y_{/f^+(U)}$ et $Y'_{/f'^+(u^+(U))}$ des topologies induites par celles de $Y$ et $Y'$, respectivement 
(\cite{sga4} III 1.6 et 3.3). 
Pour tout morphisme $t\colon U'\rightarrow U$ de $X$, le diagramme de foncteurs 
\begin{equation}\label{rec3i}
\xymatrix{
{Y_{/f^+(U)}}\ar[r]^-(0.5){v_U^+}\ar[d]_{j_t^+}&{Y'_{/f'^+(u^+(U))}}\ar[d]^{j'^+_{u^+(t)}}\\
{Y_{/f^+(U')}}\ar[r]^-(0.5){v_{U'}^+}&{Y'_{/f'^+(u^+(U'))}}}
\end{equation}
où $j_t^+$ et $j'^+_{u^+(t)}$ sont les foncteurs de changement de base par $f^+(t)$ et $f'^+(u^+(t))$,
respectivement, est commutatif à isomorphisme canonique près. 

Pour tout préfaisceau $F=\{U'\mapsto F_{U'}\}$ sur $\cD_{f'^+}$, on a 
\begin{equation}\label{rec3h}
F\circ \Phi^+=\{U\mapsto F_{u^+(U)}\circ v^+_U\}.
\end{equation}
Si $F$ est un faisceau sur $\cD_{f'^+}$, $F\circ \Phi^+$ est un faisceau sur $\cD_{f^+}$ en vertu de (\cite{agt} VI.4.4). 
Par suite, $\Phi^+$ définit un morphisme de topos 
\begin{equation}\label{rec3j}
\Phi\colon \tX'\gtimes_{\tX'}\tY'\rightarrow \tX\gtimes_{\tX}\tY.
\end{equation}

Les carrés du diagramme de morphismes de topos
\begin{equation}\label{rec3k}
\xymatrix{
{\tX'}\ar[d]_{u}&{\tX'\gtimes_{\tX'}\tY'}\ar[l]_-(0.5){\p'_1}\ar[r]^-(0.5){\p'_2}\ar[d]^{\Phi}&{\tY'}\ar[d]^{v}\\
{\tX}&{\tX\gtimes_{\tX}\tY}\ar[l]_-(0.5){\p_1}\ar[r]^-(0.5){\p_2}&{\tY}}
\end{equation}
où les flèches horizontales sont les projections canoniques, sont commutatifs à des isomorphismes canoniques près~:
\begin{eqnarray}
\gamma\colon u \p'_1&\stackrel{\sim}{\rightarrow}& \p_1 \Phi,\label{rec3l}\\
\delta\colon v  \p'_2&\stackrel{\sim}{\rightarrow} &\p_2 \Phi,\label{rec3m}
\end{eqnarray}
le premier étant induit par $a$ \eqref{rec3b} et le second tautologique. On vérifie que le diagramme de 2-morphismes
\begin{equation}\label{rec3n}
\xymatrix{
{u f' \p'_2}\ar[r]^-(0.4){\alpha*\p'_2}\ar[d]_{u*\tau'}&{fv\p'_2}\ar[r]^-(0.4){f*\delta}&{f\p_2\Phi}\ar[d]^{\tau*\Phi}\\
{u\p'_1}\ar[rr]^-(0.4)\gamma&&{\p_1\Phi}}
\end{equation}
où $\tau$ et $\tau'$ sont les $2$-morphismes canoniques \eqref{rec1e}, est commutatif. Il résulte alors de la
propriété universelle de $\tX\gtimes_{\tX}\tY$ que $\Phi$ est uniquement déterminé par 
le quadruplet $(u,v,\gamma,\delta)$. 

Le diagramme 
\begin{equation}\label{rec3o}
\xymatrix{
{\tY'}\ar[r]^-(0.5){\Psi'}\ar[d]_v&{\tX'\gtimes_{\tX'}\tY'}\ar[d]^\Phi\\
{\tY}\ar[r]^-(0.5)\Psi&{\tX\gtimes_{\tX}\tY}}
\end{equation}
où $\Psi$ et $\Psi'$ sont les morphismes des cycles co-proches \eqref{rec1i}, 
est clairement commutatif à isomorphisme canonique près.

\begin{rema}\label{rec2}
Soient $X$ et $Y$ deux $\mU$-sites dans lesquels les 
limites projectives finies sont représentables, $f^+\colon X\rightarrow Y$ un foncteur continu et exact à gauche. 
On note $\tX$ et $\tY$ les topos des faisceaux de $\mU$-ensembles sur $X$ et $Y$, respectivement, 
$\varepsilon_X\colon X\rightarrow \tX$ et $\varepsilon_Y\colon Y\rightarrow \tY$ les foncteurs canoniques et 
$f\colon \tY\rightarrow \tX$ le morphisme de topos associé à $f^+$.  On notera que le diagramme 
\begin{equation}
\xymatrix{
X\ar[r]^{f^+}\ar[d]_{\varepsilon_X}&Y\ar[d]^{\varepsilon_Y}\\
\tX\ar[r]^{f^*}&\tY}
\end{equation}
est commutatif à isomorphisme canonique près. 
On considère $\tX$ et $\tY$ comme des $\mU$-sites
munis des topologies canoniques. Les foncteurs $\varepsilon_X$ et $\varepsilon_Y$  
sont continus et exacts à gauche (\cite{sga4} III 3.5). Il en est évidemment de même du foncteur $f^*$. 
Les topos des faisceaux de $\mU$-ensembles sur $\tX$ et $\tY$ s'identifient à $\tX$ et $\tY$
par les  morphismes de topos associés à $\varepsilon_X$ et $\varepsilon_Y$, respectivement (\cite{sga4} IV 1.2). 
Le morphisme de topos associé à $f^*$ s'identifie à $f$. 
D'après \ref{rec3}, le foncteur \eqref{rec3e}
\begin{equation}
\Phi^+\colon \cD_{f^+}\rightarrow \cD_{f^*}, \ \ \ (V\rightarrow U)\mapsto (\varepsilon_Y(V)\rightarrow \varepsilon_X(U))
\end{equation}
est continu et exact à gauche, et le morphisme de topos associé $\Phi\colon \tcD_{f^*}\rightarrow \tcD_{f^+}$ est 
une équivalence de topos. On retrouve le fait que  le topos
co-évanescent $\tX\gtimes_{\tX}\tY$ ne dépend que des topos $\tX$ et $\tY$
et du morphisme $f\colon \tY\rightarrow \tX$. 
\end{rema}

\begin{prop}\label{rec22}
Reprenons les hypothèses et notations de \ref{rec1}; supposons de plus, que tous les objets de $X$ et de $Y$ soient quasi-compacts. Alors,
\begin{itemize}
\item[{\rm (i)}] Les topos $\tX$ et $\tY$ sont cohérents et le morphisme $f\colon \tY\rightarrow \tX$ est cohérent.
\item[{\rm (ii)}] Pour tout objet $(V\rightarrow U)$ de $\cD_{f^+}$, le faisceau associé $(V\rightarrow U)^a$ est un objet cohérent de $\tX\gtimes_{\tX}\tY$. 
\item[{\rm (iii)}] Le topos $\tX\gtimes_{\tX}\tY$ est cohérent. 
\end{itemize}
\end{prop}

(i) Cela résulte de (\cite{sga4} VI 2.4.5 et 3.3).

(ii) \& (iii) Pour tout objet $V$ de $Y$, notant $V^a$ le faisceau associé de $\tY$, le topos $\tY_{/V^a}$ est cohérent (\cite{sga4} VI 2.4.2 et 2.4.5).
Les propositions résultent donc de (\cite{agt} VI.5.27 et VI.5.5). 

\begin{cor}[\cite{agt} VI.5.30]\label{rec9}
Soient $\tX$, $\tY$ deux $\mU$-topos cohérents, $f\colon \tY\rightarrow \tX$ un morphisme cohérent. 
Alors le produit orienté $\tX\gtimes_{\tX}\tY$ est cohérent. En particulier, il a suffisamment de points. 
\end{cor}

\subsection{}\label{rec5}
Soient $X$ et $Y$ deux $\mU$-sites dans lesquels les 
limites projectives finies sont représentables, $f^+\colon X\rightarrow Y$ un foncteur continu et exact à gauche. 
On note $\tX$ et $\tY$ les topos des faisceaux de $\mU$-ensembles sur $X$ et $Y$, respectivement, 
$f\colon \tY\rightarrow \tX$ le morphisme de topos associé à $f^+$ et 
$e_X$ et $e_Y$ des objets finaux de $X$ et $Y$, respectivement. 
Soient $(B\rightarrow A)$ un objet de $\cD_{f^+}$, $A^\ra$ 
le faisceau de $\tX$ associé à $A$, $B^\ra$ le faisceau de $\tY$ associé à $B$.
On désigne par $j_{A} \colon X_{/A}\rightarrow X$ et $j_B\colon Y_{/B}\rightarrow Y$ les foncteurs canoniques et 
aussi abusivement par $j_{A} \colon \tX_{/A^\ra}\rightarrow \tX$ et $j_{B}\colon \tY_{/B^\ra}\rightarrow \tY$ les morphismes de 
localisation de $\tX$ en $A^\ra$ et de $\tY$ en de $B^\ra$, respectivement (cf. \ref{notconv6}). On munit $X_{/A}$
et $Y_{/B}$ des topologies induites par celles de $X$ et $Y$, respectivement. Le foncteur 
\begin{equation}\label{rec5b}
f'^+\colon  X_{/A}\rightarrow Y_{/B}, \ \ \ U\mapsto f^+(j_A(U))\times_{f^+(A)}B
\end{equation}
est exact à gauche et continu en vertu de (\cite{sga4} III 1.6 et 3.3). 
Le morphisme de topos 
\begin{equation}\label{rec5c}
f'\colon \tY_{/B^\ra}\rightarrow \tX_{/A^\ra}
\end{equation}
associé à $f'^+$ s'identifie au morphisme composé 
\[
\xymatrix{
{\tY_{/B^\ra}}\ar[r]&{\tY_{/f^*(A^\ra)}}\ar[r]^-(0.5){f_{/A^\ra}}&{\tX_{/A^\ra}}},
\]
où la première flèche est le morphisme de localisation associé à $B^\ra\rightarrow f^*(A^\ra)$ 
(\cite{sga4} IV 5.5) et la seconde flèche est le morphisme déduit de $f$ (\cite{sga4} IV 5.10).

On a un foncteur 
\begin{equation}\label{rec5d}
j_{(B\rightarrow A)}\colon \cD_{f'^+}\rightarrow \cD_{f^+}
\end{equation}
qui à tout objet $(V\rightarrow U)$ de $\cD_{f'^+}$ associe l'objet $(j_B(V)\rightarrow j_A(U))$
de $\cD_{f^+}$ défini par le composé
\begin{equation}
j_B(V)\rightarrow j_B(f'^+(U))\stackrel{\sim}{\rightarrow} f^+(j_A(U))\times_{f^+(A)}B\rightarrow f^+(j_A(U)),
\end{equation}
où la seconde flèche est l'isomorphisme canonique \eqref{rec5b} et la dernière flèche est la projection canonique. 
Celui-ci se factorise à travers une équivalence canonique de catégories 
\begin{equation}\label{rec5e}
n\colon \cD_{f'^+}\stackrel{\sim}{\rightarrow}(\cD_{f^+})_{/(B\rightarrow A)}.
\end{equation}

D'après (\cite{agt} VI.4.18), la topologie co-évanescente de $\cD_{f'^+}$ est induite par la topologie co-évanescente de 
$\cD_{f^+}$ au moyen du foncteur $j_{(B\rightarrow A)}$; en particulier, $n$ 
induit une équivalence de topos 
\begin{equation}\label{rec5f}
m\colon (\tX\gtimes_{\tX}\tY)_{/(B\rightarrow A)^\ra}\stackrel{\sim}{\rightarrow} 
\tX_{/A^\ra}\gtimes_{\tX_{/A^\ra}}\tY_{/B^\ra}.
\end{equation}

Par ailleurs, le diagramme 
\begin{equation}\label{rec5g}
\xymatrix{
X\ar[r]^{f^+}\ar[d]_{j_A^+}&Y\ar[d]^{j_B^+}\\
{X_{/A}}\ar[r]^{f'^+}&{Y_{/B}}}
\end{equation}
où $j_A^+$ et $j_B^+$ sont les foncteurs de changement de base par $A\rightarrow e_X$ et $B\rightarrow e_Y$, 
respectivement, est commutatif à isomorphisme canonique près. Le foncteur induit par ce diagramme \eqref{rec3e}
\begin{equation}\label{rec5h}
\Phi^+\colon \cD_{f^+}\rightarrow \cD_{f'^+}, \ \ \ (V\rightarrow U)\mapsto (j_B^+(V)\rightarrow j_A^+(U))
\end{equation}
est exact à gauche et continu d'après \ref{rec3} et \ref{notconv6}. Il définit donc un morphisme de topos 
\begin{equation}\label{rec5i}
\Phi\colon \tX_{/A^\ra}\gtimes_{\tX_{/A^\ra}}\tY_{/B^\ra}\rightarrow \tX\gtimes_{\tX}\tY.
\end{equation}
On notera que le composé $n\circ \Phi^+$ n'est autre que le foncteur de changement de base par $(B\rightarrow A)$. 
Par suite, $\Phi m$ est le morphisme de localisation 
\begin{equation}\label{rec5j}
j_{(B\rightarrow A)}\colon (\tX\gtimes_{\tX}\tY)_{/(B\rightarrow A)^\ra}\rightarrow \tX\gtimes_{\tX}\tY.
\end{equation}

\begin{lem}\label{rec7}
Soient $Y$ un $\mU$-topos, $B$ un ouvert de $Y$ ({\em i.e.}, un  sous-objet de l'objet final de $Y$), $V$ un objet de $Y$. 
On note $j_B\colon Y_{/B}\rightarrow Y$ et $j_V\colon Y_{/V}\rightarrow Y$
les morphismes de localisation de $Y$ en $B$ et $V$, respectivement, $Z$ le sous-topos
fermé de $Y$ complémentaire de l'ouvert $B$  {\rm (\cite{sga4} IV 9.3.5)} et 
$i\colon Z\rightarrow Y$ le morphisme d'inclusion. Alors, 
\begin{itemize}
\item[{\rm (i)}] Le morphisme
\begin{equation}\label{rec7a}
j_{B/V}\colon Y_{/B\times V}\rightarrow Y_{/V}
\end{equation}
induit par $j_B$ {\rm (\cite{sga4} IV 5.10)} est un plongement ouvert, d'ouvert de $Y_{/V}$ associé $B\times V$ 
{\rm (\cite{sga4} IV 9.2)}.
\item[{\rm (ii)}] Le morphisme
\begin{equation}\label{rec7b}
i_{/V}\colon Z_{/i^*(V)}\rightarrow Y_{/V}
\end{equation}
induit par $i$ {\rm (\cite{sga4} IV 5.10)} est un plongement, et le foncteur $(i_{/V})_*$
identifie $Z_{/i^*(V)}$ au sous-topos fermé de $Y_{/V}$ complémentaire de l'ouvert $B\times V$. 
\end{itemize}
\end{lem}

(i) En effet, on a $B\times V=j_V^*(B)$ (\cite{sga4} III 5.4); c'est donc un ouvert de $Y_{/V}$. 
Le foncteur $j_{B/V}^*$, induit par $j_B^*$, est le changement de base par le morphisme 
$B\times V\rightarrow V$. Par suite, $j_{B/V}$ est isomorphe au morphisme de localisation de $Y_{/V}$
en $B\times V$. 

(ii) Notons $Z'$ le sous-topos fermé de $Y_{/V}$
complémentaire de $B\times V$ et $i'\colon Z'\rightarrow Y_{/V}$ le morphisme d'inclusion. D'après (\cite{sga4} IV 9.4.3), 
il existe un morphisme $j'\colon Z'\rightarrow Z$, unique à isomorphisme près, tel que le diagramme 
\begin{equation}\label{rec7d}
\xymatrix{
Z'\ar[r]^-(0.5){i'}\ar[d]_{j'}&{Y_{/V}}\ar[d]^{j_V}\\
Z\ar[r]^-(0.5)i&Y}
\end{equation}
soit $2$-cartésien. Par ailleurs, le diagramme 
\begin{equation}\label{rec7c}
\xymatrix{
Z_{/i^*(V)}\ar[r]^-(0.5){i_{/V}}\ar[d]_{j_{i^*(V)}}&{Y_{/V}}\ar[d]^{j_V}\\
Z\ar[r]^-(0.5)i&Y}
\end{equation}
où $j_{i^*(V)}$ est le morphisme de localisation de $Z$ en $i^*(V)$, est $2$-cartésien en vertu de (\cite{sga4} IV 5.11);
d'où la proposition.

\begin{prop}\label{rec6}
Soient $f\colon Y\rightarrow X$ un morphisme de $\mU$-topos, $e_X$ un objet final de $X$,
$B$ un ouvert de $Y$ ({\em i.e.}, un sous-objet de l'objet final de $Y$).
On note $j_B\colon Y_{/B}\rightarrow Y$ le morphisme de localisation de $Y$ en $B$, $Z$ le sous-topos
fermé de $Y$ complémentaire de l'ouvert $B$  {\rm (\cite{sga4} IV 9.3.5)} et 
$i\colon Z\rightarrow Y$ le morphisme d'inclusion.
Pour tout $V\in \ob(Y)$, on désigne par $j_{B/V}\colon Y_{/B\times V}\rightarrow Y_{/V}$
le morphisme de localisation de $Y_{/V}$ en $B\times V$ et par
$i_{/V}\colon Z_{/i^*(V)}\rightarrow Y_{/V}$ le morphisme induit par $i$ \eqref{rec7}. 
Notons $X\gtimes_XY$, $X\gtimes_XY_{/B}$ et $X\gtimes_XZ$ les topos co-évanescents
de $f$, $fj_B$ et $fi$, respectivement. Alors, 
\begin{itemize}
\item[{\rm (i)}] Le morphisme
\begin{equation}\label{rec6a}
J\colon X\gtimes_XY_{/B}\rightarrow X\gtimes_XY
\end{equation}
induit par $j_B$ \eqref{rec3j} est un plongement ouvert, d'ouvert de $X\gtimes_XY$ associé 
$(B\rightarrow e_X)^\ra$  {\rm (\cite{sga4} IV 9.2)}.
Pour tout objet $G=\{U\mapsto G_U\}$ de $X\gtimes_XY_{/B}$, on a un isomorphisme canonique fonctoriel
\begin{equation}\label{rec6bb}
J_*(G)\stackrel{\sim}{\rightarrow}\{U\mapsto (j_{B/f^*(U)})_*(G_U)\}.
\end{equation}
Pour tout objet $F=\{U\mapsto F_U\}$ de $X\gtimes_XY$, on a un isomorphisme canonique fonctoriel
\begin{equation}\label{rec6b}
J^*(F)\stackrel{\sim}{\rightarrow}\{U\mapsto F_U\times B\},
\end{equation}
où $F_U\times B$ désigne abusivement l'objet $F_U\times_{f^*(U)}(f^*(U)\times B)$ de $Y_{/f^*(U)\times B}$.
\item[{\rm (ii)}] Le morphisme
\begin{equation}\label{rec6c}
I\colon X\gtimes_XZ\rightarrow X\gtimes_XY
\end{equation}
induit par $i$ \eqref{rec3j} est un plongement. Le foncteur $I_*$ identifie $X\gtimes_XZ$ au sous-topos fermé de 
$X\gtimes_XY$ complémentaire de l'ouvert $(B\rightarrow e_X)^\ra$. 
Pour tout objet $G=\{U\mapsto G_U\}$ de $X\gtimes_XZ$, on a un isomorphisme canonique fonctoriel
\begin{equation}\label{rec6d}
I_*(G)\stackrel{\sim}{\rightarrow}\{U\mapsto (i_{/f^*(U)})_*(G_U)\}.
\end{equation}
\item[{\rm (iii)}] Pour tout objet $F=\{U\mapsto F_U\}$ de $X\gtimes_XY$, $I^*(F)$ est le faisceau de $X\gtimes_XZ$
associé au préfaisceau $\{U\mapsto (i_{/f^*(U)})^*(F_U)\}$ sur le site $\cD_{i^*f^*}$ \eqref{rec2}.
\end{itemize}
\end{prop}

(i) En effet, d'après \ref{rec5}, on a une équivalence canonique de topos \eqref{rec5f}
\begin{equation}\label{rec6e}
m\colon (X\gtimes_{X}Y)_{/(B\rightarrow e_X)^\ra}\stackrel{\sim}{\rightarrow} X\gtimes_{X}Y_{/B},
\end{equation}
et $Jm$ est le morphisme de localisation de $X\gtimes_XY$ en $(B\rightarrow e_X)^\ra$. 
On a $(B\rightarrow e_X)^\ra=\p_2^*(B)$; c'est donc un ouvert de $X\gtimes_XY$. 
L'isomorphisme \eqref{rec6bb} résulte aussitôt des définitions \eqref{rec3h}.
Pour tout objet $F=\{U\mapsto F_U\}$ de $X\gtimes_XY$, on a un isomorphisme canonique fonctoriel 
\begin{equation}\label{rec6f}
m^*(J^*(F))\stackrel{\sim}{\rightarrow}F\times \p_2^*(B).
\end{equation}
En vertu de (\cite{agt} VI.4.12), $\p_2^*(B)$ est le faisceau $\{U\mapsto f^*(U)\times B\}$. 
Compte tenu de (\cite{sga4} II 4.1(3)), on en déduit un isomorphisme canonique fonctoriel 
\begin{equation}\label{rec6g}
m^*(J^*(F))\stackrel{\sim}{\rightarrow}\{U\mapsto F_U\times_{f^*(U)}(f^*(U)\times B)\};
\end{equation}
d'où l'isomorphisme \eqref{rec6b}.

(ii)  Pour tout objet $G=\{U\mapsto G_U\}$ de $X\gtimes_XZ$, on a un isomorphisme canonique fonctoriel \eqref{rec3h} 
\begin{equation}\label{rec6h}
I_*(G)\stackrel{\sim}{\rightarrow}\{U\mapsto (i_{/f^*(U)})_*(G_U)\}.
\end{equation}
On observera que pour tout morphisme $t\colon U'\rightarrow U$ de $X$, le diagramme de morphismes de topos
\begin{equation}\label{rec6i}
\xymatrix{
{Z_{/i^*(f^*(U'))}}\ar[r]^-(0.5){i_{/f^*(U')}}\ar[d]&{Y_{/f^*(U')}}\ar[d]\\
{Z_{/i^*(f^*(U))}}\ar[r]^-(0.5){i_{/f^*(U)}}&{Y_{/f^*(U)}}}
\end{equation}
où les flèches verticales sont les morphismes de localisation associés aux morphismes $f^*(t)$ de $Y$ et $i^*(f^*(t))$ de $Z$,
est commutatif à isomorphisme canonique près. Pour tout $U\in \ob(X)$, le foncteur $(i_{/f^*(U)})_*$
est pleinement fidèle d'après \ref{rec7}. Il en est donc de même du foncteur $I_*$ \eqref{rec4h}, 
de sorte que $I$ est un plongement. Pour tout objet $G$ de $X\gtimes_XZ$, 
$J^*(I_*(G))$ est un objet final de $X\gtimes_XY_{/B}$ en vertu de \ref{rec7} et \eqref{rec6b}.
Inversement, soit $\{U\mapsto F_U\}$ un objet de $X\gtimes_XY$ tel que $J^*(F)$ soit un objet final de $X\gtimes_XY_{/B}$.
Pour tout $U\in \ob(X)$, $F_U\times B$ est donc un objet final de $Y_{/f^*(U)\times B}$ \eqref{rec6b}. 
D'après \ref{rec7}, il existe un objet $G_U$ de $Z_{/i^*(f^*(U))}$ tel que $F_U=(i_{/f^*(U)})_*(G_U)$. 
Il résulte aussitôt de \eqref{rec6i}  et (\cite{agt} VI.4.4) que la collection $G=\{U\mapsto G_U\}$ 
est naturellement un objet de $X\gtimes_XZ$. On a $I_*(G)\simeq F$ \eqref{rec6h}. Par suite, 
le foncteur $I_*$ identifie $X\gtimes_XZ$ au sous-topos fermé de 
$X\gtimes_XY$ complémentaire de l'ouvert $(B\rightarrow e_X)^\ra$. 

(iii) On notera d'abord que $\{U\mapsto (i_{/f^*(U)})^*(F_U)\}$ est naturellement un préfaisceau sur $\cD_{i^*f^*}$,
compte tenu de \eqref{rec6i}. D'après (\cite{sga4} I 5.1 et III 1.3), 
$I^*(F)$ est le faisceau sur $\cD_{i^*f^*}$ associé au préfaisceau $P$ défini pour tout 
$(W\rightarrow U)\in \ob(\cD_{i^*f^*})$ par 
\begin{equation}
P(W\rightarrow U)=\underset{\underset{(V'\rightarrow U')\in C^\circ_{(W\rightarrow U)}}{\longrightarrow}}\lim\ F(V'\rightarrow U'),
\end{equation}
où $C_{(V\rightarrow U)}$ est la catégorie des couples $((V'\rightarrow U'),t)$ formés d'un objet $(V'\rightarrow U')$
de $\cD_{f^*}$ et d'un morphisme $t\colon (W\rightarrow U)\rightarrow (i^*(V')\rightarrow U')$ de $\cD_{i^*f^*}$.
Notons $B_{(W\rightarrow U)}$ la catégorie des couples $(V,k)$ formés d'un objet $V$ de $Y_{/f^*(U)}$ et 
d'un morphisme $k\colon W\rightarrow (i_{/f^*(U)})^*(V)$ de $Z_{/i^*(f^*(U))}$. 
Les catégories $B_{(W\rightarrow U)}$ et $C_{(W\rightarrow U)}$ sont clairement cofiltrantes. 
On a un foncteur pleinement fidèle 
\begin{equation}
\varphi\colon B_{(W\rightarrow U)}\rightarrow C_{(W\rightarrow U)}, \ \ \ V\mapsto (V\rightarrow U).
\end{equation} 
Pour tout objet $((V'\rightarrow U'),t)$ de $C_{(W\rightarrow U)}$, $t$ induit un morphisme $u\colon U\rightarrow U'$ de $X$
et un morphisme 
\[
k\colon W\rightarrow (i_{/f^*(U)})^*(V'\times_{f^*(U')}f^*(U))
\] 
de $Z_{/i^*(f^*(U))}$, de sorte que $(V'\times_{f^*(U')}f^*(U),k)$ est un objet de $B_{(W\rightarrow U)}$.
On a un morphisme 
\[
\varphi(V'\times_{f^*(U')}f^*(U),k)\rightarrow ((V'\rightarrow U'),t)
\]
de $C_{(W\rightarrow U)}$. Le foncteur $\varphi^\circ$ est donc cofinal d'après (\cite{sga4} I 8.1.3(c)). 
Par suite, on a un un isomorphisme canonique 
\begin{equation}
P(W\rightarrow U)\stackrel{\sim}{\rightarrow}
\underset{\underset{V\in B^\circ_{(W\rightarrow U)}}{\longrightarrow}}\lim\ F(V\rightarrow U).
\end{equation}
La proposition résulte alors de (\cite{agt} VI.5.17 et VI.5.5) et (\cite{sga4} I 5.1 et III 1.3).

\subsection{}\label{rec8}
Soient $X$ et $Y$ deux $\mU$-sites dans lesquels les 
limites projectives finies sont représentables, $f^+\colon X\rightarrow Y$ un foncteur continu et exact à gauche. 
On note $\tX$ et $\tY$ les topos des faisceaux de $\mU$-ensembles sur $X$ et $Y$, respectivement,
$f\colon \tY\rightarrow \tX$ le morphisme de topos associé à $f^+$ et $\tX\gtimes_\tX\tY$ le topos co-évanescent de $f$.
La donnée d'un point de $\tX\gtimes_\tX\tY$ est équivalente à 
la donnée d'une paire de points $x\colon \Pt\rightarrow \tX$ et $y\colon \Pt\rightarrow \tY$ 
et d'un $2$-morphisme $u\colon f(y)\rightarrow x$ (\cite{agt} VI.3.11). Un tel point sera noté $(y\rightarrow x)$,  
ou encore $(u\colon y\rightarrow x)$. 
Pour tous $F\in \ob(\tX)$ et $G\in \ob(\tY)$, on a des isomorphismes canoniques fonctoriels
\begin{eqnarray}
(\rp_1^*F)_{(y\rightarrow x)} &\stackrel{\sim}{\rightarrow}& F_x, \label{rec8a}\\
(\rp_2^*G)_{(y\rightarrow x)} &\stackrel{\sim}{\rightarrow}& G_y.\label{rec8b}
\end{eqnarray}
L'application 
\begin{equation}\label{rec8c}
(\rp_1^*F)_{(y\rightarrow x)} \rightarrow (\rp_2^*(f^*F))_{(y\rightarrow x)} 
\end{equation}
induite par $\tau$ \eqref{rec1e}, s'identifie canoniquement au morphisme 
de spécialisation $F_x\rightarrow F_{f(y)}$ défini par $u$ (cf. \cite{agt} VI.4.20). 
On a un isomorphisme canonique fonctoriel (\cite{agt} (VI.4.20.4))
\begin{equation}\label{rec8d}
(\Psi_*G)_{(y\rightarrow x)} \stackrel{\sim}{\rightarrow} G_y.
\end{equation}
D'après la propriété universelle de $\tX\gtimes_\tX\tY$, 
on a un isomorphisme canonique fonctoriel de points de $\tX\gtimes_\tX\tY$
\begin{equation}\label{rec8e}
\Psi(y) \stackrel{\sim}{\rightarrow} (y\rightarrow f(y)).
\end{equation}
Pour tout $(V\rightarrow U)\in \ob(\cD_{f^+})$, on a un isomorphisme canonique fonctoriel (\cite{agt} (VI.4.20.6))
\begin{equation}\label{rec8f}
(V\rightarrow U)^\ra_{(y\rightarrow x)} \stackrel{\sim}{\rightarrow} U^\ra_x\times_{U^\ra_{f(y)}}V^\ra_y,
\end{equation}
où l'exposant $^\ra$ désigne les faisceaux associés, l'application $V^\ra_y\rightarrow U^\ra_{f(y)}$ 
est induite par le morphisme structural $V\rightarrow f^+(U)$ 
et l'application $U^\ra_x\rightarrow U^\ra_{f(y)}$ est le morphisme de spécialisation défini par~$u$. 

On désigne par $\cP_{(y\rightarrow x)}$ la catégorie des voisinages de $(y\rightarrow x)$ dans $\cD_{f^+}$ 
(\cite{sga4} IV 6.8.2), autrement dit, la catégorie des triplets $((V\rightarrow U),\xi,\zeta)$ formés d'un objet 
$(V\rightarrow U)$ de $\cD_{f^+}$ et d'un élément $(\xi,\zeta)\in U^\ra_x\times_{U^\ra_{f(y)}}V^\ra_y$ \eqref{rec8f}.
Pour tout préfaisceau $F=\{U\mapsto F_U\}$ sur $\cD_{f^+}$, on a un isomorphisme canonique fonctoriel (\cite{sga4} IV 6.8.4)
\begin{equation}\label{rec8g}
(F^a)_{(y\mapsto x)}\stackrel{\sim}{\rightarrow}\underset{\underset{((V\rightarrow U),\xi,\zeta)\in \cP^\circ_{(y\rightarrow x)}}{\longrightarrow}}{\lim}F_U(V).
\end{equation}

\subsection{}\label{rec10}
Soient $X$ un $\mU$-topos local de centre $x\colon \Pt\rightarrow X$ (\cite{sga4} VI 8.4.6), 
$f\colon Y\rightarrow X$ un morphisme de $\mU$-topos, $X\gtimes_XY$ le topos co-évanescent associé. 
On note $\phi\colon X\rightarrow \Pt$ l'unique morphisme de topos (\cite{sga4} IV 4.3). 
Par définition, on a un isomorphisme $x^*\stackrel{\sim}{\rightarrow}\phi_*$. On en déduit un isomorphisme
\begin{equation}\label{rec10a}
(x\phi)^*\stackrel{\sim}{\rightarrow}\phi^*\phi_*.
\end{equation}
Composant avec le morphisme d'adjonction $\phi^*\phi_*\rightarrow \id_X$, on obtient un morphisme de foncteurs
$(x\phi)^*\rightarrow \id_X$ et par suite un $2$-morphisme 
\begin{equation}\label{rec10b}
\id_X\rightarrow x\phi.
\end{equation}
Les morphisme $x\phi f\colon Y\rightarrow X$ et $\id_Y\colon Y\rightarrow Y$ et le $2$-morphisme 
$f\rightarrow x\phi f$ déduit de \eqref{rec10b} induisent un morphisme 
\begin{equation}\label{rec10c}
\gamma\colon Y\rightarrow X\gtimes_XY
\end{equation}
tel que $\rp_1\gamma=x\phi f$ et $\rp_2\gamma=\id_Y$ \eqref{rec1}~:
\begin{equation}\label{rec10d}
\xymatrix{
X\ar[d]_{x\phi}&Y\ar[l]_f\ar[rd]^{\id_Y}\ar[d]_\gamma&\\
X\ar[rd]_{\id_X}&{X\gtimes_XY}\ar[r]^-(0.4){\rp_2}\ar[l]_-(0.4){\rp_1}&Y\ar[dl]^f\\
&X&}
\end{equation}
De la relation $\rp_2\gamma=\id_Y$ on obtient un morphisme de changement de base 
\begin{equation}\label{rec10e}
\rp_{2*}\rightarrow \gamma^*,
\end{equation}
composé de $\rp_{2*}\rightarrow \rp_{2*}\gamma_*\gamma^*\stackrel{\sim}{\rightarrow}\gamma^*$,
où la première flèche est induite par le morphisme d'adjonction $\id\rightarrow \gamma_*\gamma^*$,
et la seconde flèche par la relation $p_2\gamma=\id_Y$.

\begin{prop}[Gabber, \cite{travaux-gabber} XI 2.3]\label{rec11}
Sous les hypothèses de \ref{rec10}, pour tout objet $F$ de $X\gtimes_XY$ et tout point $y$ de $Y$, l'application 
\begin{equation}\label{rec11a}
(\rp_{2*}F)_y\rightarrow (\gamma^*F)_y
\end{equation}
déduite du morphisme de changement de base \eqref{rec10e} est bijective. 
\end{prop}

En effet, le composé $x\phi f y\colon \Pt\rightarrow X$ est canoniquement isomorphe à $x$. 
Le $2$-morphisme \eqref{rec10b} induit donc un $2$-morphisme $u\colon f(y)\rightarrow x$ définisant 
un point de $X\gtimes_XY$ \eqref{rec8}, qui n'est autre que $\gamma(y)$. 
Notons $\cP_{\gamma(y)}$ la catégorie des voisinages de $\gamma(y)$ dans $\cD_{f^*}$ \eqref{rec8}
et $\cC_y$ la catégorie des voisinages de $y$ (\cite{sga4} IV 6.8), 
autrement dit la catégorie des couples $(V,\zeta)$ formés d'un objet de $Y$ et d'un élément $\zeta\in V_y$.
Soient $e_X$ un objet final de $X$, $\xi_e$ l'unique élément de $x^*(e_X)$. 
On a un foncteur
\begin{equation}\label{rec11b}
\iota\colon \cC_y\rightarrow \cP_{\gamma(y)}, \ \ \ (V,\zeta)\mapsto ((V\rightarrow e_X),\xi_e,\zeta).
\end{equation}
Les catégories $\cP_{\gamma(y)}$ et $\cC_y$ sont cofiltrantes (\cite{sga4} IV 6.8.2) et le foncteur 
$\iota$ est pleinement fidèle. On vérifie aussitôt que \eqref{rec11a} est l'application
\begin{equation}
\underset{\underset{(V,\zeta)\in \cC^\circ_y}{\longrightarrow}}{\lim}F_{e_X}(V)\longrightarrow 
\underset{\underset{((V\rightarrow U),\xi,\zeta)\in \cP^\circ_{\gamma(y)}}{\longrightarrow}}{\lim}F_U(V)
\end{equation}
induite par $\iota^\circ$. Il suffit donc de montrer que $\iota^\circ$ est cofinal. 
Soit $((V\rightarrow U),\zeta,\xi)\in \ob(\cP_{\gamma(y)})$.
L'isomorphisme $x^*\stackrel{\sim}{\rightarrow}\phi_*$ permet d'identifier $\xi\in U_x$ à une section
que l'on note encore $\xi\in U(e_X)$.
On observera que l'image canonique $\xi_{f(y)}$ de cette section dans $U_{f(y)}$ n'est autre que l'image de $\xi$
par le morphisme de spécialisation $U_x\rightarrow U_{f(y)}$ défini par $u$ \eqref{rec10b}. 
Considérons alors l'objet $W=V\times_{f^*(U)}f^*(e_X)$ de $Y$ défini par $\xi\in U(e_X)$. 
D'après ce qui précède, $W_y$ s'identifie à l'image inverse de $\xi_{f(y)}\in U_{f(y)}$ par l'application 
$V_y\rightarrow U_{f(y)}$ définie par le morphisme structural $V\rightarrow f^*(U)$. 
En particulier, $\zeta$ est un élément de $W_y$, 
de sorte que $(W,\zeta)$ est un objet de $\cC_y$.  Par ailleurs, $\xi$ induit un morphisme 
\begin{equation}
\iota(W,\zeta)=((W\rightarrow e_X),\xi_e,\zeta)\longrightarrow ((V\rightarrow U),\xi,\zeta)
\end{equation}
de $\cP_{\gamma(y)}$. On en déduit que $\iota^\circ$ est cofinal en vertu de (\cite{sga4} IV 8.1.3(c)); d'où la proposition.

\begin{cor}\label{rec12}
Sous les hypothèses de \ref{rec10}, si, de plus, $Y$ a assez de points, alors
le morphisme de changement de base $\rp_{2*}\rightarrow \gamma^*$ \eqref{rec10e}
est un isomorphisme~; en particulier, le foncteur $\rp_{2*}$ est exact. 
\end{cor}

\begin{cor}\label{rec13}
Sous les hypothèses de \ref{rec10}, si, de plus, $Y$ est local de centre $y$, alors 
$X\gtimes_XY$ est local de centre $\gamma(y)$.
\end{cor}

En effet, pour tout objet $F$ de $X\gtimes_XY$, le diagramme 
\begin{equation}
\xymatrix{
{\Gamma(Y,\rp_{2*}F)}\ar[r]\ar@{=}[d]&{(\rp_{2*}F)_y}\ar[d]\\ 
{\Gamma(X\gtimes_XY,F)}\ar[r]&{F_{\gamma(y)}}}
\end{equation}
où les flèches horizontales sont les applications canoniques et la flèche verticale de droite est l'application \eqref{rec11a}
est commutatif. La proposition résulte donc de \ref{rec11}. 

\begin{cor}\label{rec14}
Conservons les hypothèses de \ref{rec10}, supposons, de plus, que $Y$ ait assez de points. Alors~:
\begin{itemize}
\item[{\rm (i)}] Pour tout faisceau $F$ de $X\gtimes_XY$, l'application canonique 
\begin{equation}\label{rec14a}
\Gamma(X\gtimes_XY,F)\rightarrow \Gamma(Y, \gamma^*F)
\end{equation}
est bijective. 
\item[{\rm (ii)}] Pour tout faisceau abélien $F$ de $X\gtimes_XY$, l'application canonique 
\begin{equation}\label{rec14b}
\rH^i(X\gtimes_XY,F)\rightarrow \rH^i(Y, \gamma^*F)
\end{equation}
est bijective pour tout $i\geq 0$.
\end{itemize}
\end{cor}

Cela résulte aussitôt de \ref{rec11} (cf. la preuve de l'énoncé analogue \cite{agt} VI.10.28). 

\subsection{}\label{rec15}
Soit $f\colon Y\rightarrow X$ un morphisme de schémas. On désigne par $\cD_{f^+}$ le site co-évanescent  
du foncteur $f^+\colon \Et_{/X}\rightarrow \Et_{/Y}$ de changement de base par $f$,
et par $X_\et\gtimes_{X_\et}Y_\et$ le topos co-évanescent du morphisme de topos associé $f\colon Y_\et\rightarrow X_\et$ (cf. \ref{notconv10} et \ref{rec1}). 
D'après  \ref{rec8} et (\cite{sga4} VIII 7.9), la donnée d'un point de $X_\et\gtimes_{X_\et}Y_\et$ 
est équivalente à la donnée d'une paire de points géométriques $\ox$ de $X$ et $\oy$ de $Y$
et d'une flèche de spécialisation $u$ de $f(\oy)$ vers $\ox$, c'est-à-dire, d'un $X$-morphisme 
$u\colon \oy\rightarrow X_{(\ox)}$, où $X_{(\ox)}$ désigne le localisé strict de $X$ en $\ox$. 
Un tel point sera noté $(\oy\rightsquigarrow \ox)$ ou encore $(u\colon \oy\rightsquigarrow \ox)$. 

On désigne par $\cP_{(\oy \rightsquigarrow \ox)}$ la catégorie des objets $(\oy \rightsquigarrow \ox)$-pointés de $\cD_{f^+}$,
définie comme suit. Les objets de $\cP_{(\oy \rightsquigarrow \ox)}$ sont les triplets $((V\rightarrow U), \xi,\zeta)$ 
formés d'un objet $(V\rightarrow U)$ de $\cD_{f^+}$, d'un $X$-morphisme $\xi\colon \ox\rightarrow U$ et 
d'un $Y$-morphisme $\zeta\colon \oy\rightarrow V$ tels que, notant encore 
$\xi\colon X_{(\ox)}\rightarrow U$ le $X$-morphisme induit par $\xi$ (\cite{sga4} VIII 7.3), le diagramme
\begin{equation}\label{rec15a}
\xymatrix{
\oy\ar[r]^-(0.5)u\ar[d]_{\zeta}&{X_{(\ox)}}\ar[d]^{\xi}\\
V\ar[r]&U}
\end{equation}
soit commutatif. 
Soient $((V\rightarrow U), \xi,\zeta)$, $((V'\rightarrow U'), \xi',\zeta')$ deux objets de $\cP_{(\oy \rightsquigarrow \ox)}$.
Un morphisme de $((V'\rightarrow U'), \xi',\zeta')$ dans $((V\rightarrow U), \xi,\zeta)$ est la donnée d'un 
morphisme $(g\colon U'\rightarrow U, h\colon V'\rightarrow V)$ de $\cD_{f^+}$ \eqref{rec1} 
tel que $g\circ \xi'=\xi$ et $h\circ \zeta'=\zeta$. 
Il résulte de \eqref{rec8f} et (\cite{agt} VI.10.19(i)) que $\cP_{(\oy \rightsquigarrow \ox)}$ est canoniquement équivalente à 
la catégorie des voisinages de 
$(\oy \rightsquigarrow \ox)$ dans $\cD_{f^+}$ (\cite{sga4} IV 6.8.2). Elle est donc cofiltrante et pour tout préfaisceau 
$F=\{U\mapsto F_U\}$ sur $\cD_{f^+}$, on a un isomorphisme canonique fonctoriel (\cite{sga4} IV (6.8.4))
\begin{equation}\label{rec15b}
(F^a)_{(\oy \rightsquigarrow \ox)}\stackrel{\sim}{\rightarrow} \underset{\underset{((V\rightarrow U), \xi,\zeta)\in 
\cP^\circ_{(\oy \rightsquigarrow \ox)}}{\longrightarrow}}\lim\ F_U(V).
\end{equation}

\subsection{}\label{rec16}
Soient $f\colon Y\rightarrow X$ un morphisme de schémas, 
$\ox$ un point géométrique de $X$, $\uX$ le localisé strict de $X$ en $\ox$, 
$\uY=Y\times_X\uX$, $\uf\colon \uY\rightarrow \uX$ la projection canonique. 
On désigne par $X_\et\gtimes_{X_\et}Y_\et$ et $\uX_\et\gtimes_{\uX_\et}\uY_\et$ les topos co-évanescents de 
$f$ et $\uf$, respectivement, et par 
\begin{equation}\label{rec16a}
\gamma\colon \uY_\et\rightarrow \uX_\et\gtimes_{\uX_\et}\uY_\et
\end{equation}
le morphisme défini dans \eqref{rec10c}. 
Le morphisme canonique $\uX\rightarrow X$ induit par fonctorialité un morphisme
\eqref{rec3j}
\begin{equation}\label{rec16b}
\Phi\colon \uX_\et\gtimes_{\uX_\et}\uY_\et\rightarrow X_\et\gtimes_{X_\et}Y_\et.
\end{equation}
On note 
\begin{equation}\label{rec16c}
\varphi_\ox\colon X_\et\gtimes_{X_\et}Y_\et\rightarrow \uY_\et
\end{equation}
le foncteur composé $\gamma^*\circ \Phi^*$.

\begin{prop}\label{rec17}
Conservons les hypothèses de \ref{rec16}. Alors,
\begin{itemize}
\item[{\rm (i)}] Pour tout faisceau $F$ de $X_\et\gtimes_{X_\et}Y_\et$, on a un isomorphisme canonique fonctoriel 
\begin{equation}\label{rec17a}
\rp_{1*}(F)_\ox\stackrel{\sim}{\rightarrow}\Gamma(\uY_\et,\varphi_\ox(F)).
\end{equation}
\item[{\rm (ii)}] Pour tout faisceau abélien $F$ de $X_\et\gtimes_{X_\et}Y_\et$ 
et tout entier $i\geq 0$, on a un isomorphisme canonique fonctoriel 
\begin{equation}\label{rec17b}
\rR^i\rp_{1*}(F)_\ox\stackrel{\sim}{\rightarrow}\rH^i(\uY_\et,\varphi_\ox(F)).
\end{equation}
\end{itemize}
\end{prop}

Cette proposition sera démontrée dans \ref{lptce7}. 

\begin{rema}\label{rec18}
Conservons les hypothèses de \ref{rec16}, soient, de plus, $\oy$ un point géométrique de $Y$, 
$u\colon \oy\rightarrow \uX$ un $X$-morphisme, de sorte que   
$(\oy \rightsquigarrow \ox)$ est un point de $X_\et\gtimes_{X_\et}Y_\et$.
Notons $\tx$ le point fermé de $\uX$, $v\colon \oy\rightarrow \uY$ le morphisme induit par $u$,
$\ty$ le point géométrique de $\uY$ défini par $v$ et $\nu_\ty\colon \uY_\et\rightarrow \Ens$ 
le foncteur fibre associé à $\ty$. Le foncteur composé 
\begin{equation}\label{rec18a}
\nu_\ty\circ \varphi_\ox \colon X_\et\gtimes_{X_\et}Y_\et\rightarrow \Ens
\end{equation}
est alors canoniquement isomorphe au foncteur fibre associé au point $(\oy \rightsquigarrow \ox)$. 
\end{rema}

\subsection{}\label{rec19}
Reprenons les hypothèses et notations de \ref{rec16}. On note  $\cC_\ox$ la catégorie des $X$-schémas étales 
$\ox$-pointés (\cite{sga4} VIII 3.9) que l'on identifie à la catégorie des voisinages de $\ox$ dans le site 
$\Et_{/X}$ (\cite{sga4} IV 6.8.2). C'est une catégorie cofiltrante. 
Pour tout objet $(U,\xi\colon \ox\rightarrow U)$ de $\cC_\ox$, 
on désigne encore par $\xi\colon \uX\rightarrow U$ le morphisme déduit de
$\xi$ (\cite{sga4} VIII 7.3) et par
\begin{equation}\label{rec19a}
\xi_Y\colon \uY\rightarrow U_Y
\end{equation}
son changement de base par $f$. Soit $F=\{U\mapsto F_U\}$ un préfaisceau sur $\cD_{f^+}$ 
tel que pour tout $U\in \ob(\Et_{/X})$, $F_U$ soit un faisceau de $(U_Y)_\et$. On a alors un foncteur
\begin{equation}\label{rec19b}
\cC_\ox^\circ \rightarrow \uY_\et, \ \ \ (U,\xi)\mapsto \xi_Y^*(F_U),
\end{equation} 
qui à tout morphisme $t\colon (U',\xi')\rightarrow (U,\xi)$ de $\cC_\ox$ fait correspondre le morphisme composé
\[
\xi_Y^*(F_U)\stackrel{\sim}{\rightarrow} \xi'^*_Y(t_Y^*(F_U))\rightarrow \xi'^*_Y(F_{U'}),
\]
où la première flèche est induite par la relation $\xi=t\circ \xi'$ et la seconde flèche provient du
morphisme de transition $F_U\rightarrow t_{Y*}(F_{U'})$ de $F$. 

Pour tout $(U,\xi)\in \ob(\cC_\ox)$, le topos $(X_\et\gtimes_{X_\et}Y_\et)_{/(U_Y\rightarrow U)^a}$ est canoniquement équivalent 
au topos co-évanescent $U_\et\gtimes_{U_\et}(U_Y)_\et$ du morphisme $f_U\colon U_Y\rightarrow U$ \eqref{rec5}. Notons 
\begin{equation}\label{rec19c}
\jmath_{U}\colon (X_\et\gtimes_{X_\et}Y_\et)_{/(U_Y\rightarrow U)^a}\rightarrow X_\et\gtimes_{X_\et}Y_\et
\end{equation}
le morphisme de localisation de $X_\et\gtimes_{X_\et}Y_\et$ en $(U_Y\rightarrow U)^a$. 
Le morphisme $\xi\colon \uX\rightarrow U$ induit donc par fonctorialité un morphisme
de topos \eqref{rec3j}
\begin{equation}\label{rec19e}
\Phi_{\xi}\colon \uX_\et\gtimes_{\uX_\et}\uY_\et\rightarrow (X_\et\gtimes_{X_\et}Y_\et)_{/(U_Y\rightarrow U)^a}.
\end{equation}
D'après \eqref{rec3k}, le diagramme 
\begin{equation}\label{rec19f}
\xymatrix{
{\uX_\et\gtimes_{\uX_\et}\uY_\et}\ar[r]^-(0.5){\Phi_\xi}\ar[d]_{\urp_2}&{(X_\et\gtimes_{X_\et}Y_\et)_{/(U_Y\rightarrow U)^a}}
\ar[d]^{\rp_2^U}\\
{\uY_\et}\ar[r]^-(0.5){\xi_Y}&{(U_Y)_\et}}
\end{equation} 
où $\urp_2$ et $\rp_2^U$ sont les projections canoniques \eqref{rec1d} est commutatif à isomorphisme canonique près. 
On en déduit un morphisme de changement de base 
\begin{equation}\label{rec19g}
\xi_Y^*\rp^U_{2*}\rightarrow \urp_{2*}\Phi_\xi^*.
\end{equation}
D'après (\cite{egr1} 1.2.4(i)), le composé 
\begin{equation}\label{rec19h}
\xi_Y^*\rp^U_{2*}\rightarrow \urp_{2*}\Phi_\xi^* \rightarrow \gamma^*\Phi_\xi^*,
\end{equation}
où la seconde flèche est induite par \eqref{rec10e}, est le morphisme de changement de base induit par 
l'isomorphisme 
\begin{equation}\label{rec19i}
\rp_2^U\circ \Phi_\xi\circ \gamma\stackrel{\sim}{\rightarrow}\xi_Y
\end{equation}
déduit de la relation $\urp_2\gamma=\id_Y$. 

Soit $F=\{U\mapsto F_U\}$ un objet de $X_\et\gtimes_{X_\et}Y_\et$. 
D'après (\cite{sga4} III 5.3), on a un isomorphisme canonique 
\begin{equation}\label{rec19j}
\rp^U_{2*}(\jmath_U^*(F))\stackrel{\sim}{\rightarrow} F_U. 
\end{equation}
Le morphisme \eqref{rec19h} induit donc un morphisme fonctoriel 
\begin{equation}\label{rec19k}
\xi_Y^*(F_U)\rightarrow \varphi_\ox(F).
\end{equation}
D'après (\cite{egr1} 1.2.4(i)), c'est un morphisme de systèmes projectifs sur la catégorie $\cC_\ox^\circ$ \eqref{rec19b}. 
On en déduit un morphisme fonctoriel 
\begin{equation}\label{rec19l}
\underset{\underset{(U,\xi)\in \cC_\ox^\circ}{\longrightarrow}}{\lim}\ \xi_Y^*(F_U) \rightarrow
\varphi_\ox(F).
\end{equation}

\begin{remas}\label{rec20}
Conservons les hypothèses de \ref{rec19}. 
\begin{itemize}
\item[(i)] D'après (\cite{sga4} XVII 2.1.3), le morphisme \eqref{rec19h} est égal au composé 
\begin{equation}\label{rec20a}
\xi_Y^*\rp^U_{2*}\stackrel{\sim}{\rightarrow}
\gamma^*\Phi_\xi^*\rp_2^{U*}\rp^U_{2*}\rightarrow \gamma^*\Phi_\xi^*,
\end{equation}
où la première flèche est induite par \eqref{rec19i} et la seconde flèche est induite 
par le morphisme d'adjonction $\rp^{U*}_2\rp^U_{2*}\rightarrow\id$.
\item[(ii)] Choisissons un objet affine $(X_0,\xi_0)$ de $\cC_\ox$ et notons $I$ la catégorie
des $X_0$-schémas étales $\xi_0$-pointés qui sont affines au-dessus de $X_0$. 
Le foncteur canonique $I\rightarrow \cC_\ox$ est alors cofinal (\cite{sga4} VIII 4.5). 
On peut donc remplacer dans la limite inductive dans \eqref{rec19l} la catégorie $\cC_\ox$ par $I$. 
\end{itemize}
\end{remas}

\begin{prop}\label{rec21}
Les hypothèses étant celles de \ref{rec19}, soient, de plus $F=\{U\mapsto F_U\}$ un préfaisceau sur $\cD_{f^+}$, 
$F^a$ le faisceau de $X_\et\gtimes_{X_\et}Y_\et$ associé à $F$, et pour tout $U\in \ob(\Et_{/X})$, 
$F_U^a$ le faisceau de $(U_Y)_\et$ associé à $F_U$. Alors, on a un isomorphisme canonique fonctoriel
\begin{equation}\label{rec21a}
\underset{\underset{(U,\xi)\in \cC_\ox^\circ}{\longrightarrow}}{\lim}\ \xi_Y^*(F^a_U) 
\stackrel{\sim}{\rightarrow} \varphi_\ox(F^a),
\end{equation}
qui coïncide avec le morphisme \eqref{rec19l} lorsque $F$ est un faisceau sur $\cD_{f^+}$. 
\end{prop}

On notera d'abord que $\{U\mapsto F_U^a\}$ est naturellement un préfaisceau sur $\cD_{f^+}$ et que le morphisme 
canonique $\{U\mapsto F_U\}\rightarrow \{U\mapsto F_U^a\}$ induit un isomorphisme entre les faisceaux
associés, d'après (\cite{agt} VI.5.17 et VI.5.5). 
On peut donc se borner au cas où $F_U$ est un objet de $(U_Y)_\et$ pour tout $U\in \ob(\Et_{/X})$, de sorte 
que le système inductif $\cC_\ox^\circ \rightarrow \uY_\et$, $(U,\xi)\mapsto \xi_Y^*(F_U)$
\eqref{rec19b} est bien défini. Posant $F^a=\{U\mapsto G_U\}$, on a des morphismes canoniques fonctoriels 
\begin{equation}\label{rec21b}
\underset{\underset{(U,\xi)\in \cC_\ox^\circ}{\longrightarrow}}{\lim}\ (\xi_Y)_\fet^*(F_U) \stackrel{u}{\rightarrow}
\underset{\underset{(U,\xi)\in \cC_\ox^\circ}{\longrightarrow}}{\lim}\ (\xi_Y)_\fet^*(G_U) \stackrel{v}{\rightarrow} 
\varphi_\ox(F^a),
\end{equation}
le second étant \eqref{rec19l}. Montrons que $v\circ u$ est un isomorphisme. 
Il suffit de montrer que pour tout point géométrique $\ty$ de $\uY$, 
notant $\nu_{\ty}\colon \uY_\et\rightarrow \Ens$  le foncteur fibre associé, $\nu_\ty(v\circ u)$ est un isomorphisme. 
Soient $\pr_Y\colon \uY\rightarrow Y$ la projection canonique, $\oy=\pr_Y(\ty)$. 
Le point $\ty$ définit une flèche de spécialisation de $f(\oy)$ dans $\ox$, 
{\em i.e.}, $(\oy \rightsquigarrow \ox)$ est un point de $X_\et\gtimes_{X_\et}Y_\et$.
Notons $\cP_{(\oy \rightsquigarrow \ox)}$ la catégorie des objets $(\oy \rightsquigarrow \ox)$-pointés de 
$\cD_{f^+}$. On a un foncteur 
\begin{equation}
\pi\colon \cP_{(\oy \rightsquigarrow \ox)}\rightarrow \cC_\ox, \ \ \ ((V\rightarrow U),\xi,\zeta)\mapsto (U,\xi).
\end{equation}
Pour tout $(U,\xi)\in \ob(\cC_\ox)$, la fibre de $\pi$ au-dessus de $(U,\xi)$ est canoniquement équivalente 
à la catégorie $\cD^\ty_{(U,\xi)}$ des $U_Y$-schémas étales, $\xi_Y(\ty)$-pointés \eqref{rec19a}. 
D'après \ref{rec18}, $\nu_\ty\circ \varphi_\ox$ est le foncteur fibre de $X_\et\gtimes_{X_\et}Y_\et$ 
associé au point $(\oy \rightsquigarrow \ox)$. Compte tenu de (\cite{sga4} IV (6.8.3)), 
on peut alors considérer $\nu_\ty(v\circ u)$ comme une application 
\begin{equation}
\nu_\ty(v\circ u)\colon \underset{\underset{(U,\xi)\in \cC_\ox^\circ}{\longrightarrow}}{\lim}\ 
\underset{\underset{(V,\zeta)\in (\cD^\ty_{(U,\xi)})^\circ}{\longrightarrow}}{\lim}\ F_U(V)
\rightarrow (F^a)_{(\oy\rightsquigarrow \ox)}.
\end{equation}
Il suffit de montrer que cette application est l'inverse de l'isomorphisme canonique \eqref{rec15b}, ou encore
que pour tout objet $((V\rightarrow U),\xi,\zeta)$ de $\cP_{(\oy\rightsquigarrow \ox)}$, l'application canonique 
\begin{equation}
F_U(V)\rightarrow (F^a)_{(\oy\rightsquigarrow \ox)}
\end{equation}
est le composé
\begin{equation}
F_U(V)\rightarrow \nu_\ty(\xi_Y^*(F_U))\rightarrow \nu_\ty(\varphi_\ox(F^a)),
\end{equation}
où la première flèche est définie par l'objet $(V,\zeta)$ de $\cD^\ty_{(U,\xi)}$ et la seconde flèche est l'image de \eqref{rec19k}
par $\nu_\ty$. Par localisation \eqref{rec5}, on peut supposer $U=X$. D'après \ref{rec20}(i), 
l'image par $\nu_\ty$ du morphisme $(\pr_Y)^*_\fet(F_X)\rightarrow \varphi_\ox(F)$ \eqref{rec19k} 
coïncide avec l'application 
\begin{equation}
(F_X)_{\oy}\rightarrow F_{(\oy \rightsquigarrow \ox)}
\end{equation}
induite par le morphisme d'adjonction $\rp_2^*(\rp_{2*}(F))\rightarrow F$ et l'isomorphisme \eqref{rec8b}, d'où l'assertion recherchée.

\begin{prop}\label{rec23}
Soient $X$, $Y$ deux schémas cohérents, $f\colon Y\rightarrow X$ un morphisme. 
Pour tout point géométrique $\ox$ de $X$, soient $X_{(\ox)}$ le localisé strict de $X$ en $\ox$, 
$Y_{(\ox)}=Y\times_XX_{(\ox)}$,
\begin{equation}\label{rec23a}
\varphi_\ox\colon X_\et\gtimes_{X_\et}Y_\et\rightarrow Y_{(\ox),\et}
\end{equation}
le foncteur défini dans \eqref{rec16c}. Alors, la famille des foncteurs $(\varphi_\ox)$,
lorsque $\ox$ décrit l'ensemble des points géométriques de $X$, est conservative. 
\end{prop}

En effet, la famille des foncteurs fibres de $X_\et\gtimes_{X_\et}Y_\et$ 
associés aux points $(\oy \rightsquigarrow \ox)$ est conservative d'après \ref{rec9}. 
La proposition s'ensuit compte tenu de \ref{rec18}.

\section{Limite projective de topos co-évanescents} \label{lptce}

\subsection{}\label{lptce1} 
On désigne par $\fM$ la catégorie des morphismes de schémas et par $\fC$ la catégorie des morphismes 
de $\fM$. Les objets de $\fC$ sont donc des diagrammes commutatifs de morphismes de schémas 
\begin{equation}\label{lptce1a} 
\xymatrix{
V\ar[r]\ar[d]&U\ar[d]\\
Y\ar[r]&X}
\end{equation} 
où l'on considère les flèches horizontales comme des objets de $\fM$
et les flèches verticales comme des morphismes de $\fM$; un tel objet sera noté $(V,U,Y,X)$.   
On désigne par $\fD$ la sous-catégorie pleine de $\fC$ formée des objets $(V,U,Y,X)$ tels que les morphismes
$U\rightarrow X$ et $V\rightarrow Y$ soient étales de présentation finie. Le ``foncteur but''
\begin{equation}\label{lptce1b} 
\fD\rightarrow \fM, \ \ \ (V,U,Y,X)\mapsto (Y\rightarrow X),
\end{equation}
fait de $\fD$ une catégorie fibrée, clivée et normalisée. Avec les notations de \ref{notconv10}, tout morphisme de schémas $f\colon Y\rightarrow X$ induit un foncteur
\begin{equation}
f^+_\coh\colon \Et_{\coh/X}\rightarrow \Et_{\coh/Y}.
\end{equation}
La catégorie fibre $\fD_f$ de \eqref{lptce1b} au-dessus de $f$ est la catégorie sous-jacente au site co-évanescent $\cD_{f^+_\coh}$ associé 
au foncteur $f^+_\coh$ \eqref{rec1}. Pour tout diagramme commutatif de morphismes de schémas
\begin{equation}\label{lptce1c}
\xymatrix{
Y'\ar[r]^{f'}\ar[d]_{g'}&X'\ar[d]^g\\
Y\ar[r]^f&X}
\end{equation}
le foncteur image inverse de \eqref{lptce1b} associé au morphisme $(g',g)$ de $\fM$ 
est le foncteur 
\begin{equation}\label{lptce1d}
\Phi^+\colon \fD_f\rightarrow \fD_{f'}, \ \ \ (V\rightarrow U)\mapsto (V\times_YY'\rightarrow U\times_XX').
\end{equation}
Celui-ci est continu et exact à gauche d'après \ref{rec3}. 
Munissant chaque fibre de $\fD/\fM$ de la topologie co-évanescente \eqref{rec1}, 
$\fD$ devient un $\mU$-site fibré (\cite{sga4} VI 7.2.4). On désigne par 
\begin{equation}\label{lptce1e}
\fF\rightarrow \fM
\end{equation}
le $\mU$-topos fibré associé à $\fD/\fM$ (\cite{sga4} VI 7.2.6).  
La catégorie fibre de $\fF$ au-dessus d'un objet $f\colon Y\rightarrow X$
de $\fM$ est le topos $\tfD_f$ des faisceaux de $\mU$-ensembles sur le site co-évanescent $\fD_f$, 
et le foncteur image inverse relatif au morphisme défini par 
le diagramme \eqref{lptce1c} est le foncteur $\Phi^*\colon \tfD_{f}\rightarrow \tfD_{f'}$ image inverse par le morphisme 
de topos $\Phi\colon \tfD_{f'}\rightarrow \tfD_{f}$ associé au morphisme de sites $\Phi^+$ \eqref{lptce1d}. 
On note 
\begin{equation}\label{lptce1f}
\fF^\vee\rightarrow \fM^\circ
\end{equation}
la catégorie fibrée obtenue en associant à tout objet $f\colon Y\rightarrow X$
de $\fM$ la catégorie $\fF_{f}=\tfD_{f}$, et à tout morphisme défini par 
un diagramme \eqref{lptce1c} le foncteur $\Phi_*\colon \tfD_{f'}\rightarrow \tfD_{f}$ image directe par le morphisme 
de topos $\Phi\colon \tfD_{f'}\rightarrow \tfD_{f}$. 

\subsection{}\label{lptce2}
Soient $I$ une catégorie cofiltrante essentiellement petite (\cite{sga4} I 2.7 et 8.1.8), 
\begin{equation}\label{lptce2a} 
\varphi\colon I\rightarrow \fM, \ \ \ i\mapsto (f_i\colon Y_i\rightarrow X_i)
\end{equation}
un foncteur tel que pour tout morphisme $j\rightarrow i$ de $I$, 
les morphismes $Y_j\rightarrow Y_i$ et $X_j\rightarrow X_i$ soient affines. 
On suppose qu'il existe $i_0\in \ob(I)$ tel que $X_{i_0}$ et $Y_{i_0}$ soient cohérents.
On désigne par 
\begin{eqnarray}
\fD_\varphi&\rightarrow& I\label{lptce2b}\\
\fF_\varphi&\rightarrow& I\label{lptce2c}\\
\fF^\vee_\varphi&\rightarrow& I^\circ\label{lptce2d}
\end{eqnarray}
les site, topos et catégorie fibrés déduits de $\fD$ \eqref{lptce1b}, $\fF$ \eqref{lptce1e} et $\fF^\vee$ \eqref{lptce1f}, respectivement,
par changement de base par le foncteur $\varphi$. On notera que $\fF_\varphi$ est le topos fibré associé à $\fD_\varphi$
(\cite{sga4} VI 7.2.6.8). D'après (\cite{ega4} 8.2.3), les limites projectives 
\begin{equation}\label{lptce2e}
X=\underset{\underset{i\in \ob(I)}{\longleftarrow}}{\lim}\ X_i \ \ {\rm et}\ \ 
Y=\underset{\underset{i\in \ob(I)}{\longleftarrow}}{\lim}\ Y_i
\end{equation}
sont représentables dans la catégorie des schémas. Les morphismes $(f_i)_{i\in I}$ induisent 
un morphisme $f\colon Y\rightarrow X$, qui représente la limite projective du foncteur \eqref{lptce2a}. 

Pour tout $i\in \ob(I)$, on a un diagramme commutatif canonique
\begin{equation}\label{lptce2f}
\xymatrix{
Y\ar[r]^f\ar[d]&X\ar[d]\\
Y_i\ar[r]^{f_i}&X_i}
\end{equation}
Il lui correspond un foncteur image inverse \eqref{lptce1d}
\begin{equation}\label{lptce2g}
\Phi_i^+\colon \fD_{f_i}\rightarrow \fD_f,
\end{equation}
qui est continu et exact à gauche, et par suite un morphisme de topos 
\begin{equation}\label{lptce2h}
\Phi_i\colon \tfD_{f}\rightarrow \tfD_{f_i}.
\end{equation}

On a un foncteur naturel 
\begin{equation}\label{lptce2i}
\fD_\varphi \rightarrow \fD_f,
\end{equation}
dont la restriction à la fibre au-dessus de tout $i\in \ob(I)$ est le foncteur $\Phi_i^+$ \eqref{lptce2g}. 
Ce foncteur transforme morphisme cartésien en isomorphisme. Il se factorise donc de façon unique à travers 
un foncteur (\cite{sga4} VI 6.3)
\begin{equation}\label{lptce2j}
\underset{\underset{I^\circ}{\longrightarrow}}{\lim}\ \fD_\varphi \rightarrow \fD_{f}.
\end{equation}
Le $I$-foncteur $\fD_\varphi\rightarrow \fD_f\times I$ déduit de \eqref{lptce2i} est un morphisme cartésien
de sites fibrés (\cite{sga4} VI 7.2.2). Il induit donc un morphisme cartésien de topos fibrés (\cite{sga4} VI 7.2.7)
\begin{equation}\label{lptce2k}
\tfD_f\times I\rightarrow \fF_\varphi.
\end{equation}

\begin{prop}\label{lptce3}
Le couple formé du topos $\tfD_f$ et du morphisme \eqref{lptce2k} 
est une limite projective du topos fibré $\fF_\varphi/I$ {\rm (\cite{sga4} VI 8.1.1)}.
\end{prop}

On notera d'abord que le foncteur \eqref{lptce2j} est une équivalence de catégories en vertu de 
(\cite{ega4} 8.8.2 et 17.7.8). Soient $T$ un $\mU$-topos, 
\begin{equation}
h\colon T\times I\rightarrow \fF_\varphi
\end{equation}
un morphisme cartésien de topos fibrés au-dessus de $I$. Notons $\varepsilon_I\colon \fD_\varphi\rightarrow \fF_\varphi$
le foncteur cartésien canonique (\cite{sga4} VI (7.2.6.7)), et posons 
\begin{equation}
h^+= h^*\circ \varepsilon_I\colon \fD_\varphi\rightarrow T\times I.
\end{equation}
Pour tout $i\in \ob(I)$, on désigne par
\begin{equation}
h_i^+\colon \fD_{f_i}\rightarrow T
\end{equation}
la restriction de $h^+$ aux fibres au-dessus de $i$. 
Compte tenu de l'équivalence de catégories \eqref{lptce2j} et de  (\cite{sga4} VI 6.2), il existe un et un essentiellement unique foncteur 
\begin{equation}
g^+\colon \fD_{f}\rightarrow T
\end{equation}
tel que $h^+$ soit isomorphe au composé 
\begin{equation}
\xymatrix{
{\fD_\varphi}\ar[r]&{\fD_{f}\times I}\ar[rr]^{g^+\times \id_I}&&{T\times I}},
\end{equation}
où la première flèche est le foncteur déduit de \eqref{lptce2i}. Montrons que $g^+$ est un morphisme de sites. 
Pour tout objet $(V\rightarrow U)$ de $\fD_f$, il existe $i\in \ob(I)$, un objet $(V_i\rightarrow U_i)$ de $\fD_{f_i}$ et un isomorphisme de $\fD_f$
\begin{equation}
(V\rightarrow U)\stackrel{\sim}{\rightarrow} \Phi_i^+(V_i\rightarrow U_i).
\end{equation}
Comme les foncteurs $h_i^+$ et $\Phi_i^+$ sont exacts à gauche, on en déduit que $g^+$ est exact à gauche. 
D'autre part, tout recouvrement {\em fini} de type $(\alpha)$ (resp. $(\beta)$) de $\fD_f$ \eqref{rec1}  
est l'image inverse d'un recouvrement de type $(\alpha)$ (resp. $(\beta)$) de $\fD_{f_i}$ pour un objet $i\in I$,  
en vertu de (\cite{ega4} 8.10.5(vi)). Comme les schémas $X$ et $Y$ sont cohérents, on en déduit 
que $g^+$ transforme les recouvrements de type $(\alpha)$ (resp. $(\beta)$) de $\fD_f$ en familles épimorphiques de $T$. 
Par suite, $g^+$ est continu en vertu de (\cite{agt} VI.4.4).  Il définit donc un morphisme de topos 
\begin{equation}
g\colon T\rightarrow \tfD_f
\end{equation}
tel que $h$ soit isomorphe au composé 
\begin{equation}
\xymatrix{
{T\times I}\ar[rr]^{g\times \id_I}&&{\tfD_f\times I}\ar[r]&{\fF_\varphi}},
\end{equation}
où la seconde flèche est le morphisme \eqref{lptce2k}. Un tel morphisme $g$ est essentiellement unique car la ``restriction'' 
$g^+\colon \fD_f\rightarrow T$ du foncteur $g^*$ est essentiellement unique d'après ce qui précède, d'où la proposition.

\subsection{}\label{lptce4}
Munissons $\fD_\varphi$ de la topologie totale (\cite{sga4} VI 7.4.1) 
et notons $\Top(\fD_\varphi)$ le topos des faisceaux de $\mU$-ensembles sur $\fD_\varphi$. 
D'après (\cite{sga4} VI 7.4.7), on a une équivalence canonique de catégories \eqref{notconv4}
\begin{equation}\label{lptce4a}
\Top(\fD_\varphi)\stackrel{\sim}{\rightarrow}\bHom_{I^\circ}(I^\circ, \fF^\vee_\varphi). 
\end{equation}
D'autre part, le foncteur naturel $\fD_\varphi\rightarrow \fD_f$ \eqref{lptce2i} est un morphisme de sites (\cite{sga4} VI 7.4.4)
et définit donc un morphisme de topos 
\begin{equation}\label{lptce4b}
\varpi\colon \tfD_f\rightarrow \Top(\fD_\varphi). 
\end{equation}
En vertu de \ref{lptce3} et (\cite{sga4} VI 8.2.9), il existe une équivalence de catégories $\Theta$
qui s'insère dans un diagramme commutatif 
\begin{equation}\label{lptce4c}
\xymatrix{
{\tfD_f}\ar[r]^-(0.5){\Theta}_-(0.5)\sim\ar[d]_{\varpi_*}&{\bHom_{\cart/I^\circ}(I^\circ, \fF^\vee_\varphi)}\ar@{^(->}[d]\\
{\Top(\fD_\varphi)}\ar[r]^-(0.5)\sim&{\bHom_{I^\circ}(I^\circ, \fF^\vee_\varphi)}}
\end{equation}
où la flèche horizontale inférieure est l'équivalence de catégories \eqref{lptce4a} et 
la flèche verticale de droite est l'injection canonique. 

Pour tout objet $F$ de $\Top(\fD_\varphi)$, si 
$\{i\mapsto F_i\}$ est la section correspondante de $\bHom_{I^\circ}(I^\circ, \fF^\vee_\varphi)$, on a un 
isomorphisme canonique fonctoriel (\cite{sga4} VI 8.5.2)
\begin{equation}\label{lptce4d}
\varpi^*(F)\stackrel{\sim}{\rightarrow} \underset{\underset{i\in I^\circ}{\longrightarrow}}{\lim}\ \Phi_i^*(F_i).
\end{equation}

\begin{cor}\label{lptce5}
Soient $F$ un faisceau de $\Top(\fD_\varphi)$, 
$\{i\mapsto F_i\}$ la section de $\bHom_{I^\circ}(I^\circ, \fF^\vee_\varphi)$ qui lui est 
associée par l'équivalence de catégories \eqref{lptce4a}.
Alors on a un isomorphisme canonique fonctoriel
\begin{equation}\label{lptce5a}
\underset{\underset{i\in I^\circ}{\longrightarrow}}{\lim}\ \Gamma(\tfD_{f_i},F_i)\stackrel{\sim}{\rightarrow}
\Gamma(\tfD_{f},\underset{\underset{i\in I^\circ}{\longrightarrow}}{\lim}\ \Phi_i^*(F_i)).
\end{equation}
\end{cor}

\begin{cor}\label{lptce6}
Soit $F$ un faisceau abélien de $\Top(\fD_\varphi)$ et soit
$\{i\mapsto F_i\}$ la section de $\bHom_{I^\circ}(I^\circ, \fF^\vee_\varphi)$ qui lui est 
associée par l'équivalence de catégories \eqref{lptce4a}.
Alors pour tout entier $q\geq 0$, on a un isomorphisme canonique fonctoriel
\begin{equation}\label{lptce6a}
\underset{\underset{i\in I^\circ}{\longrightarrow}}{\lim}\ \rH^q(\tfD_{f_i},F_i)\stackrel{\sim}{\rightarrow}
\rH^q(\tfD_{f},\underset{\underset{i\in I^\circ}{\longrightarrow}}{\lim}\ \Phi_i^*(F_i)).
\end{equation}
\end{cor}

Les corollaires \ref{lptce5} et \ref{lptce6} résultent de \ref{lptce3} et (\cite{sga4} VI 8.7.7). 
On notera que les conditions requises dans 
(\cite{sga4} VI 8.7.1 et 8.7.7) sont satisfaites en vertu de \ref{rec22} et (\cite{sga4} VI 3.3, 5.1 et 5.2).

\subsection{}\label{lptce7}
Nous pouvons maintenant démontrer la proposition \ref{rec17}. Choisissons un voisinage étale affine 
$X_0$ de $\ox$ dans $X$. Notons $I$ la catégorie des $X_0$-schémas étales $\ox$-pointés
qui sont affines au-dessus de $X_0$ (cf. \cite{sga4} VIII 3.9 et 4.5), et $\varphi \colon I\rightarrow \fM$ 
le foncteur qui à un objet $U$ de $I$, associe la projection canonique $f_U\colon U_Y\rightarrow U$. 
Alors, $\uf$ s'identifie canoniquement à la limite projective du foncteur $\varphi$. 
Pour tout $U\in \ob(I)$, le topos co-évanescent $U_\et\gtimes_{U_\et}(U_Y)_\et$ est canoniquement équivalent au topos $\tfD_{f_U}$ \eqref{lptce1}.
Notant $(U_Y\rightarrow U)^a$ le faisceau de $X_\et\gtimes_{X_\et}Y_\et$ associé à l'objet $f_U$ de $\fD_f$, 
le topos $U_\et\gtimes_{U_\et}(U_Y)_\et$ est aussi canoniquement équivalent au topos $(X_\et\gtimes_{X_\et}Y_\et)_{/(U_Y\rightarrow U)^a}$, 
en vertu de \ref{rec5}. Par suite, avec les notations de cette section, pour tout faisceau $F$ de $X_\et\gtimes_{X_\et}Y_\et$, 
$\{U\mapsto F|(U_Y\rightarrow U)^a\}$ est naturellement une section de $\bHom_{I^\circ}(I^\circ,\fF_\varphi^\vee)$.
Elle définit donc un faisceau de $\Top(\fD_\varphi)$ \eqref{lptce4a}. On a un isomorphisme canonique fonctoriel \eqref{lptce4d}
\begin{equation}
\Phi^*(F)\stackrel{\sim}{\rightarrow} \varpi^*(\{U\mapsto F|(U_Y\rightarrow U)^a\}).
\end{equation}

(i) D'après (\cite{sga4} VIII 3.9 et 4.5), on a un isomorphisme canonique 
\begin{equation}\label{lptce7a}
\rp_{1*}(F)_\ox\stackrel{\sim}{\rightarrow} \underset{\underset{U\in I^\circ}{\longrightarrow}}{\lim}\ 
\Gamma((U_Y\rightarrow U)^a,F).
\end{equation}
Celui-ci induit un isomorphisme fonctoriel
\begin{equation}\label{lptce7b}
\rp_{1*}(F)_\ox\stackrel{\sim}{\rightarrow} \Gamma(\uX_\et\gtimes_{\uX_\et}\uY_\et,\Phi^*(F)),
\end{equation}
en vertu de \ref{lptce5}. La proposition s'en déduit compte tenu de \ref{rec14}(i).

(ii) Cela résulte, comme pour (i), de \ref{rec14}(ii), \ref{lptce6} et de l'isomorphisme canonique (\cite{sga4} V 5.1(1))
\begin{equation}
\rR^i\rp_{1*}(F)_\ox\stackrel{\sim}{\rightarrow} \underset{\underset{U\in I^\circ}{\longrightarrow}}{\lim}\ 
\rH^i((U_Y\rightarrow U)^a,F).
\end{equation}

\section{Topos de Faltings}\label{tf}

On renvoie à (\cite{agt} VI.10) pour la définition et les principales propriétés du topos de Faltings. 

\subsection{}\label{tf1}
Dans cette section, $f\colon Y\rightarrow X$ désigne un morphisme de schémas et
\begin{equation}\label{tf1a}
\pi\colon E\rightarrow \Et_{/X}
\end{equation}
le site fibré de Faltings associé (\cite{agt} VI.10.1). On rappelle que $E$ est la catégorie des morphismes 
de schémas $V\rightarrow U$ au-dessus de $f$ tels que le morphisme
$U\rightarrow X$ soit étale et que le morphisme $V\rightarrow U_Y=U\times_XY$ soit étale fini. 
Pour tout $(V\rightarrow U)\in \ob(E)$, on a 
\begin{equation}\label{tf1b}
\pi(V\rightarrow U)=U.
\end{equation}
La fibre de $E$ au-dessus de $U$ s'identifie canoniquement au site fini étale de $U_Y$ \eqref{notconv10}. On note 
\begin{equation}\label{tf1c}
\iota_{U!}\colon \Et_{\rf/U_Y}\rightarrow E, \ \ \ V\mapsto (V\rightarrow U)
\end{equation}
le foncteur canonique (\cite{agt} (VI.5.1.2)). 

On désigne par 
\begin{equation}\label{tf1d}
\fF\rightarrow  \Et_{/X}
\end{equation}
le $\mU$-topos fibré associé à $\pi$. La catégorie fibre de $\fF$ au-dessus de tout $U\in \ob(\Et_{/X})$ est 
canoniquement équivalente au topos fini étale $(U_Y)_\fet$ de $U_Y$ \eqref{notconv10}
et le foncteur image inverse 
pour tout morphisme $t\colon U'\rightarrow U$ de $\Et_{/X}$ s'identifie au foncteur 
$(t_Y)_{\fet}^*\colon (U_Y)_\fet\rightarrow (U'_Y)_\fet$ image inverse par le morphisme 
de topos $(t_Y)_\fet\colon (U'_Y)_\fet\rightarrow (U_Y)_\fet$ (\cite{agt} VI.9.3). On désigne par
\begin{equation}\label{tf1e}
\fF^\vee\rightarrow (\Et_{/X})^\circ
\end{equation}
la catégorie fibrée obtenue en associant à tout $U\in \ob(\Et_{/X})$ la catégorie $(U_Y)_\fet$, et à tout morphisme 
$t\colon U'\rightarrow U$ de $\Et_{/X}$ le foncteur 
$(t_Y)_{\fet*}\colon (U'_Y)_\fet\rightarrow (U_Y)_\fet$ 
image directe par le morphisme de topos $(t_Y)_\fet$. On désigne par
\begin{equation}\label{tf1f}
\cP^\vee\rightarrow (\Et_{/X})^\circ
\end{equation}
la catégorie fibrée obtenue en associant à tout $U\in \ob(\Et_{/X})$ la catégorie $(\Et_{\rf/U_Y})^\wedge$ 
des préfaisceaux de $\mU$-ensembles sur $\Et_{\rf/U_Y}$, et à tout morphisme $t\colon U'\rightarrow U$ de $\Et_{/X}$ 
le foncteur 
\begin{equation}\label{tf1g}
(t_Y)_{\fet*}\colon (\Et_{\rf/U'_Y})^\wedge\rightarrow (\Et_{\rf/U_Y})^\wedge
\end{equation} 
obtenu en composant avec le foncteur image inverse $t^+_Y\colon \Et_{\rf/U_Y}\rightarrow \Et_{\rf/U'_Y}$.  

On note $\hE$ la catégorie des préfaisceaux de $\mU$-ensembles sur $E$. On a alors
une équivalence de catégories (\cite{agt} VI.5.2)
\begin{eqnarray}\label{tf1h}
\hE&\rightarrow& \bHom_{(\Et_{/X})^\circ}((\Et_{/X})^\circ,\cP^\vee)\\
F&\mapsto &\{U\mapsto F\circ \iota_{U!}\}.\nonumber
\end{eqnarray}
On identifiera dans la suite $F$ à la section $\{U\mapsto F\circ \iota_{U!}\}$ qui lui est associée par cette équivalence.

On munit $E$ de la {\em topologie co-évanescente} définie par $\pi$ (\cite{agt} VI.5.3) 
et on note $\tE$ le topos des faisceaux de $\mU$-ensembles sur $E$. 
Les site et topos ainsi définis sont appelés {\em site et topos de Faltings}  
associés à $f$ (\cite{agt} VI.10.1). Si $F$ est un préfaisceau sur $E$, on note $F^a$ le faisceau associé. 
D'après (\cite{agt} VI.5.11), le foncteur \eqref{tf1h} induit un foncteur pleinement fidèle 
\begin{equation}\label{tf1i}
\tE\rightarrow \bHom_{(\Et_{/X})^\circ}((\Et_{/X})^\circ,\fF^\vee),
\end{equation}
d'image essentielle les sections $\{U\mapsto F_U\}$ vérifiant une condition de recollement. 

Les foncteurs 
\begin{eqnarray}
\sigma^+\colon \Et_{/X}\rightarrow E,&& U\mapsto (U_Y\rightarrow U),\label{tf1jj}\\
\alpha_{X!}\colon \Et_{\rf/Y}\rightarrow E,&& V\mapsto (V\rightarrow X).\label{tf1j}
\end{eqnarray}
sont exacts à gauche et continus (\cite{agt} VI.10.6). Ils définissent donc deux morphismes de topos 
\begin{eqnarray}
\sigma\colon \tE\rightarrow X_\et,\label{tf1kk}\\
\beta\colon \tE\rightarrow Y_\fet.\label{tf1k}
\end{eqnarray} 
Pour tout faisceau $F=\{U\mapsto F_U\}$ sur $E$, on a $\beta_*(F)=F_X$. 

Le foncteur 
\begin{equation}\label{tf1l}
\psi^+\colon E\rightarrow \Et_{/Y},\ \ \ (V\rightarrow U)\mapsto V
\end{equation}
est continu et exact à gauche (\cite{agt} VI.10.7). Il définit donc un morphisme de topos 
\begin{equation}\label{tf1m}
\psi\colon Y_\et\rightarrow \tE.
\end{equation}
Nous changeons ici de notations par rapport à {\em loc. cit.}, 
en réservant la notation $\Psi$ aux foncteurs des cycles co-proches dans le sens 
strict \eqref{rec1i}. Pour tout faisceau $F$ de $Y_\et$, on a un isomorphisme canonique de $\tE$
\begin{equation}\label{tf1n}
\psi_*(F)\stackrel{\sim}{\rightarrow} \{U\mapsto \rho_{U_Y*}(F|U_Y)\},
\end{equation}
où pour tout objet $U$ de $\Et_{/X}$, $\rho_{U_Y}\colon (U_Y)_\et\rightarrow (U_Y)_\fet$ est le morphisme canonique \eqref{notconv10a}.

\subsection{}\label{tf2}
On désigne par $\Et_{\coh/X}$ la sous-catégorie pleine de $\Et_{/X}$ formée des schémas étales 
de présentation finie sur $X$, munie de la topologie induite par celle de $\Et_{/X}$ \eqref{notconv10}, et par 
\begin{equation}
\pi_\coh\colon E_\coh\rightarrow \Et_{\coh/X}
\end{equation}
le site fibré déduit de $\pi$ \eqref{tf1a} par changement de base par le foncteur d'injection canonique 
\begin{equation}
\Et_{\coh/X}\rightarrow \Et_{/X}. 
\end{equation} 
On munit $E_\coh$ de la topologie co-évanescente définie par $\pi_\coh$ (\cite{agt} VI.5.3).
D'après (\cite{agt} VI.10.4), si $X$ est quasi-séparé, la projection canonique $E_\coh\rightarrow E$ induit par restriction une équivalence entre le topos 
$\tE$ et le topos des faisceaux de $\mU$-ensembles sur $E_\coh$. De plus, sous la même hypothèse,
la topologie co-évanescente de $E_\coh$ est induite par celle de $E$.

\subsection{}\label{tf3}
On désigne par $\cD_{f^+}$ le site co-évanescent du foncteur $f^+\colon \Et_{/X}\rightarrow \Et_{/Y}$ de changement de base par $f$, 
et par $X_\et\gtimes_{X_\et}Y_\et$ le topos co-évanescent du morphisme $f_\et\colon Y_\et\rightarrow X_\et$ (cf. \ref{rec1}). 
Tout objet de $E$ est naturellement un objet de $\cD_{f^+}$. 
On définit ainsi un foncteur pleinement fidèle et exact à gauche
\begin{equation}\label{tf3a}
\rho^+\colon E\rightarrow \cD_{f^+}.
\end{equation} 
Celui-ci est exact à gauche et continu d'après (\cite{agt} VI.10.15). Il définit donc un morphisme de topos 
\begin{equation}\label{tf3b}
\rho\colon X_\et\gtimes_{X_\et}Y_\et\rightarrow \tE.
\end{equation}
Pour tout faisceau $F=\{U\mapsto F_U\}$ de $X_\et\gtimes_{X_\et}Y_\et$ \eqref{rec4h}, on a un isomorphisme canonique
\begin{equation}
\rho_*(F)=\{U\mapsto \rho_{U_Y*}(F_U)\}.
\end{equation}
D'autre part, pour tout faisceau $G=\{U\mapsto G_U\}$ de $\tE$, $\rho^*(G)$
est le faisceau associé au préfaisceau sur $\cD_{f^+}$ défini par $\{U\mapsto \rho_{U_Y}^*(G_U)\}$ 
en vertu de (\cite{agt} VI.6.5).

Il résulte aussitôt des définitions qu'avec les notations de \ref{rec1}, les carrés du diagramme
\begin{equation}\label{tf3c}
\xymatrix{
{X_\et}\ar@{=}[d]&{X_\et\gtimes_{X_\et}Y_\et}\ar[l]_-(0.5){\rp_1}\ar[d]^{\rho}\ar[r]^-(0.5){\rp_2}&
{Y_\et}\ar[d]^{\rho_Y}\\
{X_\et}&{\tE}\ar[l]_{\sigma}\ar[r]^{\beta}&{Y_\fet}}
\end{equation}
et le diagramme 
\begin{equation}\label{tf3d}
\xymatrix{
{Y_\et}\ar[r]^-(0.5){\Psi}\ar[dr]_{\psi}&{X_\et\gtimes_{X_\et}Y_\et}\ar[d]^{\rho}\\
&{\tE}}
\end{equation}
sont commutatifs à isomorphismes canoniques près.

\begin{prop}\label{tf4}
Supposons $X$ et $Y$ cohérents. Alors,
\begin{itemize}
\item[{\rm (i)}] Pour tout objet $(V\rightarrow U)$ de $E_{\coh}$ \eqref{tf2}, le faisceau associé $(V\rightarrow U)^a$ est un objet cohérent de $\tE$. 
\item[{\rm (ii)}] Le topos $\tE$ est cohérent; en particulier, il a suffisamment de points. 
\item[{\rm (iii)}] Le morphisme $\rho$ \eqref{tf3b} est cohérent. 
\end{itemize}
\end{prop}

Les propositions (i) et (ii) sont mentionnées à titre de rappel (\cite{agt} VI.10.5). Montrons la proposition (iii). 
Avec les notations de \ref{notconv10}, le morphisme $f$ induit un foncteur 
\begin{equation}
f^+_\coh\colon \Et_{\coh/X}\rightarrow \Et_{\coh/Y}.
\end{equation}
On note $\cD_{f^+_\coh}$ le site co-évanescent associé à $f^+_\coh$ \eqref{rec1}. 
D'après \ref{rec3}, le foncteur d'injection canonique $\cD_{f^+_\coh}\rightarrow \cD_{f^+}$ \eqref{tf3} est continu et exact à gauche. 
Il induit une équivalence entre les topos associés. Pour tout objet $(V\rightarrow U)$ de $E_\coh$, $\rho^+(V\rightarrow U)$ est un objet de $\cD_{f^+_\coh}$.
La proposition (iii) résulte alors de (i), \ref{rec22}(ii) et (\cite{sga4} VI 3.2) appliqué à la famille topologiquement génératrice $E_\coh$ de $E$.

\begin{prop}[\cite{agt} VI.10.21]\label{tf5} Si les schémas $X$ et $Y$ sont cohérents,
lorsque $(\oy \rightsquigarrow \ox)$ décrit la famille des points de $X_\et\gtimes_{X_\et}Y_\et$ \eqref{rec15}, 
la famille des foncteurs fibres de $\tE$ associés aux points $\rho(\oy \rightsquigarrow \ox)$ est conservative
\eqref{tf3b} {\rm (\cite{sga4} IV 6.4.0)}. 
\end{prop}

On notera que sous les mêmes hypothèses, la famille des foncteurs fibres de $X_\et\gtimes_{X_\et}Y_\et$ 
associés aux points $(\oy \rightsquigarrow \ox)$ est conservative en vertu de \ref{rec9}.

\subsection{}\label{tf6}
Soient $\uX$ un $X$-schéma, $\uY=Y\times_X\uX$, $\uf\colon \uY\rightarrow \uX$ la projection canonique,
$\uX_\et\gtimes_{\uX_\et}\uY_\et$ le topos co-évanescent et $\tuE$ le topos de Faltings associés à $\uf$. On vérifie aussitôt que le diagramme
\begin{equation}\label{tf6a}
\xymatrix{
{\uX_\et\gtimes_{\uX_\et}\uY_\et}\ar[d]_{\urho}\ar[r]^{\Xi}&{X_\et\gtimes_{X_\et}Y_\et}\ar[d]^\rho\\
{\tuE}\ar[r]^{\Phi}&\tE}
\end{equation}
où les morphismes $\Xi$ et $\Phi$ sont définis par fonctorialité \eqref{rec3j} et (\cite{agt} (VI.10.12.5)), 
et $\rho$ et $\urho$ sont les morphismes canoniques \eqref{tf3b}, est commutatif à isomorphisme canonique près. 
On en déduit un morphisme de changement de base 
\begin{equation}\label{tf6b}
\fc\colon \Phi^*\rho_*\rightarrow \urho_*\Xi^*.
\end{equation}
En restreignant aux faisceaux abéliens, on déduit un morphisme pour tout entier $q\geq 0$ (cf. \cite{egr1} 1.2.2)
\begin{equation}\label{tf6c}
\fc^q\colon \Phi^*\rR^q\rho_*\rightarrow \rR^q\urho_*\Xi^*.
\end{equation}

\begin{prop}\label{tf7}
Conservons les hypothèses et notations de \ref{tf6}, supposons, de plus, que l'une des deux hypothèses suivantes soit satisfaite~:
\begin{itemize}
\item[{\rm (a)}] $\uX$ est étale au-dessus de $X$;
\item[{\rm (b)}] le morphisme $f\colon Y\rightarrow X$ est cohérent, et $\uX$ est le localisé strict de $X$ en un point géométrique $\ox$ de $X$. 
\end{itemize} 
Alors, 
\begin{itemize}
\item[{\rm (i)}] Pour tout faisceau $F$ de $X_\et\gtimes_{X_\et}Y_\et$, le morphisme de changement de base \eqref{tf6b}
\begin{equation}\label{tf7a}
\fc(F)\colon \Phi^*(\rho_*F)\rightarrow \urho_*(\Xi^*F)
\end{equation}
est un isomorphisme. 
\item[{\rm (ii)}] Pour tout faisceau abélien $F$ de $X_\et\gtimes_{X_\et}Y_\et$ et tout entier $q\geq 0$, le morphisme de changement de base \eqref{tf6c}
\begin{equation}\label{tf7b}
\fc^q(F)\colon \Phi^*(\rR^q\rho_*F)\rightarrow \rR^q\urho_*(\Xi^*F)
\end{equation}
est un isomorphisme. 
\end{itemize}
\end{prop}

Considérons d'abord le cas (a). Notons $(\uY\rightarrow \uX)^a$ le faisceau de $\tE$ associé à 
l'objet $(\uX\rightarrow \uY)$ de $E$. D'après (\cite{agt} VI.10.14), on a une équivalence canonique
\begin{equation}\label{tf7c}
\tuE\stackrel{\sim}{\rightarrow} \tE_{/(\uY\rightarrow \uX)^a},
\end{equation}
dont le composé avec le morphisme de localisation de $\tE$ en  $(\uY\rightarrow \uX)^a$ est égal à $\Phi$. D'après \ref{rec5}, 
on a une équivalence canonique
\begin{equation}\label{tf7d}
\uX_\et\gtimes_{\uX_\et}\uY_\et\stackrel{\sim}{\rightarrow} (X_\et\gtimes_{X_\et}Y_\et)_{/\rho^*((\uY\rightarrow \uX)^a)},
\end{equation}
dont le composé avec le morphisme de localisation de $X_\et\gtimes_{X_\et}Y_\et$ en  $\rho^*((\uY\rightarrow \uX)^a)$ est égal à $\Xi$.
De plus, $\urho$ s'identifie à la restriction de $\rho$ au-dessus de $(\uY\rightarrow \uX)^a$ (\cite{sga4} IV (5.10.3)). La proposition (i) s'ensuit aussitôt,
et la proposition (ii) est une conséquence de (\cite{sga4} V 5.1(3)).   

Considérons ensuite le cas (b). 
Choisissons un voisinage étale affine $X'_0$ de $\ox$ dans $X$ et notons $\fV_\ox$ la catégorie des $X'_0$-schémas étales $\ox$-pointés
qui sont affines au-dessus de $X'_0$ (cf. \cite{sga4} VIII 3.9 et 4.5). Tout objet $X'$ de $\fV_\ox$ donne lieu à un diagramme à carrés 
commutatifs à isomorphismes canoniques près
\begin{equation}\label{tf7e}
\xymatrix{
{\uX_\et\gtimes_{\uX_\et}\uY_\et}\ar[d]_{\urho}\ar[r]^-(0.5){\uXi_{X'}}&{X'_\et\gtimes_{X'_\et}(Y\times_XX')_\et}\ar[d]^{\rho_{X'}}\ar[r]^-(0.5){\Xi_{X'}}&
{X_\et\gtimes_{X_\et}Y_\et}\ar[d]^\rho\\
{\tuE}\ar[r]^{\uPhi_{X'}}&{\tE_{X'}}\ar[r]^{\Phi_{X'}}&\tE}
\end{equation}
où $X'_\et\gtimes_{X'_\et}(Y\times_XX')_\et$ est le topos co-évanescent et $\tE_{X'}$ est le topos de Faltings associés à la projection canonique 
$f'\colon Y\times_XX'\rightarrow X'$, $\rho_{X'}$ est le morphisme canonique \eqref{tf3b} et les flèches horizontales sont les morphismes de fonctorialité. 

Suivant \ref{lptce1}, on désigne par 
\begin{equation}\label{tf7f}
\fF\rightarrow \fV_\ox
\end{equation}
le topos fibré obtenu en associant à tout objet $X'$ de $\fV_\ox$ le topos $X'_\et\gtimes_{X'_\et}(Y\times_XX')_\et$, 
et à tout morphisme $X''\rightarrow X'$ de $\fV_\ox$ le foncteur 
\[
X'_\et\gtimes_{X'_\et}(Y\times_XX')_\et\rightarrow X''_\et\gtimes_{X''_\et}(Y\times_XX'')_\et
\] 
image inverse par le morphisme de fonctorialité. 
En vertu de \ref{lptce3}, les morphismes $\uXi_{X'}$ identifient le topos $\uX_\et\gtimes_{\uX_\et}\uY_\et$ à la limite projective du topos fibré $\fF$.

Suivant (\cite{agt} VI.11.1), on désigne par  
\begin{equation}\label{tf7g}
\fG\rightarrow \fV_\ox
\end{equation}
le topos fibré obtenu en associant à tout objet $X'$ de $\fV_\ox$ le topos $\tE_{X'}$, et à tout morphisme $X''\rightarrow X'$ de $\fV_\ox$ le foncteur 
$\tE_{X'}\rightarrow \tE_{X''}$ image inverse par le morphisme de fonctorialité. 
En vertu de (\cite{agt} VI.11.3), les morphismes $\uPhi_{X'}$ identifient le topos $\tuE$ à la limite projective du topos fibré $\fG$.

Les morphismes $\rho_{X'}$, pour $X'\in \ob(\fV_\ox)$, définissent un morphisme de topos fibrés (\cite{sga4} VI 7.1.6)
\begin{equation}\label{tf7h}
\varrho\colon \fF\rightarrow \fG.
\end{equation}
Le morphisme $\urho$ se déduit de $\varrho$ par passage à la limite projective (\cite{sga4} VI 8.1.4). 

D'après (\cite{sga4} VI 8.5.10), pour tout objet $F$ de $X_\et\gtimes_{X_\et}Y_\et$, on a un isomorphisme canonique
\begin{equation}\label{tf7i}
\urho_*(\Xi^*(F))\stackrel{\sim}{\rightarrow} \underset{\underset{X'\in \ob(\fV_\ox)}{\longrightarrow}}{\lim}\ \uPhi_{X'}^*(\rho_{X'*}(\Xi_{X'}^*(F))).
\end{equation}
On notera que les conditions requises dans 
(\cite{sga4} VI 8.5.10) sont satisfaites en vertu de \ref{rec22}, \ref{tf4} et (\cite{sga4} VI 3.3 et 5.1).
Par ailleurs, compte tenu de (\cite{egr1} 1.2.4(i)), le morphisme $\fc(F)$ \eqref{tf7a} s'identifie à la limite inductive des morphismes 
\begin{equation}\label{tf7j}
\uPhi_{X'}^*(\Phi_{X'}^*(\rho_*(F)))\rightarrow \uPhi_{X'}^*(\rho_{X'*}(\Xi_{X'}^*(F)))
\end{equation}
déduits des morphismes de changement de base relativement au carré de droite de \eqref{tf7e}. 
Ces derniers sont des isomorphismes d'après le cas (a); d'où la proposition (i) dans le cas (b).

On démontre de même la proposition (ii) en utilisant (\cite{sga4} VI 8.7.5).

\subsection{}\label{tf8}
Supposons $X$ strictement local. D'après (\cite{agt} VI.10.23 et VI.10.24), on a un morphisme canonique de topos 
\begin{equation}\label{tf8a}
\theta\colon Y_\fet\rightarrow \tE,
\end{equation}
qui est une section de $\beta$ \eqref{tf1k}, {\em i.e.}, on a un isomorphisme canonique 
\begin{equation}\label{tf8g}
\beta\theta\stackrel{\sim}{\rightarrow} \id_{Y_\fet}.
\end{equation} 
On obtient un morphisme de changement de base (\cite{agt} (VI.10.24.4))
\begin{equation}\label{tf8b}
\beta_*\rightarrow \theta^*.
\end{equation}

En vertu de (\cite{agt} VI.10.25), le diagramme 
\begin{equation}\label{tf8c}
\xymatrix{
{Y_\et}\ar[r]^-(0.5){\gamma}\ar[d]_{\rho_Y}&{X_\et\gtimes_{X_\et}Y_\et}\ar[d]^\rho\\
{Y_\fet}\ar[r]^{\theta}&{\tE}}
\end{equation}
où $\rho_Y$ est le morphisme canonique \eqref{notconv10a} et $\gamma$ est le morphisme défini dans \eqref{rec10c}, 
est commutatif à isomorphisme canonique près. 
Il induit donc un morphisme de changement de base
\begin{equation}\label{tf8d}
\theta^*\rho_*\rightarrow \rho_{Y*}\gamma^*.
\end{equation}
En restreignant aux faisceaux abéliens, on déduit un morphisme pour tout entier $q\geq 0$ (cf. \cite{egr1} 1.2.2)
\begin{equation}\label{tf8e}
\theta^*\rR^q\rho_*\rightarrow \rR^q\rho_{Y*}\gamma^*.
\end{equation}

\begin{prop}\label{tf9}
Conservons les hypothèses et notations de \ref{tf8}. Alors, 
\begin{itemize}
\item[{\rm (i)}] Le morphisme de changement de base $\beta_*\rightarrow \theta^*$ \eqref{tf8b} est un isomorphisme~; en particulier, 
le foncteur $\beta_*$ est exact.   
\item[{\rm (ii)}] Pour tout faisceau $F$ de $X_\et\gtimes_{X_\et}Y_\et$, le morphisme de changement de base \eqref{tf8d}
\begin{equation}\label{tf9a}
\theta^*(\rho_*(F))\rightarrow \rho_{Y*}(\gamma^*(F))
\end{equation}
est un isomorphisme. 
\item[{\rm (iii)}] Pour tout faisceau abélien $F$ de $X_\et\gtimes_{X_\et}Y_\et$ et tout entier $q\geq 0$, le morphisme de changement de base \eqref{tf8e}
\begin{equation}\label{tf9b}
\theta^*(\rR^q\rho_*(F))\rightarrow \rR^q\rho_{Y*}(\gamma^*(F))
\end{equation}
est un isomorphisme. 
\end{itemize}
\end{prop}

La proposition (i) est mentionnée à titre de rappel (\cite{agt} VI.10.27). Considérons le diagramme 
\begin{equation}\label{tf9f}
\xymatrix{
{Y_\et}\ar[r]^-(0.5){\gamma}\ar[d]_{\rho_Y}&{X_\et\gtimes_{X_\et}Y_\et}\ar[d]^\rho\ar[r]^-(0.5){\rp_2}&{Y_\et}\ar[d]^{\rho_Y}\\
{Y_\fet}\ar[r]^{\theta}&{\tE}\ar[r]^-(0.5){\beta}&{Y_\fet}}
\end{equation}
dont les carrés sont commutatifs à isomorphismes canoniques près \eqref{tf3c} et \eqref{tf8c}.

Le diagramme de morphismes de foncteurs
\begin{equation}\label{tf9c}
\xymatrix{
{\theta^*\rho_*}\ar[r]^-(0.5)c&{\rho_{Y*}\gamma^*}\\
{\beta_*\rho_*}\ar[u]^a\ar[r]^-(0.5)d&{\rho_{Y*}\rp_{2*}}\ar[u]_b}
\end{equation}
où $a$ est induit par \eqref{tf8b}, $b$ est induit par \eqref{rec10e}, $c$ est le morphisme de changement de base \eqref{tf8d} 
et $d$ est l'isomorphisme de commutativité du carré de droite de \eqref{tf9f}, est commutatif. 
En effet, d'après (\cite{egr1} 1.2.4(iii)), $c\circ a$ est le morphisme de changement de base relativement au carré commutatif 
\begin{equation}\label{tf9d}
\xymatrix{
{Y_\et}\ar[r]^-(0.5){\gamma}\ar[d]_{\rho_Y}&{X_\et\gtimes_{X_\et}Y_\et}\ar[d]^{\beta\rho}\\
{Y_\fet}\ar@{=}[r]&{Y_\fet}}
\end{equation}
La proposition (ii) résulte de \eqref{tf9c}, (i) et \ref{rec12}. 

Restreignant aux faisceaux abéliens, pour tout entier $q\geq 0$, le diagramme de morphismes de foncteurs
\begin{equation}\label{tf9e}
\xymatrix{
{\theta^*\rR^q\rho_*}\ar[r]^-(0.5){c^q}&{\rR^q\rho_{Y*}\gamma^*}\\
{\beta_*\rR^q\rho_*}\ar[u]^{a^q}\ar[r]^-(0.5){d^q}&{\rR^q\rho_{Y*}\rp_{2*}}\ar[u]_{b^q}}
\end{equation}
où $a^q$ est induit par \eqref{tf8b}, $b^q$ est induit par \eqref{rec10e}, $c^q$ est le morphisme de changement de base \eqref{tf8e} 
et $d^q$ est l'isomorphisme déduit de la relation $\beta\rho=\rho_Y\rp_2$ \eqref{tf9f} et de l'exactitude des foncteurs $\beta_*$ et $\rp_{2*}$ (cf. (i) et \ref{rec12}), 
est commutatif. En effet, d'après (\cite{egr1} 1.2.4(v)), $c^q\circ a^q$ est le morphisme de changement de base relativement au carré commutatif \eqref{tf9d}.
La proposition (iii) résulte de \eqref{tf9e}, (i) et \ref{rec12}.

\subsection{}\label{tf10}
Soient $\ox$ un point géométrique de $X$, $\uX$ le localisé strict de $X$ en $\ox$, 
$\uY=Y\times_X\uX$, $\uf\colon \uY\rightarrow \uX$ la projection canonique,
$\uX_\et\gtimes_{\uX_\et}\uY_\et$ le topos co-évanescent et $\tuE$ 
le topos de Faltings associés à $\uf$. Les deux carrés du diagramme
\begin{equation}\label{tf10a}
\xymatrix{
{\uY_\et}\ar[r]^-(0.5)\gamma\ar[d]_{\rho_\uY}&{\uX_\et\gtimes_{\uX_\et}\uY_\et}\ar[d]^\urho\ar[r]^\Xi&
{X_\et\gtimes_{X_\et}Y_\et}\ar[d]^\rho\\
{\uY_\fet}\ar[r]^-(0.5)\theta&{\tuE}\ar[r]^\Phi&\tE}
\end{equation}
où $\gamma$ est le morphisme \eqref{rec10c}, $\theta$ est le morphisme \eqref{tf8a},
$\Phi$ et $\Xi$ sont induits par fonctorialité par le morphisme canonique $\uX\rightarrow X$, 
$\urho$ et $\rho$ sont les morphismes canoniques \eqref{tf3b}
et $\rho_\uY$ est le morphisme canonique \eqref{notconv10a}, 
sont commutatifs à isomorphismes canoniques près \eqref{tf6a} et \eqref{tf8c}. On pose
\begin{eqnarray}
\phi_\ox=\gamma^*\circ \Xi^* \colon X_\et\gtimes_{X_\et}Y_\et&\rightarrow& \uY_\et,\label{tf10c}\\
\varphi_\ox=\theta^*\circ \Phi^*\colon \tE&\rightarrow& \uY_\fet.\label{tf10b}
\end{eqnarray}
On note 
\begin{equation}\label{tf10d}
\fc\colon \varphi_\ox\rho_*\rightarrow \rho_{\uY*}\phi_\ox
\end{equation}
le morphisme de changement de base relativement au rectangle extérieur du diagramme \eqref{tf10a}. 
En restreignant aux faisceaux abéliens, on déduit un morphisme pour tout entier $q\geq 0$ (cf. \cite{egr1} 1.2.2)
\begin{equation}\label{tf10e}
\fc^q\colon \varphi_\ox\rR^q\rho_*\rightarrow \rR^q\rho_{\uY*}\phi_\ox.
\end{equation}

\begin{prop}\label{tf11}
Conservons les hypothèses et notations de \ref{tf10}; supposons de plus que le morphisme $f\colon Y\rightarrow X$ soit cohérent. Alors,
\begin{itemize}
\item[{\rm (i)}] Le diagramme 
\begin{equation}\label{tf11b}
\xymatrix{
{\varphi_\ox\rho_*\rho^*} \ar[r]^-(0.5){\fc*\rho^*}&{\rho_{\uY*}\phi_\ox\rho^*}\ar@{=}[d]\\ 
{\varphi_\ox}\ar[u]^a\ar[r]^-(0.5)b&{\rho_{\uY*}\rho_\uY^*\varphi_\ox }}
\end{equation}
où $a$ et $b$ sont induits par les morphismes d'adjonction $\id\rightarrow \rho_*\rho^*$
et $\id\rightarrow \rho_{\uY*}\rho^*_\uY$, est commutatif. 

\item[{\rm (ii)}]  Pour tout faisceau $F$ de $X_\et\gtimes_{X_\et}Y_\et$, le morphisme de changement de base \eqref{tf10d}
\begin{equation}\label{tf11a}
\fc(F)\colon \varphi_\ox(\rho_*(F))\rightarrow \rho_{\uY*}(\phi_\ox (F))
\end{equation}
est un isomorphisme.

\item[{\rm (iii)}] Pour tout faisceau abélien $F$ de $X_\et\gtimes_{X_\et}Y_\et$ et tout entier $q\geq 0$, le morphisme de changement de base \eqref{tf10e}
\begin{equation}\label{tf11c}
\fc^q\colon \varphi_\ox(\rR^q\rho_*(F))\rightarrow \rR^q\rho_{\uY*}(\phi_\ox (F))
\end{equation}
est un isomorphisme. 
\end{itemize}
\end{prop}

(i) En effet, les triangles du diagramme
\begin{equation}\label{tf11g}
\xymatrix{
{\rho_\uY^*\varphi_\ox\rho_*\rho^*}\ar@{=}[r]\ar[rd]^{\fc'}&{\phi_\ox\rho^*\rho_*\rho^*}\ar[d]^d\\
{\rho_\uY^*\varphi_\ox}\ar@{=}[r]\ar[u]^{\rho_\uY^**a}&{\phi_\ox\rho^*}}
\end{equation}
où $\fc'$ est le morphisme adjoint de $\fc*\rho^*$ et $d$ est induit par le morphisme d'adjonction $\rho^*\rho_*\rightarrow \id$,
sont commutatifs en vertu de (\cite{sga4} XVII 2.1.3). 

(ii) Notons 
\begin{eqnarray}
\fc_1\colon \Phi^*\rho_*&\rightarrow& \urho_*\Xi^*,\label{tf11d}\\
\fc_2\colon \theta^*\urho_*&\rightarrow& \rho_{\uY*}\gamma^*,\label{tf11e}
\end{eqnarray}
les morphismes de changement de base relativement au carré de droite et au carré de gauche du diagramme \eqref{tf10a}, respectivement. 
D'après (\cite{egr1} 1.2.4(i)), on a 
\begin{equation}\label{tf11f}
\fc(F)= \fc_2(\Xi^*(F))\circ \theta^*(\fc_1(F)).
\end{equation}
Comme $\fc_1(F)$ est un isomorphisme en vertu de \ref{tf7}(i), et que $\fc_2(\Xi^*(F))$ est un isomorphisme en vertu de \ref{tf9}(ii), $\fc(F)$ est un isomorphisme. 

(iii) Il suffit de calquer la preuve de (ii) en utilisant \ref{tf7}(ii), \ref{tf9}(iii) et (\cite{egr1} 1.2.4(ii)).

\begin{cor}\label{tf12}
Supposons que les schémas $X$ et $Y$ soient cohérents, et que pour tout point géométrique $\ox$ de $X$,
notant $X_{(\ox)}$ le localisé strict de $X$ en $\ox$, le schéma $Y_{(\ox)}=Y\times_XX_{(\ox)}$ ait un nombre
fini de composantes connexes.
Alors, le morphisme d'adjonction $\id\rightarrow \rho_*\rho^*$ est un isomorphisme
\eqref{tf3b}~; en particulier, le foncteur 
\begin{equation}\label{tf12a}
\rho^*\colon \tE\rightarrow X_\et\gtimes_{X_\et}Y_\et
\end{equation}
est pleinement fidèle. 
\end{cor}

En effet, soient $\ox$ un point géométrique de $X$, $\rho_{Y_{(\ox)}}\colon Y_{(\ox),\et}\rightarrow Y_{(\ox),\fet}$ le foncteur canonique \eqref{notconv10a}, 
$\varphi_\ox\colon \tE\rightarrow Y_{(\ox),\fet}$ le foncteur défini dans \eqref{tf10b}. 
En vertu de (\cite{agt} VI.9.18), le morphisme d'adjonction $\id\rightarrow \rho_{Y_{(\ox)}*}\rho_{Y_{(\ox)}}^*$ est un isomorphisme. 
On en déduit, d'après \ref{tf11}(i)-(ii), que le morphisme 
$\varphi_\ox\rightarrow \varphi_\ox \rho_*\rho^*$ induit par le morphisme d'adjonction $\id\rightarrow \rho_*\rho^*$ est un isomorphisme. 
La proposition s'ensuit puisque la famille des foncteurs $\varphi_\ox$ lorsque $\ox$ décrit l'ensemble des points géométriques de $X$, 
est conservative d'après (\cite{agt} VI.10.32).

\begin{cor}\label{tf13}
Soient $\mP$ un ensemble de nombres premiers, $F$ un faisceau abélien de $\mP$-torsion de $\tE$, $q$ un entier $\geq 1$.  
Supposons que  les schémas $X$ et $Y$ soient cohérents, et que pour tout point géométrique $\ox$ de $X$,
notant $X_{(\ox)}$ le localisé strict de $X$ en $\ox$, le schéma $Y_{(\ox)}=Y\times_XX_{(\ox)}$ soit $K(\pi,1)$ pour les faisceaux abéliens de $\mP$-torsion  \eqref{Kpun2}.
Alors, on a 
\begin{equation}\label{tf13a}
\rR^q\rho_*(\rho^*(F))=0.
\end{equation}
\end{cor}

En effet, soient $\ox$ un point géométrique de $X$, $\rho_{Y_{(\ox)}}\colon Y_{(\ox),\et}\rightarrow Y_{(\ox),\fet}$ le foncteur canonique \eqref{notconv10a}, 
\begin{eqnarray}
\phi_\ox \colon X_\et\gtimes_{X_\et}Y_\et&\rightarrow& Y_{(\ox),\et},\label{tf13b}\\
\varphi_\ox \colon \tE&\rightarrow& Y_{(\ox),\fet},\label{tf13c}
\end{eqnarray}
les foncteurs définis dans \eqref{tf10c} et \eqref{tf10b}. En vertu de \ref{tf11}(iii), on a un isomorphisme canoniqe
\begin{equation}
\varphi_\ox(\rR^q\rho_*(\rho^*(F)))\stackrel{\sim}{\rightarrow} \rR^q\rho_{Y_{(\ox)}*}(\phi_\ox (\rho^*(F))).
\end{equation}
Par ailleurs, on a un isomorphisme canonique $\phi_\ox (\rho^*(F))\stackrel{\sim}{\rightarrow}\rho_{Y_{(\ox)}}^*(\varphi_\ox (F))$ \eqref{tf10a}.
On en déduit que $\varphi_\ox(\rR^q\rho_*(\rho^*(F)))=0$. La proposition s'ensuit compte tenu de (\cite{agt} VI.10.32).

\chapter[Théorie de Hodge $p$-adique]{\texorpdfstring{Théorie de Hodge $p$-adique}
{Théorie de Hodge p-adique}}\label{thpadique}

\section{Hypothèses et notations; schéma logarithmique adéquat}\label{TFA}

On renvoie à (\cite{agt} II.5) pour un lexique de géométrie logarithmique.

\subsection{}\label{TFA1}
Dans ce chapitre, $K$ désigne un corps de valuation discrète complet de 
caractéristique $0$, à corps résiduel {\em algébriquement clos} $k$ de caractéristique $p>0$,  
$\co_K$ l'anneau de valuation de $K$, $\oK$ une clôture algébrique de $K$, $\co_\oK$ la clôture intégrale de $\co_K$ dans $\oK$,
$\fm_\oK$ l'idéal maximal de $\co_\oK$ et $G_K$ le groupe de Galois de $\oK$ sur $K$.
On note $\co_C$ le séparé complété $p$-adique de $\co_\oK$, $\fm_C$ son idéal maximal,
$C$ son corps des fractions et $v$ sa valuation, normalisée par $v(p)=1$. 
On désigne par par $\hmZ(1)$ et $\mZ_p(1)$ les $\mZ[G_K]$-modules 
\begin{eqnarray}
\hmZ(1)&=&\underset{\underset{n\geq 1}{\longleftarrow}}{\lim}\ \mu_{n}(\co_\oK),\label{TFA1a}\\
\mZ_p(1)&=&\underset{\underset{n\geq 0}{\longleftarrow}}{\lim}\ \mu_{p^n}(\co_\oK),\label{TFA1aa}
\end{eqnarray}  
où $\mu_n(\co_\oK)$ désigne le sous-groupe des racines $n$-ièmes de l'unité dans $\co_\oK$. 
Pour tout $\mZ_p[G_K]$-module $M$ et tout entier $n$, on pose 
$M(n)=M\otimes_{\mZ_p}\mZ_p(1)^{\otimes n}$.

On pose $S=\Spec(\co_K)$ et $\oS=\Spec(\co_\oK)$ et 
on note $s$ (resp.  $\eta$, resp. $\oeta$) le point fermé de $S$ (resp.  générique de $S$, resp. générique de $\oS$).
Pour tout entier $n\geq 1$, on pose $S_n=\Spec(\co_K/p^n\co_K)$. Pour tout $S$-schéma $X$, on pose 
\begin{equation}\label{TFA1b}
\oX=X\times_S\oS \ \ \ {\rm et}\ \ \ X_n=X\times_SS_n.
\end{equation} 

On munit $S$ de la structure logarithmique $\cM_S$ définie par son point fermé, 
autrement dit, $\cM_S=u_*(\co_\eta^\times)\cap \co_S$, où $u\colon \eta\rightarrow S$ est l'injection canonique.

\subsection{}\label{TFA100}
Comme $\co_\oK$ satisfait les conditions requises dans \ref{mptf1}, 
il est loisible de considérer les notions introduites dans les sections \ref{alpha}--\ref{mptf} pour cet anneau.   
On choisit un système compatible $(\beta_n)_{n>0}$ 
de racines $n$-ièmes de $p$ dans $\co_\oK$. Pour tout nombre rationnel $\varepsilon>0$, 
on pose $p^\varepsilon=(\beta_n)^{\varepsilon n}$, où $n$ est un entier $>0$ tel que $\varepsilon n$ soit entier.

\subsection{}\label{TFA3}
Dans la suite de ce chapitre, $f\colon (X,\cM_X)\rightarrow (S,\cM_S)$ 
désigne un morphisme {\em adéquat} de schémas logarithmiques (\cite{agt} III.4.7). 
On note $d=\dim(X/S)$ la dimension relative de $X$ sur $S$.
On désigne par $X^\circ$ le sous-schéma ouvert maximal de $X$
où la structure logarithmique $\cM_X$ est triviale~; c'est un sous-schéma ouvert de $X_\eta$.
On note $j\colon X^\circ\rightarrow X$ l'injection canonique. D'après (\cite{agt} III.4.2(iv)), $j$ 
est schématiquement dominant et on a un isomorphisme canonique
\begin{equation}\label{TFA3e}
\cM_X\stackrel{\sim}{\rightarrow}j_*(\co_{X^\circ}^\times)\cap \co_X.
\end{equation} 
En particulier, l'homomorphisme canonique $\cM_X\rightarrow \co_X$ est un monomorphisme.
Pour tout $X$-schéma $U$, on pose  
\begin{equation}\label{TFA3a}
U^\circ=U\times_XX^\circ.
\end{equation} 
On note $\hbar\colon \oX\rightarrow X$ et $h\colon \oX^\circ\rightarrow X$ les morphismes canoniques \eqref{TFA1b}, de sorte que 
l'on a $h=\hbar\circ j_\oX$. Pour alléger les notations, on pose
\begin{equation}\label{TFA3b}
\tOmega^1_{X/S}=\Omega^1_{(X,\cM_X)/(S,\cM_S)},
\end{equation}
que l'on considère comme un faisceau de $X_\zar$ ou $X_\et$, selon le contexte (cf. \ref{notconv12}). 
La dérivation logarithmique canonique
\begin{equation}\label{TFA3c}
d\log \colon \cM_X\rightarrow \tOmega^1_{X/S}
\end{equation}
induit un morphisme $\co_X$-linéaire et surjectif (\cite{ogus} IV 1.2.10)
\begin{equation}\label{TFA3d}
\cM_X^\gp\otimes_\mZ\co_X\rightarrow \tOmega^1_{X/S}.
\end{equation}

\subsection{}\label{TFA5}
Pour tout entier $n\geq 1$, on note $a\colon X_s\rightarrow X$, $a_n\colon X_s\rightarrow X_n$, 
$\iota_n\colon X_n\rightarrow X$ et $\oiota_n\colon \oX_n\rightarrow \oX$ les injections canoniques \eqref{TFA1b}. 
Le corps résiduel de $\co_K$ étant algébriquement clos, il existe un unique $S$-morphisme $s\rightarrow \oS$. 
Celui-ci induit des immersions fermées $\oa\colon X_s\rightarrow \oX$ et $\oa_n\colon X_s\rightarrow \oX_n$
qui relèvent $a$ et $a_n$, respectivement. 
\begin{equation}\label{TFA5a}
\xymatrix{
{X_s}\ar[r]_{\oa_n}\ar@{=}[d]\ar@/^1pc/[rr]^{\oa}&{\oX_n}\ar[r]_{\oiota_n}\ar[d]^{\hbar_n}&\oX\ar[d]^\hbar\\
{X_s}\ar[r]^{a_n}\ar@/_1pc/[rr]_{a}&{X_n}\ar[r]^{\iota_n}&X}
\end{equation}
Comme $\hbar$ est entier et que $\hbar_n$ est un homéomorphisme
universel, pour tout faisceau $\cF$ de $\oX_\et$, le morphisme de changement de base 
\begin{equation}\label{TFA5b}
a^*(\hbar_*(\cF))\rightarrow \oa^*(\cF)
\end{equation}
est un isomorphisme (\cite{sga4} VIII 5.6). Par ailleurs, $\oa_n$ étant un homéomorphisme universel, on peut considérer $\co_{\oX_n}$
comme un faisceau de $X_{s,\zar}$ ou $X_{s,\et}$, selon le contexte (cf. \ref{notconv12}). 

On pose 
\begin{equation}\label{TFA5c}
\tOmega^1_{\oX_n/\oS_n}=\tOmega^1_{X/S}\otimes_{\co_X}\co_{\oX_n},
\end{equation}
que l'on considère aussi un faisceau de $X_{s,\zar}$ ou $X_{s,\et}$, selon le contexte.

\subsection{}\label{TFA6}
On désigne par
\begin{equation}\label{TFA6a}
\pi\colon E\rightarrow \Et_{/X}
\end{equation}
le $\mU$-site fibré de Faltings associé au morphisme $h\colon \oX^\circ\rightarrow X$ (cf. \ref{tf1}).
Pour tout $U\in \ob(\Et_{/X})$, on note 
\begin{equation}\label{TFA6F}
\iota_{U!}\colon \Et_{\rf/\oU^\circ}\rightarrow E
\end{equation} 
le foncteur canonique \eqref{tf1c}. On munit $E$ de la topologie co-évanescente définie par $\pi$ (\cite{agt} VI.5.3)  
et on note $\tE$ le topos des faisceaux de $\mU$-ensembles sur $E$. On désigne par 
\begin{eqnarray}
\sigma\colon \tE \rightarrow X_\et,\label{TFA6c}\\
\beta\colon \tE \rightarrow \oX^\circ_\fet,\label{TFA6b}\\
\psi\colon \oX^\circ_\et\rightarrow \tE,\label{TFA6d}\\
\rho\colon X_\et\gtimes_{X_\et}\oX^\circ_\et\rightarrow \tE,\label{TFA6e}
\end{eqnarray}
les morphismes canoniques \eqref{tf1kk}, \eqref{tf1k}, \eqref{tf1m} et \eqref{tf3b}. 

\subsection{}\label{TFA66}
On note $\tE_s$ le sous-topos fermé de $\tE$ complémentaire de l'ouvert $\sigma^*(X_\eta)$ (\cite{agt} III.9.3), 
\begin{equation}\label{TFA66a}
\delta\colon \tE_s\rightarrow \tE
\end{equation} 
le plongement canonique et 
\begin{equation}\label{TFA66b}
\sigma_s\colon \tE_s\rightarrow X_{s,\et}
\end{equation} 
le morphisme de topos induit par $\sigma$ \eqref{TFA6c} (\cite{agt} (III.9.8.3)). Le diagramme de morphismes de topos 
\begin{equation}\label{TFA66c}
\xymatrix{
{\tE_s}\ar[r]^{\sigma_s}\ar[d]_{\delta}&{X_{s,\et}}\ar[d]^{a}\\
{\tE}\ar[r]^\sigma&{X_\et}}
\end{equation}
est commutatif à isomorphisme près. Les foncteurs $a_*$ et $\delta_*$ étant exacts, 
pour tout groupe abélien $F$ de $\tE_s$ et tout entier $i\geq 0$, on a un isomorphisme canonique 
\begin{equation}\label{TFA66d}
a_*(\rR^i\sigma_{s*}(F))\stackrel{\sim}{\rightarrow}\rR^i\sigma_*(\delta_*F). 
\end{equation}

Notons $\Pt(\tE)$ et $\Pt(\tE_s)$ les catégories des points de $\tE$ et $\tE_s$, respectivement. Le foncteur
\begin{equation}\label{TFA66e}
\Pt(\tE_s)\rightarrow \Pt(\tE), \ \ \ q\mapsto \delta \circ q,
\end{equation}
est pleinement fidèle (\cite{sga4} IV 9.7.2).

\begin{lem}[\cite{agt} III.9.5]\label{TFA7}
{\rm (i)}\ Soit $(\oy\rightsquigarrow \ox)$ un point de $X_\et\gtimes_{X_\et}\oX^\circ_\et$ \eqref{rec15}. 
Pour que $\rho(\oy\rightsquigarrow \ox)$ appartienne à l'image essentielle du foncteur \eqref{TFA66e}, 
il faut et il suffit que $\ox$ soit au-dessus de~$s$. 

{\rm (ii)}\ La famille des points de $\tE_s$ définie par la famille des points 
$\rho(\oy\rightsquigarrow \ox)$ de $\tE$ tels que $\ox$ soit au-dessus de $s$ est conservative. 
\end{lem}

\subsection{}\label{TFA2}
D'après (\cite{agt} III.4.2(iii)), $\oX$ est normal et localement irréductible (\cite{agt} III.3.1). 
Par ailleurs, l'immersion $j\colon X^\circ\rightarrow X$ est quasi-compacte puisque $X$ est noethérien.  
Pour tout objet $(V\rightarrow U)$ de $E$, on note $\oU^V$ la fermeture intégrale de $\oU$ dans $V$. 
Pour tout morphisme $(V'\rightarrow U')\rightarrow (V\rightarrow U)$ de $E$, on a un morphisme canonique 
$\oU'^{V'}\rightarrow \oU^V$ qui s'insère dans un diagramme commutatif 
\begin{equation}\label{TFA2a}
\xymatrix{
V'\ar[r]\ar[d]&{\oU'^{V'}}\ar[d]\ar[r]&\oU'\ar[r]\ar[d]&U'\ar[d]\\
V\ar[r]&{\oU^V}\ar[r]&\oU\ar[r]&U}
\end{equation} 
On désigne par $\ocB$ le préfaisceau sur $E$ défini pour tout $(V\rightarrow U)\in \ob(E)$, par 
\begin{equation}\label{TFA2b}
\ocB((V\rightarrow U))=\Gamma(\oU^V,\co_{\oU^V}).
\end{equation}
C'est un faisceau pour la topologie co-évanescente de $E$ en vertu de (\cite{agt} III.8.16). 
Pour tout $U\in \ob(\Et_{/X})$, on pose 
\begin{equation}\label{TFA2d}
\ocB_{U}=\ocB\circ \iota_{U!}.
\end{equation} 

D'après (\cite{agt} III.8.17), on a un homomorphisme canonique 
\begin{equation}\label{TFA2c}
\sigma^*(\hbar_*(\co_\oX))\rightarrow \ocB.
\end{equation}
Sauf mention explicite du contraire, on considère $\sigma$ \eqref{TFA6c}
comme un morphisme de topos annelés
\begin{equation}\label{TFA2e}
\sigma\colon (\tE,\ocB)\rightarrow (X_\et,\hbar_*(\co_\oX)).
\end{equation}

Notant encore $\co_\oK$ le faisceau constant de $\oX^\circ_\fet$ de valeur $\co_\oK$. 
Sauf mention explicite du contraire, on considère $\beta$ \eqref{TFA6b} comme un morphisme de topos annelés
\begin{equation}\label{TFA2g}
\beta\colon (\tE,\ocB)\rightarrow (\oX^\circ_\fet,\co_\oK).
\end{equation}

\subsection{}\label{TFA8}
Pour tout entier $n\geq 0$ et tout $U\in \ob(\Et_{/X})$, on pose 
\begin{eqnarray}
\ocB_n&=&\ocB/p^n\ocB,\label{TFA8a}\\
\ocB_{U,n}&=&\ocB_U/p^n\ocB_U.\label{TFA8b}
\end{eqnarray}
Les correspondances $\{U\mapsto p^n\ocB_U\}$ et $\{U\mapsto \ocB_{U,n}\}$
forment naturellement des préfaisceaux sur $E$ \eqref{tf1h}, et les morphismes canoniques 
\begin{eqnarray}
\{U\mapsto p^n\ocB_U\}^\tta&\rightarrow&p^n\ocB,\label{TFA8c}\\
\{U\mapsto \ocB_{U,n}\}^\tta&\rightarrow&\ocB_n,\label{TFA8d}
\end{eqnarray}
où  les termes de gauche désignent les faisceaux associés dans $\tE$, sont des isomorphismes
en vertu de (\cite{agt} VI.8.2 et VI.8.9). 

D'après (\cite{agt} III.9.7), $\ocB_n$ est un anneau de $\tE_s$. Si $n\geq 1$, on note
\begin{equation}\label{TFA8e}
\sigma_n\colon (\tE_s,\ocB_n)\rightarrow (X_{s,\et},\co_{\oX_n})
\end{equation}
le morphisme de topos annelés induit par $\sigma$ \eqref{TFA2e} (cf. \cite{agt} (III.9.9.4)). Pour tout entier $q\geq 0$, 
on désigne par 
\begin{equation}\label{TFA8f}
\rR^q\sigma_{n*}\colon \bMod(\ocB_n,\tE_s)\rightarrow \bMod(\co_{\oX_n},X_{s,\et})
\end{equation}
le $q$-ième foncteur dérivé droit de $\sigma_{n*}$ (cf. \ref{notconv9} pour les notations).

Notant encore $\co_\oK/p^n\co_\oK$ le faisceau constant de $\oX^\circ_\fet$ de valeur $\co_\oK/p^n\co_\oK$, on désigne par 
\begin{equation}\label{TFA8g}
\beta_n\colon (\tE_s,\ocB_n)\rightarrow (\oX^\circ_\fet,\co_\oK/p^n\co_\oK).
\end{equation}
le morphisme induit par $\beta$ \eqref{TFA2g}.

\begin{rema}\label{TFA81}
Sous les hypothèses de \ref{TFA8}, on notera que l'homomorphisme canonique \eqref{TFA6F}
\begin{equation}\label{TFA81a}
\ocB_{U,n}\rightarrow \ocB_n\circ \iota_{U!}
\end{equation} 
n'est pas en général un isomorphisme~; c'est pourquoi nous n'utiliserons pas la notation $\ocB_{n,U}$. 
Toutefois, sous certaines hypothèses \eqref{cg21}, nous montrerons dans \ref{cg28} que pour 
tout schéma affine $U$ de $\Et_{/X}$, l'homomorphisme \eqref{TFA81a} est un $\alpha$-isomorphisme.
\end{rema}

\subsection{}\label{TFA9}
Soient $U$ un objet de $\Et_{/X}$, $\oy$ un point géométrique de $\oU^\circ$. 
Le schéma $\oU$ étant localement irréductible d'après (\cite{agt} III.3.3 et III.4.2(iii)),  
il est la somme des schémas induits sur ses composantes irréductibles. On note $\oU^\star$
la composante irréductible de $\oU$ (ou ce qui revient au même, sa composante connexe) contenant $\oy$. 
De même, $\oU^\circ$ est la somme des schémas induits sur ses composantes irréductibles
et $\oU^{\star \circ}=\oU^\star\times_{X}X^\circ$ est la composante irréductible de $\oU^\circ$ contenant $\oy$. 
On note $\bB_{\pi_1(\oU^{\star \circ},\oy)}$ le topos classifiant du groupe profini $\pi_1(\oU^{\star \circ},\oy)$ et
\begin{equation}\label{TFA9a}
\nu_\oy\colon \oU^{\star \circ}_\fet \stackrel{\sim}{\rightarrow}\bB_{\pi_1(\oU^{\star \circ},\oy)}
\end{equation}
le foncteur fibre  de $\oU^{\star \circ}_\fet$ en $\oy$ \eqref{notconv11c}. On pose
\begin{equation}\label{TFA9b}
\oR^\oy_U=\nu_\oy(\ocB_U|\oU^{\star \circ}).
\end{equation}
Explicitement, soit $(V_i)_{i\in I}$ le revêtement universel normalisé de $\oU^{\star \circ}$ en $\oy$ \eqref{notconv11}.
Pour chaque $i\in I$, $(V_i\rightarrow U)$ est naturellement un objet de $E$. 
On note $\oU^{V_i}$ la fermeture intégrale de $\oU$ dans $V_i$.
Les schémas $(\oU^{V_i})_{i\in I}$ forment alors un système projectif filtrant, et on a   
\begin{equation}\label{TFA9c}
\oR^\oy_U=\underset{\underset{i\in I}{\longrightarrow}}{\lim}\ \Gamma(\oU^{V_i},\co_{\oU^{V_i}}).
\end{equation} 

\begin{rema}\label{TFA10}
Sous les hypothèses de \ref{TFA9}, si, de plus, $U$ est affine, 
l'anneau $\oR^\oy_U$ est intègre et normal. En effet, le schéma $\oU^\star$ est affine, intègre et normal. 
Pour tout $i\in I$, $V_i$ étant connexe, il est alors intègre et normal. 
Par suite, $\oU^{V_i}$ est intègre, normal et entier sur $\oU^\star$ (\cite{ega2} 6.3.7). 
Par ailleurs, pour tous $(i,j)\in I^2$ avec  
$i\geq j$, le morphisme canonique $\oU^{V_i}\rightarrow \oU^{V_j}$ est entier et dominant. 
L'assertion recherchée s'ensuit d'après (\cite{ega1n} 0.6.1.6(i) et 0.6.5.12(ii)). 
\end{rema}

\subsection{}\label{TFA14}
Soient $\ox$ un point géométrique de $X$, $X'$ le localisé strict de $X$ en $\ox$.
On désigne par $\fV_\ox$ la catégorie des $X$-schémas étales $\ox$-pointés, ou ce qui revient au même,
la catégorie des voisinages du point de $X_\et$ associé à $\ox$ dans le site $\Et_{/X}$ (\cite{sga4} IV 6.8.2 et VIII 3.9), et par
$\fW_\ox$ la sous-catégorie pleine de $\fV_\ox$ formée des objets $(U,\fp\colon \ox\rightarrow U)$ 
tels que le schéma $U$ soit affine.
Ce sont des catégories cofiltrantes, et le foncteur d'injection canonique
$\fW^\circ_\ox\rightarrow \fV^\circ_\ox$ est cofinal (\cite{sga4} I 8.1.3(c)). 
Pour tout objet $(U,\fp\colon \ox\rightarrow U)$ de $\fV_\ox$, 
on note encore $\fp\colon X'\rightarrow U$ le morphisme déduit de 
$\fp$ (\cite{sga4} VIII 7.3) et on pose
\begin{equation}\label{TFA14e}
\ofp^\circ=\fp\times_X\oX^\circ\colon \oX'^\circ \rightarrow \oU^\circ.
\end{equation}

On note
\begin{equation}\label{TFA14a}
\varphi_\ox\colon \tE\rightarrow \oX'^\circ_\fet
\end{equation}
le foncteur canonique défini dans \eqref{tf10b}.

Le foncteur composé $\varphi_\ox\circ\beta^*$ est canoniquement isomorphe au foncteur image inverse par 
le morphisme canonique $\oX'_\fet\rightarrow \oX_\fet$, d'après \eqref{tf8g} et (\cite{agt} (VI.10.12.6)). 
Pour tout objet $F$ de $X_\et$, on a un isomorphisme canonique et fonctoriel 
\begin{equation}\label{TFA14b}
\varphi_\ox(\sigma^*(F))\stackrel{\sim}{\rightarrow} F_\ox,
\end{equation}
de but le faisceau constant de $\oX'_\fet$ de valeur $F_\ox$, en vertu de (\cite{agt} VI.10.24 et (VI.10.12.6)). 

Pour tout objet $F=\{U\mapsto F_U\}$ de $\hE$,
on note $F^\tta$ le faisceau de $\tE$ associé à $F$, et pour tout $U\in \ob(\Et_{/X})$, 
$F_U^\tta$ le faisceau de $\oU^\circ_\fet$ associé à $F_U$. 
D'après (\cite{agt} VI.10.37),  on a un isomorphisme canonique et fonctoriel
\begin{equation}\label{TFA14c}
\varphi_\ox(F^\tta) \stackrel{\sim}{\rightarrow}\underset{\underset{(U,\fp)\in \fV_\ox^\circ}{\longrightarrow}}{\lim}\ (\ofp^\circ)_\fet^*(F^\tta_U).
\end{equation}

En vertu de (\cite{agt} VI.10.30), pour tout groupe abélien $F$ de $\tE$
et tout entier $q\geq 0$, on a un isomorphisme canonique et fonctoriel
\begin{equation}\label{TFA14d}
\rR^q\sigma_*(F)_\ox\stackrel{\sim}{\rightarrow}\rH^q(\oX'^\circ_\fet,\varphi_\ox(F)). 
\end{equation}

\subsection{}\label{TFA11}
Soient $(\oy\rightsquigarrow \ox)$ un point de $X_\et\gtimes_{X_\et}\oX^\circ_\et$ \eqref{rec15}, 
$X'$ le localisé strict de $X$ en $\ox$. Le $X$-morphisme $u\colon \oy\rightarrow X'$
définissant $(\oy\rightsquigarrow \ox)$ se relève en un $\oX^\circ$-morphisme $v\colon \oy\rightarrow \oX'^\circ$ et 
induit donc un point géométrique de $\oX'^\circ$ que l'on note aussi (abusivement) $\oy$.
On désigne par $\fV_\ox$ la catégorie des $X$-schémas étales $\ox$-pointés. 
Pour tout objet $(U,\fp\colon \ox\rightarrow U)$ de $\fV_\ox$, 
on note encore $\fp\colon X'\rightarrow U$ le morphisme déduit de $\fp$ (\cite{sga4} VIII 7.3) et on pose
\begin{equation}\label{TFA11b}
\ofp^\circ=\fp\times_X\oX^\circ\colon \oX'^\circ \rightarrow \oU^\circ.
\end{equation}
On note aussi (abusivement) $\oy$ le point géométrique $\ofp^\circ(v(\oy))$ de $\oU^\circ$.

Pour tout objet $F=\{U\mapsto F_U\}$ de $\hE$ \eqref{tf1h}, on note
$F^\tta$ le faisceau de $\tE$ associé à $F$, et pour tout $U\in \ob(\Et_{/X})$, 
$F_U^\tta$ le faisceau de $\oU^\circ_\fet$ associé à $F_U$. 
D'après (\cite{agt} VI.10.36 et (VI.9.3.4)), on a un isomorphisme canonique et fonctoriel 
\begin{equation}\label{TFA11a}
(F^\tta)_{\rho(\oy\rightsquigarrow \ox)} \stackrel{\sim}{\rightarrow} 
\underset{\underset{(U,\fp)\in \fV_\ox^\circ}{\longrightarrow}}{\lim}\ (F^\tta_U)_{\rho_{\oU^\circ}(\oy)},
\end{equation}
où $\rho$ est le morphisme \eqref{TFA6e} et $\rho_{\oU^\circ}\colon \oU^\circ_\et\rightarrow \oU^\circ_\fet$ 
est le morphisme canonique \eqref{notconv10a}. 
Compte tenu de \eqref{TFA9b} et (\cite{agt} VI.9.9), on en déduit un isomorphisme canonique de $\Gamma(\oX',\co_{\oX'})$-algèbres
\begin{equation}\label{TFA11c}
\ocB_{\rho(\oy\rightsquigarrow \ox)}\stackrel{\sim}{\rightarrow} 
\underset{\underset{(U,\fp)\in \fV_\ox^\circ}{\longrightarrow}}{\lim}\ \oR^{\oy}_U,
\end{equation}
où $\oR^{\oy}_U$ est la $\Gamma(\oU,\co_\oU)$-algèbre définie dans \eqref{TFA9c}.

\subsection{}\label{TFA12}
Conservons les hypothèses et notations de \ref{TFA11}; supposons, de plus, que {\em $\ox$ soit au-dessus de $s$}. 
D'après (\cite{agt} III.3.7), $\oX'$ est normal et strictement local (et en particulier intègre). 
Pour tout objet $(U,\fp\colon \ox\rightarrow U)$ de $\fV_\ox$, on désigne par $\oU^\star$  la composante irréductible de $\oU$ contenant $\oy$ 
et on pose $\oU^{\star\circ}=\oU^\star\times_XX^\circ$ qui est la composante irréductible de $\oU^\circ$ contenant $\oy$ (cf. \ref{TFA9}).
Le morphisme $\ofp^\circ \colon \oX'^\circ\rightarrow \oU^\circ$
\eqref{TFA11b} se factorise donc à travers $\oU^{\star \circ}$. On désigne par $\bB_{\pi_1(\oX'^{\circ},\oy)}$ le topos classifiant du groupe profini 
$\pi_1(\oX'^{\circ},\oy)$ et par
\begin{equation}\label{TFA12d}
\nu_{\oy}\colon \oX'^\circ_\fet \stackrel{\sim}{\rightarrow}\bB_{\pi_1(\oX'^\circ,\oy)}
\end{equation}
le foncteur fibre de $\oX'^\circ_\fet$ en $\oy$ \eqref{notconv11c}. 
D'après (\cite{agt} VI.10.31 et VI.9.9), le foncteur composé
\begin{equation}\label{TFA12e}
\xymatrix{
{\tE}\ar[r]^-(0.5){\varphi_\ox}&{\oX'^\circ_\fet}\ar[r]^-(0.5){\nu_\oy}&{\bB_{\pi_1(\oX'^\circ,\oy)}}\ar[r]&\Ens},
\end{equation}
où $\varphi_\ox$ est le foncteur canonique \eqref{TFA14a} et la dernière flèche est le foncteur d'oubli de l'action de $\pi_1(\oX'^\circ,\oy)$, 
est canoniquement isomorphe au foncteur fibre associé au point $\rho(\oy\rightsquigarrow \ox)$ de $\tE$.  
On observera (\cite{agt} (VI.10.25.6))
que pour tout objet $F$ de $X_\et$, l'isomorphisme \eqref{TFA14b} induit l'isomorphisme canonique (\cite{agt} (VI.10.18.1))
\begin{equation}
(\sigma^*F)_{\rho(\oy\rightsquigarrow \ox)} \stackrel{\sim}{\rightarrow} F_{\ox}.
\end{equation}

On définit la $\Gamma(\oX',\co_{\oX'})$-algèbre $\oR_{X'}^\oy$ de $\bB_{\pi_1(\oX'^\circ,\oy)}$ par la formule
\begin{equation}\label{TFA12f}
\oR_{X'}^\oy=\underset{\underset{(U,\fp)\in \fV_\ox^\circ}{\longrightarrow}}{\lim}\ \oR^{\oy}_U, 
\end{equation}
où l'on considère $\oR^{\oy}_U$ comme une $\Gamma(\oU,\co_\oU)$-algèbre de $\bB_{\pi_1(\oU^{\star\circ},\oy)}$ \eqref{TFA9c}.  
On note $\hoR^\oy_{X'}$ son séparé complété $p$-adique. 
L'isomorphisme \eqref{TFA14c} induit un isomorphisme de $\Gamma(\oX',\co_{\oX'})$-algèbres de 
$\bB_{\pi_1(\oX'^\circ,\oy)}$
\begin{equation}\label{TFA12g}
\nu_{\oy}(\varphi_\ox(\ocB))\stackrel{\sim}{\rightarrow} \oR_{X'}^\oy,
\end{equation} 
dont l'isomorphisme de $\Gamma(\oX',\co_{\oX'})$-algèbres sous-jacent est \eqref{TFA11c}.

\begin{remas}\label{TFA13}
Soient $\ox$ un point géométrique de $X$ au-dessus de $s$, $X'$ le localisé strict de $X$ en $\ox$, 
$\varphi_\ox\colon \tE\rightarrow \oX'^\circ_\fet$ le foncteur canonique \eqref{TFA14a}. Alors, 
\begin{itemize}
\item[(i)] Il existe un point $(\oy\rightsquigarrow \ox)$ de $X_\et\gtimes_{X_\et}\oX^\circ_\et$ \eqref{rec15}. 
En effet, d'après (\cite{agt} III.3.7), $\oX'$ est normal et strictement local (et en particulier intègre). 
Comme $X^\circ$ est schématiquement dense dans $X$ \eqref{TFA3}, $\oX'^\circ$ est intègre et non-vide 
(\cite{ega4} 11.10.5). Soit $v\colon \oy\rightarrow \oX'^\circ$ un point géométrique de $\oX'^\circ$.
On note encore $\oy$ le point géométrique de $\oX^\circ$ et 
$u\colon \oy\rightarrow X'$ le $X$-morphisme induits par~$v$. On obtient ainsi 
un point $(\oy\rightsquigarrow \ox)$ de $X_\et\gtimes_{X_\et}\oX^\circ_\et$.  
\item[(ii)] Le morphisme d'adjonction $\id\rightarrow \delta_*\delta^*$ \eqref{TFA66a} induit 
un isomorphisme de foncteurs 
\begin{equation}\label{TFA13a}
\varphi_\ox\stackrel{\sim}{\rightarrow}\varphi_\ox\delta_*\delta^*.
\end{equation}
En effet, d'après \ref{TFA7}(i), pour tout objet $F$ de $\tE$, le morphisme d'adjonction 
\begin{equation}\label{TFA13b}
F_{\rho(\oy\rightsquigarrow \ox)}\rightarrow(\delta_*(\delta^*(F)))_{\rho(\oy\rightsquigarrow \ox)}
\end{equation}
est un isomorphisme. 
Comme $\oX'^\circ$ est intègre (\cite{agt} III.3.7), la proposition s'ensuit compte tenu de la description \eqref{TFA12e} 
du foncteur fibre en $\rho(\oy\rightsquigarrow \ox)$.
\end{itemize}
\end{remas}

\section{Tour de revêtements associée à une carte logarithmique adéquate}\label{cad}

\subsection{}\label{cad1}
On suppose dans cette section que le morphisme $f\colon (X,\cM_X)\rightarrow (S,\cM_S)$ \eqref{TFA3} admet une carte adéquate 
$((P,\gamma),(\mN,\iota),\vartheta)$ (\cite{agt} III.4.4), autrement dit qu'il existe une carte $(P,\gamma)$ pour  $(X,\cM_X)$ (\cite{agt} II.5.13), 
une carte $(\mN,\iota)$ pour $(S,\cM_S)$ et un homomorphisme de monoïdes $\vartheta\colon \mN\rightarrow P$ 
tels que les conditions suivantes soient remplies~:  
\begin{itemize}
\item[(i)] Le diagramme d'homomorphismes de monoïdes
\begin{equation}\label{cad1a}
\xymatrix{
P\ar[r]^-(0.5)\gamma&{\Gamma(X,\cM_X)}\\
\mN\ar[r]^-(0.5){\iota}\ar[u]^\vartheta&{\Gamma(S,\cM_S)}\ar[u]_{f^\flat}}
\end{equation}
est commutatif, ou ce qui revient au même (avec les notations de \ref{notconv2}),
le diagramme associé de morphismes de schémas logarithmiques 
\begin{equation}\label{cad1b}
\xymatrix{
{(X,\cM_X)}\ar[r]^-(0.5){\gamma^a}\ar[d]_f&{\bA_P}\ar[d]^{\bA_\vartheta}\\
{(S,\cM_S)}\ar[r]^-(0.5){\iota^a}&{\bA_\mN}}
\end{equation}
est commutatif.
\item[(ii)] Le monoïde $P$ est torique, {\em i.e.}, $P$ est fin et saturé et $P^\gp$ est un $\mZ$-module libre (\cite{agt} II.5.1).
\item[(iii)] L'homomorphisme $\vartheta$ est saturé (\cite{agt} II.5.2).
\item[(iv)] L'homomorphisme $\vartheta^\gp\colon \mZ\rightarrow P^\gp$ est injectif, 
le sous-groupe de torsion de $\coker(\vartheta^\gp)$ est d'ordre premier à $p$ et le morphisme de schémas usuels
\begin{equation}\label{cad1c}
X\rightarrow S\times_{\bA_\mN}\bA_P
\end{equation}
déduit de \eqref{cad1b} est étale.  
\item[(v)] Posons $\lambda=\vartheta(1)\in P$, 
\begin{eqnarray}
L&=&\Hom_{\mZ}(P^\gp,\mZ),\label{cad1d}\\
\rH(P)&=&\Hom(P,\mN).\label{cad1e}
\end{eqnarray} 
On notera que $\rH(P)$ est un monoïde fin, saturé et affûté et que l'homomorphisme canonique 
$\rH(P)^\gp\rightarrow \Hom((P^\sharp)^\gp,\mZ)$ est un isomorphisme (\cite{ogus} I 2.2.3). 
On suppose qu'il existe $h_1,\dots,h_r\in \rH(P)$, qui sont $\mZ$-linéairement indépendants dans $L$, tels que  
\begin{equation}\label{cad1f}
\ker(\lambda)\cap \rH(P)=\{\sum_{i=1}^ra_ih_i | \ (a_1,\dots,a_r)\in \mN^r\},
\end{equation}
où l'on considère $\lambda$ comme un homomorphisme $L\rightarrow \mZ$. 
\end{itemize}

\vspace{2mm}
On pose $\pi=\iota(1)$ qui est une uniformisante de $\co_K$. On a alors \eqref{notconv2}
\begin{equation}\label{cad1h}
S\times_{\bA_\mN}\bA_P=\Spec(\co_K[P]/(\pi-e^\lambda)).
\end{equation}

\subsection{}\label{cad2}
Pour tout entier $n\geq 1$, on pose
\begin{eqnarray}\label{cad2a}
\co_{K_n}=\co_K[\zeta]/(\zeta^{n}-\pi),
\end{eqnarray}
qui est un anneau de valuation discrète. On note $K_n$ le corps des fractions de $\co_{K_n}$
et $\pi_n$ la classe de $\zeta$ dans $\co_{K_n}$, qui est une uniformisante de $\co_{K_n}$.  
On pose $S^{(n)}=\Spec(\co_{K_n})$
que l'on munit de la structure logarithmique $\cM_{S^{(n)}}$ définie par son point fermé. 
On désigne par $\iota_n\colon \mN\rightarrow \Gamma(S^{(n)},\cM_{S^{(n)}})$
l'homomorphisme défini par $\iota_n(1)=\pi_n$; c'est une carte pour $(S^{(n)},\cM_{S^{(n)}})$.  

Considérons le système inductif de monoïdes $(\mN^{(n)})_{n\geq 1}$, 
indexé par l'ensemble $\mZ_{\geq 1}$ ordonné par la relation de divisibilité, 
défini par $\mN^{(n)}=\mN$ pour tout $n\geq 1$ et dont l'homomorphisme de transition
$\mN^{(n)}\rightarrow \mN^{(mn)}$ (pour $m,n\geq 1$) est l'endomorphisme de Frobenius d'ordre $m$ 
de $\mN$ ({\em i.e.}, l'élévation à la puissance $m$-ième). On notera $\mN^{(1)}$ simplement $\mN$. Les  schémas logarithmiques
$(S^{(n)},\cM_{S^{(n)}})_{n\geq 1}$ forment naturellement un système projectif.
Pour tous entiers $m, n\geq 1$, avec les notations de  \ref{notconv2}, on a un diagramme cartésien de morphismes de schémas logarithmiques
\begin{equation}\label{cad2b}
\xymatrix{
{(S^{(mn)},\cM_{S^{(mn)}})}\ar[r]^-(0.5){\iota^a_{mn}}\ar[d]&{\bA_{\mN^{(mn)}}}\ar[d]\\
{(S^{(n)},\cM_{S^{(n)}})}\ar[r]^-(0.5){\iota^a_n}&{\bA_{\mN^{(n)}}}}
\end{equation}
où $\iota^a_n$ (resp. $\iota^a_{mn}$) est le morphisme associé à $\iota_n$ (resp. $\iota_{mn}$) (\cite{agt} II.5.13).

\subsection{}\label{cad3}
Considérons le système inductif de monoïdes $(P^{(n)})_{n\geq 1}$, 
indexé par l'ensemble $\mZ_{\geq 1}$ ordonné par la relation de divisibilité, 
défini par $P^{(n)}=P$ pour tout $n\geq 1$ et dont l'homomorphisme de transition
$i_{n,mn}\colon P^{(n)}\rightarrow P^{(mn)}$ (pour $m, n\geq 1$) est l'endomorphisme 
de Frobenius d'ordre $m$ de $P$ ({\em i.e.}, l'élévation à la puissance $m$-ième). Pour tout $n\geq 1$, on note
\begin{equation}\label{cad3a}
P\stackrel{\sim}{\rightarrow}P^{(n)}, \ \ \ t\mapsto t^{(n)},
\end{equation} 
l'isomorphisme canonique. Pour tout $t\in P$ et tous $m, n\geq 1$, on a donc
\begin{equation}\label{cad3b}
i_{n,mn}(t^{(n)})=(t^{(mn)})^m. 
\end{equation}
On notera $P^{(1)}$ simplement $P$. 

Pour tout entier $n\geq 1$, on pose (avec les notations de \ref{notconv2})
\begin{equation}\label{cad3c}
(X^{(n)},\cM_{X^{(n)}})=(X,\cM_{X})\times_{\bA_P}\bA_{P^{(n)}}.
\end{equation}
On notera que la projection canonique $(X^{(n)},\cM_{X^{(n)}})\rightarrow \bA_{P^{(n)}}$ est stricte.

Comme le diagramme \eqref{cad2b} est cartésien, il existe un unique morphisme
\begin{equation}\label{cad3d}
f^{(n)}\colon (X^{(n)},\cM_{X^{(n)}})\rightarrow (S^{(n)},\cM_{S^{(n)}}),
\end{equation} 
qui s'insère dans le diagramme commutatif
\begin{equation}\label{cad3e}
\xymatrix{
{(X^{(n)},\cM_{X^{(n)}})}\ar[rrr]\ar[ddd]\ar[rd]_{f^{(n)}}\ar@{}[drrr]|*+[o][F-]{1}&&&
{\bA_{P^{(n)}}}\ar[ld]^{\bA_{\vartheta}}\ar[ddd]\\
&{(S^{(n)},\cM_{S^{(n)}})}\ar[r]^-(0.5){\iota^a_n}\ar[d]&{\bA_{\mN^{(n)}}}\ar[d]&\\
&{(S,\cM_{S})}\ar[r]^-(0.5){\iota^a}&{\bA_\mN}&\\
{(X,\cM_{X})}\ar[rrr]\ar[ru]^f&&&{\bA_P}\ar[lu]_{\bA_\vartheta}}
\end{equation}

\begin{prop}[\cite{agt} II.6.6]\label{cad4}
Soient $n$ un entier $\geq 1$, $h\colon \oS\rightarrow S^{(n)}$ un $S$-morphisme. Alors,
\begin{itemize}
\item[{\rm (i)}] La face $\xymatrix{\ar@{}|*+[o][F-]{1}}$ du diagramme \eqref{cad3e} est une carte adéquate pour $f^{(n)}$ 
{\rm (\cite{agt} III.4.4)}; en particulier, $f^{(n)}$ est lisse et saturé. 
\item[{\rm (ii)}] Le schéma $X^{(n)}$ est intègre, normal, Cohen-Macaulay et $S^{(n)}$-plat. 
\item[{\rm (iii)}] Le schéma $X^{(n)}\times_{S^{(n)}}\oS$ est normal et localement irréductible {\rm (\cite{agt} III.3.1)}; 
il est donc la somme des schémas induits sur ses composantes irréductibles. 
\item[{\rm (iv)}] Si $n$ est une puissance de $p$, l'image inverse de toute composante irréductible de $\oX$ 
dans $X^{(n)}\times_{S^{(n)}}\oS$ est intègre.
\item[{\rm (v)}]  Le sous-schéma ouvert maximal de $X^{(n)}$ 
où la structure logarithmique $\cM_{X^{(n)}}$ est triviale est égal à $X^{(n)\circ}=X^{(n)}\times_XX^\circ$ \eqref{TFA3a}. 
\item[{\rm (vi)}] Le morphisme $X^{(n)\circ}\rightarrow X^\circ$ est étale, fini et surjectif~; 
plus précisément, $X^{(n)\circ}$ est un espace principal homogène pour la topologie étale 
au-dessus de $X^\circ$ sous le groupe 
\[
\Hom_\mZ(P^\gp,\mu_n(\oK)).
\]
\end{itemize}
\end{prop}

\subsection{}\label{cad5}
Pour tout entier $n\geq 1$ et tout $S$-morphisme $\tau\colon \oS\rightarrow S^{(n)}$, on pose \eqref{TFA3a}
\begin{equation}
Y^{(n)}_\tau=X^{(n)}\times_{S^{(n)}}\oeta \ \ \ {\rm et}\ \ \ Y^{(n)\circ}_\tau=Y^{(n)}_\tau\times_XX^\circ.
\end{equation}
On laissera tomber l'indice $\tau$ lorsqu'il n'y a aucun risque d'ambiguïté. 
On a $Y^{(1)}=X_\oeta$ et $Y^{(1)\circ}=\oX^\circ$. 
On note $g\colon \oX^\circ\rightarrow X_\oeta$ et $g_n\colon Y^{(n)\circ}\rightarrow Y^{(n)}$ les injections canoniques, et
$\varphi_n\colon Y^{(n)}\rightarrow X_\oeta$ et $\varphi_n^\circ\colon Y^{(n)\circ}\rightarrow \oX^\circ$ les morphismes canoniques.
Pour tout faisceau $\cF$ de $\oX^\circ_\et$ et tout entier $q\geq 0$, 
on a un morphisme de changement de base
\begin{equation}\label{cad5c}
\varphi_n^*(\rR^q g_*\cF)\rightarrow \rR^q g_{n*}(\varphi_n^{\circ *}(\cF)).
\end{equation}

\begin{lem}\label{cad6}
Il existe un entier $r\geq 0$ et un morphisme lisse $t\colon X_\eta\rightarrow \mA_\eta^r$ vérifiant les propriétés suivantes~:
\begin{itemize}
\item[{\rm (i)}] La structure logarithmique $\cM_X|X_\eta$ est définie par le diviseur à croisement normaux $D$ sur $X_\eta$ image inverse 
du diviseur des coordonnées sur $\mA_\eta^r$. En particulier, 
$X^\circ=t^{-1}(\mG_\eta^r)$, où $\mG_\eta^r$ est l'ouvert de $\mA_\eta^r$ où les coordonnées ne s'annulent pas.
\item[{\rm (ii)}] Pour tout entier $n\geq 1$, notons $Z^{(n)}$ le schéma défini par le diagramme cartésien
\begin{equation}
\xymatrix{
Z^{(n)}\ar[r]\ar[d]&X_\eta\ar[d]^t\\
{\mA_\eta^r}\ar[r]&{\mA_\eta^r}}
\end{equation}
où la flèche horizontale inférieure est définie par l'élévation à la puissance $n$-ième des coordonnées de $\mA_\eta^r$. 
Alors, pour tout $S$-morphisme $\tau\colon \oS\rightarrow S^{(n)}$, avec les notations de \ref{cad5}, il existe un $X_\oeta$-morphisme étale et fini 
\begin{equation}
Y^{(n)}_\tau\rightarrow Z^{(n)}_\oeta.
\end{equation}
\end{itemize}
\end{lem}

D'après (\cite{agt} II.6.3(v)), le schéma (usuel) $X_\eta$ est lisse sur $\eta$, $X^\circ$ est l'ouvert 
complémentaire dans $X_\eta$ d'un diviseur à croisements normaux stricts $D$ et  $\cM_X|X_\eta$ 
est la structure logarithmique sur $X_\eta$ définie par $D$. Rappelons brièvement la démonstration de {\em loc. cit.} 
Soit $F$ la face de $P$ engendrée par $\lambda$, c'est-à-dire l'ensemble des éléments 
$x\in P$ tels qu'il existe $y\in P$ et $m\in \mN$ tels que $x+y=m\lambda$ (cf. \cite{ogus} I 1.4.2).
On note $F^{-1}P$ la localisation de $P$ par $F$ (\cite{ogus} I 1.4.4), $Q$ le sous-groupe des unités de $F^{-1}P$,
$(F^{-1}P)^\sharp$ le quotient de $F^{-1}P$ par $Q$ (\cite{agt} II.5.1) et $P/F$ le conoyau dans la catégorie des monoïdes de l'injection canonique 
$F\rightarrow P$ (\cite{ogus} I 1.1.5). Il résulte aussitôt des propriétés universelles des localisations de monoïdes et d'anneaux que 
l'homomorphisme canonique $\mZ[P]\rightarrow \mZ[F^{-1}P]$ induit un isomorphisme
\begin{equation}\label{cad6b}
\mZ[P]_{\lambda}\stackrel{\sim}{\rightarrow} \mZ[F^{-1}P],
\end{equation}
de sorte qu'on a un diagramme cartésien
\begin{equation}\label{cad6c}
\xymatrix{
{X_\eta}\ar[r]\ar[d]&{\bA_{F^{-1}P}}\ar[d]\\
X\ar[r]&{\bA_P}}
\end{equation} 
Par ailleurs, on a un isomorphisme canonique 
\begin{equation}\label{cad6d}
P/F\stackrel{\sim}{\rightarrow}(F^{-1}P)^\sharp.
\end{equation} 
La preuve de (\cite{agt} II.6.3(v)) montre que $P/F$ est un monoïde libre de type fini. On en déduit un isomorphisme 
\begin{equation}\label{cad6e}
Q\times (P/F) \stackrel{\sim}{\rightarrow} F^{-1}P.
\end{equation}

(i) Comme le monoïde $P/F$ est libre de type fini, il existe un entier $r\geq 0$ et un isomorphisme $\Spec(K[P/F])\simeq \mA_\eta^r$.
Les diagrammes \eqref{cad1b} et \eqref{cad6c} et l'isomorphisme \eqref{cad6e} induisent donc un morphisme 
\begin{equation}\label{cad6i}
X_\eta\rightarrow \mA_\eta^r\times_\eta\Spec(K[Q]/(\pi-e^\lambda)),
\end{equation}
et par suite un morphisme $t\colon X_\eta\rightarrow \mA_\eta^r$. La structure logarithmique $\cM_X|X_\eta$ sur $X_\eta$ est l'image inverse par $t$
de la structure logarithmique sur $\mA_\eta^r$ définie par les axes des coordonnées (\cite{agt} II.5.10).
Compte tenu de \ref{cad1}(iv), le morphisme \eqref{cad6i} est étale et le morphisme $t$ est lisse; d'où la proposition.

(ii) Notant $u$ l'endomorphisme de Frobenius d'ordre $n$ de $P$, le diagramme \eqref{cad3e}  induit un diagramme commutatif à carrés cartésiens
\begin{equation}
\xymatrix{
{X^{(n)}}\ar[r]\ar[d]&{S^{(n)}\times_{\bA_\mN}\bA_P}\ar[r]\ar[d]&{\bA_P}\ar[d]^{\bA_u}\\
X\ar[r]&{S\times_{\bA_\mN}\bA_P}\ar[r]&{\bA_P}}
\end{equation}
On en déduit un diagramme cartésien 
\begin{equation}
\xymatrix{
{X^{(n)}\times_{S^{(n)},\tau}\oS}\ar[r]\ar[d]&{\Spec(\co_{\oK}[P]/(\pi_n-e^\lambda))}\ar[d]\\
{\oX}\ar[r]&{\Spec(\co_{\oK}[P]/(\pi-e^\lambda))}}
\end{equation}
et par suite, compte tenu de \eqref{cad6b}, un diagramme cartésien
\begin{equation}
\xymatrix{
{Y^{(n)}_\tau}\ar[r]\ar[d]_{\varphi_n}&{\Spec(\oK[F^{-1}P]/(\pi_n-e^\lambda))}\ar[d]\\
{X_\oeta}\ar[r]&{\Spec(\oK[F^{-1}P]/(\pi-e^\lambda))}}
\end{equation}
où les flèches verticales de droite sont induites par les homomorphismes de Frobenius d'ordre $n$ de $P$ et $F^{-1}P$. 
On notera que le morphisme 
\begin{equation}
\Spec(\oK[Q]/(\pi_n-e^\lambda))\rightarrow \Spec(\oK[Q]/(\pi-e^\lambda))
\end{equation}
induit par l'homomorphisme de Frobenius d'ordre $n$ de $Q$, est étale et fini. 
Compte tenu de \eqref{cad6e} et des définitions, on en déduit un $X_\oeta$-morphisme étale et fini $Y^{(n)}_\tau\rightarrow Z^{(n)}_\oeta$.

\begin{lem}\label{cad7}
Soient $n$ un entier $\geq 1$, $F$ un faisceau de $(\mZ/n\mZ)$-modules 
localement constant et constructible de $X_{\oeta,\et}$, $\tau\colon \oS\rightarrow S^{(n)}$ un $S$-morphisme, $q$ un entier $\geq 1$. 
Alors, avec les notations de \ref{cad5}, le morphisme de changement de base \eqref{cad5c}
\begin{equation}
\varphi_n^*(\rR^q g_*(g^*F))\rightarrow \rR^q g_{n*}(g_n^*(\varphi_n^*(F)))
\end{equation}
est nul.
\end{lem}

Cela résulte de \ref{cad6} et \ref{Kpun14}.

\begin{lem}\label{cad8}
Pour tout faisceau abélien de torsion localement constant et constructible $F$ de $\oX^\circ_\et$, il existe un entier $n\geq 1$ 
tel que pour tout $S$-morphisme $\tau\colon \oS\rightarrow S^{(n)}$, avec les notations de \ref{cad5}, $\varphi^{\circ*}_n(F)$ se prolonge en un 
faisceau abélien de torsion localement constant et constructible de $Y^{(n)}_\et$. 
\end{lem}

Par descente, $F$ est représentable par un schéma en groupes abéliens étale et fini $V$ au-dessus de $\oX^\circ$.  
Notons $T$ la fermeture intégrale de $X_\oeta$ dans $V$. Soit $n$ un entier $\geq 1$ et divisible par tous les indices de ramification
du morphisme $T\rightarrow X_\oeta$ en les points de codimension un de $T$. Reprenons les notations de \ref{cad6} et notons 
$T^{(n)}$ la fermeture intégrale de $Z_\oeta^{(n)}$ dans $V\times_{X_\oeta}Z^{(n)}_\oeta$. On observera que $Z^{(n)}$ est lisse sur $\eta$.  
D'après le lemme d'Abhyankar (\cite{sga1} X 3.6), $T^{(n)}\rightarrow Z_\oeta^{(n)}$ est étale au-dessus des points de codimension un de $Z_\oeta^{(n)}$. 
Par suite, le morphisme canonique $T^{(n)}\rightarrow Z_\oeta^{(n)}$ est étale en vertu du théorème de pureté de Zariski-Nagata (\cite{sga2} X 3.4). 
Par suite, $T^{(n)}\times_{Z_\oeta^{(n)}}T^{(n)}$ est la fermeture intégrale de $Z_\oeta^{(n)}$ dans 
$(V\times_{\oX^\circ}V)\times_{X_\oeta}Z^{(n)}_\oeta$. On en déduit une structure de schéma en groupes abéliens sur $T^{(n)}$ au-dessus de $Z_\oeta^{(n)}$
qui prolonge celle de $V\times_{X_\oeta}Z^{(n)}_\oeta$ au-dessus de $\oX^\circ\times_{X_\oeta}Z^{(n)}_\oeta$. Le faisceau 
$F|(\oX^\circ\times_{X_\oeta}Z^{(n)}_\oeta)$ se prolonge donc en un faisceau abélien de torsion localement constant et constructible de $(Z^{(n)}_\oeta)_\et$.
La proposition s'ensuit compte tenu de \ref{cad6}(ii).

\begin{prop}\label{cad9}
Soient $F$ un faisceau abélien de torsion, localement constant et constructible de $\oX^\circ_\et$, $V\rightarrow \oX^\circ$ un revêtement étale, 
$q$ un entier $\geq 1$, $\xi\in \rH^q(V,F)$, $\ox$ un point géométrique de $X$.
Alors, il existe un voisinage étale $U$ de $\ox$ dans $X$ et un revêtement étale surjectif $W\rightarrow V\times_{\oX^\circ}\oU^\circ$ tel que l'image canonique
de $\xi$ dans $\rH^q(W,F)$ soit nulle. 
\end{prop}

Par la preuve de l'implication (iv)$\Rightarrow$(v) de \ref{Kpun1}, on peut se borner au cas où $V=\oX^\circ$. 
En vertu de \ref{cad8}, il existe un entier $n\geq 1$ tel que pour tout $S$-morphisme $\tau\colon \oS\rightarrow S^{(n)}$, 
avec les notations de \ref{cad5}, $F|Y^{(n)\circ}_\tau$ se prolonge en un faisceau abélien de torsion localement constant et constructible de $Y^{(n)}_{\tau,\et}$.  
D'après  \ref{cad4}(i), le morphisme $f_n\colon (X^{(n)},\cM_{X^{(n)}})\rightarrow (S^{(n)},\cM_{S^{(n)}})$ est adéquat, 
et il vérifie les conditions de \ref{cad1}. Supposons que pour tout $S$-morphisme $\tau\colon \oS\rightarrow S^{(n)}$, 
la proposition soit démontrée pour l'image de $\xi$ dans  $\rH^q(Y^{(n)\circ}_\tau,F)$ et tout point
géométrique de $X^{(n)}$ au-dessus de $\ox$. Il existe donc un morphisme étale $U^{(n)}_\tau\rightarrow X^{(n)}$ tel que  le morphisme
\begin{equation}
U^{(n)}_\tau\otimes_X\kappa(\ox)\rightarrow X^{(n)}\otimes_X\kappa(\ox)
\end{equation}
soit surjectif, et un revêtement étale surjectif $W_\tau\rightarrow U^{(n)\circ}_\tau\times_{S^{(n)},\tau}\oS$ tel que l'image canonique de $\xi$ dans 
$\rH^q(W_\tau,F)$ soit nulle.
D'après \ref{Kpun15}, il existe un $X$-schéma étale $\ox$-pointé $U$ et pour tout $S$-morphisme $\tau\colon \oS\rightarrow S^{(n)}$,
un $X^{(n)}$-morphisme $U\times_XX^{(n)}\rightarrow U^{(n)}_\tau$. Considérons le diagramme commutatif à carrés cartésiens 
\begin{equation}
\xymatrix{
{W_\tau}\ar[d]&{W_\tau\times_{U^{(n)}_\tau}(U\times_XX^{(n)})}\ar[d]^{v_\tau}\ar[l]&\\
{U^{(n)\circ}_\tau\times_{S^{(n)},\tau}\oS}\ar[d]&{(U^\circ\times_XX^{(n)})\times_{S^{(n)},\tau}\oS}\ar[l]\ar[r]^-(0.5){u_\tau}\ar[d]&{\oU^\circ}\\
{U^{(n)}_\tau}&{U\times_XX^{(n)}}\ar[l]&}
\end{equation}
Le couple formé de $U$ et du revêtement étale surjectif de $\oU^\circ$ somme disjointe des $u_\tau\circ v_\tau$ pour tous les $\tau$, 
répond alors à la question. On peut donc se borner au cas où $F$ se prolonge en un faisceau abélien de torsion, 
localement constant et constructible $G$ de $X_{\oeta,\et}$. L'isomorphisme $g^*(G)\stackrel{\sim}{\rightarrow} F$ \eqref{cad5}
induit par adjonction un morphisme $G\rightarrow g_*(F)$ qui est en fait un isomorphisme; en particulier, $g_*(F)$ est localement constant et constructible.
En effet, la question étant locale pour la topologie étale de $X$, on peut 
supposer $G$ constant, auquel cas l'assertion recherchée résulte de (\cite{sga4} IX 2.14.1). 

On désigne par $g\colon \oX^\circ\rightarrow X_\oeta$ l'injection canonique. Considérons la suite spectrale de Cartan-Leray
\begin{equation}
E_2^{a,b}=\rH^a(X_\oeta,\rR^bg_*(F))\Rightarrow \rH^{a+b}(\oX^\circ,F),
\end{equation} 
et notons $(E^q_i)_{0\leq i\leq q}$ la filtration aboutissement sur $\rH^q(\oX^\circ,F)$, de sorte que 
\begin{equation}
E_i^q/E_{i+1}^q=E_\infty^{i,q-i}
\end{equation} 
pour tout $0\leq i\leq q$. Il existe un entier $0\leq i_\xi\leq q$ tel que $\xi\in E_{i_\xi}^q\backslash E_{i_\xi+1}^q$. 
Soit $n$ un entier $\geq 1$ qui annule $F$. Reprenons de nouveau les notations de \ref{cad5}.  
Il résulte de (\cite{sga4} XII 4.4) qu'on a un morphisme de suites spectrales de Cartan-Leray
\begin{equation}
\xymatrix{
{\rH^a(X_\oeta,\rR^bg_*(F))}\ar[d]_{u_{a,b}}\ar@{=>}[r]&{\rH^{a+b}(\oX^\circ,F)}\ar[d]^{u}\\
{\rH^a(Y^{(n)},\rR^b g_{n*}(\varphi_n^{\circ*}(F)))}\ar@{=>}[r]&{\rH^{a+b}(Y^{(n)\circ},\varphi_n^{\circ*}(F))}}
\end{equation}
où $u_{a,b}$ est induit par le morphisme de changement de base \eqref{cad5c} et 
$u$ est le morphisme canonique. Pour tout $b\geq 1$, le morphisme $u_{a,b}$ est nul, en vertu de \ref{cad7}. Il en est alors de même du morphisme
induit par $u$ entre les gradués des filtrations aboutissement.
Procédant comme plus haut, on peut donc se réduire par une récurrence finie, au cas où $i_\xi=q$. 
Par suite, $\xi$ est l'image canonique d'un élément $\zeta\in \rH^q(X_\oeta,g_*(F))$. En vertu de \ref{Kpun21} et \ref{Kpun11}, 
il existe un voisinage étale $U$ de $\ox$ dans $X$ et un revêtement étale surjectif $W'\rightarrow U_\oeta$ tel que l'image canonique
de $\zeta$ dans $\rH^q(W',g_*(F))$ soit nulle. Par suite, l'image canonique de $\xi$ dans $\rH^q(W'^\circ,F)$ est nulle.

\begin{cor}[\cite{achinger} 9.5]\label{cad10}
Soient $\ox$ un point géométrique de $X$ au-dessus de $s$, 
$X'$ le localisé strict de $X$ en $\ox$. Alors, $\oX'^\circ$ est un schéma $K(\pi,1)$ \eqref{Kpun2}.
\end{cor}

En effet, $\oX'$ étant normal et strictement local (et en particulier intègre) d'après (\cite{agt} III.3.7), 
la proposition résulte de \ref{Kpun18} et \ref{cad9}. 

\begin{rema} 
La proposition \ref{cad9} lorsque le schéma (usuel) $X$ est lisse sur $S$, est due à Faltings (\cite{faltings1} Lemme 2.3 page 281). 
Le cas général est dû à Achinger (\cite{achinger} 9.5). 
On notera que le cas où $\ox$ est un point géométrique de $X_\eta$ a déjà été démontré dans \ref{Kpun16}.  
\end{rema}

\section{Conséquences du théorème de pureté de Faltings}\label{cg}

\subsection{}\label{cg21}
On suppose dans cette section que le morphisme 
$f\colon (X,\cM_X)\rightarrow (S,\cM_S)$ \eqref{TFA3} remplit les conditions de \ref{cad1}, que $X=\Spec(R)$ est affine et connexe
et que $X_s$ est non-vide. Ces conditions correspondent à celles fixées dans (\cite{agt} II.6.2). On reprend les notations de §~\ref{cad}. 
On pose \eqref{TFA3b}
\begin{equation}\label{cg21d}
\tOmega^1_{R/\co_K}=\tOmega^1_{X/S}(X).
\end{equation}
D'après \ref{cad1}(iv) et (\cite{kato1} 1.8), on a un isomorphisme $R$-linéaire canonique
\begin{equation}\label{cg21e}
\tOmega^1_{R/\co_K}\stackrel{\sim}{\rightarrow} (P^\gp/\lambda\mZ)\otimes_\mZ R.
\end{equation}
En particulier, $\tOmega^1_{R/\co_K}$ est un $R$-module libre de rang fini.

Pour tout entier $n\geq 1$, le schéma $X^{(n)}$ est affine d'anneau 
\begin{equation}\label{cg21c}
A_n=R\otimes_{\mZ[P]}\mZ[P^{(n)}].
\end{equation}
On désigne par $A_n^\circ$ l'anneau du schéma affine $X^{(n)\circ}=X^{(n)}\times_XX^\circ$ \eqref{TFA3a}. 
Les schémas logarithmiques $(X^{(n)},\cM_{X^{(n)}})_{n\geq 1}$ forment naturellement un système projectif au-dessus 
de $(\bA_{P^{(n)}})_{n\geq 1}$, indexé par l'ensemble $\mZ_{\geq 1}$ ordonné par la relation de divisibilité. On pose 
\begin{eqnarray}
A_\infty&=&\underset{\underset{n\geq 1}{\longrightarrow}}{\lim}\ A_n,\label{cg21a}\\
A_\infty^\circ&=&\underset{\underset{n\geq 1}{\longrightarrow}}{\lim}\ A_n^\circ.\label{cg21b}
\end{eqnarray}

\begin{teo}[Faltings, \cite{faltings2} §~2b]\label{cg22}
Pour toute $A^\circ_\infty$-algèbre étale finie $B'$, la fermeture intégrale  $B$  de $A_\infty$ dans $B'$
est une extension $\alpha$-étale de $A_\infty$ \eqref{aet2}.
\end{teo}

Signalons ici que Scholze dispose, dans le cadre de sa théorie des perfectoïdes,  d'une généralisation de ce  résultat  
(\cite{scholze} 1.10 et 7.9).

\subsection{}\label{cg2}
Soit $\oy$ un point géométrique de $\oX^\circ$ que l'on suppose fixé dans la suite de cette section. 
Le schéma $\oX$ étant localement irréductible d'après \ref{cad4}(iii),  
il est la somme des schémas induits sur ses composantes irréductibles. On note $\oX^\star$
la composante irréductible de $\oX$ contenant $\oy$. 
De même, $\oX^\circ$ est la somme des schémas induits sur ses composantes irréductibles
et $\oX^{\star \circ}=\oX^\star\times_{X}X^\circ$ est la composante irréductible de $\oX^\circ$ contenant $\oy$. 
On désigne par $\Delta$ le groupe profini $\pi_1(\oX^{\star \circ},\oy)$ 
et pour alléger les notations, par $\oR$ la représentation discrète $\oR^\oy_X$ de $\Delta$ définie dans \eqref{TFA9b}. 
On notera que $\oR$ est intègre et normal \eqref{TFA10}.

\subsection{}\label{cg6}
Pour tous entiers $m,n\geq 1$, le morphisme canonique $X^{(mn)}\rightarrow X^{(n)}$ 
est fini et surjectif. D'après (\cite{ega4} 8.3.8(i)), il existe alors un $X$-morphisme \eqref{cg2}
\begin{equation}\label{cg6a}
\oy\rightarrow \underset{\underset{n\geq 1}{\longleftarrow}}{\lim}\ X^{(n)},
\end{equation} 
la limite inductive étant indexée par l'ensemble $\mZ_{\geq 1}$ ordonné par la relation de divisibilité.
On fixe un tel morphisme dans toute la suite de cette section. Celui-ci induit un $S$-morphisme
\begin{equation}\label{cg6b}
\oS\rightarrow \underset{\underset{n\geq 1}{\longleftarrow}}{\lim}\ S^{(n)}.
\end{equation}
Pour tout entier $n\geq 1$, on pose
\begin{equation}\label{cg6c}
\oX^{(n)}= X^{(n)}\times_{S^{(n)}}\oS \ \ \ {\rm et}\ \ \ \oX^{(n)\circ}=\oX^{(n)}\times_XX^\circ.
\end{equation} 
On en déduit un $\oX$-morphisme 
\begin{equation}\label{cg6d}
\oy\rightarrow \underset{\underset{n\geq 1}{\longleftarrow}}{\lim}\ \oX^{(n)}.
\end{equation} 

Pour tout entier $n\geq 1$, le schéma $\oX^{(n)}$ étant localement irréductible d'après \ref{cad4}(iii), 
il est la somme des schémas induits sur ses composantes irréductibles. 
On note $\oX^{(n)\star}$ la composante irréductible de $\oX^{(n)}$ contenant l'image de $\oy$ \eqref{cg6d}.
De même, $\oX^{(n)\circ}$ est la somme des schémas induits sur ses composantes irréductibles
et $\oX^{(n)\star\circ}=\oX^{(n)\star}\times_XX^\circ$  est la composante irréductible de $\oX^{(n)\circ}$ 
contenant l'image de $\oy$. On notera que $\oX^{(n)}$ étant fini sur $\oX$ \eqref{cad3e},  
$\oX^{(n)\star}$ est la fermeture intégrale de $\oX^{\star}$ dans $\oX^{(n)\star\circ}$ d'après \ref{cad4}(iii). 
On pose 
\begin{equation}\label{cg6f}
R_n=\Gamma(\oX^{(n)\star},\co_{\oX^{(n)}}).
\end{equation}
D'après \ref{cad4}(vi), le morphisme $\oX^{(n)\star\circ}\rightarrow \oX^{\star\circ}$ 
est étale fini. Il résulte de la preuve de (\cite{agt} II.6.8(iv)) que $\oX^{(n)\star \circ}$ 
est en fait un revêtement étale fini et galoisien de $\oX^{\star \circ}$ de groupe $\Delta_n$ 
canoniquement isomorphe à un sous-groupe de $\Hom_\mZ(P^\gp/\mZ\lambda,\mu_{n}(\oK))$.  
Le groupe $\Delta_n$ agit naturellement sur $R_n$. 
Le morphisme \eqref{cg6d} induit un $\oX^\star$-morphisme 
\begin{equation}\label{cg6e}
\Spec(\oR)\rightarrow \underset{\underset{n\geq 1}{\longleftarrow}}{\lim}\ \oX^{(n)\star}. 
\end{equation} 
Si $n$ est une puissance de $p$, on a $\Delta_n\simeq \Hom_\mZ(P^\gp/\mZ\lambda,\mu_{n}(\oK))$ 
et $\oX^{(n)\star}\simeq\oX^{(n)}\times_\oX\oX^\star$ en vertu de \ref{cad4}(iv).

Les anneaux $(R_n)_{n\geq 1}$ forment naturellement un système inductif. On pose 
\begin{eqnarray}
R_\infty&=&\underset{\underset{n\geq 1}{\longrightarrow}}{\lim}\ R_n,\label{cg6g}\\
R_{p^\infty}&=&\underset{\underset{n\geq 0}{\longrightarrow}}{\lim}\ R_{p^n}.\label{cg6gg}
\end{eqnarray}
Ce sont des anneaux normaux et intègres d'après (\cite{ega1n} 0.6.1.6(i) et 0.6.5.12(ii)).
Le morphisme \eqref{cg6e} induit des homomorphismes injectifs de $R_1$-algèbres 
\begin{equation}\label{cg6h}
R_{p^\infty}\rightarrow R_\infty\rightarrow \oR.
\end{equation}

Les groupes $(\Delta_n)_{n\geq 1}$ forment naturellement un système projectif. On pose 
\begin{eqnarray}
\Delta_\infty&=&\underset{\underset{n\geq 1}{\longleftarrow}}{\lim}\ \Delta_n,\label{cg6k}\\
\Delta_{p^\infty}&=&\underset{\underset{n\geq 0}{\longleftarrow}}{\lim}\ \Delta_{p^n}.\label{cg6kk}
\end{eqnarray}
On a des homomorphismes canoniques 
\begin{equation}\label{cg6l}
\xymatrix{
{\Delta_{\infty}}\ar@{->>}[d]\ar@{^(->}[r]&{\Hom_{\mZ}(P^\gp/\mZ\lambda,\hmZ(1))}\ar[d]\\
{\Delta_{p^\infty}}\ar[r]^-(0.5)\sim&{\Hom_{\mZ}(P^\gp/\mZ \lambda,\mZ_p(1))}}
\end{equation}
Le noyau $\Sigma_0$ de l'homomorphisme canonique $\Delta_\infty\rightarrow \Delta_{p^\infty}$
est un groupe profini d'ordre premier à $p$. Par ailleurs, le morphisme \eqref{cg6a} détermine un homomorphisme 
surjectif $\Delta\rightarrow \Delta_\infty$. On note $\Sigma$ son noyau. Les homomorphismes \eqref{cg6h} sont alors $\Delta$-équivariants.

\begin{equation}\label{cg6m}
\xymatrix{
R_1\ar[rr]^{\Delta_{p^\infty}}\ar@/^2pc/[rrrr]|{\Delta_\infty}\ar@/_2pc/[rrrrrr]|{\Delta}&&
{R_{p^\infty}}\ar[rr]^{\Sigma_0}&&{R_\infty}\ar[rr]^\Sigma&&\oR}
\end{equation}
On note $\oF$ (resp. $F_\infty$, resp. $F_{p^\infty}$) le corps des fractions de $\oR$ (resp. $R_\infty$, resp. $R_{p^\infty}$).

Comme le sous-groupe de torsion de $P^\gp/\mZ\lambda$ est d'ordre premier à $p$,  on a un isomorphisme canonique
\begin{equation}\label{cg6i}
(P^\gp/\mZ\lambda)\otimes_\mZ\mZ_p(-1)\stackrel{\sim}{\rightarrow}
\Hom_{\mZ_p}(\Delta_{p^\infty}, \mZ_p).
\end{equation}
Compte tenu de \eqref{cg21e}, on en déduit, pour tout élément non nul $a$ de $\co_\oK$, un isomorphisme $R_1$-linéaire
\begin{equation}\label{cg6j}
\tOmega^1_{R/\co_K}\otimes_R(R_1/aR_1)(-1)\stackrel{\sim}{\rightarrow} 
\Hom_{\mZ}(\Delta_{p^\infty},R_1/aR_1).
\end{equation}

\begin{prop}[\cite{agt} II.6.17] \label{cg35}
L'extension $\oF$ de $F_\infty$ est la réunion d'un système inductif filtrant de 
sous-extensions galoisiennes finies $N$ de $F_\infty$ telles que la clôture intégrale
de $R_\infty$ dans $N$ soit $\alpha$-étale sur $R_\infty$ \eqref{aet2}.
\end{prop}

\begin{prop}[\cite{agt} II.8.17]\label{cg9}
Pour tout élément non nul $a$ de $\co_\oK$, il existe un unique homomorphisme de $R_1$-algèbres graduées 
\begin{equation}\label{cg9a}
\wedge(\tOmega^1_{R/\co_K}\otimes_R(R_1/aR_1)(-1))\rightarrow \rH^*(\Delta,\oR/a\oR)
\end{equation}
dont la composante en degré un est induite par \eqref{cg6j}. 
Celui-ci est $\alpha$-injectif et son conoyau est annulé par $p^{\frac{1}{p-1}}\fm_\oK$.   
\end{prop}

\begin{prop}[\cite{agt} II.9.10]\label{cg67} 
L'endomorphisme de Frobenius absolu de $\oR/p\oR$ est surjectif.
\end{prop}

\begin{lem}\label{cg36}
Soient $m$ un entier $\geq 1$, $M$ un $(R_1/p^mR_1)$-module de présentation finie, 
muni d'une action linéaire et continue de $\Delta_{p^\infty}$ (resp. $\Delta_{\infty}$). Alors, 
\begin{itemize}
\item[{\rm (i)}] Pour tout entier $i\geq 0$, $\rH^i(\Delta_{p^\infty},M)$ (resp. $\rH^i(\Delta_{\infty},M)$)
est un $R_1$-module de présentation finie, et est nul pour tout $i\geq d+1=\dim(X/S)+1$.
\item[{\rm (ii)}] Pour tout $\gamma\in \fm_\oK$, il existe un entier $n_0\geq 1$ tel que pour tout entier $n\geq n_0$, 
tout homomorphisme {\em surjectif} $\nu\colon \Delta_{p^\infty}\rightarrow \mu_{p^n}(\co_\oK)$ et tout entier $i\geq 0$, 
notant $\co_\oK(\nu)$ le $(\co_\oK)$-$\Delta_{p^\infty}$-module topologique $\co_\oK$, 
muni de la topologie $p$-adique et de l'action de $\Delta_{p^\infty}$ définie par $\nu$,
$\rH^i(\Delta_{p^\infty},M\otimes_{\co_\oK}\co_\oK(\nu))$ soit annulé par $\gamma$.
\end{itemize}
\end{lem}

(i) Supposons d'abord que $M$ soit une représentation de $\Delta_\infty$. Soient $n$ un entier tel que l'action de $\Delta_\infty$ sur $M$ se factorise à travers 
$\Delta_n$, $p^t$ la plus grande puissance de $p$ divisant $n$, $H_n$ le noyau de l'homomorphisme canonique $\Delta_n\rightarrow \Delta_{p^t}$.
On a donc un homomorphisme canonique surjectif $\Sigma_0\rightarrow H_n$ \eqref{cg6l}. 
L'anneau $R_1/p^mR_1$ étant cohérent (\cite{egr1} 1.10.3), $M$ est un $(R_1/p^mR_1)$-module cohérent.
On en déduit que le $(R_1/p^mR_1)$-module $M^{\Sigma_0}=M^{H_n}$ est cohérent. 
Comme $\Sigma_0$ est un groupe profini d'ordre premier à $p$ \eqref{cg6l}, pour tout $i\geq 0$,
on a un isomorphisme canonique $\rH^i(\Delta_{p^\infty},M^{\Sigma_0})\stackrel{\sim}{\rightarrow} \rH^i(\Delta_\infty,M)$. 
On peut donc se borner au cas où $M$ est une représentation de $\Delta_{p^\infty}$. 
D'après (\cite{agt} II.3.23), $\rR\Gamma(\Delta_{p^\infty},M)$ est quasi-isomorphe à un complexe borné, en degrés compris entre $0$ et $d$, 
de $(R_1/p^mR_1)$-modules cohérents~; d'où la proposition.

(ii) Soit $e_1,\dots,e_d$ une $\mZ_p$-base de $\Delta_{p^\infty}$. 
Il existe un entier $t\geq 1$ tel que pour tout $1\leq j\leq d$, $e_j^{p^t}$ agisse trivialement sur $M$. 
Soient $n$ un entier $\geq t$, $\nu\colon \Delta_{p^\infty}\rightarrow \mu_{p^n}(\co_\oK)$ un homomorphisme surjectif.
Il existe $1\leq j\leq d$ tel que que $\zeta_n=\nu(e_j)$ soit une racine primitive $p^n$-ième de l'unité.
Notons $G_j$ le sous-groupe de $\Delta_{p^\infty}$ engendré par $e_j$. 
D'après (\cite{agt} II.3.23), $\rH^i(G_j,M\otimes_{\co_\oK}\co_\oK(\nu))$ est la cohomologie du complexe, concentré en degrés $0, 1$, suivant 
\begin{equation}
\zeta_n e_j-\id_M\colon M\rightarrow M.
\end{equation}
La relation 
\begin{equation}
(\zeta_n e_j-\id_M)((\zeta_n e_j)^{p^t-1}+(\zeta_n e_j)^{p^t-2}+\dots+\id_M)=(\zeta_n^{p^t}-1)\id_M
\end{equation}
montre que le noyau et le conoyau de $\zeta_n e_j -\id_M$ sont annulés par $\zeta_n^{p^t}-1$. Par suite, compte tenu de la suite spectrale (\cite{agt} II.3.4)
\begin{equation}
E_2^{a,b}=\rH^a(\Delta_{p^\infty}/G_j,\rH^b(G_j,M))\Rightarrow \rH^{a+b}(\Delta_{p^\infty}, M),
\end{equation}
pour tout $i\geq 0$, $\rH^i(\Delta_{p^\infty}, M)$ est annulé par $\zeta_n^{p^t}-1$.
Pour tout $\gamma\in \fm_\oK$, il existe $n_0\geq t$ tel que si $n\geq n_0$, alors $\gamma\in (\zeta_n^{p^t}-1) \co_\oK$ ; d'où la proposition.

\begin{lem}\label{cg68}
Soient $m, n$ deux entiers $\geq 1$. Alors, 
\begin{itemize}
\item[{\rm (i)}] Pour tout $(R_1/p^mR_1)$-module de présentation $\alpha$-finie $M$, 
muni d'une action linéaire de $\Delta_{p^n}$ et tout entier $i\geq 0$, le $R_1$-module $\rH^i(\Delta_{p^n},M)$ est de présentation $\alpha$-finie.
\item[{\rm (ii)}] Pour tout $(R_{p^n}/p^mR_{p^n})$-module de présentation finie $M$, 
muni d'une action semi-linéaire de $\Delta_{p^n}$ et tout entier $i\geq 0$, 
le $R_1$-module $\rH^i(\Delta_{p^\infty},M\otimes_{R_{p^n}}R_{p^\infty})$ est de présentation $\alpha$-finie, 
et est nul pour tout $i\geq d+1$.
\item[{\rm (iii)}] Pour tout $(R_{n}/p^mR_{n})$-module de présentation finie $M$, 
muni d'une action semi-linéaire de $\Delta_n$ et tout entier $i\geq 0$, 
le $R_1$-module $\rH^i(\Delta_{\infty},M\otimes_{R_n}R_\infty)$ est de présentation $\alpha$-finie, 
et est nul pour tout $i\geq d+1$.
\end{itemize}
\end{lem}

On notera d'abord que $R_1/p^mR_1$ est un anneau cohérent (\cite{egr1} 1.10.3).

(i) D'après \ref{afini6}, 
le $(R_1/p^mR_1)$-module $M$ est $\alpha$-cohérent et il s'agit de montrer qu'il en est de même de $\rH^i(\Delta_{p^n},M)$ pour tout $i\geq 0$. 
Compte tenu de \ref{finita17}, \ref{finita18} et de la suite spectrale de Hochschild-Serre (\cite{agt} II.3.4), on peut se borner au cas où $d=1$, de sorte que
le groupe $\Delta_{p^n}$ est cyclique d'ordre $p^n$. Soient $\sigma$ un générateur de $\Delta_{p^n}$, $\Tr$ l'endomorphisme de $M$ induit par 
$\sum_{\sigma\in \Delta_{p^n}}\sigma$. Alors, $\rR\Gamma(\Delta_{p^n},M)$ est quasi-isomorphe au complexe de $(R_1/p^mR_1)$-modules, 
placé en degrés $\geq 0$, 
\begin{equation}
M\stackrel{\sigma-1}{\longrightarrow} M\stackrel{\Tr}{\longrightarrow} M\stackrel{\sigma-1}{\longrightarrow}
M\stackrel{\Tr}{\longrightarrow}M\stackrel{\sigma-1}{\longrightarrow}  \dots,
\end{equation}
d'où la proposition \eqref{finita18}.

(ii) En effet, la $p$-dimension cohomologique de $\Delta_{p^\infty}$ est égale à $d$ (\cite{agt} II.3.24). 
Il suffit donc de montrer que pour tout $i\geq 0$, le $(R_1/p^mR_1)$-module $\rH^i(\Delta_{p^\infty},M\otimes_{R_{p^n}}R_{p^\infty})$ est 
de présentation $\alpha$-finie. Comme $R_{p^n}/p^mR_{p^n}$ est une $(R_1/p^mR_1)$-algèbre cohérente, 
tout $(R_{p^n}/p^mR_{p^n})$-module de présentation $\alpha$-finie est aussi un 
$(R_1/p^mR_1)$-module de présentation $\alpha$-finie. D'après (i) et la suite spectrale de Hochschild-Serre (\cite{agt} II.3.4),
en tenant compte de \ref{cad4}(i), il est alors loisible de remplacer $f$ par $f^{(p^n)}$ \eqref{cad3e}.
On peut donc supposer $n=0$. 

Pour tout entier $j\geq 0$, posons 
\begin{eqnarray}
\Xi_{p^j}&=&\Hom(\Delta_{p^\infty},\mu_{p^j}(\co_\oK)),\\
\Xi_{p^\infty}&=&\Hom(\Delta_{p^\infty},\mu_{p^\infty}(\co_\oK)).
\end{eqnarray}
On identifie $\Xi_{p^j}$ à un sous-groupe de $\Xi_{p^\infty}$. 
Pour tout $\nu\in \Xi_{p^\infty}$, on note $\co_\oK(\nu)$ le $(\co_\oK)$-$\Delta_{p^\infty}$-module topologique $\co_\oK$, 
muni de la topologie $p$-adique et de l'action de $\Delta_{p^\infty}$ définie par $\nu$.
D'après (\cite{agt} II.8.9), il existe une décomposition canonique de $R_{p^\infty}$ en somme directe de $R_1$-modules de présentation finie,
stables sous l'action de $\Delta_{p^\infty}$,
\begin{equation}
R_{p^\infty}=\bigoplus_{\nu\in \Xi_{p^\infty}}R_{p^\infty}^{(\nu)}\otimes_{\co_\oK}\co_\oK(\nu),
\end{equation}
où $\Delta_{p^\infty}$ agit trivialement sur $R_{p^\infty}^{(\nu)}$. De plus, pour tout $j\geq 0$, on a 
\begin{equation}
R_{p^j}=\bigoplus_{\nu\in \Xi_{p^j}}R_{p^\infty}^{(\nu)}. 
\end{equation}
En particulier, pour tout $\nu\in \Xi_{p^\infty}$, $R_{p^\infty}^{(\nu)}$ est un $R_1$-module de présentation finie. 

D'après ce qui précède, pour tout entier $i\geq 0$, il existe une décomposition canonique 
\begin{equation}
\rH^i(\Delta_{p^\infty},M\otimes_{R_1}R_{p^\infty})=\bigoplus_{\nu\in \Xi_{p^\infty}}\rH^i(\Delta_{p^\infty},M\otimes_{R_1}R_{p^\infty}^{(\nu)}\otimes_{\co_\oK}\co_\oK(\nu)).
\end{equation}
Pour tout $\gamma\in \fm_\oK$, il existe un entier $j\geq 1$ tel que pour tout $\nu\in \Xi_{p^\infty}-\Xi_{p^j}$, le $R_1$-module
$\rH^i(\Delta_{p^\infty},M\otimes_{R_1}R_{p^\infty}^{(\nu)}\otimes_{\co_\oK}\co_\oK(\nu))$ soit annulé par $\gamma$, en vertu de \ref{cg36}(ii). 
On en déduit, d'après \ref{cg36}(i), que le $R_1$-module $\rH^i(\Delta_{p^\infty},M\otimes_{R_1}R_{p^\infty})$ est de $\alpha$-présentation finie. 

(iii) Soit $n'$ le plus grand diviseur premier à $p$ de $n$.   Comme $R_{n'}/p^mR_{n'}$ est une $(R_1/p^mR_1)$-algèbre cohérente, 
tout $(R_{n'}/p^mR_{n'})$-module de présentation $\alpha$-finie est $\alpha$-cohérent d'après \ref{afini6}, 
et est aussi un $(R_1/p^mR_1)$-module $\alpha$-cohérent. Le foncteur $\Gamma(\Delta_{n'}, -)$ est exact,
et il transforme les $(R_{n'}/p^mR_{n'})$-modules $\alpha$-cohérents en $(R_1/p^mR_1)$-modules $\alpha$-cohérents \eqref{finita18}. 
Il est alors loisible de remplacer $f$ par $f^{(n')}$ \eqref{cad3e}, d'après \ref{cad4}(i).
On peut donc supposer $n$ une puissance de $p$. Comme $\Sigma_0$ est un groupe profini d'ordre premier à $p$ \eqref{cg6l}, pour tout $i\geq 0$,
on a un isomorphisme canonique 
\begin{equation}
\rH^i(\Delta_{p^\infty},(M\otimes_{R_n}R_\infty)^{\Sigma_0})\stackrel{\sim}{\rightarrow} \rH^i(\Delta_\infty,M\otimes_{R_n}R_\infty).
\end{equation} 
Par ailleurs, considérant une présentation finie de $M$ sur $R_n$, on déduit de (\cite{agt} II.6.13) que le morphisme canonique
\begin{equation}
M\otimes_{R_n}R_{p^\infty}\rightarrow (M\otimes_{R_n}R_\infty)^{\Sigma_0}
\end{equation}
est un isomorphisme. La proposition résulte donc de (ii).

\begin{prop}\label{cg37}
Soient $m$ un entier $\geq 1$, $\mL$ un $(\co_\oK/p^m\co_\oK)$-module libre de type fini, muni d'une action linéaire continue de $\Delta$. 
Alors, 
\begin{itemize}
\item[{\rm (i)}] Le $(R_\infty/p^mR_\infty)$-module $(\mL\otimes_{\co_\oK}\oR)^\Sigma$ est $\alpha$-projectif de type $\alpha$-fini \eqref{aet1}, 
et le morphisme canonique 
\begin{equation}\label{cg37a}
(\mL\otimes_{\co_\oK}\oR)^\Sigma\otimes_{R_\infty}\oR\rightarrow \mL\otimes_{\co_\oK}\oR
\end{equation}
est un $\alpha$-isomorphisme. 
\item[{\rm (ii)}]  Pour tout $\gamma\in \fm_\oK$, il existe un entier $n\geq 1$, un $(R_n/p^mR_n)$-module de présentation finie $M_n$, 
muni d'une action semi-linéaire de $\Delta_n$, et un morphisme $R_{\infty}$-linéaire et $\Delta_{\infty}$-équivariant 
\begin{equation}\label{cg37b}
M_n\otimes_{R_n}R_{\infty}\rightarrow (\mL\otimes_{\co_\oK}\oR)^\Sigma,
\end{equation}
dont le noyau et le conoyau sont annulés par $\gamma$. 
\item[{\rm (iii)}] Pour tout entier $i\geq 0$, le $R_1$-module $\rH^i(\Delta,\mL\otimes_{\co_\oK}\oR)$ est de présentation $\alpha$-finie, 
et est $\alpha$-nul pour tout $i\geq d+1$.
\end{itemize}
\end{prop}

(i) Soient $N$ une extension galoisienne finie de $F_\infty$ contenue dans $\oF$ \eqref{cg6m}, 
$D$ la clôture intégrale de $R_\infty$ dans $N$, $G=\Gal(N/F_\infty)$, 
$\Tr_G$ l'endomorphisme $R_\infty$-linéaire de $D$ (ou de $D/p^m D$) induit par $\sum_{\sigma\in G}\sigma$.
Comme on a $D=\oR\cap N$ et $R_\infty=D\cap F_\infty$ d'après \ref{cad4}(iii), 
les homomorphismes $R_\infty/p^m R_\infty\rightarrow D/p^m D\rightarrow \oR/p^m \oR$ sont injectifs. 
Supposons que $D$ soit une extension $\alpha$-étale de $R_\infty$.  
Alors, $D$ est un $\alpha$-$G$-torseur sur $R_\infty$ en vertu de \ref{aet8}. Par suite, le quotient
\[
\frac{(D/p^m D)^G}{\Tr_G(D/p^m D)}
\]
est $\alpha$-nul en vertu de \ref{aet7}. Comme $\Tr_G(D)\subset R_\infty$, l'homomorphisme 
$R_\infty/p^m R_\infty\rightarrow (D/p^m D)^G$ est un $\alpha$-isomorphisme. 
Compte tenu de \ref{cg35}, on en déduit, par passage à la limite inductive, que l'homomorphisme 
$R_\infty/p^mR_\infty\rightarrow (\oR/p^m\oR)^\Sigma$ est un $\alpha$-isomorphisme.
On démontre par le même argument, en tenant compte de (\cite{agt} V.7.8), que l'homomorphisme
\begin{equation}\label{cg37c}
D/p^mD\rightarrow (\oR/p^m\oR)^{\Gal(\oF/N)}
\end{equation}
est un $\alpha$-isomorphisme. 

Choisissons $N$ de sorte que l'action de $\Sigma$ sur $\mL$ se factorise à travers $G$ \eqref{cg35}.
On en déduit que le morphisme canonique
\begin{equation}\label{cg37d}
(\mL\otimes_{\co_\oK} D)^G\rightarrow (\mL\otimes_{\co_\oK}\oR)^\Sigma
\end{equation}
est un $\alpha$-isomorphisme. Par ailleurs, il résulte de \ref{aet5} que le morphisme canonique 
\begin{equation}\label{cg37e}
(\mL\otimes_{\co_\oK} D)^G\otimes_{R_\infty}D\rightarrow \mL\otimes_{\co_\oK} D
\end{equation}
est un $\alpha$-isomorphisme~; il en est donc de même de \eqref{cg37a}. Comme $D$ est $\alpha$-fidèlement plat sur $R_\infty$ 
d'après (\cite{agt} V.6.7), 
on en déduit que $(\mL\otimes_{\co_\oK}\oR)^\Sigma$ est un $(R_\infty/p^mR_\infty)$-module $\alpha$-projectif de type $\alpha$-fini 
en vertu de (\cite{agt} V.8.7). 

(ii) Il suffit de montrer que pour tout $\gamma\in \fm_\oK$, il existe une suite de $R_\infty$-représentations continues de $\Delta_\infty$  
\begin{equation}
(R_\infty/p^mR_\infty)^a\stackrel{v}{\rightarrow}(R_\infty/p^mR_\infty)^b\stackrel{u}{\rightarrow} (\mL\otimes_{\co_\oK}\oR)^\Sigma\rightarrow 0,
\end{equation}
où $a$ et $b$ sont deux entiers $\geq 0$, telle que la suite des $R_\infty$-modules sous-jacents soit $\gamma$-exacte \eqref{finita3}. 
En effet, il existerait alors un entier $n\geq 1$ tel que $v$ se déduise par extension des scalaires d'un morphisme de $R_n$-représentations 
de $\Delta_n$, $w\colon (R_n/p^mR_n)^a\rightarrow (R_n/p^mR_n)^b$, et on prendrait pour $M_n$ le conoyau de $w$.

D'après (i), il existe deux entiers $r,q\geq 1$ et un morphisme $R_\infty$-linéaire 
\begin{equation}\label{cg37f}
R_\infty^r\rightarrow (\mL\otimes_{\co_\oK}\oR)^\Sigma
\end{equation}
dont le conoyau est annulé par $\gamma^{1/12}$ et tels que les images de la base canonique de $R_\infty^r$ soient invariants par le noyau 
de l'homomorphisme $\Delta_\infty\rightarrow \Delta_q$. Notons $R_q\langle \Delta_q\rangle$ l'anneau
non-commutatif de groupe sous-jacent le $R_q$-module libre de base $(e_\sigma)_{\sigma\in \Delta_q}$,
et dont la multiplication est donnée par $(be_\sigma)(b'e_{\sigma'})=b\sigma(b')e_{\sigma\sigma'}$. 
La catégorie des $R_q$-modules munis d'une action semi-linéaire de $\Delta_q$ 
est équivalente à la catégorie des $R_q\langle \Delta_q\rangle$-modules à gauche. 
Le morphisme \eqref{cg37f} induit donc un morphisme $R_\infty$-linéaire et $\Delta_\infty$-équivariant
\begin{equation}\label{cg37g}
u\colon R_q\langle \Delta_q\rangle^r\otimes_{R_q}(R_\infty/p^mR_\infty)\rightarrow (\mL\otimes_{\co_\oK}\oR)^\Sigma.
\end{equation}
D'après \ref{finita12}(iii), le noyau de $u$ est type $\gamma$-fini. Recommençant la même construction pour $\ker(u)$, on obtient 
l'assertion recherchée.

(iii) En effet, pour tout $i\geq 0$, le morphisme canonique 
\begin{equation}
\rH^i(\Delta_\infty,(\mL\otimes_{\co_\oK} \oR)^\Sigma)\rightarrow \rH^i(\Delta,\mL\otimes_{\co_\oK} \oR)
\end{equation} 
est un $\alpha$-isomorphisme d'après (\cite{agt} II.6.20). La proposition résulte alors de (ii) et \ref{cg68}(iii).

\subsection{}\label{cg10}
Soit $n$ un entier $\geq 1$. D'après \ref{cad4}(vi), le morphisme $\oX^{(n)\circ}\rightarrow \oX^\circ$ est étale fini
\eqref{cg6c}, de sorte que $(\oX^{(n)\circ}\rightarrow X)$ est un objet de $E$ \eqref{TFA6a}. 
On note $E^{(n)}\rightarrow \Et_{/X}$ (resp. $\tE^{(n)}$) le site fibré (resp. topos) de Faltings associé 
au morphisme $\oX^{(n)\circ}\rightarrow X$ \eqref{tf1}. Tout objet de $E^{(n)}$
est naturellement un objet de $E$. On définit ainsi un foncteur 
\begin{equation}\label{cg10a}
\Phi_n\colon E^{(n)}\rightarrow E.
\end{equation}
On vérifie aussitôt que $\Phi_n$ se factorise à travers une équivalence de catégories 
\begin{equation}\label{cg10b}
E^{(n)}\stackrel{\sim}{\rightarrow} E_{/(\oX^{(n)\circ}\rightarrow X)}.
\end{equation}
Il résulte alors de (\cite{agt} VI.5.38) que la topologie co-évanescente de $E^{(n)}$ 
est induite par celle de $E$ au moyen du foncteur $\Phi_n$. 
Par suite, $\Phi_n$ est continu et cocontinu (\cite{sga4} III 5.2). 
Il définit donc une suite de trois foncteurs adjoints~:
\begin{equation}\label{cg10c}
(\Phi_{n})_!\colon \tE^{(n)}\rightarrow \tE, \ \ \ \Phi^*_n\colon \tE\rightarrow \tE^{(n)}, 
\ \ \ \Phi_{n*}\colon \tE^{(n)}\rightarrow \tE,
\end{equation}
dans le sens que pour deux foncteurs consécutifs de la suite, celui de droite est
adjoint à droite de l'autre. D'après (\cite{sga4} III 5.4), le foncteur $(\Phi_{n})_!$ se factorise à travers 
une équivalence de catégories 
\begin{equation}\label{cg10d}
\tE^{(n)}\stackrel{\sim}{\rightarrow} \tE_{/(\oX^{(n)\circ}\rightarrow X)^\tta},
\end{equation}
où $(\oX^{(n)\circ}\rightarrow X)^\tta$ est l'objet de $\tE$ associé à $(\oX^{(n)\circ}\rightarrow X)$.
Le couple de foncteurs $(\Phi^*_n,\Phi_{n*})$ définit le morphisme de localisation de $\tE$ en 
$(\oX^{(n)\circ}\rightarrow X)^\tta$ que l'on note aussi
\begin{equation}\label{cg10e}
\Phi_n\colon \tE^{(n)}\rightarrow \tE.
\end{equation}
Le foncteur \eqref{cg10a} est un adjoint à gauche du foncteur 
\begin{equation}\label{cg10f}
\Phi_n^+\colon E\rightarrow E^{(n)}, \ \ \ (V\rightarrow U)\mapsto (V\times_{\oX^\circ}\oX^{(n)\circ}\rightarrow U).
\end{equation}
D'après (\cite{sga4} III 2.5), le morphisme \eqref{cg10e} s'identifie donc au morphisme  
défini par fonctorialité (\cite{agt} VI.10.12) à partir du diagramme 
\begin{equation}\label{cg10g}
\xymatrix{
{\oX^{(n)\circ}}\ar[r]\ar[d]&X\ar@{=}[d]\\
{\oX^\circ}\ar[r]&X}
\end{equation}

\subsection{}\label{cg11}
Soit $n$ un entier $\geq 1$. 
Pour tout objet $(V\rightarrow U)$ de $E$, on pose $\oU^{(n)}=U\times_X\oX^{(n)}$ et on note 
$\oU^{V(n)}$ la fermeture intégrale de $\oU^{(n)}$ dans $V\times_{\oU}\oU^{(n)}$. 
Pour tout morphisme $(V'\rightarrow U')\rightarrow (V\rightarrow U)$ de $E$, on a un morphisme canonique 
$\oU'^{V'(n)}\rightarrow \oU^{V(n)}$ qui s'insère dans un diagramme commutatif 
\begin{equation}\label{cg11a}
\xymatrix{
{V'\times_{\oU'}\oU'^{(n)}}\ar[r]\ar[d]&{\oU'^{V'(n)}}\ar[d]\ar[r]&{\oU'^{(n)}}\ar[r]\ar[d]&U'\ar[d]\\
{V\times_{\oU}\oU^{(n)}}\ar[r]&{\oU^{V(n)}}\ar[r]&{\oU^{(n)}}\ar[r]&U}
\end{equation} 
Comme l'injection canonique $X^{(n)\circ}\rightarrow X^{(n)}$ est schématiquement dominante d'après \ref{cad4}(i) et (\cite{agt} III.4.2(iv)), 
$\oU^{(n)\circ}$ est un ouvert schématiquement dominant de $\oU^{(n)}$ (\cite{ega4} 11.10.5).
On a donc $\oU^{\oU^\circ(n)}=\oU^{(n)}$ en vertu de \ref{cad4}(iii).

On désigne par $\ocB^{(n)}$ le préfaisceau sur $E$ défini pour tout $(V\rightarrow U)\in \ob(E)$, par 
\begin{equation}\label{cg11b}
\ocB^{(n)}(V\rightarrow U)=\Gamma(\oU^{V(n)},\co_{\oU^{V(n)}}).
\end{equation}
Comme $\oX^{(n)}$ est fini sur $\oX$, on a 
\begin{equation}\label{cg11c}
\ocB^{(n)}(V\rightarrow U)=\ocB(V\times_{\oX^\circ}\oX^{(n)\circ}\rightarrow U).
\end{equation}
Par suite, $\ocB^{(n)}$ est un faisceau sur $E$, et 
on a un isomorphisme canonique d'anneaux sur $\tE$
\begin{equation}\label{cg11d}
\ocB^{(n)}\stackrel{\sim}{\rightarrow} \Phi_{n*}(\Phi_n^*\ocB),
\end{equation}
où $\Phi_n$ est le morphisme de localisation \eqref{cg10e}. Pour tout $U\in \ob(\Et_{/X})$, on pose 
\begin{equation}\label{cg11e}
\ocB^{(n)}_{U}=\ocB^{(n)}\circ \iota_{U!},
\end{equation} 
où $\iota_{U!}$ est le foncteur \eqref{TFA6F}. 

Les $(\ocB^{(n)})_{n\geq 1}$ forment naturellement un système inductif d'anneaux de $\tE$, indexé par l'ensemble 
$\mZ_{\geq 1}$ ordonné par la relation de divisibilité. On pose 
\begin{equation}\label{cg11f}
\ocB^{(\infty)}=\underset{\underset{n\geq 1}{\longrightarrow}}{\lim}\ \ocB^{(n)}.
\end{equation}
Pour tout $U\in \ob(\Et_{/X})$, on pose 
\begin{equation}\label{cg11g}
\ocB^{(\infty)}_{U}=\ocB^{(\infty)}\circ \iota_{U!}.
\end{equation} 

D'après \ref{tf4}, le topos $\tE$ est cohérent. Supposons $U$ cohérent. 
Le faisceau associé à tout objet $(V\rightarrow U)$ de $E$ est alors cohérent d'après {\em loc. cit.}
Par suite, le morphisme canonique
\begin{equation}\label{cg11h}
\underset{\underset{n\geq 1}{\longrightarrow}}{\lim}\ \ocB^{(n)}_U(V)\rightarrow \ocB^{(\infty)}_U(V)
\end{equation}
est un isomorphisme en vertu de (\cite{sga4} VI 5.2). Comme le schéma $\oU^\circ$ est cohérent, on en déduit 
par (\cite{agt} VI.9.12) que le morphisme canonique
\begin{equation}\label{cg11i}
\underset{\underset{n\geq 1}{\longrightarrow}}{\lim}\ \ocB^{(n)}_U\rightarrow \ocB^{(\infty)}_U
\end{equation}
est un isomorphisme de $\oU^\circ_\fet$. 

\begin{lem}\label{cg12}
Soit $(V\rightarrow U)$ un objet de $E$. Alors,
\begin{itemize}
\item[{\rm (i)}] Pour tout entier $n\geq 1$, 
le schéma $\oU^{V(n)}$ est normal et localement irréductible et le morphisme canonique 
$V\times_{\oU}\oU^{(n)}\rightarrow \oU^{V(n)}$ est une immersion ouverte schématiquement dominante. 
\item[{\rm (ii)}] Si $U$ est quasi-compact, le nombre de composantes irréductibles de $\oU^{V(n)}$ lorsque $n$
décrit les entiers $\geq 1$, est borné. 
\end{itemize}
\end{lem}

(i) En effet, $V$ étant entier sur $\oU^\circ$, 
le morphisme canonique $V\times_{\oU}\oU^{(n)}\rightarrow \oU^\circ \times_{\oU}\oU^{V(n)}$ est un isomorphisme.
Le morphisme canonique $V\times_{\oU}\oU^{(n)}\rightarrow \oU^{V(n)}$ est donc une immersion ouverte 
schématiquement dominante.
Le schéma $\oU^{(n)}$ est normal et localement irréductible d'après \ref{cad4}(iii) et (\cite{agt} III.3.3). 
Soit $P$ un ouvert de $\oU^{(n)}$ n'ayant qu'un nombre fini de composantes irréductibles. 
Alors, $P\times_{\oU^{(n)}}\oU^{V(n)}$ est la somme finie des fermetures intégrales de $P$ dans 
les points génériques de $V\times_{\oU}\oU^{(n)}$ qui sont au-dessus de $P$, 
dont chacune est un schéma intègre et normal en vertu de (\cite{ega2} 6.3.7). 
Le schéma $\oU^{V(n)}$ est donc normal et localement irréductible.

(ii) Compte tenu de (i), comme le morphisme canonique $V\rightarrow \oU^\circ$ est étale et fini, il suffit de montrer que 
le nombre de composantes irréductibles de $\oU^{(n)\circ}$ (pour $n\geq 1$) est borné. 
Quitte à remplacer $X$ par des ouverts affines qui couvrent $U$, on peut se borner au cas où $U=X$. 
D'après \eqref{cad3e} et avec les notations de \ref{cad3}, on a un diagramme cartésien de $S$-morphismes
\begin{equation}
\xymatrix{
{X^{(n)}}\ar[r]\ar[d]&{\Spec(\co_{K_n}[P^{(n)}]/(\pi_n-e^{\lambda^{(n)}}))}\ar[d]\\
{X}\ar[r]&{\Spec(\co_K[P]/(\pi-e^\lambda))}}
\end{equation}
où les flèches horizontales sont étales. On en déduit 
un diagramme cartésien de $\oS$-morphismes
\begin{equation}
\xymatrix{
{\oX^{(n)}}\ar[r]\ar[d]&{\Spec(\co_{\oK}[P^{(n)}]/(\pi_n-e^{\lambda^{(n)}}))}\ar[d]\\
{\oX}\ar[r]&{\Spec(\co_{\oK}[P]/(\pi-e^\lambda))}}
\end{equation}
Compte tenu de \ref{cad4}(v), on en déduit un diagramme cartésien de $\oK$-morphismes
\begin{equation}
\xymatrix{
{\oX^{(n)\circ}}\ar[r]\ar[d]&{\Spec(\oK[(P^{(n)})^\gp]/(\pi_n-e^{\lambda^{(n)}}))}\ar[d]\\
{\oX^\circ}\ar[r]&{\Spec(\oK[P^\gp]/(\pi-e^\lambda))}}
\end{equation}
La flèche verticale de gauche est étale d'après \ref{cad4}(vi), ainsi que les flèches horizontales. 
Par suite, les points génériques de $\oX^{(n)\circ}$ s'envoient sur ceux de 
$\Spec(\oK[(P^{(n)})^\gp]/(\pi_n-e^{\lambda^{(n)}}))$ et sur ceux de $\oX^{\circ}$, et donc aussi sur ceux de 
$\Spec(\oK[P^\gp]/(\pi-e^\lambda))$. 

Par ailleurs, comme le $\mZ$-module $P^\gp$ est libre de type fini, pour tout entier $n\geq 1$, on a 
un isomorphisme de $\oK$-algèbres
\begin{equation}
\oK[P^\gp]/(1-e^\lambda)\stackrel{\sim}{\rightarrow} \oK[(P^{(n)})^\gp]/(\pi_n-e^{\lambda^{(n)}}).
\end{equation}
Le nombre de composantes irréductibles de $\Spec(\oK[(P^{(n)})^\gp]/(\pi_n-e^{\lambda^{(n)}}))$ 
(pour $n\geq 1$) est donc constant~; d'où la proposition.

\subsection{}\label{cg25}
On désigne par $\cC$ la sous-catégorie pleine de $\Et_{/X}$ formée des schémas affines, que l'on munit 
de la topologie induite par celle de $\Et_{/X}$. On voit aussitôt que $\cC$ est 
une famille $\mU$-petite, topologiquement génératrice du site $\Et_{/X}$ et est stable par produits fibrés. 
On désigne par 
\begin{equation}\label{cg25a}
\pi_\cC\colon E_\cC\rightarrow \cC
\end{equation} 
le site fibré déduit de $\pi$ \eqref{TFA6a}
par changement de base par le foncteur d'injection canonique $\cC\rightarrow \Et_{/X}$. 
On munit $E_\cC$ de 
la topologie co-évanescente définie par $\pi_\cC$ et on note $\tE_\cC$ le topos des faisceaux de $\mU$-ensembles 
sur $E_\cC$. D'après (\cite{agt} VI.5.21 et VI.5.22), 
la topologie de $E_\cC$ est induite par celle de $E$ au moyen du foncteur de projection canonique 
$E_\cC\rightarrow E$, et celui-ci induit par restriction une équivalence de catégories 
\begin{equation}\label{cg25b}
\tE\stackrel{\sim}{\rightarrow}\tE_\cC.
\end{equation}

\begin{lem}\label{cg16}
Soit $(V\rightarrow U)$ un objet de $E_\cC$. Alors,
\begin{itemize}
\item[{\rm (i)}] Le schéma $\Spec(\ocB^{(\infty)}_U(V))$ est normal et localement irréductible. 
Il est la fermeture intégrale de $\oU$ dans $\Spec(\ocB_U^{(\infty)}(V))\times_XX^\circ$.
\item[{\rm (ii)}] Le diagramme de morphismes canoniques
\begin{equation}\label{cg16a}
\xymatrix{
{\Spec(\ocB_U^{(\infty)}(V))\times_XX^\circ}\ar[r]\ar[d]&V\ar[d]\\
{\Spec(\ocB_U^{(\infty)}(\oU^\circ))\times_XX^\circ}\ar[r]&{\oU^\circ}}
\end{equation}
est cartésien. 
\end{itemize}
\end{lem}

(i) Pour tout entier $n\geq 1$, comme $X^{(n)}$ est fini sur $X$, $\oU^{V(n)}$ est la fermeture intégrale de 
$\oU$ dans $V\times_\oU\oU^{(n)}$.  D'après \ref{cg12}(i), $\oU^{V(n)}$ est normal et localement irréductible.
Pour tout entier $m\geq 1$, le morphisme canonique $X^{(mn)\circ}\rightarrow X^{(n)\circ}$ étant étale, fini et surjectif, 
le morphisme canonique $\oU^{V(mn)}\rightarrow \oU^{V(n)}$ est entier et surjectif,  
et toute composante irréductible de $\oU^{V(mn)}$ domine une composante irréductible de 
$\oU^{V(n)}$. De plus, en vertu de \ref{cg12}(ii), il existe un entier $n\geq 1$ tel que pour toute composante 
irréductible $P$ de $\oU^{V(n)}$ et tout entier $m\geq 1$, $P\times _{\oU^{V(n)}}\oU^{V(mn)}$ soit irréductible. 
La proposition résulte alors de (\cite{ega1n} 0.6.1.6 et 0.6.5.12(ii)).

(ii) En effet, $\Spec(\ocB_U^{(\infty)}(V))$ est la limite projective des schémas $(\oU^{V(n)})_{n\geq 1}$ (\cite{ega4} 8.2.3). 
Par suite,   $\Spec(\ocB_U^{(\infty)}(V))\times_XX^\circ$ est la limite projective des schémas 
$(V\times_\oU\oU^{(n)})_{n\geq 1}$ (\cite{ega4} 8.2.5); d'où la proposition.

\begin{prop}\label{cg17}
Soient $(V\rightarrow U)$ un objet de $E_\cC$ tel que $V$ soit connexe, $\ov$ un point géométrique de $V$.  
Alors, il existe un revêtement universel $(V_i)_{i\in I}$ de $V$ en $\ov$ tel que pour tout $i\in I$, $V_i$ soit un revêtement étale galoisien de $V$
et que $\ocB_U^{(\infty)}(V_i)$ soit une extension $\alpha$-étale de $\ocB_U^{(\infty)}(V)$ \eqref{aet2}. 
\end{prop}

En effet, on peut se borner au cas où $U=X$.  Soient $Y$ un revêtement étale fini de $X^\circ$, $W=V\times_{X^\circ}Y$.
Considérons le diagramme de morphismes canoniques 
\begin{equation}\label{cg17a}
\xymatrix{
{\Spec(\ocB^{(\infty)}_X(W))\times_XX^\circ}\ar[rd]\ar[ddd]\ar[rrr]
\ar@{}[rddd]|*+[o][F-]{1}&&&{\Spec(A_\infty)\times_XY}\ar[ddd]\ar[ld]\\
&W\ar[r]\ar[d]&Y\ar[d]&\\
&V\ar[r]&X^\circ&\\
{\Spec(\ocB^{(\infty)}_X(V))\times_XX^\circ}\ar[rrr]\ar[ru]&&&{\Spec(A_\infty)\times_XX^\circ}\ar[lu]}
\end{equation}
où $A_\infty$ est l'algèbre définie dans \eqref{cg21a}.
Notons $\Spec(B)$ la fermeture intégrale de $\Spec(A_\infty)$ dans $\Spec(A_\infty)\times_XY$. En vertu de \ref{cg22}, 
$B$ est une extension $\alpha$-étale de $A_\infty$. 
Par suite, $B\otimes_{A_\infty}\ocB^{(\infty)}_X(V)$ est une extension 
$\alpha$-étale de $\ocB^{(\infty)}_X(V)$ (\cite{agt} V.7.4(3)).

Le carré $\xymatrix{\ar@{}|*+[o][F-]{1}}$ du diagramme \eqref{cg17a} étant cartésien en vertu de 
\ref{cg16}(ii), il en est de même du grand carré extérieur. On en déduit un diagramme cartésien
\begin{equation}
\xymatrix{
{\Spec(\ocB^{(\infty)}_X(W))\times_XX^\circ}\ar[r]\ar[d]&{\Spec(B\otimes_{A_\infty}\ocB^{(\infty)}_X(V)[\frac 1 p])}\ar[d]\\
{\Spec(\ocB^{(\infty)}_X(V))\times_XX^\circ}\ar[r]&{\Spec(\ocB^{(\infty)}_X(V)[\frac 1 p])}}
\end{equation}
Comme $B[\frac 1 p]$ est une extension finie étale de $A_\infty[\frac 1 p]$ (\cite{agt} V.7.3), 
on en déduit en vertu de \ref{cg16}(i) que $\ocB^{(\infty)}_X(W)$ est la fermeture intégrale de 
$\ocB^{(\infty)}_X(V)$ dans $B\otimes_{A_\infty}\ocB^{(\infty)}_X(V)[\frac 1 p]$.
Par suite, $\ocB^{(\infty)}_X(W)$ est une extension $\alpha$-étale de $\ocB^{(\infty)}_X(V)$ d'après \ref{aet3}.

Comme $\oX$ est normal et localement irréductible (\cite{agt} III.4.2(iii)), il en est de même de $W$ 
(\cite{agt} III.3.3). Par suite, $W$ est la réunion des schémas induits sur ses composantes irréductibles.
Pour toute composante irréductible $W'$ de $W$, 
$\ocB^{(\infty)}_X(W')$ est une extension $\alpha$-étale de $\ocB^{(\infty)}_X(V)$ d'après (\cite{agt} V.7.4).
On conclut la preuve en observant que l'homomorphisme $\pi_1(V,\ov)\rightarrow \pi_1(X^\circ,\ov)$ est injectif puisque $U=X$.

\subsection{}\label{cg130}
Soient $U$ un objet de $\Et_{/X}$, $n$ un entier $\geq 1$. 
Avec les notations de \ref{cg6}, on pose 
\begin{equation}\label{cg130a}
\oU^{\star\circ}=U\times_X\oX^{\star\circ} \ \ \ {\rm et}\ \ \ \oU^{(n)\star\circ}=U\times_X\oX^{(n)\star\circ},
\end{equation}
et on note $\pi_U^{(n)\star}\colon \oU^{(n)\star\circ}\rightarrow \oU^{\star\circ}$ le morphisme canonique.
On rappelle que $\oX^{(n)\star\circ}$ est un revêtement étale fini et galoisien de $\oX^{\star\circ}$ 
de groupe $\Delta_n$. On définit le foncteur 
\begin{equation}\label{cg130b}
\oU^{\circ}_\fet\rightarrow \oU^{\star \circ}_\fet, \ \ \ F\mapsto F^{(n)\star}=(\pi_{U}^{(n)\star})_*(F|\oU^{(n)\star\circ}).
\end{equation}
Celui-ci est exact puisque le foncteur $(\pi_U^{(n)\star})_*$ admet un adjoint à gauche et est donc exact.

Soit $F$ un objet de $\oU^\circ_\fet$. D'après (\cite{agt} VI.9.4), pour tout $V\in \ob(\Et_{\rf/\oU^{\star \circ}})$, on a 
\begin{equation}\label{cg130c}
F^{(n)\star}(V)=F(V\times_{\oX^{\star\circ}}\oX^{(n)\star\circ}).
\end{equation}
Cet ensemble est donc naturellement muni d'une action de $\Delta_n$, et on a un isomorphisme canonique 
\begin{equation}\label{cg130d}
F(V) \stackrel{\sim}{\rightarrow}(F^{(n)\star}(V))^{\Delta_n}.
\end{equation}

Les faisceaux $(F^{(n)\star})_{n\geq 1}$  forment naturellement un système inductif de $\oU^{\star\circ}_\fet$, 
indexé par l'ensemble $\mZ_{\geq 1}$ ordonné par la relation de divisibilité. On pose 
\begin{equation}\label{cg130e}
F^{(\infty)\star}=\underset{\underset{n\geq 1}{\longrightarrow}}{\lim}\ F^{(n)\star}.
\end{equation}

On voit aussitôt que le morphisme canonique $\oX^{\star\circ}\rightarrow X$ est cohérent \eqref{cg2}.
Supposons $U$ cohérent, de sorte que $\oU^{\star\circ}$ est cohérent.
D'après (\cite{agt} VI.9.12) et (\cite{sga4} VI 5.2), pour tout $V\in \ob(\Et_{\rf/\oU^{\star \circ}})$, on a alors un isomorphisme
canonique
\begin{equation}\label{cg130f}
F^{(\infty)\star}(V)\stackrel{\sim}{\rightarrow}\underset{\underset{n\geq 1}{\longrightarrow}}{\lim}\ 
F^{(n)\star}(V).
\end{equation}
Par suite, l'ensemble discret $F^{(\infty)\star}(V)$ est naturellement muni d'une action continue de $\Delta_\infty$,
et on a un isomorphisme canonique
\begin{equation}\label{cg130g}
F(V) \stackrel{\sim}{\rightarrow}(F^{(\infty)\star}(V))^{\Delta_\infty}.
\end{equation}

\subsection{}\label{cg13}
Soient $U$ un objet de $\Et_{/X}$, $m, n$ deux entiers tels que $m\geq 0$ et $n\geq 1$. Appliquant le foncteur \eqref{cg130b} aux anneaux 
$\ocB_U$ \eqref{TFA2d} et $\ocB_{U,m}=\ocB_U/p^m\ocB_U$ \eqref{TFA8b} de $\oU^\circ_\fet$, 
on obtient deux anneaux $\ocB^{(n)\star}_U$ et $\ocB^{(n)\star}_{U,m}$ de $\oU^{\star\circ}_\fet$. Ce foncteur étant exact, on a 
un isomorphisme canonique 
\begin{equation}\label{cg13a}
\ocB^{(n)\star}_{U,m}\stackrel{\sim}{\rightarrow} \ocB^{(n)\star}_U/p^m\ocB^{(n)\star}_U.
\end{equation}
Compte tenu de \eqref{cg11c}, l'injection canonique $\oX^{(n)\star}\rightarrow \oX^{(n)}$ induit un homomorphisme \eqref{cg11e}
\begin{equation}\label{cg13b}
\ocB^{(n)}_U|\oU^{\star \circ}\rightarrow \ocB^{(n)\star}_U.
\end{equation}
On en déduit par passage à la limite inductive un homomorphisme \eqref{cg11g}
\begin{equation}\label{cg13c}
\ocB^{(\infty)}_U|\oU^{\star \circ}\rightarrow \ocB^{(\infty)\star}_U.
\end{equation}
Par ailleurs, les isomorphismes \eqref{cg13a} induisent un isomorphisme 
\begin{equation}\label{cg13d}
\ocB^{(\infty)\star}_{U,m}\stackrel{\sim}{\rightarrow} \ocB^{(\infty)\star}_U/p^m\ocB^{(\infty)\star}_U.
\end{equation}

\subsection{}\label{cg14}
On désigne par $E^\star_\cC$ la sous-catégorie pleine de $E_\cC$ \eqref{cg25} formée des objets $(V\rightarrow U)$ 
tels que le morphisme canonique $V\rightarrow \oU^{\circ}$ se factorise à travers  $\oU^{\star\circ}$ \eqref{cg130a}.
Le foncteur \eqref{cg25a} induit un foncteur fibrant
\begin{equation}\label{cg14a}
E^\star_\cC\rightarrow\cC.
\end{equation}

\begin{lem}\label{cg15}
Soit $(V\rightarrow U)$ un objet de $E^\star_\cC$. Alors, 
\begin{itemize}
\item[{\rm (i)}] Le morphisme canonique $\Spec(\ocB_U^{(\infty)\star}(V))\rightarrow \Spec(\ocB_U^{(\infty)}(V))$ 
\eqref{cg13b} est une immersion ouverte et fermée. 
\item[{\rm (ii)}] Le diagramme de morphismes canoniques
\begin{equation}
\xymatrix{
{\Spec(\ocB_U^{(\infty)\star}(V))}\ar[r]\ar[d]&{\Spec(\ocB_U^{(\infty)}(V))}\ar[d]\\
{\Spec(\ocB_U^{(\infty)\star}(\oU^\circ))}\ar[r]&{\Spec(\ocB_U^{(\infty)}(\oU^\circ))}}
\end{equation}
est cartésien. 
\end{itemize}
\end{lem}

(i) Pour tout entier $n\geq 1$, $\Spec(\ocB^{(n)\star}_U(V))$ est un ouvert et fermé 
de $\Spec(\ocB^{(n)}_U(V))=\oU^{V(n)}$. D'après \ref{cg12}, il existe $n\geq1$ tel que pour tout $m\geq1$, 
le diagramme canonique
\begin{equation}
\xymatrix{
{\oX^{(mn)\star}}\ar[r]\ar[d]&{\oX^{(mn)}}\ar[d]\\ 
{\oX^{(n)\star}}\ar[r]&{\oX^{(n)}}}
\end{equation}
est cartésien. Il en est donc de même du diagramme canonique
\begin{equation}
\xymatrix{
{\Spec(\ocB^{(mn)\star}_U(V))}\ar[r]\ar[d]&{\Spec(\ocB^{(mn)}_U(V))}\ar[d]\\ 
{\Spec(\ocB^{(n)\star}_U(V))}\ar[r]&{\Spec(\ocB^{(n)}_U(V))}}
\end{equation}
La proposition s'ensuit en vertu de (\cite{ega4} 8.3.12). 

(ii) En effet, pour tout entier $n\geq 1$, le diagramme 
\begin{equation}
\xymatrix{
{\Spec(\ocB^{(n)\star}_U(V))}\ar[r]\ar[d]&{\oU^{V(n)}}\ar[d]\\ 
{\oU^{(n)\star}}\ar[r]&{\oU^{(n)}}}
\end{equation}
est cartésien.

\begin{prop}\label{cg18}
Soient $(V\rightarrow U)$ un objet de $E^\star_\cC$ tel que $V$ soit connexe,  
$\ov$ un point géométrique de $V$.  
Alors, il existe un revêtement universel $(V_i)_{i\in I}$ de $V$ en $\ov$ tel que pour tout $i\in I$, $V_i$ soit un revêtement étale galoisien de $V$ et 
que $\ocB_U^{(\infty)\star}(V_i)$ soit une extension $\alpha$-étale de $\ocB_U^{(\infty)\star}(V)$. 
\end{prop}

Cela résulte de \ref{cg17}, \ref{cg15} et (\cite{agt} V.7.4(2) et V.7.8).

\begin{cor}\label{cg19}
Soient $(V\rightarrow U)$ un objet de $E^\star_\cC$, $\cF$ un  $(\ocB^{(\infty)\star}_U|V)$-module de $V_\fet$, $q$ un entier $\geq 1$.  
Alors, $\rH^q(V_\fet, \cF)$ est $\alpha$-nul. 
\end{cor}

Comme $\oX$ est normal et localement irréductible (\cite{agt} III.4.2(iii)), il en est de même de $V$ (\cite{agt} III.3.3).
Par ailleurs, $j\colon X^\circ\rightarrow X$ étant quasi-compact, $V$ est quasi-compact. On peut donc supposer
$V$ connexe. Soit $\ov$ un point géométrique de $V$.
D'après \ref{cg18}, il existe un revêtement universel $(V_i)_{i\in I}$ de $V$ en $\ov$ 
tel que pour tout $i\in I$, $V_i$ soit un revêtement étale galoisien de $V$ et 
que $\ocB_U^{(\infty)\star}(V_i)$ soit une extension $\alpha$-étale de $\ocB_U^{(\infty)\star}(V)$.
Pour tout $i\in I$, notons $G_i$ le groupe des $V$-automorphismes de $V_i$. Il résulte de \ref{cg16} et \ref{cg15} que
$\Spec(\ocB_U^{(\infty)\star}(V))$ est normal et localement irréductible,
de fermeture intégrale $\Spec(\ocB_U^{(\infty)\star}(V_i))$ 
dans $\Spec(\ocB_U^{(\infty)\star}(V_i))\times_XX^\circ$ et que le diagramme 
\begin{equation}\label{cg19a}
\xymatrix{
{\Spec(\ocB_U^{(\infty)\star}(V_i))\times_XX^\circ}\ar[r]\ar[d]&V_i\ar[d]\\
{\Spec(\ocB_U^{(\infty)\star}(V))\times_XX^\circ}\ar[r]&V}
\end{equation}
est cartésien. En vertu de \ref{aet8}, $\ocB_U^{(\infty)\star}(V_i)$ est donc un $\alpha$-$G_i$-torseur de $\ocB_U^{(\infty)\star}(V)$. 
Par suite, $\rH^q(G_i,\cF(V_i))$ est $\alpha$-nul d'après \ref{aet7}. 
La proposition s'ensuit par passage à la limite. En effet, notant
\begin{equation}
\nu_\ov\colon V_\fet\stackrel{\sim}{\rightarrow} \bB_{\pi_1(V,\ov)},
\ \ \ M\mapsto \underset{\underset{i\in I}{\longrightarrow}}{\lim}\ M(V_i)
\end{equation}
le foncteur fibre en $\ov$ \eqref{notconv11c}, on a des isomorphismes canoniques 
\begin{equation}
\rH^q(V_\fet,\cF)\stackrel{\sim}{\rightarrow}\rH^q(\pi_1(V,\ov),\nu_\ov(\cF))\stackrel{\sim}{\rightarrow}
\underset{\underset{i\in I}{\longrightarrow}}{\lim}\ \rH^q(G_i,\cF(V_i)).
\end{equation}

\begin{cor}\label{cg20}
Pour tout objet $(V\rightarrow U)$ de $E^\star_\cC$ et tout entier $m\geq 1$, 
l'homomorphisme canonique
\begin{equation}
\ocB^{(\infty)\star}_U(V)/p^m\ocB^{(\infty)\star}_U(V)\rightarrow \ocB^{(\infty)\star}_{U,m}(V)
\end{equation} 
est un $\alpha$-isomorphisme. 
\end{cor}

En effet, d'après \eqref{cg13d}, comme $\ocB^{(\infty)\star}_U$ est $\mZ_p$-plat, 
on a une suite exacte de groupes abéliens de $\oU^{\star \circ}_\fet$
\begin{equation}
0\longrightarrow \ocB^{(\infty)\star}_U\stackrel{\cdot p^m}{\longrightarrow} \ocB^{(\infty)\star}_U\longrightarrow \ocB^{(\infty)\star}_{U,m}\longrightarrow 0.
\end{equation}
La proposition résulte alors du fait que $\rH^1(V_\fet, \ocB^{(\infty)\star}_U|V)=0$ en vertu de \ref{cg19}.

\begin{cor}\label{cg23}
Pour tout morphisme cartésien $(V'\rightarrow U')\rightarrow (V\rightarrow U)$ de $E^\star_\cC$ \eqref{cg14a} 
et tout entier $m\geq 1$, l'homomorphisme canonique 
\begin{equation}
\ocB_{U,m}^{(\infty)\star}(V)\otimes_{\co_X(U)}\co_X(U')\rightarrow \ocB^{(\infty)\star}_{U',m}(V')
\end{equation}
est un $\alpha$-isomorphisme. 
\end{cor}
En effet, pour tout entier $n\geq 1$, l'homomorphisme canonique 
\begin{equation}
\ocB_U^{(n)\star}(V)\otimes_{\co_X(U)}\co_X(U')\rightarrow \ocB^{(n)\star}_{U'}(V')
\end{equation}
est un isomorphisme en vertu de (\cite{agt} III.8.11(iii)). 
On en déduit que l'homomorphisme canonique 
\begin{equation}
\ocB_U^{(\infty)\star}(V)\otimes_{\co_X(U)}\co_X(U')\rightarrow \ocB^{(\infty)\star}_{U'}(V')
\end{equation}
est un isomorphisme \eqref{cg130f}. La proposition est donc une conséquence de \ref{cg20}.

\begin{cor}\label{cg40}
Soient $m$ un entier $\geq 1$, $\mL$ un $(\co_\oK/p^m\co_\oK)$-module localement libre de type fini de $\oX^\circ_\fet$. 
Pour tout $X$-schéma étale $U$, on note $\mL_U$ l'image inverse de $\mL$ par le morphisme canonique $\oU^\circ\rightarrow \oX^\circ$,
on pose $\cL_U=\mL_U\otimes_{\co_\oK}\ocB_U$ et on désigne par $\cL_U^{(\infty)\star}$ 
le $(\ocB_U^{(\infty)\star})$-module de $\oU^\circ_\fet$ associé à $\cL_U$ \eqref{cg130e}. 
Alors, pour tout morphisme cartésien $(V'\rightarrow U')\rightarrow (V\rightarrow U)$ de $E^\star_\cC$, le morphisme canonique 
\begin{equation}
\cL_U^{(\infty)\star}(V)\otimes_{\co_X(U)}\co_X(U')\rightarrow \cL_{U'}^{(\infty)\star}(V')
\end{equation}
est un $\alpha$-isomorphisme. 
\end{cor} 

En effet, comme $\oX$ est normal et localement irréductible (\cite{agt} III.4.2(iii)), il en est de même de $V$ (\cite{agt} III.3.3).
Par ailleurs, $j\colon X^\circ\rightarrow X$ étant quasi-compact, $V$ est quasi-compact. 
On peut donc supposer $V$ connexe. 
Pour tout revêtement étale galoisien $V_1$ de $V$ de groupe $G$, posant $V'_1=V'\times_VV_1$, les morphismes canoniques
\begin{equation}
\cL^{(\infty)\star}_U(V)\rightarrow (\cL^{(\infty)\star}_U(V_1))^G\ \ \ {\rm et}\ \ \ \cL^{(\infty)\star}_{U'}(V')\rightarrow (\cL^{(\infty)\star}_{U'}(V'_1))^G
\end{equation} 
sont des isomorphismes. Compte tenu de (\cite{agt} II.3.15 et III.2.11), comme $\co_X(U')$ est plat sur $\co_X(U)$, 
on peut donc se réduire au cas où $\mL_U|V$ est constant de valeur un $(\co_\oK/p^m\co_\oK)$-module libre de type fini $M$.
On a un isomorphisme canonique \eqref{cg130f}
\begin{equation}
\cL_U^{(\infty)\star}(V)\stackrel{\sim}{\rightarrow}M\otimes_{\co_\oK/p^m\co_\oK}\ocB_{U,m}^{(\infty)\star}(V);
\end{equation}
et de même pour $(V'\rightarrow U')$. La proposition résulte alors de \ref{cg23}.

\begin{cor}\label{cg240}
Soient $m$ un entier $\geq 1$, $\mL$ un $(\co_\oK/p^m\co_\oK)$-module localement libre de type fini de $\oX^\circ_\fet$. 
Pour tout $X$-schéma étale $U$, on note $\mL_U$ l'image inverse de $\mL$ par le morphisme canonique $\oU^\circ\rightarrow \oX^\circ$,
et on pose $\cL_U=\mL_U\otimes_{\co_\oK}\ocB_U$. Alors, 
pour tout morphisme cartésien $(V'\rightarrow U')\rightarrow (V\rightarrow U)$ de $E_\cC$ \eqref{cg25a}, le morphisme canonique 
\begin{equation}
\cL_U(V)\otimes_{\co_X(U)}\co_X(U')\rightarrow \cL_{U'}(V')
\end{equation}
est un $\alpha$-isomorphisme. 
\end{cor}

En effet, on peut se borner au cas où $V$ est connexe (cf. la preuve de \ref{cg40}). 
Quitte à changer de point géométrique $\oy$ de $\oX^\circ$ \eqref{cg2}, 
on peut supposer de plus que $(V\rightarrow U)$ est un objet de $E^\star_\cC$. 
La proposition résulte alors de \ref{cg40} et \eqref{cg130g} puisque $\co_X(U')$ est plat sur $\co_X(U)$ (\cite{agt} II.3.15).

\begin{cor}\label{cg24}
Pour tout morphisme cartésien $(V'\rightarrow U')\rightarrow (V\rightarrow U)$ de $E_\cC$
et tout entier $m\geq 1$, l'homomorphisme canonique 
\begin{equation}
\ocB_{U,m}(V)\otimes_{\co_X(U)}\co_X(U')\rightarrow \ocB_{U',m}(V')
\end{equation}
est un $\alpha$-isomorphisme. 
\end{cor}

\begin{prop}\label{cg30}
Soient $(V\rightarrow U)$ un objet de $E^\star_\cC$, $\cF$ un  $\ocB_U$-module de $\oU^\circ_\fet$,
$\cF^{(\infty)\star}$ le $(\ocB_U^{(\infty)\star})$-module de $\oU^{\star\circ}_\fet$ associé à $\cF$ \eqref{cg130e}, 
$q$ un entier $\geq 0$. Alors, il existe un morphisme canonique
\begin{equation}\label{cg30a}
\rH^q(\Delta_\infty,\cF^{(\infty)\star}(V))\rightarrow \rH^q(V_\fet,\cF|V),
\end{equation}
qui est un $\alpha$-isomorphisme.
\end{prop}

Soit $n$ un entier $\geq 1$. Avec les notations de \ref{cg130}, on désigne par $\Gamma^{\Delta_n}$  le 
foncteur ``sous-faisceau des sections $\Delta_n$-invariantes'' sur la catégorie des $\mZ[\Delta_n]$-modules de $\oU^{\star\circ}_\fet$. 
Par descente, les foncteurs $(\pi_U^{(n)\star})^*$ et $\Gamma^{\Delta_n}\circ (\pi_U^{(n)\star})_*$ établissent des équivalences de 
catégories quasi-inverse l'une de l'autre entre la catégorie des groupes abéliens de $\oU^{\star\circ}_\fet$ et 
la catégorie des groupes abélien de $\oU^{(n)\star\circ}_\fet$ sur lesquels $\Delta_n$ agit d'une façon compatible 
avec son action sur $\oU^{(n)\star \circ}$. On a donc un isomorphisme canonique de groupes abéliens de $\oU^{\star\circ}_\fet$
\begin{equation}
\cF|\oU^{\star\circ}\stackrel{\sim}{\rightarrow}\Gamma^{\Delta_n}(\cF^{(n)\star}).
\end{equation}
Le foncteur $(\pi_U^{(n)\star})_*$ étant exact, appliquant (\cite{verdier} II §~2 prop.~3.1), on obtient un isomorphisme de groupes abéliens
\begin{equation}
\cF|\oU^{\star\circ}\stackrel{\sim}{\rightarrow}\rR\Gamma^{\Delta_n}(\cF^{(n)\star}).
\end{equation}
Le foncteur $\Gamma^{\Delta_n}$ transforme les $\mZ[\Delta_n]$-modules injectifs de $V_\fet$ en groupes abéliens injectifs
de $V_\fet$, et le foncteur $\Gamma(V_\fet,-)$ transforme les $\mZ[\Delta_n]$-modules injectifs de $V_\fet$
en $\mZ[\Delta_n]$-modules injectifs (\cite{sga4} V 0.2). D'autre part, on a  
\begin{equation}
\Gamma(V_\fet, \Gamma^{\Delta_n}(-))=\Gamma^{\Delta_n}(\Gamma(V_\fet,-)).
\end{equation}
Prenant les foncteurs dérivés des deux membres et appliquant de nouveau (\cite{verdier} II §~2 prop.~3.1), on en déduit une suite spectrale 
\begin{equation}
{^{(n)}E}_2^{a,b}=\rH^a(\Delta_n,\rH^b(V_\fet,\cF^{(n)\star}|V))
\Rightarrow \rH^{a+b}(V_\fet,\cF|V).
\end{equation}
Comme l'immersion $j\colon X^\circ\rightarrow X$ est quasi-compacte, le schéma $V$ est cohérent.
Par suite, en vertu de (\cite{sga4} VI 5.3) et (\cite{agt} VI.9.12), pour tout entier $b\geq 0$, on a un isomorphisme canonique
\begin{equation}
\underset{\underset{n\geq 1}{\longrightarrow}}{\lim}\ \rH^b(V_\fet,\cF^{(n)\star}|V) 
\stackrel{\sim}{\rightarrow} \rH^b(V_\fet,\cF^{(\infty)\star}|V).
\end{equation}
On en déduit par passage à la limite inductive une suite spectrale 
\begin{equation}\label{cg30b}
E_2^{a,b}=\rH^a(\Delta_\infty,\rH^b(V_\fet,\cF^{(\infty)\star}|V))\Rightarrow \rH^{a+b}(V_\fet,\cF|V).
\end{equation}
La proposition s'ensuit puisque $\rH^b(V_\fet, \cF^{(\infty)\star}|V)$ est $\alpha$-nul pour tout $b\geq 1$ en vertu de \ref{cg19}.

\begin{cor}\label{cg310}
Soient $m$ un entier $\geq 1$, $\mL$ un $(\co_\oK/p^m\co_\oK)$-module localement libre de type fini de $\oX^\circ_\fet$. 
Pour tout $X$-schéma étale $U$, on note $\mL_U$ l'image inverse de $\mL$ par le morphisme canonique $\oU^\circ\rightarrow \oX^\circ$,
et on pose $\cL_U=\mL_U\otimes_{\co_\oK}\ocB_U$. Alors, pour tout morphisme $U'\rightarrow U$ de $\cC$
et tout entier $q\geq 0$, le morphisme canonique
\begin{equation}\label{cg310a}
\rH^q(\oU^\circ_\fet,\cL_U)\otimes_{\co_X(U)}\co_X(U')\rightarrow \rH^q(\oU'^\circ_\fet,\cL_{U'})
\end{equation}
est un $\alpha$-isomorphisme.
\end{cor}

En effet, d'après \ref{cg30}, on a un morphisme canonique 
\begin{equation}\label{cg310b}
\rH^q(\Delta_\infty,\cL_U^{(\infty)\star}(\oU^{\star\circ}))\rightarrow \rH^q(\oU^{\star\circ}_\fet,\cL_U|\oU^{\star\circ}),
\end{equation}
qui est un $\alpha$-isomorphisme. Celui-ci est fonctoriel en $U\in \ob(\cC)$ d'après la preuve de {\em loc. cit.}
Par ailleurs, d'après \ref{cg40}, le morphisme canonique 
\begin{equation}
\cL_{U}^{(\infty)\star}(\oU^{\star\circ})\otimes_{\co_X(U)}\co_X(U')\rightarrow \cL^{(\infty)\star}_{U'}(\oU'^{\star\circ})
\end{equation}
est un $\alpha$-isomorphisme. Comme $\co_X(U')$ est plat sur $\co_X(U)$, on en déduit que le morphisme canonique
\begin{equation}\label{cg310c}
\rH^q(\oU^{\star\circ}_\fet,\cL_{U}|\oU^{\star\circ})\otimes_{\co_X(U)}\co_X(U')\rightarrow \rH^q(\oU'^{\star\circ}_\fet,\cL_{U'}|\oU'^{\star\circ})
\end{equation}
est un $\alpha$-isomorphisme (\cite{agt} II.3.15). La proposition s'ensuit en variant le point géométrique $\oy$ de $\oX^\circ$ \eqref{cg2}.

\begin{cor}\label{cg31}
Pour tout morphisme $U'\rightarrow U$ de $\cC$ et tous entiers $m\geq 1$ et $q\geq 0$, le morphisme canonique
\begin{equation}\label{cg31a}
\rH^q(\oU^\circ_\fet,\ocB_{U,m})\otimes_{\co_X(U)}\co_X(U')\rightarrow \rH^q(\oU'^\circ_\fet,\ocB_{U',m})
\end{equation}
est un $\alpha$-isomorphisme.
\end{cor}

\begin{prop}\label{cg29}
Soient $m$ un entier $\geq 1$, $\mL$ un $(\co_\oK/p^m\co_\oK)$-module localement libre de type fini de $\oX^\circ_\fet$. 
Pour tout $X$-schéma étale $U$, on note $\mL_U$ l'image inverse de $\mL$ par le morphisme canonique $\oU^\circ\rightarrow \oX^\circ$,
et on pose $\cL_U=\mL_U\otimes_{\co_\oK}\ocB_U$. Alors, pour tout entier $q\geq 0$, 
le préfaisceau de $\alpha$-$(\co_\oK)$-modules sur $\cC$ qui à tout objet $U$ de $\cC$ 
associe le $\alpha$-$(\co_\oK)$-module $\alpha(\rH^q(\oU^\circ_\fet,\cL_U))$ \eqref{alpha4} est un faisceau pour la topologie étale \eqref{alpha18}. 
\end{prop}

Soit $(U_i\rightarrow U)_{i\in I}$ un recouvrement de $\cC$. Pour tout $(i,j)\in I^2$, posons $U_{ij}=U_i\times_UU_j$. 
Par descente fidèlement plate, la suite de $\co_\oX(\oU)$-modules 
\begin{eqnarray}
0\rightarrow \rH^q(\oU^\circ_\fet,\cL_U)&\rightarrow& \prod_{i\in I}\rH^q(\oU^\circ_\fet,\cL_U)\otimes_{\co_X(U)}\co_X(U_i) \\
&\rightarrow& \prod_{(i,j)\in I^2}\rH^q(\oU^\circ_\fet,\cL_U)\otimes_{\co_X(U)}\co_X(U_{ij}),\nonumber
\end{eqnarray}
où la dernière flèche est la différence des morphismes induits par les projections de $U_{ij}$ sur les deux facteurs, 
est exacte. On en déduit par \ref{cg310} que la suite de $\co_\oK$-modules 
\begin{equation}
0\rightarrow \rH^q(\oU^\circ_\fet,\cL_U)\rightarrow \prod_{i\in I}\rH^q((\oU^\circ_i)_{\fet},\cL_{U_i})
\rightarrow \prod_{(i,j)\in I^2}\rH^q((\oU^\circ_{ij})_{\fet},\cL_{U_{ij}}),
\end{equation}
où la dernière flèche est la différence des morphismes induits par les projections de $U_{ij}$ 
sur les deux facteurs, est $\alpha$-exacte \eqref{finita3}. La proposition s'ensuit en vertu de \ref{alpha19}.

\begin{prop}\label{cg32}
Soient $m,q$ deux entiers tels que $m\geq 1$ et $q\geq 0$, $\mL$ un $(\co_\oK/p^m\co_\oK)$-module localement libre de type fini de $\oX^\circ_\fet$. 
On note $\rH^q(\mL)$ le $\co_\oX(\oX)$-module 
$\rH^q(\oX^\circ_\fet,\mL\otimes_{\co_\oK}\ocB_X)$ et $\trH^q(\mL)$ le $\hbar_*(\co_\oX)$-module associé de $X_\et$ \eqref{notconv12a}. 
Alors, on a un morphisme $\hbar_*(\co_\oX)$-linéaire canonique de $X_\et$ 
\begin{equation}\label{cg32a}
\trH^q(\mL)\rightarrow \rR^q\sigma_*(\beta^*(\mL)),
\end{equation}
qui est un $\alpha$-isomorphisme, où $\sigma$ et $\beta$ sont les morphismes de topos annelés \eqref{TFA2e} et \eqref{TFA2g}.
\end{prop}

Pour tout $X$-schéma étale $U$, on note $\mL_U$ l'image inverse de $\mL$ par le morphisme canonique $\oU^\circ\rightarrow \oX^\circ$.  
La correspondance $\cL=\{U\mapsto  \mL_U\otimes_{\co_\oK}\ocB_U\}$ définit un préfaisceau abélien sur $E$ \eqref{tf1h}. 
Notons $\cL^\tta$ le faisceau associé de $\tE$. 
D'après (\cite{agt} VI.5.34, VI.8.9 et VI.5.17), on a un isomorphisme canonique 
\begin{equation}
\beta^*(\mL)\stackrel{\sim}{\rightarrow}\cL^\tta. 
\end{equation}
Suivant (\cite{agt} VI.10.39), notons $\cH^q(\cL)$ le faisceau de $X_\et$ associé au préfaisceau sur $\Et_{/X}$ 
défini pour tout $U\in \ob(\Et_{/X})$ par le groupe abélien $\rH^q(\oU^\circ_\fet,\mL_U\otimes_{\co_\oK}\ocB_U)$. 
D'après (\cite{agt} (VI.10.39.10)), on a un morphisme canonique de $X_\et$
\begin{equation}\label{cg32b}
\cH^q(\cL)\rightarrow \rR^q\sigma_*(\beta^*(\mL)).
\end{equation}
Celui-ci est clairement $\hbar_*(\co_\oX)$-linéaire. Comme la source et le but sont annulés par $p^m$, 
\eqref{cg32b} est un isomorphisme au-dessus de $X_\eta$. 
Pour tout point géométrique $\ox$ de $X$ au-dessus de $s$, notant $X_{(\ox)}$ le localisé strict de $X$ en $\ox$,
le schéma $X_{(\ox)}\times_X\oX$ est normal et strictement local (et en particulier intègre), d'après (\cite{agt} III.3.7).
Par suite, la fibre de \eqref{cg32b} en $\ox$ est un isomorphisme en vertu de (\cite{agt} VI.10.40). 
On en déduit que \eqref{cg32b} est un isomorphisme. 

Compte tenu de (\cite{sga4} II 3.0.4) et de la définition du foncteur ``faisceau associé'' (\cite{sga4} II 3.4),   
$\cH^q(\cL)$ est le faisceau de $X_\et$ associé au préfaisceau sur $\cC$ \eqref{cg25}
défini pour tout $U\in \ob(\cC)$ par le groupe abélien $\rH^q(\oU^\circ_\fet,\mL_U\otimes_{\co_\oK}\ocB_U)$.
Il résulte alors de \ref{alpha27} et \ref{cg29}, que pour tout $U\in \ob(\cC)$, le morphisme canonique
\begin{equation}\label{cg32c}
\rH^q(\oU^\circ_\fet,\mL_U\otimes_{\co_\oK}\ocB_U)\rightarrow \cH^q(\cL)(U)
\end{equation}
est un $\alpha$-isomorphisme. On en déduit d'après \ref{cg310} que le morphisme canonique
\begin{equation}\label{cg32d}
\rH^q(\mL)\otimes_R\co_X(U)\rightarrow \cH^q(\cL)(U)
\end{equation}
est un $\alpha$-isomorphisme~; d'où la proposition.

\begin{cor}\label{cg33}
Soient $m,q$ deux entiers tels que $m\geq 1$ et $q\geq 0$, $\mL$ un $(\co_\oK/p^m\co_\oK)$-module localement libre de type fini de $\oX^\circ_\fet$. 
On note $\rM^q(\mL)$ le $\co_{\oX_m}(\oX_m)$-module $\rH^q(\oX^\circ_\fet,\mL\otimes_{\co_\oK}\ocB_X)$ 
et $\trM^q(\mL)$ le $(\co_{\oX_m})$-module associé de $(\oX_m)_\et$ \eqref{notconv12a}. 
Alors, on a un morphisme $(\co_{\oX_m})$-linéaire canonique de $(\oX_m)_\et$  
\begin{equation}
\trM^q(\mL)\rightarrow \rR^q\sigma_{m*}(\beta^*_m(\mL))
\end{equation}
qui est un $\alpha$-isomorphisme, où $\sigma_m$ et $\beta_m$ sont les morphismes de topos annelés \eqref{TFA8e} et \eqref{TFA8g}.
\end{cor}

En effet, on a $\beta^*(\mL)=\beta^{-1}(\mL)\otimes_{(\co_\oK/p^m\co_\oK)}\ocB_m$, où $\beta$ est le morphisme de topos annelés \eqref{TFA2g},
qui est donc un objet de $\tE_s$ d'après (\cite{agt} III.9.6 et III.9.7). 
Le morphisme canonique $\delta^*(\beta^*(\mL))\rightarrow \beta^*_m(\mL)$ est un isomorphisme, 
où $\delta$ est le plongement \eqref{TFA66a}. On en déduit par adjonction un isomorphisme
\begin{equation}
\beta^*(\mL)\stackrel{\sim}{\rightarrow} \delta_*(\beta^*_m(\mL)).
\end{equation} 
Il résulte alors de \eqref{TFA66d} qu'on a un isomorphisme canonique
\begin{equation}
\rR^q\sigma_{m*}(\beta^*_m(\mL))\stackrel{\sim}{\rightarrow}a^{-1}(\rR^q\sigma_*(\beta^*(\mL))),
\end{equation}
où $a^{-1}$ désigne l'image inverse au sens des faisceaux abéliens par l'injection canonique $a\colon X_s\rightarrow X$. 
Par ailleurs, avec les notations de \ref{cg32}, on a un isomorphisme canonique 
\begin{equation}
\trH^q(\mL)\stackrel{\sim}{\rightarrow} a_*(\trM^q(\mL)).
\end{equation}
On en déduit par adjonction un isomorphisme $a^{-1}(\trH^q(\mL))\stackrel{\sim}{\rightarrow} \trM^q(\mL)$. La proposition résulte donc de \ref{cg32}. 

\begin{rema}\label{cg34}
Conservons les hypothèses de \ref{cg33}. 
D'après (\cite{egr1} 2.6.18(ii)), $\co_{\oX_m}$ est un faisceau d'anneaux cohérent de $(\oX_m)_\zar$. 
Il résulte alors de \ref{afini6} et \ref{cg37}(iii) que le $(\co_{\oX_m})$-module de $(\oX_m)_\zar$ associé à $\rM^q(\mL)$ est $\alpha$-cohérent \eqref{finita9}
et est $\alpha$-nul pour tout $i\geq d+1$, où $d=\dim(X/S)$. 
\end{rema}

\subsection{}\label{cg26}
On désigne par $\hE$ (resp. $\hE_\cC$) la catégorie des préfaisceaux de $\mU$-ensembles sur $E$ (resp. $E_\cC$)
et par 
\begin{equation}\label{cg26a}
\cP^\vee\rightarrow \Et_{/X}^\circ
\end{equation}
la catégorie fibrée obtenue en associant à tout $U\in \ob(\Et_{/X})$ la catégorie $(\Et_{\rf/\oU^\circ})^\wedge$ 
des préfaisceaux de $\mU$-ensembles
sur $\Et_{\rf/\oU^\circ}$, et à tout morphisme $f\colon U'\rightarrow U$ de $\Et_{/X}$ 
le foncteur 
\begin{equation}\label{cg26b}
\of^\circ_{\fet*}\colon (\Et_{\rf/\oU'^\circ})^\wedge\rightarrow (\Et_{\rf/\oU^\circ})^\wedge
\end{equation} 
défini par la composition avec le foncteur image inverse $\Et_{\rf/\oU^\circ} \rightarrow \Et_{\rf/\oU'^\circ}$
par le morphisme $\of^\circ\colon \oU'^\circ\rightarrow \oU^\circ$ déduit de $f$. On note 
\begin{equation}\label{cg26c}
\cP^\vee_\cC\rightarrow \cC^\circ
\end{equation}
la catégorie fibrée au-dessus de $\cC^\circ$ déduite de la catégorie fibrée \eqref{cg26a} 
par changement de base par le foncteur d'injection canonique $\cC\rightarrow \Et_{/X}$.
On a alors une équivalence de catégories (\cite{agt} VI.5.2)
\begin{eqnarray}
\hE&\stackrel{\sim}{\rightarrow}& \bHom_{(\Et_{/X})^\circ}((\Et_{/X})^\circ,\cP^\vee)\\
F&\mapsto &\{U\mapsto F\circ \iota_{U!}\},\nonumber
\end{eqnarray}
où $\iota_{U!}$ est le foncteur \eqref{TFA6F}. 
On identifiera dans la suite $F$ à la section $\{U\mapsto F\circ \iota_{U!}\}$ qui lui est associée par cette équivalence.
De même, on a une équivalence de catégories
\begin{eqnarray}
\hE_\cC&\stackrel{\sim}{\rightarrow}& \bHom_{\cC^\circ}(\cC^\circ,\cP^\vee_\cC)\label{cg26d}\\
F&\mapsto &\{U\mapsto F\circ \iota_{U!}\}.\nonumber
\end{eqnarray}

On note 
\begin{equation}\label{cg26e}
\hE\rightarrow \tE, \ \ \ F\mapsto F^\tta.
\end{equation}
le foncteur ``faisceau associé'' sur $E$. 
Comme $E_\cC$ est une sous-catégorie topologiquement génératrice de $E$, le foncteur ``faisceau associé'' sur $E_\cC$
induit un foncteur que l'on note aussi 
\begin{equation}\label{cg26f}
\hE_\cC\rightarrow \tE, \ \ \ F\mapsto F^\tta.
\end{equation}

Soient $F=\{W\mapsto G_W\}$ $(W\in \ob(\Et_{/X}))$ un objet de $\hE$, $F_\cC=\{U\mapsto G_U\}$ 
$(U\in \ob(\cC))$ l'objet de $\hE_\cC$ obtenu en restreignant $F$ à $E_\cC$. 
Il résulte aussitôt de (\cite{sga4} II 3.0.4) et de la définition du foncteur ``faisceau associé'' (\cite{sga4} II 3.4)
qu'on a un isomorphisme canonique de $\tE$
\begin{equation}\label{cg26g}
(F_\cC)^\tta\stackrel{\sim}{\rightarrow} F^\tta.
\end{equation}

On désigne par $\ahE_\cC$ la catégorie de $\alpha$-$(\co_\oK)$-préfaisceaux sur $E_\cC$ \eqref{alpha14}, par 
$\atE_\cC$ la catégorie des $\alpha$-$(\co_\oK)$-faisceaux \eqref{alpha18} sur $E_\cC$ et par 
\begin{eqnarray}
\halpha\colon \bMod(\co_\oK,\hE_\cC)&\rightarrow& \ahE_\cC,\label{cg26h}\\
\talpha\colon \bMod(\co_\oK,\tE_\cC)&\rightarrow& \atE_\cC,\label{cg26i}
\end{eqnarray}
les foncteurs canoniques \eqref{alpha16a} et \eqref{alpha22a}. 

\begin{prop}\label{cg27}
Soient $m$ un entier $\geq 1$, $\mL$ un $(\co_\oK/p^m\co_\oK)$-module localement libre de type fini de $\oX^\circ_\fet$. 
Pour tout $X$-schéma étale $U$, on note $\mL_U$ l'image inverse de $\mL$ par le morphisme canonique $\oU^\circ\rightarrow \oX^\circ$,
et on pose $\cL_U=\mL_U\otimes_{\co_\oK}\ocB_U$. Alors, 
le préfaisceau de $\alpha$-$(\co_\oK)$-modules $\halpha(\{U\mapsto \cL_U\})$ $(U\in \ob(\cC))$ \eqref{cg26h} est un faisceau.
\end{prop}

Soient $(V\rightarrow U)$ un objet de $E_\cC$, $(U_i\rightarrow U)_{i\in I}$ un recouvrement de $\cC$. 
Pour tout $(i,j)\in I^2$, posons $U_{ij}=U_i\times_UU_j$, $V_i=V\times_UU_i$ et $V_{ij}=V\times_UU_{ij}$. 
Par descente fidèlement plate, la suite de $\co_X(U)$-modules 
\begin{equation}
0\rightarrow \cL_U(V)\rightarrow \prod_{i\in I}\cL_U(V)\otimes_{\co_X(U)}\co_X(U_i)
\rightarrow \prod_{(i,j)\in I^2}\cL_U(V)\otimes_{\co_X(U)}\co_X(U_{ij}),
\end{equation}
où la dernière flèche est la différence des morphismes induits par les projections de $U_{ij}$ sur les deux facteurs, 
est exacte. En en déduit par \ref{cg240} que la suite de $\co_\oK$-modules 
\begin{equation}
0\rightarrow \cL_U(V)\rightarrow \prod_{i\in I}\cL_{U_i}(V_i)
\rightarrow \prod_{(i,j)\in I^2}\cL_{U_{ij}}(V_{ij}),
\end{equation}
où la dernière flèche est la différence des morphismes induits par les projections de $U_{ij}$ 
sur les deux facteurs, est $\alpha$-exacte \eqref{finita3}. La proposition s'ensuit en vertu de \ref{alpha181}.

\begin{cor}\label{cg28}
Pour tout entier $m\geq 1$ et tout schéma affine $U$ de $\Et_{/X}$, l'homomorphisme canonique
\begin{equation}
\ocB_{U,m}\rightarrow \ocB_m\circ \iota_{U!},
\end{equation}
où $\ocB_m=\ocB/p^m\ocB$ \eqref{TFA8a} et $\iota_{U!}$ est le foncteur \eqref{TFA6F},  est un $\alpha$-isomorphisme.
\end{cor}
En effet, le préfaisceau de $\alpha$-$(\co_\oK)$-modules $\halpha(\{U\mapsto \ocB_{U,m}\})$ $(U\in \ob(\cC))$ \eqref{cg26h}
est un faisceau en vertu de \ref{cg27}.  La proposition résulte alors de \ref{alpha27}, compte tenu de \eqref{TFA8d} et \eqref{cg26g}.

\section{Le principal théorème de comparaison de Faltings}\label{TPCF}

\subsection{}\label{TPCF1}
Pour tout entier $n\geq 1$, on dispose des morphismes canoniques \eqref{TFA8e} et \eqref{notconv12j} (cf. \ref{TFA5})
\begin{eqnarray}
\sigma_n\colon (\tE_s,\ocB_n)&\rightarrow& (X_{s,\et},\co_{\oX_n}),\label{TPCF1a}\\
u_n\colon (X_{s,\et},\co_{\oX_n})&\rightarrow& (X_{s,\zar},\co_{\oX_n}).\label{TPCF1b}
\end{eqnarray}
On note le composé $u_n\circ \sigma_n$  
\begin{equation}\label{TPCF1c}
\tau_n\colon (\tE_s,\ocB_n)\rightarrow (X_{s,\zar},\co_{\oX_n}).
\end{equation}

\begin{prop}\label{TPCF2}
Soient $n,q$ deux entiers tels que $n\geq 1$ et $q\geq 0$, $\mL$ un $(\co_\oK/p^n\co_\oK)$-module localement libre de type fini de $\oX^\circ_\fet$. 
Alors, le $(\co_{\oX_n})$-module $\fm_\oK\otimes_{\co_\oK}\rR^q\sigma_{n*}(\beta_n^*(\mL))$ \eqref{TFA8g} est associé à un  $(\co_{\oX_n})$-module
quasi-cohérent et $\alpha$-cohérent \eqref{finita9} de $(\oX_n)_\zar$ \eqref{notconv12a}, 
et il est nul pour tout $q\geq d+1$, où $d=\dim(X/S)$.
\end{prop}

En effet, soit $(u_j\colon X_j\rightarrow X)_{j\in J}$ un recouvrement étale tel que $J$ soit fini et que 
pour tout $j\in J$, $(X_j,\cM_X|X_j)$ vérifie les hypothèses de \ref{cg21}. 
Soit $j\in J$. D'après (\cite{agt} (VI.10.12.6) et (III.9.11.12)), on a un diagramme de morphismes de topos annelés 
\begin{equation}\label{TPCF2f}
\xymatrix{
{((\oX_j^\circ)_\fet,\co_\oK/p^n\co_\oK)}\ar[d]_{\ou^\circ_j}&{(\tE_{j,s},\ocB_{j,n})}\ar[l]_-(0.5){\beta_{j,n}}\ar[d]^{\Phi_j}\ar[r]^-(0.5){\sigma_{j,n}}&
{((X_{j,s})_{\et},\co_{\oX_{j,n}})}\ar[d]^{\ou_{j,s}}\\
{(\oX^\circ_\fet,\co_\oK/p^n\co_\oK)}&{(\tE_s,\ocB_n)}\ar[l]_-(0.5){\beta_n}\ar[r]^-(0.5){\sigma_n}&{(X_{s,\et},\co_{\oX_n})}}
\end{equation}
où la ligne supérieure est l'analogue pour $X_j$ de la ligne inférieure pour $X$, les flèches
verticales sont définies par fonctorialité, et dont les carrés sont commutatifs à isomorphismes canoniques près. 
Par ailleurs, le morphisme $\Phi_j$ s'identifie au morphisme de localisation de $(\tE_s,\ocB_n)$ en $\sigma_s^*(X_{j,s})$ en vertu de (\cite{agt} III.9.12 et III.9.13). 
En fait, $\sigma_{j,n}$ s'identifie au morphisme $(\sigma_n)_{/X_{j,s}}$, de sorte que le carré de droite de \eqref{TPCF2f} est 2-cartésien (cf. \cite{sga4} IV 5.10).
Par suite, en vertu de \ref{cg33} et \ref{cg34}, il existe un $(\co_{\oX_{j,n}})$-module quasi-cohérent et $\alpha$-cohérent
$\cF_j$ de $(\oX_{j,n})_\zar$ et un morphisme $(\co_{\oX_{j,n}})$-linéaire
\begin{equation}\label{TPCF2a}
\iota_j(\cF_j)\rightarrow  \rR^q\sigma_{n*}(\beta_n^*(\mL))|X_{j,s},
\end{equation}
qui est un $\alpha$-isomorphisme, où $\iota_j(\cF_j)$ est le $(\co_{\oX_{j,n}})$-module de $(\oX_{j,n})_\et$ associé à $\cF_j$ \eqref{notconv12a}.
En particulier, le morphisme induit
\begin{equation}\label{TPCF2b}
\fm_\oK\otimes_{\co_\oK}\iota_j(\cF_j)\rightarrow  \fm_\oK\otimes_{\co_\oK}\rR^q\sigma_{n*}(\beta_n^*(\mL))|X_{j,s}
\end{equation}
est un isomorphisme \eqref{alpha21}. 
Par suite, $\fm_\oK\otimes_{\co_\oK}\rR^q\sigma_{n*}(\beta_n^*(\mL))$ est nul pour tout $q\geq d+1$, d'après \ref{cg37}(iii).
Compte tenu de \eqref{notconv12e}, on a un isomorphisme canonique 
\begin{equation}\label{TPCF2c}
\fm_\oK\otimes_{\co_\oK}\iota_j(\cF_j)\stackrel{\sim}{\rightarrow} \iota_j(\fm_\oK\otimes_{\co_\oK}\cF_j).
\end{equation}
Les isomorphismes \eqref{TPCF2b} définissent une donnée de descente sur $(\iota_j(\fm_\oK\otimes_{\co_\oK}\cF_j))_{j\in J}$ relativement 
au recouvrement étale $(\oX_{j,n}\rightarrow \oX_n)_{j\in J}$. Comme le foncteur \eqref{notconv12a} est pleinement fidèle 
et compte tenu de \eqref{afini12l}, 
on en déduit une donnée de descente sur $(\fm_\oK\otimes_{\co_\oK}\cF_j)_{j\in J}$ relativement au même recouvrement. 
D'après (\cite{sga1} VIII 1.1), il existe alors un $(\co_{\oX_n})$-module quasi-cohérent $\cG$ et pour tout $j\in J$, 
un isomorphisme
\begin{equation}\label{TPCF2d}
\cG\otimes_{\co_{\oX_n}}\co_{\oX_{j,n}}\stackrel{\sim}{\rightarrow} \fm_\oK\otimes_{\co_\oK}\cF_j
\end{equation} 
tels que si l'on munit $(\cG\otimes_{\co_{\oX_n}}\co_{\oX_{j,n}})_{j\in J}$ de la donnée de descente définie par $\cG$, 
les isomorphismes \eqref{TPCF2d} soient compatibles aux données de descente.
Les isomorphismes \eqref{TPCF2b} se descendent en 
un isomorphisme 
\begin{equation}\label{TPCF2e}
\iota(\cG)\rightarrow  \fm_\oK\otimes_{\co_\oK}\rR^q\sigma_{n*}(\beta_n^*(\mL)),
\end{equation}
où $\iota(\cG)$ est le $(\co_{\oX_n})$-module de $(\oX_n)_\et$ associé à $\cG$.
Par ailleurs, $\cG$ est $\alpha$-cohérent, d'après \ref{finita16} et \ref{afini12}; d'où la proposition.

\begin{cor}\label{TPCF3}
Soient $n,q$ deux entiers tels que $n\geq 1$ et $q\geq 0$, $\mL$ un $(\co_\oK/p^n\co_\oK)$-module localement libre de type fini de $\oX^\circ_\fet$. 
Alors, le $(\co_{\oX_n})$-module $\rR^q\tau_{n*}(\beta_n^*(\mL))$ est $\alpha$-quasi-cohérent \eqref{afini3} et $\alpha$-cohérent \eqref{finita9}, 
et est $\alpha$-nul pour tout $q\geq d+1$, où $d=\dim(X/S)$.
\end{cor}

En effet, pour tout $(\co_{\oX_n})$-module $\cF$ de $(\oX_n)_\et$ et pour tout entier $i\geq 0$, les morphismes canoniques \eqref{TPCF1b}
\begin{equation}\label{TPCF3a}
\fm_\oK\otimes_{\co_\oK}\rR^ju_{n*}(\cF)\rightarrow \rR^iu_{n*}(\fm_{\co_\oK}\otimes_{\co_\oK}\cF)\rightarrow 
\rR^ju_{n*}(\cF)
\end{equation}
sont des $\alpha$-isomorphismes d'après \ref{alpha3}. Par suite, 
$\rR^iu_{n*}(\rR^j\sigma_{n*}(\beta_n^*(\mL)))$ est $\alpha$-nul pour tous $i\geq 1$ et $j\geq 0$, en vertu de \ref{TPCF2} et (\cite{sga4} VII 4.3). 
La suite spectrale de Cartan-Leray
\begin{equation}\label{TPCF3b}
E_2^{i,j}=\rR^iu_{n*}(\rR^j\sigma_{n*}(\beta_n^*(\mL)))\Rightarrow \rR^{i+j}\tau_{n*}(\beta_n^*(\mL))
\end{equation}
induit donc un isomorphisme de $\alpha$-$(\co_{\oX_n})$-modules de $(\oX_n)_\zar$ \eqref{alpha20}
\begin{equation}\label{TPCF3c}
\alpha(u_{n*}(\fm_\oK\otimes_{\co_\oK}\rR^q\sigma_{n*}(\beta_n^*(\mL))))\stackrel{\sim}{\rightarrow} \alpha(\rR^q\tau_{n*}(\beta_n^*(\mL))).
\end{equation}
La proposition s'ensuit compte tenu de \ref{TPCF2} et \ref{finita16}.

\begin{cor}\label{TPCF4}
Supposons $X$ propre sur $S$, et soient $n,q$ deux entiers tels que $n\geq 1$ et $q\geq 0$, 
$\mL$ un $(\co_\oK/p^n\co_\oK)$-module localement libre de type fini de $\oX^\circ_\fet$. 
Alors, le $\co_\oK$-module $\rH^q(\tE_s,\beta_n^*(\mL))$ est de présentation $\alpha$-finie et est $\alpha$-nul pour tout $q\geq 2d+1$.
\end{cor}

On notera d'abord que l'anneau $\co_\oK/p^n\co_\oK$ est universellement cohérent (\cite{egr1} 1.4.1 et 1.12.15).
Considérons la suite spectrale de Cartan Leray
\begin{equation}\label{TPCF4a}
E_2^{i,j}=\rH^i(\oX_n,\rR^j\tau_{n*}(\beta_n^*(\mL)))\Rightarrow \rH^{i+j}(\tE_s,\beta_n^*(\mL)).
\end{equation}
En vertu de \ref{TPCF3} et \ref{afini9}, pour tous entiers $i,j\geq 0$, le $(\co_\oK/p^n\co_\oK)$-module $E_2^{i,j}$ 
est $\alpha$-cohérent, et il est $\alpha$-nul si $i\geq d+1$ ou si $j\geq d+1$.
La proposition s'ensuit compte tenu de \ref{finita18}.

\subsection{}\label{TPCF5}
Considérons l'anneau \eqref{eipo3a}
\begin{equation}\label{TPCF5a}
(\co_\oK)^\flat= \underset{\underset{x\mapsto x^p}{\longleftarrow}}{\lim}\co_\oK/p\co_\oK,
\end{equation}
qui est intègre, parfait de caractéristique $p$, et notons $\oK^\flat$ son corps des fractions, qui est algébriquement clos (cf. \cite{scholze} 3.8).
On rappelle que l'application 
\begin{equation}\label{TPCF5b}
\{(x_n)_{n\geq 0}\in (\co_C)^{\mN} \ |\ (x_{n+1})^p=x_n,\ \forall n\in \mN\}\rightarrow (\co_\oK)^\flat, \ \ \  
(x_n)_{n\geq 0}\mapsto (\ox_n)_{n\geq 0},
\end{equation}
où $\ox_n$ désigne la classe de $x_n$ dans $\co_\oK/p\co_\oK$, est un isomorphisme de monoïdes multiplicatifs. 
Pour tout $x\in (\co_\oK)^\flat$, d'image inverse $(x_n)_{n\geq 0}$ par \eqref{TPCF5b}, on pose  
\begin{equation}\label{TPCF5c}
x^\sharp=x_0.
\end{equation}
L'application 
\begin{equation}\label{TPCF5d}
(\co_\oK)^\flat\rightarrow \mQ, \ \ \ x\mapsto v(x^\sharp),
\end{equation}
où $v$ est la valuation normalisée de $\co_C$ \eqref{TFA1}, est une valuation non discrète, de hauteur $1$ de $(\co_\oK)^\flat$. 
Il est alors naturel de noter $(\co_\oK)^\flat$ aussi $\co_{\oK^\flat}$, et son idéal maximal $\fm_{\oK^\flat}$.

\subsection{}\label{TPCF15}
On fixe une suite $(p_n)_{n\geq 0}$ d'éléments de $\co_\oK$ telle que $p_0=p$ et $p_{n+1}^p=p_n$ pour tout $n\geq 0$
et on note $\varpi$ l'élément associé de $\co_{\oK^\flat}$. Pour tout entier $i\geq 0$, l'homomorphisme 
\begin{equation}\label{TPCF15a}
\pi_i \colon \co_{\oK^\flat}\rightarrow \co_\oK/p\co_\oK, \ \ \ (x_n)_{n\geq 0}\mapsto x_i
\end{equation}
est surjectif. On vérifie aussitôt que $\ker(\pi_0)=\varpi \co_{\oK^\flat}$. On note $\varphi$ l'endomorphisme de Frobenius de $\co_{\oK^\flat}$.
Comme $\pi_0=\pi_i\circ \varphi^i$, on en déduit que $\ker(\pi_i)=\varpi^{p^i}\co_{\oK^\flat}$. En particulier, $\pi_i$ induit un isomorphisme
\begin{equation}\label{TPCF15b}
\co_{\oK^\flat}/\varpi^{p^i}\co_{\oK^\flat}\stackrel{\sim}{\rightarrow}\co_\oK/p\co_\oK.
\end{equation}
Par ailleurs, le diagramme 
\begin{equation}\label{TPCF15c}
\xymatrix{
{\co_{\oK^\flat}/\varpi^{p^{i+1}}\co_{\oK^\flat}}\ar[r]^-(0.5)\sim\ar[d]&{\co_\oK/p\co_\oK}\ar[d]\\
{\co_{\oK^\flat}/\varpi^{p^i}\co_{\oK^\flat}}\ar[r]^-(0.5)\sim&{\co_\oK/p\co_\oK}}
\end{equation}
où la flèche verticale de gauche est le morphisme canonique et la flèche verticale de droite est l'endomorphisme de Frobenius, est commutatif.
On en déduit que $\co_{\oK^\flat}$ est complet et séparé pour la topologie $\varpi$-adique, et que cette dernière 
coïncide avec la topologie limite projective des topologies discrètes de $\co_\oK/p\co_\oK$.

\subsection{}\label{TPCF6}
On désigne par $\tE_s^{\mN^\circ}$ le topos des systèmes projectifs de $\tE_s$, indexés par l'ensemble ordonné $\mN$ des entiers naturels 
(\cite{agt} III.7.7). On rappelle qu'on a un morphisme de topos 
\begin{equation}\label{TPCF6a}
\lambda\colon \tE_s^{\mN^\circ}\rightarrow \tE_s,
\end{equation}
dont le foncteur image inverse $\lambda^*$ associe à tout objet $F$ de $\tE_s$ le système projectif constant de valeur $F$,
et dont le foncteur image directe $\lambda_*$ associe à tout système projectif sa limite projective (\cite{agt} III.7.4).

Pour tout objet $F=(F_n)_{n\geq 0}$ de $\tE_s^{\mN^\circ}$ et tout entier $i\geq 0$, on note $\varepsilon_{\leq i}(F)$ l'objet $(F_n^i)_{n\geq 0}$  
de $\tE_s^{\mN^\circ}$, défini par 
\begin{equation}\label{TPCF6b}
F_n^i=\left\{
\begin{array}{clcr}
F_n&{\rm si}\ 0\leq n\leq i,\\
F_i&{\rm si}\ n\geq i,
\end{array}
\right.
\end{equation}
le morphisme de transition $F_{n+1}^i\rightarrow F_n^i$ étant égal au morphisme de transition $F_{n+1}\rightarrow F_n$ de $F$ si $0\leq n<i$ et à l'identité 
si $n\geq i$. La correspondance $F\mapsto \varepsilon_{\leq i}(F)$ est clairement fonctorielle, et on a un morphisme fonctoriel canonique
\begin{equation}\label{TPCF6c}
F\rightarrow  \varepsilon_{\leq i}(F).
\end{equation}
Si $F$ est un anneau (resp. module), il en est de même de $\varepsilon_{\leq i}(F)$ et le morphisme \eqref{TPCF6c} est un homomorphisme. 

Si $A$ est un anneau, on note encore $A$ l'anneau constant de valeur $A$ de $\tE_s$ ou de $\tE_s^{\mN^\circ}$.
On a donc un isomorphisme canonique $\lambda^*(A)\stackrel{\sim}{\rightarrow} A$ de $\tE_s^{\mN^\circ}$.
On note $(\co_\oK)^\wp$ l'anneau $(\co_\oK/p\co_\oK)_{\mN}$ de $\tE_s^{\mN^\circ}$ dont les morphismes de transition sont les itérés
de l'endomorphisme de Frobenius de $\co_\oK/p\co_\oK$. On le considère comme  une $\co_{\oK^\flat}$-algèbre 
via les homomorphismes \eqref{TPCF15b}.

On désigne par $\cP$ la catégorie abélienne 
des $\mF_p$-modules de $\tE_s^{\mN^\circ}$ et par $\cP_\AR$ son quotient par la sous-catégorie 
épaisse des $\mF_p$-modules Artin-Rees-nuls (AR-nuls en abrégé) (\cite{sga5} V 2.2.1). On rappelle que $\cP_\AR$
est canoniquement équivalente à la catégorie des systèmes projectifs de $\cP$ à translation près (\cite{sga5} V 2.4.2 et 2.4.4).

\subsection{}\label{TPCF7}
On note $\phi$ l'endomorphisme de Frobenius de $\ocB_1$ \eqref{TFA8a}
et on désigne par $\ocB^\wp$ la $(\co_\oK)^\wp$-algèbre $(\ocB_1)_{\mN}$ de $\tE_s^{\mN^\circ}$  
dont les morphismes de transition sont les itérés de $\phi$ . 
On considère $\ocB^\wp$ comme une $\co_{\oK^\flat}$-algèbre via les homomorphismes \eqref{TPCF15b}.
On note encore $\phi$ l'endomorphisme de Frobenius de $\ocB^\wp$.
Cette notation n'induit aucun risque d'ambiguïté puisque l'endomorphisme de Frobenius de $\ocB^\wp$ est induit par celui de $\ocB_1$. 
Pour tout entier $i\geq 0$, on considère $\varepsilon_{\leq i}(\ocB^\wp)$ comme une $\ocB^\wp$-algèbre via l'homomorphisme canonique \eqref{TPCF6c}.

\begin{lem}\label{TPCF8}
\begin{itemize}
\item[{\rm (i)}] Pour tout entier $n\geq 0$, notant $\op_n$ la classe de $p_n$ dans $\co_\oK/p\co_\oK$ \eqref{TPCF15}, la suite
\begin{equation}\label{TPCF8a}
\xymatrix{
{\ocB_1}\ar[r]^-(0.5){\cdot \op_n}&{\ocB_1}\ar[r]^{\phi^n}&{\ocB_1}\ar[r]&0}
\end{equation}
est exacte.
\item[{\rm (ii)}] L'endomorphisme de Frobenius $\phi$ de $\ocB^\wp$ induit un isomorphisme de $\cP_\AR$ \eqref{TPCF6}.
\item[{\rm (iii)}] Pour tout entier $i\geq 0$, la suite de $\ocB^\wp$-modules
\begin{equation}\label{TPCF8b}
\xymatrix{
0\ar[r]&{\ocB^\wp}\ar[r]^-(0.5){\cdot \varpi^{p^i}}&{\ocB^\wp}\ar[r]&{\varepsilon_{\leq i}(\ocB^\wp)}\ar[r]&0}
\end{equation}
où la deuxième flèche est le morphisme canonique \eqref{TPCF6c}, est exacte dans $\cP_\AR$. De plus, privée du zéro de gauche, elle est exacte dans $\cP$.
\end{itemize}
\end{lem}

(i) En effet, soient $(\oy\rightsquigarrow \ox)$ un point de $X_\et\gtimes_{X_\et}\oX^\circ_\et$ \eqref{rec15} tel que $\ox$ soit au-dessus de $s$,
$X'$ le localisé strict de $X$ en $\ox$. On désigne par $\fV_\ox$ la catégorie des $X$-schémas étales $\ox$-pointés, 
et par $\fW_\ox$ la sous-catégorie pleine de $\fV_\ox$ formée des objets $(U,\fp\colon \ox\rightarrow U)$ 
tels que la restriction $(U,\cM_X|U)\rightarrow (S,\cM_S)$ du morphisme $f$ \eqref{TFA3} vérifie les hypothèses de \ref{cg21}.
Ce sont des catégories cofiltrantes, et le foncteur d'injection canonique
$\fW^\circ_\ox\rightarrow \fV^\circ_\ox$ est cofinal (\cite{sga4} I 8.1.3(c)).
Avec les notations de \ref{TFA11}, on a donc un isomorphisme canonique \eqref{TFA11c} 
\begin{equation}
\ocB_{\rho(\oy\rightsquigarrow \ox)} \stackrel{\sim}{\rightarrow} 
\underset{\underset{(U,\fp)\in \fW^\circ_\ox}{\longrightarrow}}{\lim}\ \oR^{\oy}_U.
\end{equation} 
La proposition s'ensuit en vertu de \ref{cg67} et \ref{TFA7}(ii).

(ii) Il résulte de (i) et (\cite{agt} III.7.3(i)) que l'endomorphisme $\phi$ de $\ocB^\wp$ est surjectif, de noyau $(\op_1\ocB_1)_{\mN}$, où les morphismes
de transition sont induits par l'endomorphisme de Frobenius $\phi$ de $\ocB_1$, et sont donc nuls. 

(iii) Pour tout entier $j<0$, posons $p_j=p$ et $\op_j=0$. On notera que $\varpi^{p^i}$ correspond à l'élément $(\op_{n-i})_{n\geq 0}$ via l'égalité \eqref{TPCF5a}.
Compte tenu de (i) et (\cite{agt} III.7.3(i)), la suite privée du zéro de gauche est donc exacte dans $\cP$. 
Pour tout entier $n$, posons $q_n=p/p_n\in \co_\oK$. Comme $\ocB$ est $\co_\oK$-plat,
le noyau de la multiplication par $\varpi^{p^i}$ dans $\ocB^\wp$ est  $(q_{n-i}\ocB_1)_{n\geq 0}$, 
où les morphismes de transition sont induits par $\phi$ et sont donc nuls en degrés $\geq i$ puisque $v(q_n)=1-1/p^n$; d'où la proposition.

\begin{prop}\label{TPCF9}
Supposons $X$ propre sur $S$ et soient  $\mL$ un $\mF_p$-module de type fini de $\oX^\circ_\fet$, $q$ un entier $\geq 0$. 
Posons $\cL=\delta^*(\beta^*(\mL))$, où $\beta$ est le morphisme \eqref{TFA6b} et $\delta$ est le plongement \eqref{TFA66a},
et notons $\uphi$ l'endomorphisme $\co_{\oK^\flat}$-semi-linéaire de $\rH^q(\tE_s^{\mN^\circ},\lambda^*(\cL)\otimes_{\mF_p}\ocB^\wp)$
induit par l'endomorphisme de Frobenius $\phi$ de $\ocB^\wp$. Alors, il existe un entier $r\geq 0$ et un morphisme $\co_{\oK^\flat}$-linéaire 
\begin{equation}\label{TPCF9a}
\fm_{\oK^\flat}\otimes_{\co_{\oK^\flat}}\rH^q(\tE_s^{\mN^\circ},\lambda^*(\cL)\otimes_{\mF_p}\ocB^\wp)\rightarrow \co_{\oK^\flat}^r,
\end{equation}
compatible à $\uphi$ et à l'endomorphisme de Frobenius $\varphi$ de $\co_{\oK^\flat}$, et qui est un $\alpha$-isomorphisme.
\end{prop}

En effet, considérons le système projectif de $\co_{\oK^\flat}$-modules $M=(M_n)_{n\geq 1}$ défini pour tout entier $n\geq 1$, par 
\begin{equation}\label{TPCF9b}
M_n=\rH^q(\tE_s^{\mN^\circ},\lambda^*(\cL)\otimes_{\mF_p}(\ocB^\wp/\varpi^n\ocB^\wp)),
\end{equation} 
les morphismes de transition étant induits par les projections canoniques 
\begin{equation}
\ocB^\wp/\varpi^{n+1}\ocB^\wp\rightarrow \ocB^\wp/\varpi^n\ocB^\wp.
\end{equation}
D'après \ref{TPCF8}(iii), pour tout entier $i\geq 0$, on a un isomorphisme canonique 
\begin{equation}
\ocB^\wp/\varpi^{p^i}\ocB^\wp\stackrel{\sim}{\rightarrow}  \varepsilon_{\leq i}(\ocB^\wp).
\end{equation}
On en déduit, compte tenu de (\cite{agt} III.7.11) et (\cite{jannsen} 1.15), un isomorphisme $\co_{\oK^\flat}$-linéaire
\begin{equation}\label{TPCF9c}
M_{p^i}\stackrel{\sim}{\rightarrow} \rH^q(\tE_s,\cL\otimes_{\mF_p} \ocB_1),
\end{equation}
où le but est considéré comme un $\co_{\oK^\flat}$-module via $\pi_i$ \eqref{TPCF15a}.
De plus, le diagramme 
\begin{equation}\label{TPCF9d}
\xymatrix{
{M_{p^{i+1}}}\ar[r]^-(0.5){\sim}\ar[d]&{\rH^q(\tE_s,\cL\otimes_{\mF_p} \ocB_1)}\ar[d]\\
{M_{p^i}}\ar[r]^-(0.5){\sim}&{\rH^q(\tE_s,\cL\otimes_{\mF_p} \ocB_1)}}
\end{equation}
où la flèche verticale de gauche est le morphisme canonique et la flèche verticale de droite est induite par l'endomorphisme de Frobenius de $\ocB_1$, 
est commutatif. Les $\co_{\oK^\flat}$-modules $M_{p^i}$ sont donc de type $\alpha$-fini en vertu de \ref{TPCF4}.

Reprenons les notations de \ref{mptf10} pour $R=\co_{\oK^\flat}$. 
La multiplication par $\varpi$ dans $\ocB^\wp$ induit un morphisme $\co_{\oK^\flat}$-linéaire 
\begin{equation}
u\colon M[-1]\rightarrow M.
\end{equation}
Il est clair que le composé $M[-1]\stackrel{u}{\rightarrow}M\rightarrow M[-1]$, où la seconde flèche est le morphisme canonique, est induit par 
la multiplication par $\varpi$. Par ailleurs, il résulte de \ref{TPCF8}(iii) (pour $i=0$) que la suite
\begin{equation}
\xymatrix{
0\ar[r]&{\ocB^\wp/\varpi \ocB^\wp}\ar[r]^-(0.5){\cdot\varpi^n}&{\ocB^\wp/\varpi^{n+1}\ocB^\wp}\ar[r]&{\ocB^\wp/\varpi^n\ocB^\wp}\ar[r]&0}
\end{equation}
est exacte dans $\cP_\AR$. On en déduit que la suite 
\begin{equation}
\xymatrix{M_1\ar[r]^-(0.5){u^{\otimes n}}&M_{n+1}\ar[r]&M_n}
\end{equation}
où la seconde flèche est le morphisme canonique, est exacte au centre. 

Il résulte de \ref{TPCF8}(ii) que l'endomorphisme de Frobenius de $\ocB^\wp$ induit pour tout $n\geq 1$, un isomorphisme de $\cP_\AR$ 
\begin{equation}\label{TPCF9e}
\phi\colon \ocB^\wp/\varpi^n\ocB^\wp\rightarrow\ocB^\wp/\varpi^{pn}\ocB^\wp,
\end{equation}
et par suite un isomorphisme $\co_{\oK^\flat}$-semi-linéaire
\begin{equation}\label{TPCF9f}
M\stackrel{\sim}{\rightarrow}M^{(p)}.
\end{equation}
On notera que pour tout $i\geq 0$, le digramme 
\begin{equation}
\xymatrix{
{M_{p^i}}\ar[r]\ar[d]&{\rH^q(\tE_s,\cL\otimes_{\mF_p} \ocB_1)}\ar@{=}[d]\\
{M_{p^{i+1}}}\ar[r]&{\rH^q(\tE_s,\cL\otimes_{\mF_p} \ocB_1)}}
\end{equation}
où la flèche verticale de gauche est la composante de degré $p^i$ de 
l'isomorphisme \eqref{TPCF9f} et les flèches horizontales sont les isomorphismes \eqref{TPCF9c}, est commutatif. 

En vertu de \ref{mptf11}(i), les morphismes de transition du système projectif $M$ sont $\alpha$-surjectifs.  
Par ailleurs, d'après (\cite{agt} III.7.11), on a une suite exacte
\begin{eqnarray}
\lefteqn{0\rightarrow \rR^1\underset{\underset{\mN}{\longleftarrow}}{\lim}\ \rH^{q-1}(\tE_s,\cL\otimes_{\mF_p}\ocB_1)}\\
&&\rightarrow \rH^q(\tE_s^{\mN^\circ},\lambda^*(\cL)\otimes_{\mF_p}\ocB^\wp)\rightarrow \underset{\underset{\mN}{\longleftarrow}}{\lim}\ 
\rH^q(\tE_s,\cL\otimes_{\mF_p}\ocB_1)\rightarrow 0,\nonumber
\end{eqnarray}
où les morphismes de transition des systèmes projectifs sont induits par l'endomorphisme de Frobenius $\phi$ de $\ocB_1$.
Compte tenu de \ref{alpha10}, \eqref{TPCF9d}, (\cite{jannsen} 1.15) et (\cite{roos} théo.~1), le terme de gauche de cette suite est donc $\alpha$-nul. 
On en déduit que le morphisme canonique
\begin{equation}
\rH^q(\tE_s^{\mN^\circ},\lambda^*(\cL)\otimes_{\mF_p}\ocB^\wp)\rightarrow \underset{\longleftarrow}{\lim}\ M_n
\end{equation}
est un $\alpha$-isomorphisme. La proposition s'ensuit en vertu de \ref{mptf11}(iii).

\begin{cor}\label{TPCF90}
Sous les hypothèses de \eqref{TPCF9}, pour tout entier $i\geq 0$, il existe un morphisme canonique $\co_{\oK^\flat}$-linéaire
\begin{equation}
\rH^q(\tE_s^{\mN^\circ},\lambda^*(\cL)\otimes_{\mF_p}\ocB^\wp)\otimes_{\co_{\oK^\flat}}(\co_{\oK^\flat}/\varpi^{p^i}\co_{\oK^\flat})\rightarrow
\rH^q(\tE_s,\cL\otimes_{\mF_p}\ocB_1),
\end{equation}
où le but est considéré comme un $\co_{\oK^\flat}$-module via $\pi_i$ \eqref{TPCF15a}, qui est un $\alpha$-isomorphisme.
\end{cor}

\begin{prop}\label{TPCF10}
Pour tout élément $b\in \co_\oK$ tel que $b^p\not\in p\co_\oK$, la suite 
\begin{equation}\label{TPCF10a}
\xymatrix{
0\ar[r]&{\mF_p \cdot b\oplus (\frac{p}{b^{p-1}})\cdot \ocB_1}\ar[r]&{\ocB_1}\ar[rr]^{\phi-b^{p-1}\id}&&{\ocB_1}\ar[r]&0}
\end{equation}
est exacte.
\end{prop} 

En effet, la question étant locale, on peut se borner au cas où $X$ est affine. 
Soient $(\oy\rightsquigarrow \ox)$ un point de $X_\et\gtimes_{X_\et}\oX^\circ_\et$ \eqref{rec15} tel que $\ox$ soit au-dessus de $s$,
$X'$ le localisé strict de $X$ en $\ox$. On désigne par $\fV_\ox$ la catégorie des $X$-schémas étales $\ox$-pointés, 
et par $\fW_\ox$ la sous-catégorie pleine de $\fV_\ox$ formée des objets $(U,\fp\colon \ox\rightarrow U)$
tels que le schéma $U$ soit affine. 
Ce sont des catégories cofiltrantes, et le foncteur d'injection canonique
$\fW^\circ_\ox\rightarrow \fV^\circ_\ox$ est cofinal (\cite{sga4} I 8.1.3(c)). 
Avec les notations de \ref{TFA11}, on a un isomorphisme canonique \eqref{TFA11c} 
\begin{equation}\label{TPCF10b}
\ocB_{\rho(\oy\rightsquigarrow \ox)} \stackrel{\sim}{\rightarrow} 
\underset{\underset{(U,\fp)\in \fW^\circ_\ox}{\longrightarrow}}{\lim}\ \oR^{\oy}_U.
\end{equation} 
Montrons que pour tout $z\in \ocB_{\rho(\oy\rightsquigarrow \ox)}$, il existe $t\in \ocB_{\rho(\oy\rightsquigarrow \ox)}$ tel que 
$t^p-b^{p-1}t=z$. On peut supposer que $z\in \oR^\oy_X$. 
Comme l'anneau $\ocB_{\rho(\oy\rightsquigarrow \ox)}$ est local (\cite{agt} III.10.10(i)), d'idéal maximal contenant $p$, 
\begin{equation}\label{TPCF10c}
(\frac{p}{b^{p-1}})^pz^{p-1}-(1-p)^{p-1}
\end{equation} 
est une unité de $\ocB_{\rho(\oy\rightsquigarrow \ox)}$. 
Il existe alors un objet $(U,\fp\colon \ox\rightarrow U)$ de $\fW_\ox$ tel que l'élément \eqref{TPCF10c} soit une unité de $\oR^\oy_U$. 
Il suffit de montrer que l'extension 
\begin{equation}\label{TPCF10d}
\oR^\oy_U[T]/(T^p-b^{p-1}T-z)
\end{equation}
est étale au-dessus de $\oR^\oy_U[\frac 1 p]$, ou encore que pour tout corps $L$ et tout homomorphisme $u\colon \oR^\oy_U[\frac 1 p]\rightarrow L$, 
les équations $T^p-b^{p-1} T -u(z)$ et $p T^{p-1}-b^{p-1}$ n'ont aucune solution commune dans $L$. 
Si une telle solution $y\in L$ existait, on aurait $y^{p-1}=\frac{b^{p-1}}{p}$ et $u(z)=y^p-b^{p-1} y\in L$, et par suite,
\begin{equation}\label{TPCF10e}
(\frac{p}{b^{p-1}})^pu(z)^{p-1}=(1-p)^{p-1},
\end{equation}
ce qui est absurde puisque $(\frac{p}{b^{p-1}})^pz^{p-1}-(1-p)^{p-1}$ est une unité de $\oR^\oy_U$.

Considérons ensuite, pour tout $t\in \ocB_{\rho(\oy\rightsquigarrow \ox)}$, l'équation
\begin{equation}\label{TPCF10f}
t^p-b^{p-1}t\in p\ocB_{\rho(\oy\rightsquigarrow \ox)},
\end{equation}
qui est équivalente, puisque $\ocB_{\rho(\oy\rightsquigarrow \ox)}$ est $\co_\oK$-plat, à l'équation
\begin{equation}\label{TPCF10g}
(\frac{t}{b})^p-\frac{t}{b}\in \frac{p}{b^p}\ocB_{\rho(\oy\rightsquigarrow \ox)}.
\end{equation}
Comme l'anneau $\ocB_{\rho(\oy\rightsquigarrow \ox)}$ est normal et strictement local (\cite{agt} III.10.10(i)), d'idéal maximal contenant $p$, 
on en déduit d'abord que $\frac t b\in \ocB_{\rho(\oy\rightsquigarrow \ox)}$, puis que 
\begin{equation}
t\in b\cdot (\mZ_p+\frac{p}{b^{p}}\ocB_{\rho(\oy\rightsquigarrow \ox)}).
\end{equation}

Il résulte de ce qui précède et de \ref{TFA7}(ii) que la suite \eqref{TPCF10a} est exacte.

\begin{cor}\label{TPCF11}
Pour tout élément non-nul $b$ de $\co_{\oK^\flat}$, la suite  
\begin{equation}\label{TPCF11a}
\xymatrix{
0\ar[r]&{\mF_p}\ar[r]^-(0.4){\cdot b}&{\ocB^\wp}\ar[rr]^{\phi-b^{p-1}\id}&&{\ocB^\wp}\ar[r]&0}
\end{equation}
est exacte dans $\cP_\AR$.
\end{cor}

Soit $(b_n)_{n\geq 0}\in (\co_C)^{\mN}$ l'image inverse de $b$ par l'isomorphisme \eqref{TPCF5b}. 
Il existe $n_0\geq 0$ tel que pour tout $n\geq n_0$, $b_n^p\not\in p\co_C$. En vertu de \ref{TPCF10}, la suite 
\begin{equation}\label{TPCF11b}
\xymatrix{
0\ar[r]&{\mF_p \cdot b_n\oplus (\frac{p}{b_n^{p-1}})\cdot \ocB_1}\ar[r]&{\ocB_1}\ar[rr]^{\phi-b_n^{p-1}\id}&&{\ocB_1}\ar[r]&0}
\end{equation}
est exacte. Comme $v(b_n)<\frac 1 p$, on a $pv(\frac{p}{b_n^{p-1}})>1$, d'où la proposition.

\begin{prop}\label{TPCF12}
Supposons $X$ propre sur $S$, et soient  $\mL$ un $\mF_p$-module de type fini de $\oX^\circ_\fet$, $q$ un entier $\geq 0$, 
$b$ un élément non-nul de $\co_{\oK^\flat}$. 
Posons $\cL=\delta^*(\beta^*(\mL))$, où $\beta$ est le morphisme \eqref{TFA6b} et $\delta$ est le plongement \eqref{TFA66a},
et notons $\uphi$  l'endomorphisme $\co_{\oK^\flat}$-semi-linéaire de 
\[
\rH^q(\tE_s^{\mN^\circ},\lambda^*(\cL)\otimes_{\mF_p}\ocB^\wp)
\]
induit par l'endomorphisme de Frobenius $\phi$ de $\ocB^\wp$. Alors,
\begin{itemize}
\item[{\rm (i)}] La suite  
\begin{equation}\label{TPCF12a}
\xymatrix{ 
0\ar[r]&{\rH^q(\tE_s,\cL)}\ar[r]^-(0.5){\cdot b}&{\rH^q(\tE_s^{\mN^\circ},\lambda^*(\cL)\otimes_{\mF_p}\ocB^\wp)}\ar[rr]^-(0.5){\uphi-b^{p-1}\id}&&
{\rH^q(\tE_s^{\mN^\circ},\lambda^*(\cL)\otimes_{\mF_p}\ocB^\wp)}\ar[r]& 0}
\end{equation}
est exacte.
\item[{\rm (ii)}] Le morphisme canonique
\begin{equation}\label{TPCF12b}
\ker(\uphi-b^{p-1}\id|\rH^q(\tE_s^{\mN^\circ},\lambda^*(\cL)\otimes_{\mF_p}\ocB^\wp))\rightarrow 
\ker(\uphi-b^{p-1}\id|\rH^q(\tE_s^{\mN^\circ},\lambda^*(\cL)\otimes_{\mF_p}\ocB^\wp)\otimes_{\co_{\oK^\flat}}\oK^\flat)
\end{equation}
est un isomorphisme.
\end{itemize}
\end{prop}
En effet, d'après (\cite{agt} III.7.11) et (\cite{jannsen} 1.15), la suite \eqref{TPCF11a} induit une suite exacte longue de cohomologie
\begin{equation}
\xymatrix{
\rM^{q-1}\ar[r]&{\rH^q(\tE_s,\cL)}\ar[r]&{\rM^q}\ar[rr]^{\uphi-b^{p-1}\id}&&{\rM^q}\ar[r]&{\rH^{q+1}(\tE_s,\cL)}},
\end{equation}
où on a posé $\rM^j=\rH^j(\tE_s^{\mN^\circ},\lambda^*(\cL)\otimes_{\mF_p}\ocB^\wp)$ pour $j\geq 0$ et $\rM^{-1}=0$. 
En vertu de \ref{TPCF9}, il existe un entier $r\geq 0$ et un morphisme $\co_{\oK^\flat}$-linéaire
\begin{equation}\label{TPCF12c}
u\colon \fm_{\oK^\flat}\otimes_{\co_{\oK^\flat}}\rM^q\rightarrow \co_{\oK^\flat}^r 
\end{equation}
compatible à $\uphi$ et à l'endomorphisme de Frobenius $\varphi$ de $\co_{\oK^\flat}$, et qui est un $\alpha$-isomorphisme. 
Comme $\oK^\flat$ est algébriquement clos, pour tout $c\in \co_{\oK^\flat}$, l'endomorphisme $\varphi-c^{p-1}\id$ de $\co_{\oK^\flat}$ est surjectif.
Pour tous $y\in \rM^q$ et $t_1\in \fm_{\oK^\flat}$, il existe $z\in \co_{\oK^\flat}^r$ tel que 
\begin{equation}\label{TPCF12d}
u(t_1^p\otimes y)=\varphi(z)-c^{p-1}z,
\end{equation}
où on a encore noté $\varphi$ l'endomorphisme de $\co_{\oK^\flat}^r$ induit par l'endomorphisme $\varphi$ de $\co_{\oK^\flat}$. 
Pour tout $t_2\in \fm_{\oK^\flat}$, il existe $x\in \rM^q$ et $t\in \fm_{\oK^\flat}$ tels que $u(t\otimes x)=t_2z$. On en déduit que 
\begin{equation}\label{TPCF12e}
u((t_1t_2)^p\otimes y)=u(t^p\otimes \uphi(x)-(ct_2)^{p-1} t\otimes x)
\end{equation}
Pour tout $t_3\in \fm_{\oK^\flat}$, on a donc 
\begin{equation}\label{TPCF12f}
(t_1t_2t_3)^p y=\uphi(tt_3x)-(ct_2t_3)^{p-1} (tt_3 x).
\end{equation}
Par ailleurs, pour tout élément non nul $b'$ de $\co_{\oK^\flat}$ on a un diagramme commutatif
\begin{equation}\label{TPCF12g}
\xymatrix{
{\ker(\uphi-b^{p-1}\id|\rM^q)}\ar[r]\ar[d]_{\cdot b'}&{\rM^q}\ar[rr]^{\uphi-b^{p-1}\id}\ar[d]_{\cdot b'}&&{\rM^q}\ar[r]\ar[d]^{\cdot b'^p}&
{\coker(\uphi-b^{p-1}\id|\rM^q)}\ar[d]^{\cdot b'^p}\\
{\ker(\uphi-(bb')^{p-1}\id|\rM^q)}\ar[r]&{\rM^q}\ar[rr]^{\uphi-(bb')^{p-1}\id}&&{\rM^q}\ar[r]&{\coker(\uphi-(bb')^{p-1}\id|\rM^q)}}
\end{equation}
On déduit de ce qui précède (en prenant $c=bt_1 $ et $b'=t_1t_2t_3$) 
que pour tout $b'\in \fm_{\oK^\flat}$, le dernier morphisme vertical du diagramme \eqref{TPCF12g} est nul.

D'après \ref{TPCF11}, pour tout $b'\in \fm_{\oK^\flat}$, on a un diagramme commutatif à lignes exactes
\begin{equation}\label{TPCF12h}
\xymatrix{
0\ar[r]&{\mF_p}\ar@{=}[d]\ar[r]^-(0.5){\cdot b}&{\ocB^\wp}\ar[d]^{\cdot b'}\ar[rr]^-(0.5){\phi-b^{p-1}\id}&&{\ocB^\wp}\ar[d]^{\cdot b'^p}\ar[r]&0\\
0\ar[r]&{\mF_p}\ar[r]^-(0.5){\cdot bb'}&{\ocB^\wp}\ar[rr]^-(0.5){\phi-(bb')^{p-1}\id}&&{\ocB^\wp}\ar[r]&0}
\end{equation}
Il induit un diagramme commutatif à lignes exactes
\begin{equation}\label{TPCF12i}
\xymatrix{
0\ar[r]&{\coker(\uphi'-b^{p-1}\id|\rM^{q-1})}\ar[r]\ar[d]_{\cdot b'^p}&{\rH^q(\tE_s,\cL)}\ar@{=}[d]\ar[r]&{\ker(\uphi-b^{p-1}\id|\rM^q)}\ar[d]^{\cdot b'}\ar[r]&0\\
0\ar[r]&{\coker(\uphi'-(bb')^{p-1}\id|\rM^{q-1})}\ar[r]&{\rH^q(\tE_s,\cL)}\ar[r]&{\ker(\uphi-(bb')^{p-1}\id|\rM^q)}\ar[r]&0}
\end{equation}
où on a noté $\uphi'$ l'endomorphisme de $M^{q-1}$ induit par $\phi$. 
Comme le morphisme vertical de gauche est nul, on en déduit que l'endomorphisme $\uphi'-b^{p-1}\id$ de $\rM^{q-1}$ est surjectif, 
et il en est donc de même de l'endomorphisme $\uphi-b^{p-1}\id$ de $\rM^q$; 
d'où l'exactitude de la suite \eqref{TPCF12a}.
On en déduit  aussi que la flèche verticale de droite est un isomorphisme~:
\begin{equation}\label{TPCF12j}
\cdot b'\colon \ker(\uphi-b^{p-1}\id|\rM^q)\stackrel{\sim}{\rightarrow}\ker(\uphi-(bb')^{p-1}\id|\rM^q).
\end{equation}
Il s'ensuit que l'application \eqref{TPCF12b} est injective. 

Soit $x\in \rM^q\otimes_{\co_{\oK^\flat}}\oK^\flat$ tel que 
$\uphi(x)=b^{p-1}x$. Il existe $z\in \rM^q$ et $b'\in \fm_{\oK^\flat}$ tels que $z=b'x\in \rM^q\otimes_{\co_{\oK^\flat}}\oK^\flat$.
Par suite, $\uphi(z)-(bb')^{p-1}z$ est de torsion. Quitte à remplacer $b'$ par un multiple, on peut supposer que $\uphi(z)=(bb')^{p-1}z$.
D'après l'isomorphisme \eqref{TPCF12j}, il existe $x'\in \ker(\uphi-b^{p-1}\id|\rM^q)$ tel que $z=b'x'$. On a alors $x=x'$ dans 
$\rM^q\otimes_{\co_{\oK^\flat}}\oK^\flat$; d'où la surjectivité de l'application \eqref{TPCF12b}. 

\begin{cor}\label{TPCF13}
Supposons $X$ propre sur $S$, et soient  $\mL$ un $\mF_p$-module de type fini de $\oX^\circ_\fet$, $\cL=\delta^*(\beta^*(\mL))$.
Alors, pour tout entier $q\geq 0$, $\rH^q(\tE_s,\cL)$ est un $\mF_p$-espace vectoriel de dimension finie,
et le morphisme canonique
\begin{equation}\label{TPCF13a}
\rH^q(\tE_s,\cL)\otimes_{\mF_p}\co_{\oK^\flat}\rightarrow \rH^q(\tE_s^{\mN^\circ},\lambda^*(\cL)\otimes_{\mF_p}\ocB^\wp)
\end{equation}
est un $\alpha$-isomorphisme.
\end{cor}

En effet, en vertu de \ref{TPCF12} et avec les notations de sa preuve, on a des isomorphismes canoniques 
\begin{equation}
\rH^q(\tE_s,\cL)\stackrel{\sim}{\rightarrow}\ker(\uphi-\id |\rM^q)\stackrel{\sim}{\rightarrow}
\ker(\uphi-\id |\rM^q\otimes_{\co_{\oK^\flat}}\oK^\flat).
\end{equation}
D'après \ref{TPCF9}, $\rM^q\otimes_{\co_{\oK^\flat}}\oK^\flat$ est un $\oK^\flat$-espace vectoriel de dimension finie,
et $\uphi$ est un isomorphisme $\co_{\oK^\flat}$-semi-linéaire de $\rM^q$. 
Il résulte alors de (\cite{katz} 4.1.1) que $\rH^q(\tE_s,\cL)$ est un $\mF_p$-espace vectoriel de dimension finie, et que le morphisme $\oK^\flat$-linéaire
\begin{equation}
\rH^q(\tE_s,\cL)\otimes_{\mF_p}\oK^\flat\rightarrow \rM^q\otimes_{\co_{\oK^\flat}}\oK^\flat
\end{equation}
est bijectif. Par suite, le morphisme \eqref{TPCF13a} est injectif. Montrons qu'il est $\alpha$-surjectif.
Considérons le morphisme $u\colon \fm_{\oK^\flat}\otimes_{\co_{\oK^\flat}}\rM^q\rightarrow \co_{\oK^\flat}^r$ \eqref{TPCF12c} 
et notons $e_1,\dots,e_r$ la base canonique de $\co_{\oK^\flat}^r$. Comme $u$ est $\alpha$-injectif, 
il suffit de montrer que pour tout $1\leq i\leq r$ et tout $\gamma\in \fm_{\oK^\flat}$, il existe $x\in \rM^q$ tel que 
$\uphi(x)=x$ et $u(\gamma\otimes x)=\gamma e_i$. Pour tout $t\in \fm_{\oK^\flat}$,  
il existe $y\in \rM^q$ et $\beta\in \fm_{\oK^\flat}$ tels que $u(\beta\otimes y)=te_i$. 
On a $u(\beta^p\otimes \uphi(y))=t^pe_i$ et donc $u(\beta^p\otimes \uphi(y)-t^{p-1}\beta\otimes y)=0$. Par suite, pour 
tout $t'\in \fm_{\oK^\flat}$, on a 
\begin{equation}
\uphi(t'\beta y)=(tt')^{p-1}t'\beta y.
\end{equation}
Compte tenu de l'isomorphisme \eqref{TPCF12j} (pour $b=1$ et $b'=tt'$), il existe $x\in \rM^q$ tel que $\uphi(x)=x$ et $t'\beta y=tt'x$. 
Par suite, pour tout $t''\in \fm_{\oK^\flat}$, on a $t't''\beta\otimes y=tt't''\otimes x$ dans $\fm_{\oK^\flat}\otimes_{\co_{\oK^\flat}}\rM^q$
et donc $u(tt't''\otimes x)=tt't''e_i$, d'où l'assertion recherchée. 

\begin{cor}\label{TPCF14}
Supposons $X$ propre sur $S$, et soient  $n$ un entier $\geq 1$, $\mL$ un $(\mZ/p^n\mZ)$-module localement libre de type fini de 
$\oX^\circ_\fet$, $\cL=\delta^*(\beta^*(\mL))$. Alors, pour tout entier $q\geq 0$, le morphisme canonique
\begin{equation}\label{TPCF14a}
\rH^q(\tE_s,\cL)\otimes_{\mZ_p}\co_\oK\rightarrow \rH^q(\tE_s,\cL\otimes_{\mZ_p}\ocB_n)
\end{equation}
est un $\alpha$-isomorphisme.
\end{cor}

En effet, par dévissage, on peut se réduire au cas où $n=1$. 
La proposition résulte alors de \ref{TPCF13}, \ref{TPCF90} et \eqref{TPCF15b}

\section{Acyclicité locale}\label{acycloc}

\subsection{}\label{acycloc1}
On rappelle que l'on dispose des morphismes canoniques
\begin{eqnarray}
\psi\colon \oX^\circ_\et\rightarrow \tE,\label{acycloc1a}\\
\rho\colon X_\et\gtimes_{X_\et}\oX^\circ_\et\rightarrow \tE,\label{acycloc1b}
\end{eqnarray}
\eqref{tf1m} et \eqref{tf3b}, et du plongement canonique $\delta\colon \tE_s\rightarrow \tE$ \eqref{TFA66a}.

\begin{prop}\label{acycloc2}
Pour tout faisceau abélien de torsion, localement constant et constructible $F$ de $\oX^\circ_\et$ et 
tout entier $q\geq 1$, on a $\rR^q\psi_*(F)=0$.
\end{prop}

En effet, d'après (\cite{sga4} V 5.1), $\rR^q\psi_*(F)$ est le faisceau de $\tE$ associé au préfaisceau sur $E$ défini par 
\begin{equation}
(V\rightarrow U)\mapsto \rH^q(V,F).
\end{equation}
Soient $(V\rightarrow U)$ un objet de $E$, $\xi\in \rH^q(V,F)$, $\ox$ un point géométrique de $U$. 
D'après \ref{cad9}, il existe un voisinage étale $U_1$ de $\ox$ dans $X$
et un objet $(V_1\rightarrow U_1)$ de $E$ au-dessus de $(V\rightarrow U)$ tels que le morphisme canonique
$V_1\rightarrow V\times_{\oU^\circ}\oU^\circ_1$ soit surjectif et que l'image canonique de $\xi$ dans $\rH^q(V_1,F)$ soit nulle. 
Par suite, il existe un recouvrement 
$((V_\alpha\rightarrow U_\alpha)\rightarrow (V\rightarrow U))_{\alpha\in \Sigma}$ de $E$ pour la topologie co-évanescente \eqref{tf1} 
tel que pour tout $\alpha\in \Sigma$, l'image canonique de $\xi$ dans $\rH^q(V_\alpha,F)$ soit nulle; d'où la proposition. 

\begin{prop}
\begin{itemize}
\item[{\rm (i)}] Le morphisme d'adjonction $\id\rightarrow \rho_*\rho^*$ induit un isomorphisme $\delta_*\stackrel{\sim}{\rightarrow} \rho_*\rho^*\delta_*$;
en particulier, le foncteur composé
\begin{equation}
\rho^* \delta_*\colon \tE_s\rightarrow X_\et\gtimes_{X_\et}\oX^\circ_\et
\end{equation}
est pleinement fidèle. 
\item[{\rm (ii)}] Pour tout faisceau abélien $F$ de $\tE_s$ et tout entier $q\geq 1$, on a $\rR^q\rho_*(\rho^*(\delta_*(F)))=0$. 
\end{itemize}
\end{prop}

Soient $\ox$ un point géométrique de $X$ au-dessus de $s$,  $X'$ le localisé strict de $X$ en $\ox$. 
Notons $\rho_{\oX'^\circ}\colon \oX'^\circ_\et\rightarrow \oX'^\circ_\fet$ et 
$\varphi_\ox\colon \tE\rightarrow \oX'^\circ_\fet$ les foncteurs canoniques \eqref{notconv10a} et \eqref{TFA14a}. 
D'après (\cite{agt} III.3.7), $\oX'$ est normal et strictement local (et en particulier intègre). 
La preuve de \ref{tf12} montre alors que le morphisme 
d'adjonction $\id\rightarrow \rho_*\rho^*$ induit un isomorphisme $\varphi_\ox \stackrel{\sim}{\rightarrow} \varphi_\ox \rho_*\rho^*$. 
En vertu de \ref{cad10}, $\oX'^\circ$ est un schéma $K(\pi,1)$. 
La preuve de \ref{tf13} établit alors que pour tout faisceau abélien $G$ de $\tE$ et tout entier $q\geq 1$, on a $\varphi_\ox(\rR^q\rho_*(\rho^*(G)))=0$. 

Par ailleurs, pour tout faisceau $F$ de $\tE_s$, $\rho_*(\rho^*(\delta_*(F)))|\sigma^*(X_\eta)$ est l'objet final de $\tE_{/\sigma^*(X_\eta)}$. 
Pour tout faisceau abélien $G$ de $\tE_s$ et tout entier $q\geq 0$, on a $\rR^q\rho_*(\rho^*(\delta_*(G)))|\sigma^*(X_\eta)=0$. La proposition s'ensuit 
puisque la famille des foncteurs $\varphi_\ox\circ \delta_*$ lorsque $\ox$ décrit l'ensemble des points géométriques de $X$ au-dessus de $s$, est conservative en vertu de \ref{TFA7}(ii) et (\cite{agt} VI.10.31).

\section{Le torseur des déformations}\label{deform}

\subsection{}\label{deform1}
On pose $\coS=\Spec(\co_C)$, que l'on munit de la structure logarithmique $\cM_\coS$ image inverse de $\cM_S$ \eqref{TFA1}. 
Pour tout $S$-schéma $Y$, on pose 
\begin{equation}\label{deform1d}
\coY=Y\times_S\coS. 
\end{equation}

On fixe une suite $(p_n)_{n\geq 0}$ d'éléments de $\co_\oK$ telle que $p_0=p$ et $p_{n+1}^p=p_n$ pour tout $n\geq 0$, 
et on note $\varpi$ l'élément associé de $\co_{\oK^\flat}$ \eqref{TPCF5}. On pose  
\begin{equation}\label{deform1a}
\xi=[\varpi]-p \in \rW(\co_{\oK^\flat}),
\end{equation}
où $[\ ]$ est le représentant multiplicatif. D'après \ref{eip4} (avec les notations de \ref{eipo3}), la suite 
\begin{equation}\label{deform1b}
0\longrightarrow \rW(\co_{\oK^\flat})\stackrel{\cdot \xi}{\longrightarrow} \rW(\co_{\oK^\flat})
\stackrel{\theta}{\longrightarrow} \co_C \longrightarrow 0
\end{equation}
est exacte. Elle induit une suite exacte \eqref{eipo3e}
\begin{equation}\label{deform1c}
0\longrightarrow \co_C\stackrel{\cdot \xi}{\longrightarrow} \cA_2(\co_\oK)
\stackrel{\theta}{\longrightarrow} \co_C \longrightarrow 0,
\end{equation}
où on a encore noté $\cdot \xi$ le morphisme induit par la multiplication par $\xi$ dans $\cA_2(\co_\oK)$. 
L'idéal $\ker(\theta)$ de $\cA_2(\co_\oK)$ est de carré nul. 
C'est un $\co_C$-module libre de base $\xi$. Il sera noté $\xi\co_C$. 
On observera que contrairement à $\xi$, ce module ne dépend pas du choix de la suite $(p_n)_{n\geq 0}$. 
On note $\xi^{-1}\co_C$ le $\co_C$-module dual de $\xi\co_C$. 
Pour tout $\co_C$-module $M$, on désigne les $\co_C$-modules $M\otimes_{\co_C}(\xi \co_C)$ 
et $M\otimes_{\co_C}(\xi^{-1} \co_C)$ simplement par $\xi M$ et $\xi^{-1} M$, respectivement.

\subsection{}\label{deform2}
Considérons le système projectif de monoïdes multiplicatifs $(\co_\oK)_{n\in \mN}$, 
où les morphismes de transition sont tous égaux à l'élévation à la puissance $p$-ième.
On note $Q_K$ le produit fibré du diagramme d'homomorphismes de monoïdes 
\begin{equation}\label{deform2a}
\xymatrix{
&{\co_K-\{0\}}\ar[d]\\
{\underset{\underset{x\mapsto x^p}{\longleftarrow}}{\lim}\ \co_\oK}\ar[r]&\co_\oK}
\end{equation}
où  la flèche verticale est l'homomorphisme canonique
et la flèche horizontale est la projection sur la première composante ({\em i.e.}, d'indice $0$).
On désigne par $q_K$ l'homomorphisme composé
\begin{equation}\label{deform2b}
Q_K\longrightarrow
\underset{\underset{x\mapsto x^p}{\longleftarrow}}{\lim}\ \co_\oK \longrightarrow \co_{\oK^\flat} \stackrel{[\ ]}{\longrightarrow} 
\rW(\co_{\oK^\flat}) \rightarrow \cA_2(\co_\oK),
\end{equation} 
où $[\ ]$ est le représentant multiplicatif et les autres flèches sont les morphismes canoniques. 
Il résulte aussitôt des définitions que le diagramme 
\begin{equation}\label{deform2c}
\xymatrix{
Q_K\ar[r]\ar[d]_{q_K}&{\co_K-\{0\}}\ar[d]\\
{\cA_2(\co_\oK)}\ar[r]^-(0.4)\theta&{\co_C}}
\end{equation}
où les flèches non libellées sont les morphismes canoniques, est commutatif.

On pose 
\begin{equation}\label{deform2d}
\cA_2(\oS)=\Spec(\cA_2(\co_\oK))
\end{equation} 
que l'on munit de la structure logarithmique $\cM_{\cA_2(\oS)}$
associée à la structure pré-logarithmique définie par l'homomorphisme $q_K$. D'après \ref{eip6}, 
le schéma logarithmique $(\cA_2(\oS),\cM_{\cA_2(\oS)})$ est fin et saturé, 
et $\theta$ induit une immersion fermée exacte 
\begin{equation}\label{deform2e}
i_K\colon (\coS,\cM_\coS)\rightarrow (\cA_2(\oS),\cM_{\cA_2(\oS)}).
\end{equation}

Le groupe de Galois $G_K$ agit naturellement sur $\rW(\co_{\oK^\flat})$ par des automorphismes d'anneaux,
et l'homomorphisme $\theta$ est $G_K$-équivariant. On en déduit une action de $G_K$ 
sur $\cA_2(\co_{\oK})$ par des automorphismes d'anneaux tel que l'homomorphisme $\theta$ \eqref{deform1c} soit 
$G_K$-équivariant. Par ailleurs, $G_K$ agit naturellement sur le monoïde $Q_K$,
et l'homomorphisme $q_K$ est $G_K$-équivariant. 
On en déduit une action à gauche de $G_K$ sur le schéma logarithmique $(\cA_2(\oS),\cM_{\cA_2(\oS)})$. 
Le morphisme $i_K$ est $G_K$-équivariant.

\subsection{}\label{THT1}
On suppose dans la suite de cette section que le schéma $X=\Spec(R)$ est affine et connexe \eqref{TFA3},  que $X_s$ est non-vide, 
que le morphisme $f\colon (X,\cM_X)\rightarrow (S,\cM_S)$ admet une carte adéquate \eqref{cad1} et qu'il existe une carte fine et saturée 
$M\rightarrow \Gamma(X,\cM_X)$ pour $(X,\cM_X)$ induisant un isomorphisme 
\begin{equation}\label{THT1a}
M\stackrel{\sim}{\rightarrow} \Gamma(X,\cM_X)/\Gamma(X,\co^\times_X).
\end{equation}
Cette dernière est a priori indépendante de la carte adéquate. On pose \eqref{TFA3b}
\begin{equation}\label{THT1c}
\tOmega^1_{R/\co_K}=\tOmega^1_{X/S}(X).
\end{equation}
Soit $\oy$ un point géométrique de $\oX^\circ$. Le schéma $\oX$ étant localement irréductible d'après \ref{cad4}(iii),  
il est la somme des schémas induits sur ses composantes irréductibles. On note $\oX^\star$
la composante irréductible de $\oX$ contenant $\oy$. 
De même, $\oX^\circ$ est la somme des schémas induits sur ses composantes irréductibles
et $\oX^{\star \circ}=\oX^\star\times_{X}X^\circ$ est la composante irréductible de $\oX^\circ$ contenant $\oy$. 
On pose 
\begin{equation}\label{THT1b}
R_1=\Gamma(\oX^\star,\co_\oX).
\end{equation}
On désigne par $\Delta$ le groupe profini $\pi_1(\oX^{\star \circ},\oy)$ 
et pour alléger les notations, par $\oR$ la représentation discrète $\oR^\oy_X$ de $\Delta$ définie dans \eqref{TFA9b} 
et par $\hoR$ son séparé complété $p$-adique.  On notera que $\oR$ est intègre et normal \eqref{TFA10}. On pose 
\begin{equation}\label{THT1e}
Y=\Spec(\oR)\ \ \ {\rm et} \ \ \ \hY=\Spec(\hoR)
\end{equation}
que l'on munit des structures logarithmiques images inverses de $\cM_X$, 
notées respectivement $\cM_Y$ et $\cM_\hY$. 
Les actions de $\Delta$ sur $\oR$ et $\hoR$ induisent des actions à gauche sur 
les schémas logarithmiques $(Y,\cM_Y)$ et $(\hY,\cM_\hY)$, respectivement. 

Munissant $(\coX,\cM_\coX)$ de l'action triviale de $\Delta$ \eqref{deform3c}, on a un morphisme canonique $\Delta$-équivariant 
\begin{equation}\label{THT1d}
(\hY,\cM_\hY)\rightarrow (\coX,\cM_\coX).
\end{equation}

\subsection{}\label{deform3}
On munit  $\coX=X\times_S\coS$ de la structure logarithmique $\cM_\coX$ image inverse de $\cM_X$. 
On a alors un isomorphisme canonique
\begin{equation}\label{deform3c}
(\coX,\cM_{\coX})\stackrel{\sim}{\rightarrow}(X,\cM_X)\times_{(S,\cM_S)}(\coS,\cM_\coS),
\end{equation} 
le produit étant indifféremment pris dans la catégorie des schémas logarithmiques ou 
dans celle des schémas logarithmiques fins. 
Une {\em $(\cA_2(\oS),\cM_{\cA_2(\oS)})$-déformation lisse} de $(\coX,\cM_{\coX})$
est la donnée d'un morphisme lisse de schémas logarithmiques fins $(\tX,\cM_\tX)\rightarrow (\cA_2(\oS), \cM_{\cA_2(\oS)})$
et d'un $(\coS,\cM_{\coS})$-isomorphisme 
\begin{equation}\label{deform3d}
(\coX,\cM_{\coX})\stackrel{\sim}{\rightarrow}
(\tX,\cM_\tX)\times_{(\cA_2(\oS), \cM_{\cA_2(\oS)})}(\coS,\cM_{\coS}).
\end{equation}
Comme $(\coX,\cM_{\coX})\rightarrow (\coS,\cM_\coS)$ est lisse et que $\coX$ est affine,
une telle déformation existe et est unique à isomorphisme près en vertu de (\cite{kato1} 3.14). 
Son groupe d'automorphismes est isomorphe à 
\begin{equation}\label{deform3e}
\Hom_{\co_{\coX}}(\tOmega^1_{X/S}\otimes_{\co_X}\co_{\coX},\xi \co_{\coX}).
\end{equation}
On se donne une telle déformation $(\tX,\cM_\tX)$ {\em que l'on suppose fixée dans cette section}.

\subsection{}\label{THT2}
D'après \ref{eip4} et \ref{cg67} (avec les notations de \ref{eipo3}), la suite 
\begin{equation}\label{THT2a}
0\longrightarrow \rW(\oR^\flat)\stackrel{\cdot \xi}{\longrightarrow} \rW(\oR^\flat)
\stackrel{\theta}{\longrightarrow} \hoR \longrightarrow 0,
\end{equation}
où $\xi$ est l'élément \eqref{deform1a}, est exacte. Elle induit une suite exacte \eqref{eipo3e}
\begin{equation}\label{THT2b}
0\longrightarrow \hoR\stackrel{\cdot \xi}{\longrightarrow} \cA_2(\oR)
\stackrel{\theta}{\longrightarrow} \hoR \longrightarrow 0,
\end{equation}
où on a encore noté $\cdot \xi$ le morphisme induit par la multiplication par $\xi$ dans $\cA_2(\oR)$. 
L'idéal $\ker(\theta)$ de $\cA_2(\oR)$ est de carré nul. C'est un $\hoR$-module libre de base $\xi$, 
canoniquement isomorphe à $\xi\hoR$ \eqref{deform1}. 

On désigne par $Q$ le monoïde et par $q\colon Q\rightarrow \rW(\oR^\flat)$
l'homomorphisme définis dans \ref{eip5} en prenant pour schéma logarithmique affine $(Y,\cM_Y)$ et pour 
$u\colon \Gamma(X,\cM_X)\rightarrow \Gamma(Y,\cM_Y)$  l'homomorphisme canonique, de sorte qu'on a un diagramme cartésien 
\begin{equation}\label{THT2d}
\xymatrix{
Q\ar[r]\ar[d]&{\Gamma(X,\cM_X)}\ar[d]\\
{\underset{\underset{x\mapsto x^p}{\longleftarrow}}{\lim}\ \oR}\ar[r]&\oR}
\end{equation} 
On note encore $q\colon Q\rightarrow \cA_2(\oR)$
l'homomorphisme induit par $q$. On a un homomorphisme canonique 
$Q_K\rightarrow Q$ \eqref{deform2} qui s'insère dans un diagramme commutatif 
\begin{equation}\label{THT2c}
\xymatrix{
Q_K\ar[r]\ar[d]_{q_K}&Q\ar[d]^{q}\\
{\cA_2(\co_\oK)}\ar[r]&{\cA_2(\oR)}}
\end{equation}

Le groupe $\Delta$ agit naturellement sur $\rW(\oR^\flat)$ par des automorphismes d'anneaux,
et l'homomorphisme $\theta$ \eqref{THT2a} est $\Delta$-équivariant. On en déduit une action de $\Delta$ 
sur $\cA_2(\oR)$ par des automorphismes d'anneaux telle que l'homomorphisme $\theta$ \eqref{THT2b} soit 
$\Delta$-équivariant. Par ailleurs, $\Delta$ agit naturellement sur le monoïde $Q$,
et l'homomorphisme $q$ est $\Delta$-équivariant.  

\subsection{}\label{THT4}
On pose 
\begin{equation}\label{THT4a}
\cA_2(Y)=\Spec(\cA_2(\oR))
\end{equation} 
que l'on munit de la structure logarithmique $\cM_{\cA_2(Y)}$ associée à la structure pré-logarithmique définie 
par l'homomorphisme $q\colon Q\rightarrow \cA_2(\oR)$ \eqref{THT2}. 
D'après \ref{eip6} et les hypothèses de \ref{THT1}, le schéma logarithmique $(\cA_2(Y),\cM_{\cA_2(Y)})$ est fin et saturé, et 
l'homomorphisme $\theta$ induit une immersion fermée exacte 
\begin{equation}\label{THT4b}
i\colon (\hY,\cM_\hY)\rightarrow (\cA_2(Y),\cM_{\cA_2(Y)}).
\end{equation}
Les actions de $\Delta$ sur $\cA_2(\oR)$ et $Q$ induisent une action à gauche sur
le schéma logarithmique $(\cA_2(Y),\cM_{\cA_2(Y)})$. L'immersion fermée $i$ est $\Delta$-équivariante. 

On a un diagramme commutatif canonique \eqref{deform2e}
\begin{equation}\label{THT4c}
\xymatrix{
{(\hY,\cM_\hY)}\ar[r]^-(0.5){i}\ar[d]&{(\cA_2(Y),\cM_{\cA_2(Y)})}\ar[d]\\
{(\coS,\cM_\coS)}\ar[r]^-(0.5){i_K}&{(\cA_2(\oS),\cM_{\cA_2(\oS)})}}
\end{equation}

 \subsection{}\label{THT5}
On pose \eqref{THT1c}
\begin{equation}\label{THT5a}
\rT=\Hom_{\hoR}(\tOmega^1_{R/\co_K}\otimes_R\hoR,\xi\hoR).
\end{equation} 
On identifie le $\hoR$-module dual à $\xi^{-1}\tOmega^1_{R/\co_K}\otimes_R\hoR$ \eqref{deform1}
et on note $\cS$ la $\hoR$-algèbre symétrique associée \eqref{notconv9}
\begin{equation}\label{THT5b}
\cS=\oplus_{n\geq 0}\cS^{n}=\rS_{\hoR}(\xi^{-1}\tOmega^1_{R/\co_K}\otimes_R\hoR).
\end{equation}
On désigne par $\hY_\zar$ le topos de Zariski de $\hY=\Spec(\hoR)$ \eqref{THT1e}, 
par $\trT$ le $\co_\hY$-module associé à $\rT$
et par $\bT$ le $\hY$-fibré vectoriel associé à son dual, autrement dit,  
\begin{equation}\label{THT5c}
\bT=\Spec(\cS).
\end{equation}
 
Soient $U$ un ouvert de Zariski de $\hY$, $\tU$ l'ouvert correspondant de $\cA_2(Y)$ (cf. \ref{THT4}). 
On désigne par $\cL(U)$ 
l'ensemble des morphismes représentés par des flèches pointillées qui complètent  le diagramme canonique 
\begin{equation}\label{THT5d}
\xymatrix{
{(U,\cM_\hY|U)}\ar[r]^-(0.5){i|U}\ar[d]&{(\tU,\cM_{\cA_2(Y)}|\tU)}\ar@{.>}[d]\ar@/^2pc/[dd]\\
{(\coX,\cM_{\coX})}\ar[r]\ar[d]&{(\tX,\cM_\tX)}\ar[d]\\
{(\coS,\cM_\coS)}\ar[r]^-(0.5){i_K}&{(\cA_2(\oS),\cM_{\cA_2(\oS)})}}
\end{equation}
de façon à le laisser commutatif (cf. \ref{deform3}). D'après (\cite{agt} II.5.23),
le foncteur $U\mapsto \cL(U)$ est un $\trT$-torseur de $\hY_\zar$.  
On désigne par $\cF$ le $\hoR$-module des fonctions affines sur $\cL$ (cf. \cite{agt} II.4.9). 
Celui-ci s'insère dans une suite exacte canonique 
\begin{equation}\label{THT5e}
0\rightarrow \hoR\rightarrow \cF\rightarrow \xi^{-1}\tOmega^1_{R/\co_K} \otimes_R \hoR\rightarrow 0.
\end{equation} 
D'après (\cite{illusie1} I 4.3.1.7), cette suite induit pour tout entier $n\geq 1$, une suite exacte \eqref{notconv9}
\begin{equation}\label{THT5f}
0\rightarrow \rS^{n-1}_{\hoR}(\cF)\rightarrow \rS^{n}_{\hoR}(\cF)\rightarrow \rS^n_{\hoR}(\xi^{-1}\tOmega^1_{R/\co_K}
\otimes_R\hoR)\rightarrow 0.
\end{equation}
Les $\hoR$-modules $(\rS^{n}_{\hoR}(\cF))_{n\in \mN}$ forment donc un système inductif filtrant, 
dont la limite inductive 
\begin{equation}\label{THT5g}
\cC=\underset{\underset{n\geq 0}{\longrightarrow}}\lim\ \rS^n_{\hoR}(\cF)
\end{equation}
est naturellement munie d'une structure de $\hoR$-algèbre. D'après (\cite{agt} II.4.10), le $\hY$-schéma 
\begin{equation}\label{THT5h}
\bL=\Spec(\cC)
\end{equation}
est naturellement un $\bT$-fibré principal homogène sur $\hY$ qui représente canoniquement $\cL$. 
On prendra garde que $\cL$, $\cF$, $\cC$ et $\bL$ dépendent de $(\tX,\cM_\tX)$.

\subsection{}\label{THT6}
On munit $\hY$ de l'action naturelle à gauche de $\Delta$~; pour tout $g\in \Delta$, 
l'automorphisme de $\hY$ défini par $g$, que l'on note aussi $g$, est induit par l'automorphisme $g^{-1}$ de $\hoR$. 
On considère $\trT$ comme un $\co_\hY$-module $\Delta$-équivariant au moyen 
de la donnée de descente correspondant au $\hRun$-module 
$\Hom_{\hRun}(\tOmega^1_{R/\co_K}\otimes_R\hRun,\xi\hRun)$ (cf. \cite{agt} II.4.18). Pour tout $g\in \Delta$, on a donc 
un isomorphisme canonique de $\co_\hY$-modules
\begin{equation}\label{THT6a}
\tau_g^\trT\colon \trT\stackrel{\sim}{\rightarrow} g^*(\trT).
\end{equation}
Celui-ci induit un isomorphisme de $\hY$-schémas en groupes
\begin{equation}\label{THT6b}
\tau_g^\bT\colon \bT\stackrel{\sim}{\rightarrow} g^\bullet(\bT),
\end{equation}
où $g^\bullet$ désigne le foncteur de changement de base par l'automorphisme $g$ de $\hY$.
On obtient ainsi une structure $\Delta$-équivariante sur le $\hY$-schéma en groupes $\bT$ (cf.  \cite{agt} II.4.17)
et par suite une action à gauche de $\Delta$ sur $\bT$ compatible avec son action sur $\hY$; 
l'automorphisme de $\bT$ défini par un élément $g$ de $\Delta$ est le composé de $\tau^\bT_g$ 
et de la projection canonique $g^\bullet(\bT)\rightarrow \bT$. 
On en déduit une action de $\Delta$ sur $\cS$ par des automorphismes 
d'anneaux, compatible avec son action sur $\hoR$, que l'on appelle {\em action canonique}.
Cette dernière est concrètement  induite par l'action triviale sur 
$\rS_{\hRun}(\xi^{-1}\tOmega^1_{R/\co_K}\otimes_R\hRun)$.

L'action naturelle à gauche de $\Delta$ sur le schéma logarithmique $(\cA_2(Y),\cM_{\cA_2(Y)})$ \eqref{THT2}
induit sur le $\trT$-torseur $\cL$ une structure $\Delta$-équivariante (cf. \cite{agt} II.4.18), 
autrement dit, elle induit pour tout $g\in \Delta$, un isomorphisme $\tau_g^\trT$-équivariant
\begin{equation}\label{THT6c}
\tau^{\cL}_g\colon \cL\stackrel{\sim}{\rightarrow} g^*(\cL);
\end{equation}
ces isomorphismes étant soumis à des relations de compatibilité (cf. \cite{agt} II.4.16). 
En effet, pour tout ouvert de Zariski $U$ de $\hY$, on prend pour 
\begin{equation}\label{THT6d}
\tau^{\cL}_g(U)\colon \cL(U)\stackrel{\sim}{\rightarrow} \cL(g(U))
\end{equation}
l'isomorphisme défini de la façon suivante. 
Soient $\tU$ l'ouvert de $\cA_2(Y)$ correspondant à $U$, $\mu\in \cL(U)$ que l'on considère comme un morphisme
\begin{equation}\label{THT6e}
\mu\colon (\tU,\cM_{\cA_2(Y)}|\tU)\rightarrow (\tX,\cM_\tX).
\end{equation}
Comme $i$ \eqref{THT4b} et le morphisme \eqref{THT1d} sont $\Delta$-équivariants, le morphisme composé
\begin{equation}\label{THT6f}
(g(\tU),\cM_{\cA_2(Y)}|g(\tU))\stackrel{g^{-1}}{\longrightarrow} (\tU,\cM_{\cA_2(Y)}|\tU)
\stackrel{\mu}{\longrightarrow} (\tX,\cM_\tX)
\end{equation}
prolonge le morphisme canonique $(g(U),\cM_\hY|g(U))\rightarrow (\tX,\cM_\tX)$. 
Il correspond à l'image de $\mu$ par $\tau^{\cL}_g(U)$. On vérifie aussitôt que 
le morphisme $\tau^{\cL}_g$ ainsi défini est un isomorphisme $\tau_g^\trT$-équivariant et que 
ces isomorphismes vérifient les relations de compatibilité requises dans (\cite{agt} II.4.16).

D'après (\cite{agt} II.4.21), les structures $\Delta$-équivariantes sur $\trT$ et $\cL$ induisent une structure
$\Delta$-équi\-variante sur le $\co_\hY$-module associé à $\cF$, ou, ce qui revient au même,
une action $\hoR$-semi-linéaire de $\Delta$ sur $\cF$, telle que les morphismes de la suite \eqref{THT5e} soient 
$\Delta$-équivariants. 
On en déduit sur $\bL$ une structure de $\bT$-fibré principal homogène $\Delta$-équivariant 
sur $\hY$ (cf. \cite{agt} II.4.20). Pour tout $g\in \Delta$, on a donc un isomorphisme $\tau_g^\bT$-équivariant
\begin{equation}\label{THT6g}
\tau_g^{\bL}\colon \bL\stackrel{\sim}{\rightarrow} g^\bullet(\bL).
\end{equation}
Cette structure détermine une action à gauche de $\Delta$ sur $\bL$ compatible avec son action sur $\hY$; 
l'automorphisme de $\bL$ défini par un élément $g$ de $\Delta$ est le composé de $\tau^{\bL}_g$ 
et de la projection canonique $g^\bullet(\bL)\rightarrow \bL$. On obtient ainsi une action de $\Delta$ sur 
$\cC$ par des automorphismes d'anneaux, compatible avec son action sur $\hoR$, 
que l'on appelle {\em action canonique}. Cette dernière est concrètement induite par l'action de $\Delta$ sur $\cF$. 

Pour tout $g\in \Delta$, on désigne par 
\begin{equation}\label{THT6h}
\bL(\hY)\stackrel{\sim}{\rightarrow}\bL(\hY), \ \ \ \mu\mapsto {^g\mu} 
\end{equation}
le composé des isomorphismes 
\begin{eqnarray}
\tau_g^{\bL}\colon \bL(\hY)&\stackrel{\sim}{\rightarrow}& g^\bullet(\bL)(\hY),\\
g^\bullet(\bL)(\hY)&\stackrel{\sim}{\rightarrow}&\bL(\hY), \ \ \ \mu\mapsto \pr\circ \mu\circ g^{-1}, 
\end{eqnarray}
où $g^{-1}$ agit sur $\hY$ et $\pr\colon g^\bullet(\bL)\rightarrow \bL$ est la projection canonique,
de sorte que le diagramme 
\begin{equation}\label{THT6i}
\xymatrix{
{\bL}\ar[r]^g&{\bL}\\
\hY\ar[r]^g\ar[u]^{\mu}&\hY\ar[u]_{{^g\mu}}}
\end{equation}
est commutatif. En particulier, pour tout $\beta\in \cF$, on a 
\begin{equation}\label{THT6j}
(g^{-1}(\beta))(\mu)=g^{-1}(\beta({^g\mu})).
\end{equation}
Par ailleurs, ${^g\mu}$ est défini par le morphisme composé 
\begin{equation}\label{THT6k}
(\cA_2(Y),\cM_{\cA_2(Y)})\stackrel{g^{-1}}{\longrightarrow} (\cA_2(Y),\cM_{\cA_2(Y)})
\stackrel{\mu}{\longrightarrow} (\tX,\cM_\tX).
\end{equation}

\begin{defi}[\cite{agt} II.10.5]\label{THT7}
La $\hoR$-algèbre $\cC$ \eqref{THT5g}, munie de l'action canonique de $\Delta$ 
\eqref{THT6}, est appelée {\em l'algèbre de Higgs-Tate} associée à $(\tX,\cM_\tX)$. 
La $\hoR$-représentation $\cF$ \eqref{THT5e} de $\Delta$ est appelée 
l'{\em extension de Higgs-Tate} associée à $(\tX,\cM_\tX)$.
\end{defi}

D'après (\cite{agt} II.12.4), les actions canoniques de $\Delta$ sur $\cF$ et $\cC$ sont continues pour les topologies $p$-adiques.

\subsection{}\label{THT3}
On a un homomorphisme canonique \eqref{TFA1a}
\begin{equation}\label{THT3a}
\mZ_p(1)\rightarrow (\oR^\flat)^\times.
\end{equation} 
Pour tout $\zeta\in \mZ_p(1)$, on note encore $\zeta$ son image dans $(\oR^\flat)^\times$.  
Comme $\theta([\zeta]-1)=0$, on obtient un homomorphisme de groupes
\begin{equation}\label{THT3b}
\mZ_p(1)\rightarrow \cA_2(\oR),\ \ \ 
\zeta\mapsto\log([\zeta])=[\zeta]-1,
\end{equation}
dont l'image est contenue dans $\ker(\theta)=\xi\hoR$. Celui-ci est clairement $\mZ_p$-linéaire. D'après (\cite{agt} II.9.18),
son image engendre l'idéal $p^{\frac{1}{p-1}}\xi\hoR$ de $\cA_2(\oR)$, et le morphisme $\hoR$-linéaire induit
\begin{equation}\label{THT3c}
\hoR(1)\rightarrow p^{\frac{1}{p-1}}\xi \hoR
\end{equation}
est un isomorphisme. 

L'homomorphisme de monoïdes $\mZ_p(1)\rightarrow (\oR^\flat,\times)$ déduit de \eqref{THT3a} et l'homomorphisme
trivial $\mZ_p(1)\rightarrow \Gamma(X,\cM_X)$  (de valeur $1$) induisent un homomorphisme \eqref{THT2d}
\begin{equation}\label{THT3d}
\mZ_p(1)\rightarrow Q.
\end{equation}

\subsection{}\label{THT8}
On désigne par $\tQ$ le produit fibré du diagramme d'homomorphismes de monoïdes
\begin{equation}\label{THT8a}
\xymatrix{
&&{\Gamma(\tX,\cM_\tX)}\ar[d]\\
Q\ar[r]&{\Gamma(X,\cM_X)}\ar[r]&{\Gamma(\coX,\cM_\coX)}&}
\end{equation}
L'action naturelle de $\Delta$ sur $Q$ et son action triviale sur $\Gamma(\tX,\cM_\tX)$ induise une action sur $\tQ$.
L'homomorphisme canonique $\mZ_p(1)\rightarrow Q$ \eqref{THT3d} et l'homomorphisme
trivial $\mZ_p(1)\rightarrow \Gamma(\tX,\cM_\tX)$  (de valeur $1$) induisent un homomorphisme 
\begin{equation}\label{THT8b}
\mZ_p(1)\rightarrow \tQ.
\end{equation}

Soit $T'$ un élément de $\tQ$, d'images canoniques $T\in Q$, $t\in \Gamma(X,\cM_X)$ et $t'\in \Gamma(\tX,\cM_\tX)$. 
Tout élément $\mu\in \cL(\hY)$ détermine un morphisme de schémas logarithmiques que l'on note encore 
\begin{equation}\label{THT8c}
\mu\colon (\cA_2(Y),\cM_{\cA_2(Y)})\rightarrow (\tX,\cM_{\tX}),
\end{equation}
et en particulier un homomorphisme de monoïdes $\mu^\flat\colon \mu^{-1}(\cM_{\tX})\rightarrow \cM_{\cA_2(Y)}$.
D'après (\cite{ogus} IV 2.1.2),  il existe un et un unique élément $b_\mu \in  \xi\hoR$ tel que 
\begin{equation}\label{THT8d}
\mu^{\flat}(\mu^{-1}(t'))=(1+b_\mu) T \in \Gamma(\cA_2(Y),\cM_{\cA_2(Y)}).
\end{equation}
L'application $\mu\mapsto b_\mu$ est une fonction affine sur $\cL$ à valeurs dans $\xi \hoR$ 
de terme linéaire $d\log(t)\in \tOmega^1_{R/\co_K}\otimes_R\hoR$ (cf. \cite{agt} II.4.7). Elle définit donc naturellement 
un élément de $\xi\cF$, que l'on note $d\log(T')\in \xi\cF$, au-dessus de $d\log(t)$ \eqref{THT5e}.
On voit aussitôt que l'application 
\begin{equation}\label{THT8e}
d\log\colon \tQ\rightarrow \xi\cF
\end{equation}
ainsi définie est un homomorphisme de monoïdes. Il est $\Delta$-équivariant d'après \eqref{THT6j}.
Le diagramme
\begin{equation}\label{THT8f}
\xymatrix{
&\mZ_p(1)\ar[r]\ar[d]_{-\log([\ ])}&{\tQ}\ar[r]\ar[d]^{d\log}&{\Gamma(X,\cM_X)}\ar[d]^{d\log}&\\
0\ar[r]&{\xi\hoR}\ar[r]&{\xi\cF}\ar[r]&{\tOmega^1_{R/\co_K}\otimes_R\hoR}\ar[r]&0}
\end{equation}
est clairement commutatif.

\section{Théorie de Kummer}\label{TFSLA}
 
\subsection{}\label{TFA4}
Pour toute extension finie $L$ de $K$, on note $\co_L$ la fermeture intégrale
de $\co_K$ dans $L$ et on pose $S_L=\Spec(\co_L)$ que l'on munit de la structure logarithmique $\cM_{S_L}$
définie par son point fermé. On a un morphisme canonique $(S_L,\cM_{S_L})\rightarrow (S,\cM_S)$. On pose 
\begin{equation}\label{TFA4a}
(X_L,\cM_{X_L})=(X,\cM_X)\times_{(S,\cM_S)}(S_L,\cM_{S_L}),
\end{equation}
le produit étant pris dans la catégorie des schémas logarithmiques, de sorte que l'on a $X_L=X\times_SS_L$.
On note $f_L\colon (X_L,\cM_{X_L})\rightarrow (S_L,\cM_{S_L})$ la projection canonique. 
On voit aussitôt que $X_L^\circ=X_L\times_XX^\circ$ est le sous-schéma ouvert maximal de $X_L$
où la structure logarithmique $\cM_{X_L}$ est triviale. 
Le morphisme $f_L$ est adéquat d'après (\cite{agt} III.4.8). En vertu de (\cite{agt} III.4.2(iv)), on a donc
\begin{equation}\label{TFA4b}
\cM_{X_L}\stackrel{\sim}{\rightarrow}(j_{X_L})_*(\co_{X_L^\circ}^\times)\cap \co_{X_L}.
\end{equation} 
En particulier, l'homomorphisme canonique $\cM_{X_L}\rightarrow \co_{X_L}$ est un monomorphisme. 
On a un isomorphisme canonique \eqref{TFA3b}
\begin{equation}\label{TFA4c}
\Omega^1_{(X_L,\cM_{X_L})/(S_L,\cM_{S_L})}\stackrel{\sim}{\rightarrow}\tOmega^1_{X/S}\otimes_{\co_X}\co_{X_L}.
\end{equation}
Par suite, la dérivation logarithmique universelle induit un homomorphisme 
\begin{equation}\label{TFA4d}
d\log \colon \cM_{X_L}\rightarrow \tOmega^1_{X/S}\otimes_{\co_X}\co_{X_L}.
\end{equation}

Pour tout $K$-homomorphisme $L\rightarrow L'$ entre extensions finies de $K$, on a un diagramme cartésien canonique
dans la catégorie des schémas logarithmiques 
\begin{equation}\label{TFA4e}
\xymatrix{
{(X_{L'},\cM_{X_{L'}})}\ar[r]\ar[d]&{(X_L,\cM_{X_L})}\ar[d]\\
{(S_{L'},\cM_{S_{L'}})}\ar[r]&{(S_L,\cM_{S_L})}}
\end{equation}
Pour toute sous-$K$-extension finie $L$ de $\oK$,
on note $\hbar_L\colon \oX\rightarrow X_L$ le morphisme canonique. On pose 
\begin{equation}\label{TFA4f}
\cM_{\oX}=\underset{\underset{L\subset \oK}{\longrightarrow}}{\lim}\ h_L^*(\cM_{X_L}),
\end{equation}
la limite étant prise sur la catégorie filtrante des sous-$K$-extensions finies de $\oK$. C'est un monoïde intègre de $\oX_\et$
(\cite{ogus} I 1.3.6). On notera que l'homomorphisme canonique de $\oX_\zar$ 
\begin{equation}\label{TFA4h}
\underset{\underset{L\subset \oK}{\longrightarrow}}{\lim}\ h_L^{-1}(\co_{X_L})\rightarrow \co_\oX
\end{equation}
est un isomorphisme (\cite{ega4} 8.2.12). Il résulte de (\cite{ega1n} 0.6.1.6) et (\cite{raynaud1} I §~3 prop.~1) que 
l'homomorphisme analogue de \eqref{TFA4h} dans $\oX_\et$  est aussi un isomorphisme. 
On en déduit un monomorphisme de $\oX_\et$
\begin{equation}\label{TFA4g}
\alpha\colon \cM_\oX\rightarrow \co_\oX.
\end{equation}
C'est une structure logarithmique sur $\oX$. 

On pose 
\begin{equation}\label{TFA4j}
\tOmega^1_{\oX/\oS}=\tOmega^1_{X/S}\otimes_{\co_X}\co_{\oX},
\end{equation}
que l'on considère comme un faisceau de $\oX_{\zar}$ ou $\oX_{\et}$, selon le contexte (cf. \ref{notconv12}). Compte tenu de \eqref{notconv12e}, 
les homomorphismes \eqref{TFA4d} induisent un homomorphisme que l'on note encore 
\begin{equation}\label{TFA4k}
d\log\colon \cM_\oX\rightarrow \tOmega^1_{\oX/\oS}. 
\end{equation}
Celui-ci induit un morphisme $\co_\oX$-linéaire et surjectif 
\begin{equation}\label{TFA4i}
\cM_\oX^\gp\otimes_\mZ\co_\oX\rightarrow \tOmega^1_{\oX/\oS}.
\end{equation}

Pour tout entier $n\geq 1$, on pose 
\begin{equation}\label{TFA4l}
\tOmega^1_{\oX_n/\oS_n}=\tOmega^1_{X/S}\otimes_{\co_X}\co_{\oX_n},
\end{equation}
que l'on considère aussi un faisceau de $X_{s,\zar}$ ou $X_{s,\et}$, selon le contexte (cf. \ref{TFA5}). 
Compte tenu de \eqref{notconv12e}, le morphisme \eqref{TFA4i}
induit un morphisme $\co_{\oX_n}$-linéaire et surjectif de $X_{s,\et}$
\begin{equation}\label{TFA4m}
\oa^*(\cM_\oX^\gp)\otimes_\mZ\co_{\oX_n}\rightarrow \tOmega^1_{\oX_n/\oS_n}.
\end{equation}

\subsection{}\label{TFSLA70}
Pour tout $\mU$-topos $\fT$ et tout groupe abélien $G$, on note $G_\fT$ 
(ou simplement $G$ lorsqu'il n'y a aucun risque d'ambiguïté) le groupe constant de $\fT$ de valeur $G$. 
D'après (\cite{agt} VI.5.17), $G_\tE$ est canoniquement isomorphe au faisceau associé au préfaisceau  
$\{U\mapsto G_{\oU^\circ_\fet}\}$ sur $E$ (cf. \ref{TFA6}).
Pour tout entier $n\geq 1$, on désigne par $\mu_{n,\fT}$ 
(ou simplement $\mu_n$ lorsqu'il n'y a aucun risque d'ambiguïté) le groupe constant de $\fT$ de valeur 
le groupe $\mu_n(\co_\oK)$ des racines $n$-ièmes de l'unité dans $\co_\oK$.  
Pour tout $(\mZ/n\mZ)$-module $F$ de $\fT$ et tout entier $m$, on pose 
\begin{equation}\label{TFSLA70b}
F(m)=\left\{
\begin{array}{clcr}
F\otimes_{\mZ}\mu_{n,\fT}^{\otimes m},&{\rm si}\ m\geq 0,\\
\Hom_{\mZ}(\mu_{n,\fT}^{\otimes -m},F),&{\rm si}\ m< 0.
\end{array}
\right.
\end{equation}
Cet usage n'étant pas classique pour le topos $X_\et$, 
on prendra garde de ne pas confondre le faisceau $\mu_{n,X_\et}$ ainsi défini avec 
le faisceau des racines $n$-ièmes de l'unité de $X_\et$, c'est-à-dire le noyau de 
la puissance $n$-ième sur $\mG_{m,X}$ (qui ne sera pas utilisé dans cet article).

\begin{lem}\label{Slad4}
Soient $A$ un anneau intègre, $M$ un monoïde intègre, 
$u$ un homomorphisme injectif de $M$ dans le monoïde multiplicatif $A$, $n$ un entier 
$\geq 1$, $\phi\colon A\rightarrow A$ l'application d'élévation à la puissance $n$-ième. 
Supposons que pour tout $t\in M$, l'équation $x^n=u(t)$ admette 
$n$ racines distinctes dans $A$. Alors, $\phi^{-1}(M)$ est un monoïde intègre et l'homomorphisme  
$\phi^{-1}(M)\rightarrow M$ déduit de $\phi$ identifie $M$ au quotient de $\phi^{-1}(M)$ par le sous-monoïde $\mu_n(A)$ 
des racines $n$-ièmes de l'unité dans $A$  \eqref{notconv7}. 
\end{lem}

En effet, pour tout $t\in M$, on a $u(t)\not=0$ car $u$ est injectif et $M$ est intègre. 
Comme $A$ est intègre, $\phi^{-1}(M)$ est donc intègre. Par ailleurs, l'application canonique 
$\phi^{-1}(M)\rightarrow M$ est surjective et fait de $\phi^{-1}(M)$ un torseur sur $M$ sous le groupe $\mu_n(A)$.
Par suite, $M$ est le quotient de $\phi^{-1}(M)$ par le sous-monoïde $\mu_n(A)$  (\cite{ogus} I 1.1.6).

\subsection{}\label{TFSLA7}
Considérant l'anneau $\ocB$ de $\tE$ \eqref{TFA2} muni de sa structure multiplicative 
et le monoïde $\cM_\oX$ de $\oX_\et$ défini dans \eqref{TFA4f}, on note 
\begin{equation}\label{TFSLA7b}
\fA\colon \sigma^*(\hbar_*(\cM_\oX))\rightarrow \ocB
\end{equation}
le composé des homomorphismes canoniques \eqref{TFA4g} et \eqref{TFA2c}
\[
\sigma^*(\hbar_*(\cM_\oX))\rightarrow \sigma^*(\hbar_*(\co_\oX))\rightarrow \ocB.
\] 
Pour tout entier $n\geq 0$,  on désigne par $\cQ_n$ le monoïde de $\tE$ défini par le diagramme cartésien 
de la catégorie des monoïdes de $\tE$ \eqref{notconv7}
\begin{equation}\label{TFSLA7a}
\xymatrix{
{\cQ_n}\ar[d]\ar[r]&{\sigma^*(\hbar_*(\cM_\oX))}\ar[d]^\fA\\
{\ocB}\ar[r]^{\nu^n}&{\ocB}}
\end{equation}
où $\nu$ est l'endomorphisme d'élévation à la puissance $p$-ième de $\ocB$. 
Les monoïdes $(\cQ_n)_{n\in \mN}$ forment naturellement un système projectif.
Compte tenu de \eqref{TFA5b} et \eqref{TFA66c}, le diagramme \eqref{TFSLA7a} induit par image inverse par le plongement
$\delta$ \eqref{TFA66a} un diagramme cartésien 
de la catégorie des monoïdes de $\tE_s$
\begin{equation}\label{TFSLA7c}
\xymatrix{
{\delta^*(\cQ_n)}\ar[d]\ar[r]&{\sigma^*_s(\oa^*(\cM_\oX))}\ar[d]^{\delta^*(\fA)}\\
{\delta^*(\ocB)}\ar[r]^{\delta^*(\nu^n)}&{\delta^*(\ocB)}}
\end{equation}

\begin{prop}\label{TFSLA8}
Soit $n$ un entier $\geq 0$. 
\begin{itemize}
\item[{\rm (i)}] Le monoïde $\delta^*(\cQ_n)$ est intègre. 
\item[{\rm (ii)}] L'homomorphisme canonique $\mu_{p^n,\tE_s}\rightarrow \delta^*(\cQ_n)$  \eqref{TFSLA70} est un monomorphisme.
\item[{\rm (iii)}] La projection canonique $\delta^*(\cQ_n)\rightarrow \sigma^*_s(\oa^*(\cM_\oX))$ 
identifie $\sigma^*_s(\oa^*(\cM_\oX))$ au quotient de $\delta^*(\cQ_n)$ par le sous-monoïde $\mu_{p^n,\tE_s}$ \eqref{notconv7}.
\end{itemize}
\end{prop} 

Soient $(\oy\rightsquigarrow \ox)$ un point de $X_\et\gtimes_{X_\et}\oX^\circ_\et$ \eqref{rec15} tel que $\ox$ soit
au-dessus de $s$, $\rho(\oy\rightsquigarrow \ox)$ son image dans $\tE$ par le morphisme $\rho$ \eqref{TFA6e}. 
On note encore $\ox$ le point géométrique $\oa(\ox)$ de $\oX$ \eqref{TFA5a}. 
Compte tenu de (\cite{agt} (VI.10.18.1)) et  \eqref{TFA5b}, la fibre de $\fA$ en $\rho(\oy\rightsquigarrow \ox)$ s'identifie à un homomorphisme
\begin{equation}
\fA_{\rho(\oy\rightsquigarrow \ox)} \colon \cM_{\oX,\ox}\rightarrow \ocB_{\rho(\oy \rightsquigarrow \ox)}.
\end{equation}
Montrons que celui-ci remplit les conditions de \ref{Slad4}.
D'après (\cite{agt} III.10.10), $\ocB_{\rho(\oy \rightsquigarrow \ox)}$ 
est un anneau normal et strictement local (et en particulier intègre) et l'homomorphisme
\begin{equation}
\co_{\oX,\ox}\rightarrow \ocB_{\rho(\oy \rightsquigarrow \ox)}
\end{equation}
fibre de l'homomorphisme canonique $\sigma^*(\hbar_*(\co_\oX))\rightarrow \ocB$ \eqref{TFA2c} en $\rho(\oy \rightsquigarrow \ox)$,
est injectif et local. Comme $\alpha\colon \cM_\oX\rightarrow \co_\oX$ \eqref{TFA4g} est un monomorphisme, 
$\fA_{\rho(\oy\rightsquigarrow \ox)}$ est injectif.
D'autre part, avec les notations de \ref{TFA11}, on a un isomorphisme canonique \eqref{TFA11c} 
\begin{equation}
\ocB_{\rho(\oy\rightsquigarrow \ox)} \stackrel{\sim}{\rightarrow} 
\underset{\underset{(U,\fp)\in \fW_\ox^\circ}{\longrightarrow}}{\lim}\ \oR^{\oy}_U.
\end{equation} 
On a clairement \eqref{TFA5b}
\begin{equation}
\cM_{\oX,\ox} \stackrel{\sim}{\rightarrow} 
\underset{\underset{(U,\fp)\in \fW_\ox^\circ}{\longrightarrow}}{\lim}\ \Gamma(\oU,\cM_\oX).
\end{equation} 
De plus, l'homomorphisme $\fA_{\rho(\oy\rightsquigarrow \ox)}$ s'identifie à la limite inductive des homomorphismes composés
\[
\Gamma(\oU,\cM_\oX) \stackrel{\alpha_U}{\longrightarrow} \Gamma(\oU,\co_{\oX})\longrightarrow \oR^\oy_U,
\]
où la seconde flèche est l'homomorphisme canonique. 
Pour tout $(U,\fp)\in \ob(\fW_\ox)$ et tout $t\in \Gamma(\oU,\cM_\oX)$, comme $\alpha_\oU(t)$ est inversible sur $\oU^\circ$ \eqref{TFA4g}, le $\oU$-schéma 
\begin{equation}
\Spec(\co_\oU[T]/(T^{p^n}-\alpha_\oU(t))
\end{equation}
est fini et étale au-dessus de $\oU^\circ$. L'anneau $\oR^{\oy}_U$  étant intègre et normal \eqref{TFA10}, 
il résulte aussitôt de sa définition \eqref{TFA9c}  que l'équation $T^{p^n}=\alpha_\oU(t)$ admet une racine et donc 
$p^n$ racines distinctes dans $\oR_U^\oy$. 
On en déduit que pour tout $m\in \cM_{\oX,\ox}$, l'équation $T^{p^n}=\fA_{\rho(\oy\rightsquigarrow \ox)}(m)$ 
admet $p^n$ racines distinctes dans $\ocB_{\rho(\oy\rightsquigarrow \ox)}$. La proposition s'ensuit en vertu de 
\ref{Slad4} et \ref{TFA7}.

\begin{cor}\label{TFSLA9}
Pour tout entier $n\geq 0$, la suite d'homomorphismes canoniques
\begin{equation}\label{TFSLA9a}
0\rightarrow \mu_{p^n,\tE_s}\rightarrow \delta^*(\cQ^\gp_n)\rightarrow \sigma_s^*(\oa^*(\cM^\gp_\oX))\rightarrow 0
\end{equation}
est exacte. 
\end{cor} 

En effet, le foncteur d'injection canonique de la catégorie des groupes abéliens de $\tE_s$ 
dans celle des monoïdes de $\tE_s$ admet pour adjoint à gauche, le foncteur ``groupe associé'' 
$\cN\mapsto \cN^\gp$. Ce dernier commute donc aux limites inductives.
L'exactitude de la suite \eqref{TFSLA9a} au centre et à droite résulte alors de \ref{TFSLA8}(iii) et (\cite{agt} (II.5.7.1)). 
L'exactitude de cette suite à gauche est une conséquence de \ref{TFSLA8}(i)-(ii).

\begin{cor}\label{TFSLA18}
Soient $\ox$ un point géométrique de $X$ au-dessus de $s$, $X'$ le localisé strict de $X$ en $\ox$, 
$\varphi_\ox\colon \tE\rightarrow \oX'^\circ_\fet$ le foncteur canonique \eqref{TFA14a}, $n$ un entier $\geq 0$. 
On désigne encore par $\ox$ le point géométrique $\oa(\ox)$ de $\oX$ \eqref{TFA5a}. 
Alors, l'image du diagramme \eqref{TFSLA7a} par le foncteur $\varphi_\ox$ s'identifie à un diagramme cartésien 
de la catégorie des monoïdes de $\oX'^\circ_\fet$
\begin{equation}\label{TFSLA18a}
\xymatrix{
{\varphi_\ox(\cQ_n)}\ar[r]\ar[d]&{\cM_{\oX,\ox}}\ar[d]^{\varphi_\ox(\fA)}\\
{\varphi_\ox(\ocB)}\ar[r]^{\nu^n}&{\varphi_\ox(\ocB)}}
\end{equation}
où on a encore noté $\nu$  l'endomorphisme d'élévation à la puissance $p$-ième de $\varphi_\ox(\ocB)$.
De plus, le monoïde $\varphi_\ox(\cQ_n)$ est intègre, et  la projection canonique $\varphi_\ox(\cQ_n)\rightarrow \cM_{\oX,\ox}$
identifie $\cM_{\oX,\ox}$ au quotient de $\varphi_\ox(\cQ_n)$ par le sous-monoïde $\mu_{p^n,\oX'^\circ_\fet}$ \eqref{notconv7}. 
\end{cor}

En effet, on a un isomorphisme canonique 
$\varphi_\ox(\sigma^*(\hbar_*(\cM_\oX)))\stackrel{\sim}{\rightarrow}\cM_{\oX,\ox}$
d'après \eqref{TFA14b} et \eqref{TFA5b}; d'où la première assertion. 
Comme $\oX'^\circ$ est intègre et non-vide \eqref{TFA13},  
considérant un point $(\oy\rightsquigarrow \ox)$ de $X_\et\gtimes_{X_\et}\oX^\circ_\et$,
la seconde assertion résulte de la preuve de \ref{TFSLA8} compte tenu de la description \eqref{TFA12e} 
du foncteur fibre de $\tE$ associé au point $\rho(\oy\rightsquigarrow \ox)$. 

\subsection{}\label{TFSLA13}
Pour tout entier $n\geq 1$, on note 
\begin{equation}\label{TFSLA13a}
\partial_n\colon \oa^*(\cM^\gp_\oX)\rightarrow \rR^1\sigma_{s*}(\mu_{p^n,\tE_s})
\end{equation}
l'homomorphisme composé du morphisme d'adjonction 
$\oa^*(\cM^\gp_\oX)\rightarrow \sigma_{s*}(\sigma^*_s(\oa^*(\cM^\gp_\oX)))$
et du cobord de la suite exacte \eqref{TFSLA9a}. Celui-ci induit deux morphismes  $\co_{\oX_n}$-linéaires
\begin{eqnarray}
\oa^*(\cM^\gp_\oX)\otimes_\mZ\co_{\oX_n}&\rightarrow& \rR^1\sigma_{n*}(\ocB_n(1)),\label{TFSLA13b}\\
\oa^*(\cM^\gp_\oX)\otimes_\mZ\co_{\oX_n}&\rightarrow& \rR^1\sigma_{n*}(\xi\ocB_n),\label{TFSLA13c}
\end{eqnarray}
le second étant défini au moyen du composé des morphismes canoniques \eqref{THT3c}
\begin{equation}\label{TFSLA13d}
\co_C(1)\stackrel{\sim}{\rightarrow}p^{\frac{1}{p-1}}\xi\co_C\rightarrow \xi \co_C.
\end{equation}

\begin{teo}\label{TFSLA24}
Soit $n$ un entier $\geq 1$.  
\begin{itemize}
\item[{\rm (i)}] Le morphisme \eqref{TFSLA13c} se factorise à travers le morphisme
surjectif \eqref{TFA4m} et induit un morphisme $\co_{\oX_n}$-linéaire de $X_{s,\et}$ 
\begin{equation}\label{TFSLA24a}
\xi^{-1}\tOmega^1_{\oX_n/\oS_n}\rightarrow \rR^1\sigma_{n*}(\ocB_n).
\end{equation}
\item[{\rm (ii)}] Il existe un et un unique homomorphisme de $\co_{\oX_n}$-algèbres graduées de $X_{s,\et}$
\begin{equation}\label{TFSLA24b}
\wedge (\xi^{-1}\tOmega^1_{\oX_n/\oS_n})\rightarrow \oplus_{i\geq 0}\rR^i\sigma_{n*}(\ocB_n)
\end{equation}
dont la composante en degré un est le morphisme \eqref{TFSLA24a}. De plus, 
son noyau est annulé par $p^{\frac{2d}{p-1}}\fm_\oK$ et son conoyau est annulé par $p^{\frac{2d+1}{p-1}}\fm_\oK$, 
où $d=\dim(X/S)$ \eqref{TFA3}.  
\end{itemize}
\end{teo}

La preuve de cet énoncé sera donnée dans \ref{TFSLA27} après quelques résultats préliminaires.

\subsection{}\label{TFSLA19}
Soient $(\oy\rightsquigarrow \ox)$ un point de $X_\et\gtimes_{X_\et}\oX^\circ_\et$ \eqref{rec15} tel que $\ox$ soit
au-dessus de $s$, $X'$ le localisé strict de $X$ en $\ox$. 
On note encore $\ox$ le point géométrique $\oa(\ox)$ de $\oX$ \eqref{TFA5a}. 
D'après (\cite{agt} III.3.7), $\oX'$ est normal et strictement local (et en particulier intègre); on peut donc l'identifier 
au localisé strict de $\oX$ en $\ox$. Le $X$-morphisme $u\colon \oy\rightarrow X'$
définissant $(\oy\rightsquigarrow \ox)$ se relève en un $\oX^\circ$-morphisme $v\colon \oy\rightarrow \oX'^\circ$ et 
induit donc un point géométrique de $\oX'^\circ$ que l'on note aussi (abusivement) $\oy$.
On désigne par 
\begin{equation}\label{TFSLA19a}
\varphi_\ox\colon \tE\rightarrow \oX'^\circ_\fet
\end{equation} 
le foncteur canonique \eqref{TFA14a} et par 
\begin{equation}\label{TFSLA19b}
\nu_{\oy}\colon \oX'^\circ_\fet \stackrel{\sim}{\rightarrow}\bB_{\pi_1(\oX'^\circ,\oy)}
\end{equation}
le foncteur fibre de $\oX'^\circ_\fet$ en $\oy$ \eqref{notconv11c}. 
Compte tenu de \eqref{TFA14b}, \eqref{TFA5b} et \eqref{TFA12g}, 
l'homomorphisme $\nu_\oy(\varphi_\ox(\fA))$ \eqref{TFSLA7b} 
s'identifie à  l'homomorphisme composé
\begin{equation}\label{TFSLA19c}
\xymatrix{
{\cM_{\oX,\ox}}\ar[r]^{\alpha_\ox}&{\co_{\oX,\ox}}\ar[r]&{\oR^\oy_{X'}}},
\end{equation}
où la première flèche est la fibre de $\alpha$ \eqref{TFA4g} en $\ox$, $\oR^\oy_{X'}$ est l'algèbre définie dans \eqref{TFA12f}  et  
la seconde flèche est l'homomorphisme canonique \eqref{TFA12f}.
Pour tout entier $n\geq 0$, l'image du diagramme \eqref{TFSLA7a} par le foncteur $\nu_\oy\circ \varphi_\ox$ 
s'identifie donc à un diagramme cartésien de la catégorie des monoïdes de  $\bB_{\pi_1(\oX'^\circ,\oy)}$
\begin{equation}\label{TFSLA19d}
\xymatrix{
{\nu_\oy(\varphi_\ox(\cQ_n))}\ar[r]\ar[d]&{\cM_{\oX,\ox}}\ar[d]^{\nu_\oy(\varphi_\ox(\fA))}\\
{\oR^\oy_{X'}}\ar[r]^{\nu^n}&{\oR^\oy_{X'}}}
\end{equation}
où on a encore noté $\nu$  l'homomorphisme d'élévation à la puissance $p$-ième de $\oR^\oy_{X'}$.

Pour tout $t\in \cM_{\oX,\ox}$, la section $\alpha_\ox(t)$  est inversible sur $\oX'^\circ$. Le $\oX'$-schéma 
\begin{equation}\label{TFSLA19g}
\cT_n(t)=\Spec(\co_{\oX'}[T]/(T^{p^n}-\alpha_\ox(t))
\end{equation}
est donc fini et étale sur $\oX'^\circ$. En particulier, $\cT^\circ_n(t)=\cT_n(t)\times_XX^\circ$ induit un $\mu_{p^n}(\co_\oK)$-torseur 
de $\oX'^\circ_\fet$. On désigne par $\cT_n(t)_\oy$ la fibre de $\cT_n(t)$ au-dessus de $\oy$ et par $\cQ_n(t)$ la fibre 
de la projection canonique $\nu_\oy(\varphi_\ox(\cQ_n))\rightarrow \cM_{\oX,\ox}$ \eqref{TFSLA19d} au-dessus de $t$.
On a un isomorphisme canonique de $\mu_{p^n}(\co_\oK)$-torseurs de $\bB_{\pi_1(\oX'^\circ,\oy)}$
\begin{equation}\label{TFSLA19h}
\cT_n(t)_\oy\stackrel{\sim}{\rightarrow}\cQ_n(t).
\end{equation}

\begin{lem}\label{TFSLA17}
Sous les hypothèses de \ref{TFSLA19}, pour tout entier $n\geq 0$, l'application 
\begin{equation}\label{TFSLA17a}
\tau_n\colon \cM_{\oX,\ox}\rightarrow \rH^1(\oX'^\circ_\fet,\mu_{p^n,\oX'^\circ_\fet}), \ \ \ t\mapsto [\cT^\circ_n(t)],
\end{equation}
est un homomorphisme de monoïdes \eqref{TFSLA70}, et le diagramme 
\begin{equation}\label{TFSLA17b}
\xymatrix{
{\cM_{\oX,\ox}}\ar[r]^-(0.5){\tau_n}\ar[d]&{\rH^1(\oX'^\circ_\fet,\mu_{p^n,\oX'^\circ_\fet})}\ar[r]_-(0.5)\sim^-(0.5)v&
{\rH^1(\oX'^\circ_\fet,\varphi_\ox(\delta_*(\mu_{p^n,\tE_s})))}\\
{\cM_{\oX,\ox}^\gp}\ar[r]^-(0.5){\partial_{n,\ox}}&{\rR^1\sigma_{s*}(\mu_{p^n,\tE_s})_\ox}\ar[r]_-(0.5)\sim^-(0.5)w&
{\rR^1\sigma_*(\delta_*(\mu_{p^n,\tE_s}))_\ox}\ar[u]_u}
\end{equation}
où $\partial_{n,\ox}$ est la fibre en $\ox$ de  l'homomorphisme 
$\partial_n$ \eqref{TFSLA13a}, $u$ est l'isomorphisme \eqref{TFA14d},  
$v$ est induit par l'isomorphisme canonique $\varphi_\ox\stackrel{\sim}{\rightarrow} \varphi_\ox\delta_*\delta^*$ \eqref{TFA13a},
$w$ est induit par l'isomorphisme \eqref{TFA66d} et la flèche non libellée est l'homomorphisme canonique, est commutatif. 
\end{lem}

Considérons la suite exacte de groupes abéliens de $\tE$ 
\begin{equation}\label{TFSLA17e}
0\rightarrow \delta_*(\mu_{p^n,\tE_s})\rightarrow \delta_*(\delta^*(\cQ_n^\gp))\rightarrow 
\delta_*(\sigma^*_s(\oa^*(\cM_\oX^\gp)))\rightarrow 0
\end{equation}
déduite de la suite exacte \eqref{TFSLA9a} par image directe par le plongement $\delta$. 
D'après (\cite{agt} VI.10.30(iii)), le diagramme 
\begin{equation}
\xymatrix{
{\sigma_*(\delta_*(\sigma_s^*(\oa^*\cM^\gp_\oX)))_\ox}\ar[r]^-(0.5){u'}_-(0.5)\sim\ar[d]_a&
{\rH^0(\oX'^\circ_\fet,\varphi_\ox(\delta_*(\sigma_s^*(\oa^*\cM^\gp_\oX))))}\ar[d]^b\\
{\rR^1\sigma_*(\delta_*(\mu_{p^n,\tE_s}))_\ox}
\ar[r]^-(0.5)u_-(0.5)\sim&{\rH^1(\oX'^\circ_\fet,\varphi_\ox(\delta_*(\mu_{p^n,\tE_s})))}}
\end{equation}
où les flèches horizontales sont les isomorphismes \eqref{TFA14d} et les flèches verticales sont les bords 
des suites exactes longues de cohomologie, est commutatif. Par ailleurs, le diagramme 
\begin{equation}
\xymatrix{
{\cM^\gp_{\oX,\ox}}\ar[r]^-(0.5)i\ar[rd]_{\partial_{n,\ox}}&{\sigma_{s*}(\sigma_s^*(\oa^*\cM^\gp_\oX))_\ox}
\ar[r]^-(0.5){w'}_-(0.5)\sim\ar[d]^c&{\sigma_*(\delta_*(\sigma_s^*(\oa^*\cM^\gp_\oX)))_\ox}\ar[d]^a\\
&{\rR^1\sigma_{s*}(\mu_{p^n,\tE_s})_\ox}\ar[r]^-(0.5){w}_-(0.5)\sim&{\rR^1\sigma_*(\delta_*(\mu_{p^n,\tE_s}))_\ox}}
\end{equation}
où $i$ est induit par le morphisme d'adjonction $\id\rightarrow \sigma_{s*}\sigma_s^*$, 
$c$ est le bord de la suite exacte longue de cohomologie déduite de la suite \eqref{TFSLA9a}
et $w$ et $w'$ sont induits par l'isomorphisme \eqref{TFA66d}, est commutatif. 

D'après \eqref{TFA13a}, \eqref{TFA14b} et \eqref{TFA66c}, l'image de \eqref{TFSLA17e} par le foncteur $\varphi_\ox$ s'identifie à une suite exacte de 
groupes abéliens de $\oX'^\circ_\fet$
\begin{equation}\label{TFSLA17c}
0\rightarrow \mu_{p^n,\oX'^\circ_\fet}\rightarrow \varphi_\ox(\cQ_n^\gp)\rightarrow \cM_{\oX,\ox}^\gp\rightarrow 0. 
\end{equation}
Celle-ci induit un homomorphisme 
\begin{equation}\label{TFSLA17d}
\tau'_n\colon \cM^\gp_{\oX,\ox}\rightarrow\rH^1(\oX'^\circ_\fet,\mu_{p^n,\oX'^\circ_\fet}),
\end{equation}
qui s'identifie donc à $b$
\begin{equation}
\xymatrix{
{\cM^\gp_{\oX,\ox}}\ar[d]_{\tau'_n}\ar[r]^-(0.5){v'}_-(0.5)\sim&
{\rH^0(\oX'^\circ_\fet,\varphi_\ox(\delta_*(\sigma_s^*(\oa^*\cM^\gp_\oX))))}\ar[d]^b\\
{\rH^1(\oX'^\circ_\fet,\mu_{p^n,\oX'^\circ_\fet})}\ar[r]^-(0.5)v_-(0.5)\sim&
{\rH^1(\oX'^\circ_\fet,\varphi_\ox(\delta_*(\mu_{p^n,\tE_s})))}}
\end{equation}
On vérifie que $v'^{-1}\circ u'\circ w'\circ i$ est l'identité. 

Il résulte de \ref{TFSLA18} et (\cite{ogus} I.1.1.6) que le diagramme canonique
\begin{equation}
\xymatrix{
{\varphi_\ox(\cQ_n)}\ar[d]\ar[r]&{\cM_{\oX,\ox}}\ar[d]\\
{\varphi_\ox(\cQ^\gp_n)}\ar[r]&{\cM^\gp_{\oX,\ox}}}
\end{equation}
est cartésien. Compte tenu de \eqref{TFSLA19h}, on en déduit  
que $\tau_n$ est le composé de $\tau'_n$ et de l'homomorphisme canonique $\cM_{\oX,\ox}\rightarrow \cM_{\oX,\ox}^\gp$;
d'où la proposition.

\subsection{}\label{TFSLA14} 
Conservons les hypothèses et notations de \ref{TFSLA19}; supposons, de plus, les conditions suivantes remplies~:
\begin{itemize}
\item[(i)] Le schéma $X$ est affine et connexe et le morphisme $f\colon (X,\cM_X)\rightarrow (S,\cM_S)$ admet une carte adéquate \eqref{cad1}.
\item[(ii)] Il existe une carte fine et saturée $M\rightarrow \Gamma(X,\cM_X)$ pour $(X,\cM_X)$
induisant un isomorphisme 
\begin{equation}\label{TFSLA31a}
M\stackrel{\sim}{\rightarrow} \Gamma(X,\cM_X)/\Gamma(X,\co^\times_X).
\end{equation}
\end{itemize}
Ces conditions correspondent à celles fixées dans \ref{TFA1}. 
Notons $R$ l'anneau de $X$ et posons $\tOmega^1_{R/\co_K}=\tOmega^1_{X/S}(X)$ \eqref{TFA3b}.
Le schéma $\oX$ étant localement irréductible, 
il est la somme des schémas induits sur ses composantes irréductibles. 
On désigne par $\oX^\star$  la composante irréductible de $\oX$ contenant $\oy$. De même, $\oX^\circ$
est la somme des schémas induits sur ses composantes irréductibles, et $\oX^{\star\circ}=\oX^\star\times_XX^\circ$ 
est la composante irréductible de $\oX^\circ$ contenant $\oy$.  
On munit  $\coX=X\times_S\coS$ \eqref{deform1d} de la structure logarithmique $\cM_\coX$ image inverse de $\cM_X$
et on se donne une $(\cA_2(\oS),\cM_{\cA_2(\oS)})$-déformation lisse $(\tX,\cM_\tX)$ de $(\coX,\cM_{\coX})$ 
(cf. \ref{deform3}). On note $\hoR^\oy_X$ le complété $p$-adique de l'anneau $\oR^\oy_X$ \eqref{TFA9c} et 
\begin{equation}\label{TFSLA14b} 
0\rightarrow \hoR^\oy_X\rightarrow \cF\rightarrow \xi^{-1} \tOmega^1_{R/\co_K}\otimes_R\hoR^\oy_X\rightarrow 0
\end{equation}
l'extension de Higgs-Tate associée à $(\tX,\cM_\tX)$ \eqref{THT7};  
c'est une extension de $\hoR^\oy_X$-représentations continues de $\pi_1(\oX^{\star\circ},\oy)$. 
Pour tout entier $n\geq 0$, on désigne par  
\begin{equation}\label{TFSLA14c} 
0\rightarrow \xi(\oR^\oy_X/p^n\oR^\oy_X)\rightarrow \xi(\cF/p^n\cF)\rightarrow \tOmega^1_{R/\co_K}\otimes_R
(\oR^\oy_X/p^n\oR^\oy_X)\rightarrow 0
\end{equation}
la suite exacte déduite de \eqref{TFSLA14b}, par
\begin{equation}\label{TFSLA14e}
A_n\colon \Gamma(\oX^\star,\tOmega^1_{\oX_n/\oS_n})
\rightarrow \rH^1(\pi_1(\oX^{\star\circ},\oy),\xi(\oR^\oy_X/p^n\oR^\oy_X))
\end{equation}
le morphisme $\Gamma(\oX^\star,\co_\oX)$-linéaire induit, et par $Q_n$ le produit fibré du diagramme d'homomorphismes de monoïdes de $\bB_{\pi_1(\oX^{\star \circ},\oy)}$
\begin{equation}\label{TFSLA14d} 
\xymatrix{
&{\Gamma(X,\cM_X)}\ar[d]\\
{\oR^\oy_X}\ar[r]^{\nu^n}&{\oR^\oy_X}}
\end{equation}
où on a encore noté $\nu$  l'homomorphisme d'élévation à la puissance $p$-ième de $\oR^\oy_X$.
On observera que $A_n$ ne dépend pas du choix de $(\tX,\cM_\tX)$ (\cite{agt} II.10.10). 
L'homomorphisme canonique $\Gamma(X,\cM_X)\rightarrow \oR^\oy_X$ vérifie les hypothèses de \ref{Slad4}. 
En effet, l'homomorphisme $\Gamma(X,\cM_X)\rightarrow R$ est injectif \eqref{TFA3e}
et comme $X$ est connexe,  l'homomorphisme $R\rightarrow \oR^\oy_X$ est injectif. 
Par suite, le monoïde $Q_n$ est intègre et la projection canonique $Q_n\rightarrow \Gamma(X,\cM_X)$ 
identifie $\Gamma(X,\cM_X)$ au quotient de $Q_n$ par le sous-monoïde $\mu_{p^n}(\co_\oK)$ d'après  \ref{Slad4}. 
On a donc une suite exacte de groupes abéliens de $\bB_{\pi_1(\oX^{\star\circ},\oy)}$
\begin{equation}\label{TFSLA14f}
0\rightarrow \mu_{p^n}(\co_\oK)\rightarrow Q_n^\gp\rightarrow \Gamma(X,\cM_X)^\gp\rightarrow 0.
\end{equation}
Elle induit un homomorphisme 
\begin{equation}\label{TFSLA14g}
B_n\colon \Gamma(X,\cM_X)\rightarrow \rH^1(\pi_1(\oX^{\star\circ},\oy),\mu_{p^n}(\co_\oK)).
\end{equation}
On observera que le diagramme canonique
\begin{equation}\label{TFSLA14h}
\xymatrix{
{Q_n}\ar[d]\ar[r]&{\Gamma(X,\cM_X)}\ar[d]\\
{Q_n^\gp}\ar[r]&{\Gamma(X,\cM_X)^\gp}}
\end{equation}
est cartésien (\cite{ogus} I 1.1.6). On note encore 
\begin{equation}\label{TFSLA14i}
-\log([\ ])\colon \mu_{p^n}(\co_\oK)\rightarrow \xi(\oR^\oy_X/p^n\oR^\oy_X)
\end{equation}
l'homomorphisme induit par l'homomorphisme $-\log([\ ])$ \eqref{THT3b}.

\begin{lem}\label{TFSLA20}
Sous les hypothèses de \ref{TFSLA14}, pour tout entier $n\geq 0$, le diagramme 
\begin{equation}\label{TFSLA20a}
\xymatrix{
{\Gamma(X,\cM_X)}\ar[r]^-(0.5){B_n}\ar[d]_{d\log}&
{\rH^1(\pi_1(\oX^{\star\circ},\oy),\mu_{p^n}(\co_\oK))}\ar[d]^{-\log([\ ])}\\
{\Gamma(\oX^\star,\tOmega^1_{\oX_n/\oS_n})}\ar[r]^-(0.5){A_n}&
{\rH^1(\pi_1(\oX^{\star\circ},\oy),\xi(\oR^\oy_X/p^n\oR^\oy_X))}}
\end{equation}
où on a encore noté $d\log$ l'homomorphisme induit par l'homomorphisme $d\log$ \eqref{TFA3c}, est commutatif. 
\end{lem}

On notera d'abord que le morphisme $\coX\rightarrow \tX$ étant un homéomorphisme universel, on peut identifier 
les topos étales de $\coX$ et $\tX$. D'après (\cite{fkato} 3.6 ou \cite{ogus} IV 2.1.2), 
$\cM_{\coX}$ est le quotient de $\cM_\tX$ par le sous-monoïde $1+\xi\co_{\coX}$ (cf. \cite{ogus} I 1.1.5). 
Comme $\coX$ est affine, on en déduit que l'homomorphisme canonique 
\begin{equation}\label{TFSLA20c}
\Gamma(\tX,\cM_\tX)\rightarrow \Gamma(\coX,\cM_{\coX})
\end{equation}
est surjectif. 
Soit $t\in \Gamma(X,\cM_X)$. Notons $Q_n(t)$ l'image inverse de $t$ par la projection canonique 
$Q_n\rightarrow \Gamma(X,\cM_X)$ et $\xi\cF_n(t)$ la fibre au-dessus de $d\log(t)\otimes 1$ dans l'extension \eqref{TFSLA14c}. 
Compte tenu de \eqref{TFSLA14h}, $B_n(t)$ est la classe du $\mu_{p^n}(\co_\oK)$-torseur $Q_n(t)$ de 
$\bB_{\pi_1(\oX^{\star\circ},\oy)}$. D'autre part, $A_n(d\log(t)\otimes 1)$
est la classe du $\xi(\oR_X^\oy/p^n\oR_X^\oy)$-torseur $\xi\cF_n(t)$ de $\bB_{\pi_1(\oX^{\star\circ},\oy)}$.
Choisissant un élément $t'\in \Gamma(\tX,\cM_\tX)$ ayant même image que $t$ dans $\Gamma(\coX,\cM_\coX)$ \eqref{TFSLA20c}, l'homomorphisme \eqref{THT8e} induit une application $\pi_1(\oX^{\star\circ},\oy)$-équivariante
\begin{equation}\label{TFSLA20b} 
Q_n(t)\rightarrow \xi\cF_n(t).
\end{equation}
Cette application est aussi équivariante relativement à l'homomorphisme $-\log([\ ])$ \eqref{TFSLA14i}; d'où la proposition.

\begin{prop}\label{TFSLA21} 
Sous les hypothèses de \ref{TFSLA19}, pour tout entier $n\geq 1$, il existe un unique morphisme $\co_{\oX,\ox}$-linéaire
\begin{equation}\label{TFSLA21a} 
\phi_n\colon \tOmega^1_{\oX_n/\oS_n,\ox}\rightarrow \rH^1(\pi_1(\oX'^\circ,\oy),\xi(\oR^\oy_{X'}/p^n\oR^\oy_{X'}))
\end{equation}
qui s'insère dans le diagramme commutatif 
\begin{equation}\label{TFSLA21b} 
\xymatrix{
{\cM_{\oX,\ox}}\ar[r]^-(0.5){\tau_n}\ar[d]_{d\log}&{\rH^1(\pi_1(\oX'^\circ,\oy),\mu_{p^n}(\co_\oK))}\ar[d]\ar[d]^{-\log([\ ])}\\
{\tOmega^1_{\oX_n/\oS_n,\ox}}\ar[r]^-(0.5){\phi_n}&{\rH^1(\pi_1(\oX'^\circ,\oy),\xi(\oR^\oy_{X'}/p^n\oR^\oy_{X'}))}}
\end{equation}
où on a encore noté $d\log$ l'homomorphisme induit par l'homomorphisme $d\log$ \eqref{TFA4k}. 
\end{prop}

En effet, l'unicité de $\phi_n$ est claire \eqref{TFA4i}. Montrons son existence. 
Pour toute extension finie $L$ de $K$ contenue dans $\oK$,
on rappelle que l'on a $X_L=X\times_SS_L$ \eqref{TFA4a}.
On note encore $\ox$ l'image du point géométrique $\ox$ par la projection canonique $\oX\rightarrow X_L$. 
D'après \eqref{TFA4f} et \eqref{TFA4h}, on a des isomorphismes canoniques 
\begin{eqnarray}
\cM_{\oX,\ox}&\stackrel{\sim}{\rightarrow}&\underset{\underset{L\subset \oK}{\longrightarrow}}{\lim}\  \cM_{X_L,\ox},\\
\tOmega^1_{\oX_n/\oS_n,\ox}&\stackrel{\sim}{\rightarrow}& \underset{\underset{L\subset \oK}{\longrightarrow}}{\lim}\  \tOmega^1_{X_n/S_n,\ox}\otimes_{\co_{X,\ox}}\co_{X_L,\ox}.
\end{eqnarray}
Comme chaque morphisme $f_L$ est adéquat (\cite{agt} III.4.8), 
on est réduit à montrer que la restriction de $-\log([\ ])\circ \tau_n$ à $\cM_{X,\ox}$ se factorise à travers l'homomorphisme
\begin{equation}
\cM_{X,\ox}\rightarrow \tOmega^1_{X_n/S_n,\ox}
\end{equation}
induit par $d\log$ \eqref{TFA3c}.
En vertu de (\cite{agt} II.5.17), on peut supposer les hypothèses de \ref{TFSLA14} remplies. 
Il suffit de montrer que  la restriction de l'homomorphisme $-\log([\ ])\circ \tau_n$ à $\Gamma(X,\cM_X)$ 
se factorise à travers l'homomorphisme 
\begin{equation}
\Gamma(X,\cM_X)\rightarrow \Gamma(X,\tOmega^1_{X_n/S_n})
\end{equation}
induit par $d\log$ \eqref{TFA3c}, ce qui résulte de \ref{TFSLA20}.

\begin{cor}\label{TFSLA22}
Sous les hypothèses de \ref{TFSLA14}, pour tout entier $n\geq 1$, le diagramme 
\begin{equation}\label{TFSLA22a}
\xymatrix{
{\Gamma(\oX^\star,\tOmega^1_{\oX_n/\oS_n})}\ar[r]^-(0.5){A_n}\ar[d]&
{\rH^1(\pi_1(\oX^{\star\circ},\oy),\xi(\oR^\oy_X/p^n\oR^\oy_X))}\ar[d]\\
{\tOmega^1_{\oX_n/\oS_n,\ox}}\ar[r]^-(0.5){\phi_n}&{\rH^1(\pi_1(\oX'^\circ,\oy),\xi(\oR^\oy_{X'}/p^n\oR^\oy_{X'}))}}
\end{equation}
où $A_n$ est le morphisme \eqref{TFSLA14e} et $\phi_n$ est le morphisme \eqref{TFSLA21a}, est commutatif. 
\end{cor}

Cela résulte aussitôt de la preuve de \ref{TFSLA21}.

\begin{cor}\label{TFSLA23}
Sous les hypothèses de \ref{TFSLA19}, pour tout entier $n\geq 1$, 
il existe un et un unique homomorphisme de $\co_{\oX,\ox}$-algèbres graduées 
\begin{equation}\label{TFSLA23a}
\wedge (\xi^{-1}\tOmega^1_{\oX_n/\oS_n,\ox})\rightarrow \oplus_{i\geq 0}\rH^i(\pi_1(\oX'^\circ,\oy),\oR^\oy_{X'}/p^n\oR^\oy_{X'})
\end{equation}
dont la composante en degré un est le morphisme $\xi^{-1}\phi_n$ \eqref{TFSLA21a}. De plus, 
notant $d=\dim(X/S)$ la dimension relative de $X$ sur $S$, 
le noyau de \eqref{TFSLA23a} est annulé par $p^{\frac{2d}{p-1}}\fm_\oK$ et 
son conoyau est annulé par $p^{\frac{2d+1}{p-1}}\fm_\oK$.
\end{cor}

On désigne par $\fV_\ox$ la catégorie des $X$-schémas étales $\ox$-pointés, par
$\fW_\ox$ (resp. $\fW'_\ox$) la sous-catégorie pleine de $\fV_\ox$ formée des objets $(U,\fp\colon \ox\rightarrow U)$ 
tels que le schéma $U$ soit affine
(resp. vérifie les conditions (i) et (ii) de \ref{TFSLA14}). 
Ce sont des catégories cofiltrantes, et les foncteurs d'injection canonique
$\fW'^\circ_\ox\rightarrow \fW^\circ_\ox\rightarrow \fV^\circ_\ox$ sont cofinaux d'après (\cite{agt} II.5.17) et (\cite{sga4} I 8.1.3(c)).
Reprenons les notations de \ref{TFA12}. 
Pour tout objet $(U,\fp)$ de $\fW'_\ox$, on note 
\begin{equation}\label{TFSLA23b}
A_{(U,\fp),n}\colon \Gamma(\oU^\star,\tOmega^1_{\oX_n/\oS_n})
\rightarrow \rH^1(\pi_1(\oU^{\star\circ},\oy),\xi(\oR^\oy_U/p^n\oR^\oy_U))
\end{equation} 
le morphisme canonique \eqref{TFSLA14e}; on rappelle que celui-ci ne dépend pas de 
la déformation choisie dans \ref{TFSLA14} pour le définir (\cite{agt} II.10.10). D'après \ref{TFSLA22}, le diagramme 
\begin{equation}\label{TFSLA23c}
\xymatrix{
{\Gamma(\oU^\star,\tOmega^1_{\oX_n/\oS_n})}\ar[rr]^-(0.5){A_{(U,\fp),n}}\ar[d]&&
{\rH^1(\pi_1(\oU^{\star\circ},\oy),\xi(\oR^\oy_U/p^n\oR^\oy_U))}\ar[d]\\
{\tOmega^1_{\oX_n/\oS_n,\ox}}\ar[rr]^-(0.5){\phi_n}&&{\rH^1(\pi_1(\oX'^\circ,\oy),\xi(\oR^\oy_{X'}/p^n\oR^\oy_{X'}))}}
\end{equation}
est commutatif. D'après (\cite{agt} VI.11.10) et \eqref{TFA12f}, pour tout entier $i\geq 0$, on a un isomorphisme canonique
\begin{equation}\label{TFSLA23d}
\rH^i(\pi_1(\oX'^\circ,\oy),\oR^\oy_{X'}/p^n\oR^\oy_{X'})
\stackrel{\sim}{\rightarrow} \underset{\underset{(U,\fp)\in \fW'^\circ_\ox}{\longrightarrow}}{\lim}\
\rH^i(\pi_1(\oU^{\star\circ},\oy),\oR^\oy_U/p^n\oR^\oy_U). 
\end{equation} 
Par ailleurs, $\oX'$ étant strictement local d'après (\cite{agt} III.3.7), il s'identifie au localisé strict de $\oX$ en $\ox$.
On a donc un isomorphisme canonique
\begin{equation}\label{TFSLA23e}
\tOmega^1_{\oX_n/\oS_n,\ox}
\stackrel{\sim}{\rightarrow} \underset{\underset{(U,\fp)\in \fW'^\circ_\ox}{\longrightarrow}}{\lim}\
\Gamma(\oU^\star,\tOmega^1_{\oX_n/\oS_n}). 
\end{equation} 
Les homomorphismes $A_{(U,\fp),n}$, pour $(U,\fp)\in \ob(\fW'_\ox)$, forment un homomorphisme 
de systèmes inductifs d'après (\cite{agt} III.10.16). Leur limite inductive s'identifie donc à $\phi_n$ \eqref{TFSLA23c}. 
Il suffit alors de montrer que pour tout objet $(U,\fp)$ de $\fW'_\ox$, 
il existe un et un unique homomorphisme de $\co_{\oX_n}(\oU^\star)$-algèbres graduées 
\begin{equation}\label{TFSLA23f}
\wedge \Gamma(\oU^\star,\xi^{-1}\tOmega^1_{\oX_n/\oS_n})\rightarrow 
\oplus_{i\geq 0}\rH^i(\pi_1(\oU^{\star\circ},\oy),\oR^\oy_U/p^n\oR^\oy_U)
\end{equation}
dont la composante en degré un est le morphisme $\xi^{-1} A_{(U,\fp),n}$, que 
son noyau est annulé par $p^{\frac{2d}{p-1}}\fm_\oK$ et que son conoyau est annulé par $p^{\frac{2d+1}{p-1}}\fm_\oK$.

D'après (\cite{agt} II.10.16), on a un diagramme commutatif
\begin{equation}\label{TFSLA23g}
\xymatrix{
{\Gamma(\oU^\star,\xi^{-1}\tOmega^1_{\oX_n/\oS_n})}\ar[d]_a\ar[r]^-(0.5){\xi^{-1}A_{(U,\fp),n}}&
{\rH^1(\pi_1(\oU^{\star \circ},\oy),\oR^{\oy}_U/p^n\oR^{\oy}_U)}\\
{\rH^1(\pi_1(\oU^{\star \circ},\oy),\xi^{-1}\oR^{\oy}_U(1)/p^n\xi^{-1}\oR^{\oy}_U(1))}\ar[r]^{-b}_\sim&
{\rH^1(\pi_1(\oU^{\star \circ},\oy),p^{\frac{1}{p-1}}\oR^{\oy}_U/p^{n+\frac{1}{p-1}}\oR^{\oy}_U)}\ar[u]_c}
\end{equation}
où $a$ est induit par le morphisme défini dans \eqref{cg6j}, 
$b$ est induit par l'isomorphisme canonique $\hoR^{\oy}_U(1)\stackrel{\sim}{\rightarrow} p^{\frac{1}{p-1}}\xi \hoR^{\oy}_U$ \eqref{THT3c}
et $c$ est induit par l'injection canonique $p^{\frac{1}{p-1}}\oR^{\oy}_U \rightarrow \oR^{\oy}_U$.
En vertu de \ref{cg9}, il existe un et un unique homomorphisme de $\co_{\oX_n}(\oU^\star)$-algèbres graduées
\begin{equation}\label{TFSLA23h}
\wedge (\xi^{-1}\tOmega^1_{\oX_n/\oS_n}(\oU^\star))\rightarrow 
\oplus_{i\geq 0}\rH^i(\pi_1(\oU^{\star \circ},\oy),\xi^{-i}\oR^{\oy}_U(i)/p^n\xi^{-i}\oR^{\oy}_U(i))
\end{equation}
dont la composante en degré un est le morphisme $a$. Son noyau est $\alpha$-nul et son conoyau est annulé par 
$p^{\frac{1}{p-1}}\fm_\oK$. On en déduit qu'il existe un et un unique homomorphisme de 
$\co_{\oX_n}(\oU^\star)$-algèbres graduées
\begin{equation}\label{TFSLA23i}
\wedge (\xi^{-1}\tOmega^1_{\oX_n/\oS_n}(\oU^\star))\rightarrow 
\oplus_{i\geq 0}\rH^i(\pi_1(\oU^{\star \circ},\oy),\oR^{\oy}_U/p^n\oR^{\oy}_U)
\end{equation}
dont la composante en degré $1$ est $\xi^{-1}A_{(U,\fp),n}$. 
Une chasse au diagramme \eqref{TFSLA23g} 
montre que le noyau de \eqref{TFSLA23i} est annulé par $p^{\frac{2d}{p-1}}\fm_\oK$.
Comme $\rH^i(\pi_1(\oU^{\star \circ},\oy),\oR^{\oy}_U/p^n\oR^{\oy}_U)$ 
est $\alpha$-nul pour tout $i\geq d+1$  en vertu de \ref{cg37}(iii), 
le conoyau de \eqref{TFSLA23i} est annulé par $p^{\frac{2d+1}{p-1}}\fm_\oK$.

\subsection{}\label{TFSLA27} 
Nous pouvons maintenant démontrer le théorème \ref{TFSLA24}. La proposition étant locale pour la topologie étale de $X_s$, 
soient $\ox$ un point géométrique de $X_s$, $X'$ le localisé strict de $X$ en $\ox$. 
D'après \ref{TFA13}(i), il existe un point  $(\oy\rightsquigarrow \ox)$ de $X_\et\gtimes_{X_\et}\oX^\circ_\et$ 
\eqref{rec15}. 
En vertu de \eqref{TFA14d} et \eqref{TFA12g}, pour tout entier $i\geq 0$, on a un isomorphisme canonique 
\begin{equation}\label{TFSLA24c}
\rR^i\sigma_{n*}(\xi\ocB_n)_\ox\stackrel{\sim}{\rightarrow}\rH^i(\oX'^\circ_\fet,\xi(\oR^\oy_{X'}/p^n\oR^\oy_{X'})). 
\end{equation}
La proposition résulte alors de \ref{TFSLA17}, \ref{TFSLA21}  et \ref{TFSLA23}.

\end{document}